%% file: main.tex
\documentclass[a4paper, 12pt]{amsart}

\usepackage[utf8]{inputenc}
\usepackage[OT2, T1]{fontenc}
\DeclareSymbolFont{cyrletters}{OT2}{wncyr}{m}{n}
\DeclareMathSymbol{\Sha}{\mathalpha}{cyrletters}{"58} 

\usepackage{fancyhdr}

\usepackage{subfiles} 

\usepackage{amsthm}
\usepackage{amssymb}
\usepackage{mathtools}
 \usepackage{stmaryrd}	
\usepackage{mathrsfs}
\usepackage{bm}
\usepackage{tikz-cd}
\usepackage{array} 
\usetikzlibrary{positioning, shapes.geometric, patterns, graphs, decorations.pathmorphing}

\usepackage{paralist}

\usepackage[colorlinks]{hyperref}
\hypersetup{
	linkcolor=blue,
}

\newcommand{\Z}{\ensuremath{\mathbb{Z}}} 
\newcommand{\bZ}{\ensuremath{\breve{\Z}}} 
\newcommand{\bE}{\ensuremath{\breve{E}}} 
\newcommand{\bF}{\ensuremath{\breve{F}}} 
\newcommand{\Q}{\ensuremath{\mathbb{Q}}} 
\newcommand{\bQ}{\ensuremath{\breve{\Q}_p}} 

\newcommand{\ur}{\ensuremath{\mathrm{ur}}} 
\newcommand{\unQ}{\ensuremath{\Q^{\ur}}} 
\newcommand{\R}{\ensuremath{\mathbb{R}}} 
\newcommand{\CC}{\ensuremath{\mathbb{C}}} 
\newcommand{\A}{\ensuremath{\mathbb{A}}} 
\newcommand{\F}{\ensuremath{\mathbb{F}}} 
\newcommand{\cF}{\ensuremath{\overline{\F}}} 

\newcommand{\BL}{\operatorname{BL}} 

\newcommand{\gr}{\operatorname{gr}} 
\newcommand{\Aut}{\operatorname{Aut}} 
\newcommand{\AUT}{\operatorname{\underline{Aut}}} 

\newcommand{\Ind}{\operatorname{Ind}} 

\renewcommand{\emptyset}{\ensuremath{\varnothing}}	

\newcommand{\Gal}{\operatorname{Gal}}  

\newcommand{\lrangle}[1]{\ensuremath{\langle #1 \rangle}}
\newcommand{\lrbracket}[1]{\ensuremath{\{ #1 \}}}

\newcommand{\wdt}[1]{\ensuremath{\widetilde{#1}}} 
\newcommand{\wdh}[1]{\ensuremath{\widehat{#1}}} 
\newcommand{\ovl}[1]{\ensuremath{\overline{#1}}} 

\newcommand{\identity}{\ensuremath{\mathrm{id}}}
\newcommand{\prolim}{\ensuremath{\underleftarrow{\lim}}}

\newcommand{\Hom}{\operatorname{Hom}}

\newcommand{\rightiso}{\ensuremath{\stackrel{\sim}{\rightarrow}}}

\newcommand{\Isom}{\operatorname{\underline{\mathrm{Isom}}}}

\newcommand{\VS}{\ensuremath{\mathscr{V}}} 

\newcommand{\mT}{\ensuremath{T}} 
\newcommand{\mG}{\ensuremath{G}} 
\newcommand{\mB}{\ensuremath{\mathcal{B}}} 
\newcommand{\mQ}{\ensuremath{\mathcal{Q}}} 
\newcommand{\mH}{\ensuremath{\mathcal{H}}} 

\newcommand{\mult}{\ensuremath{\mathrm{mult}}} 
\newcommand{\abe}{\ensuremath{\mathrm{ab}}} 

\newcommand{\mhg}[1]{\ensuremath{\mathbf{H}_{\mathrm{MH}}(#1)}} 
\newcommand{\mb}[1]{\ensuremath{\mathbf{H}_{\mathrm{B}}(#1)}} 
\newcommand{\mdr}[1]{\ensuremath{\mathbf{H}_{\mathrm{dR}}(#1)}} 

\newcommand{\betti}{\operatorname{\mathrm{B}}}
\newcommand{\MH}{\operatorname{\mathrm{MH}}}

\newcommand{\Bun}{\operatorname{Bun}}

\newcommand{\Lie}{\operatorname{Lie}}

\newcommand{\ad}{\ensuremath{\mathrm{ad}}}  
\newcommand{\der}{\ensuremath{\mathrm{der}}} 

\newcommand{\Spec}{\operatorname{Spec}}
\newcommand{\SPEC}{\operatorname{\underline{Spec}}} 

\newcommand{\Gm}{\ensuremath{\mathbb{G}_\mathrm{m}}}

\newcommand{\Res}{\operatorname{Res}}

\newcommand{\Supp}{\operatorname{Supp}}

\newcommand{\Adm}{\operatorname{Adm}} 

\newcommand{\dR}{\ensuremath{\mathrm{dR}}} 
\newcommand{\et}{\textit{\'et}} 
\newcommand{\crys}{\ensuremath{\mathrm{crys}}} 

\newcommand{\an}{\ensuremath{\mathrm{an}}} 

\newcommand{\Spf}{\operatorname{Spf}}
\newcommand{\dom}{\ensuremath{\mathrm{dom}}} 

\newcommand{\Spa}{\operatorname{Spa}}
\newcommand{\Spd}{\operatorname{Spd}}

\newcommand{\lcM}{\ensuremath{M^{\loc}}} 
\newcommand{\nby}{\ensuremath{\mathrm{R}\Psi}} 

\newcommand{\Fl}{\ensuremath{\mathrm{Fl}}} 
\newcommand{\fl}{\ensuremath{\mathcal{FL}}} 

\newcommand{\Def}{\ensuremath{\mathfrak{Def}}} 
\newcommand{\loc}{\ensuremath{\mathrm{loc}}} 
\newcommand{\ZIP}{\ensuremath{\mathrm{Zip}}} 
\newcommand{\Sht}{\ensuremath{\mathrm{Sht}}} 
\newcommand{\XX}{\ensuremath{\mathcal{X}}} 
\newcommand{\XXp}{\ensuremath{\mathbb{X}}} 

\newcommand{\stein}{\ensuremath{\mathrm{st}}} 

\newcommand{\GL}{\operatorname{GL}}

\newcommand{\GSp}{\operatorname{GSp}}
\newcommand{\GSP}{\operatorname{\mathcal{GSP}}}
\newcommand{\GLL}{\operatorname{\mathcal{GL}}}

\newcommand{\PGL}{\operatorname{PGL}}

\newcommand{\Char}{\ensuremath{\operatorname{char}}} 
\newcommand{\OO}{\ensuremath{\mathcal{O}}} 

\newcommand{\DS}{\ensuremath{\mathbb{S}}} 
\newcommand{\tor}{\ensuremath{\mathrm{tor}}} 
\newcommand{\Shum}[1]{\ensuremath{\mathscr{S}_{#1}}} 
\newcommand{\shu}[1]{\ensuremath{\mathrm{Sh}_{#1}}} 
\newcommand{\shuc}[1]{\ensuremath{\mathrm{Sh}_{#1}^{\tor}}}
\newcommand{\shum}[1]{\ensuremath{\mathrm{Sh}_{#1}^{\min}}}
\newcommand{\Shumc}[1]{\ensuremath{\mathscr{S}_{#1}^{\tor}}}
\newcommand{\Shumm}[1]{\ensuremath{\mathscr{S}_{#1}^{\min}}}

\newcommand{\Cusp}{\operatorname{Cusp}} 

\newcommand{\bigsur}{\ensuremath{\star}} 

\newcommand{\imin}[1]{\ensuremath{i_{#1}^{\mathrm{min}}}}
\newcommand{\itor}[1]{\ensuremath{i_{#1}^{\tor}}}
\newcommand{\Zb}{\ensuremath{\mathrm{Z}}} 

\newcommand{\NE}{\ensuremath{\mathcal{N}}} 
\newcommand{\CE}{\ensuremath{\mathcal{C}}} 
\newcommand{\KR}{\ensuremath{\mathrm{KR}}} 
\newcommand{\EO}{\ensuremath{\mathrm{EO}}} 
\newcommand{\EKOR}{\ensuremath{\mathrm{EKOR}}} 
\newcommand{\Ig}{\ensuremath{\mathrm{Ig}}} 
\newcommand{\IGUSA}{\ensuremath{\mathfrak{Ig}}} 

\newcommand{\GG}{\ensuremath{\mathcal{G}}} 
 
\newcommand{\LL}{\ensuremath{\mathcal{L}}}

\newcommand{\LLc}{\ensuremath{\mathcal{L}^{\circ}}}

\newcommand{\GGc}{\ensuremath{\mathcal{G}^{\circ}}}
\newcommand{\GGf}{\ensuremath{\mathcal{G}^{\flat}}}

\newcommand{\KK}{\ensuremath{\Breve{K}}} 
\newcommand{\KKc}{\ensuremath{\Breve{K}^{\circ}}} 

\newcommand{\Kc}{\ensuremath{K}^{\circ}} 
\newcommand{\Kk}{\ensuremath{K}} 
\newcommand{\Kf}{\ensuremath{K}^{\flat}} 

\newcommand{\red}{\ensuremath{\textrm{red}}} 
\newcommand{\ext}{\ensuremath{\textrm{ext}}} 
\newcommand{\perf}{\ensuremath{\textrm{perf}}} 

\newcommand{\Isoc}{\operatorname{Isoc}}

\newcommand{\mM}{\ensuremath{\mathbb{M}}} 

\newcommand{\Bui}{\operatorname{\mathcal{B}}} 
\newcommand{\aff}{\ensuremath{\textrm{aff}}} 

\newcommand{\ab}{\ensuremath{\mathcal{A}}} 
\newcommand{\DD}{\operatorname{\mathbb{D}}} 

\newcommand{\ZZ}{\ensuremath{\mathfrak{Z}}} 

\newcommand{\dd}{\ensuremath{\mathrm{\ddagger}}} 

\theoremstyle{plain}
\newtheorem{proposition}{Proposition}
\newtheorem{lemma}[proposition]{Lemma}
\newtheorem{theorem}[proposition]{Theorem}
\newtheorem{corollary}[proposition]{Corollary}
\newtheorem{assumption}[proposition]{Assumption}

\theoremstyle{definition}
\newtheorem{definition}[proposition]{Definition}
\newtheorem{definition-theorem}[proposition]{Definition-Theorem}
\newtheorem{definition-proposition}[proposition]{Definition-Proposition}
\newtheorem{remark}[proposition]{Remark}

\newtheorem{assumptionaxiom}[proposition]{Assumption}

\theoremstyle{definition}

\theoremstyle{plain}

\numberwithin{equation}{section}
\numberwithin{table}{section}
\numberwithin{proposition}{section}
\numberwithin{conj}{section}	

\usepackage{amsmath}

\usepackage{geometry}
\usepackage{lmodern}

\title[Well-positioned subschemes]{Compactifications of subschemes of integral models of Hodge-type Shimura Varieties with Parahoric level structures}
\author[Shengkai Mao]{Shengkai Mao}
\address{Morningside Center of Mathematics}
\email{maoshengkai@amss.ac.cn}

\date{}

\begin{document}

\begin{abstract}

    We prove that central leaves, Igusa varieties, Newton strata, Kottwitz-Rapoport Strata, Ekedahl-Kottwitz-Oort-Rapoport strata on the special fiber of a Kisin-Pappas integral model of a Hodge-type Shimura variety with connected parahoric level structure are well-positioned, generalizing the work of Lan and Stroh \cite{lan2018compactifications}.
    
\end{abstract}

\maketitle

\tableofcontents

\section{Introduction}

\subfile{sections/introduction}

\section{Well-positioned subschemes}\label{sec: well-positioned subschemes}

\subfile{sections/well_position}

\section{Various stratification}

\subfile{sections/various_stratification}

\section{Compactifications of Hodge-type Shimura varieties}\label{sec: compactifications of HT Shimura variety}

\subfile{sections/compactifications_of_hodge_type}

\section{Newton strata and Central leaves}\label{sec: Newton strata}

\subfile{sections/newton_strata}

\section{Igusa varieties}\label{sec: Igusa varieties}

\subfile{sections/igusa_varieties}

\section{Kottwitz-Rapoport strata}\label{sec: KR strata}

\subfile{sections/KR_strata}

\section{Ekedahl-Oort strata}\label{sec: EKOR strata}

\subfile{sections/EKOR_strata}

\section{Closure of Irreducible components of EKOR strata}\label{sec: closure of irreducible components}

\subfile{sections/closure}

\section{Change-of-Parahorics}\label{sec: change of parahorics}

\subfile{sections/change_of_parahorics}


\bibliographystyle{abbrv} 
\bibliography{compactifications}

\end{document}

%% file: sections/introduction.tex
The reduction modulo $p$ of integral models of Shimura varieties has very rich geometric properties. There are several types of weak stratifications on them, specifically central leaves, Newton strata, Kottwitz-Rapoport strata, Ekedahl-Oort strata, and Ekedahl-Kottwitz-Oort-Rapoport strata. In this article, we study the toroidal and minimal compactifications of these strata.

Consider a Hodge-type Shimura datum denoted as $(G, X)$. Let $K_p=\GG(\Z_p)$ be a stablizer quasi-parahoric subgroup, here $\GG:=\GG_x$ refers to the Bruhat-Tits stabilizer group scheme associated with a point $x$ in the extended Bruhat-Tits building.

In \cite{kisin2018integral}, Kisin and Pappas constructed integral models denoted as $\Shum{K}(G, X)$ for Hodge-type Shimura varieties. The local properties of these models can be understood through the local models $\lcM$ constructed by Pappas and Zhu in \cite{pappas2013local}. The assumption that $G_{\Q_p}$ is tamely ramified over $\Q_p$ is largely relaxed in \cite{kisin2024independence} and \cite{kisin2024integral}.

In \cite{he2017stratifications}, an axiomatic approach was introduced to systematically characterize various strata on the special fiber, denoted as $\Shum{K, s}(G, X)$, of integral models for Shimura varieties with parahoric level structures. In the context of the Kisin-Pappas integral models, Rong Zhou verified all axioms except for \cite[Axiom 4(c)]{he2017stratifications}. In addition, Zhou verified \cite[Axiom 4(c)]{he2017stratifications} when the group $G$ is residually split. The remaining \cite[Axiom 3.7(c)]{he2017stratifications} was verified by Ian Gleason, Dong Gyu Lim, Yujie Xu in \cite{gleason2023connected}. In \cite{shen2021ekor}, Xu Shen, Chia-Fu Yu, and Chao Zhang further established additional geometric properties of these strata.

On the other hand, Shimura varieties and their integral models are known to have good toroidal and minimal compactifications, see \cite{faltings2013degeneration}, \cite{lan2013arithmetic}, \cite{pera2019toroidal}. It is expected that these strata on the mod-$p$ fibers also have good toroidal and minimal compactifications. In the paper \cite{boxer2015torsion}, Boxer proved that the Ekedahl-Oort stratification can be extended to the minimal compactification for PEL-type Shimura varieties of type A or C, and the Hasse-invariants can also be extended to the boundary in a compatible manner. In Hodge-type case, this was done by Goldring and Koskivirta in \cite{goldring2019strata}. 

In \cite{lan2018compactifications}, Lan and Stroh extended and generalized these concepts and ideas to a broader range of integral models for Shimura varieties and a broader types of strata. They proved that various strata, specifically $p$-rank strata, Newton strata, central leaves, Kottwitz-Rapoport strata, and Ekedahl-Oort strata are \emph{well-positioned} in the sense of \cite[Definition 2.2.1]{lan2018compactifications}, also see Definition \ref{def: well-positioned}. 

\subsection{}

In this article, we generalize some results of \cite{lan2018compactifications} to various strata on the special fiber of integral models of Hodge-type Shimura varieties with connected parahoric level structures, and study their boundary behaviors. 

 \begin{theorem}[Propositions \ref{proposition: newton strata are well-positioned}, \ref{proposition: central leaves are well--positioned}, \ref{prop: KR strata are well-positioned}, \ref{prop: EKOR are well-positioned}]\label{main theorem, KPZ}
     Let $\Shum{K}(G, X)$ be a Kisin-Pappas integral model of a Hodge-type Shimura variety with connected parahoric level structure, then Newton strata, central leaves, KR strata, EKOR strata and their connected components and corresponding pullbacks to any higher level are well-positioned.
 \end{theorem}
 For Kisin-Pappas-Zhou integral models constructed in  \cite{kisin2024independence} and \cite{kisin2024integral}, see Remark \ref{rmk: KPZ models vs KP models} 

We can also consider Pappas-Rapoport integral models with general quasi-parahoric subgroups constructed in \cite{pappas2024p} ($K_p$ being stablizer quasi-parahoric) and \cite{daniels2024conjecture} ($K_p$ being general quasi-parahoric). Such models do not rely on the conditions on $(p, G)$, they have good functoriality and are independent of the choice of Siegel embeddings. One can also define Newton strata, central leaves, KR strata on Pappas-Rapoport integral models, as in \cite[Remark 4.2.1]{pappas2024p}.

 \begin{theorem}[Propositions \ref{prop: Newton central are well-posiitoned, PR}, \ref{prop: KR strata are well-positioned}]\label{main theorem, PZ}
     Let $\Shum{K}(G, X)$ be a Pappas Rapoport integral model of a Hodge-type Shimura variety with general quasi-parahoric level structure, then Newton strata, central leaves, KR strata (if the schematic local model diagram exists) and and their connected components and corresponding pullbacks to any higher level are well-positioned.
 \end{theorem}
 
In particular, each of these strata $Y \subset \Shum{K, \bar{s}}$ ($\bar{s} = \Spec \ovl{\F}_p$) have a good theory of compactifications, they have partial toroidal compactifications $Y^{\tor}_{\Sigma} \subset \Shumc{K, \Sigma, \bar{s}}(G, X)$ and partial minimal compactifications $Y^{\min} \subset \Shumm{K, \bar{s}}$, and one can describe their boundary strata in an axiomatic way, see \cite[Theorem 2.3.2]{lan2018compactifications}, also see Proposition \ref{prop: Y same as Shum}. Roughly speaking, compactifications of well-positioned subschemes exhibit behaviors similar to the compactifications of Shimura varieties over fields of characteristic 0, and the local properties of well-positioned subschemes extend to their boundaries.

\subsection{}

Let $\varphi: \shu{K}(G, X) \to \shu{K^{\dd}}(G^{\dd}, X^{\dd})$ be the Siegel embedding introduced in \cite{kisin2018integral}, \cite{kisin2024independence}, \cite{kisin2024integral} or in \cite{pappas2024p} and \cite{daniels2024conjecture}, and $\Shum{K}(G, X)$ is defined as the relative normalization of the integral model of $\Shum{K^{\dd}}(G^{\dd}, X^{\dd})$ in $\shu{K}(G, X)$. Let us denote $\varphi: \Shum{K}(G, X) \to \Shum{K^{\dd}}(G^{\dd}, X^{\dd})$ the induced finite morphism. In Kisin-Pappas integral models, we take $K_p$ be stablizer quasi-parahoric subgroup. In Pappas-Rapoport integral models, we take $K_p$ be a general quasi-parahoric subgroup.
 
 \begin{proposition}[{Proposition \ref{proposition: newton strata are well-positioned}, \ref{proposition: central leaves are well--positioned}, \ref{prop: Newton central are well-posiitoned, PR}}]\label{prop: intro, discrete fiber}
    Under $\varphi$, the pullbacks of Newton strata, and central leaves are topologically disjoint unions of Newton strata, and central leaves, respectively.
 \end{proposition}
 Here we say a locally closed subset (resp. subscheme) $Y = \bigsqcup_{i \in I} Y_i$ of $X$ is a topological disjoint union of $\lrbracket{Y_i}_{i \in I}$ if under the subspace topology, each $Y_i$ is open and closed in $Y$.

 \subsection{}

  Note that the finite morphism $\varphi: \Shum{K} \to \Shum{K^{\dd}}$ extend to a finite morphism $\varphi^{\min}: \Shumm{K} \to \Shumm{K^{\dd}}$ and a finite morphism $\Shumc{K, \Sigma} \to \Shumc{K^{\dd}, \Sigma^{\dd}}$ (when $\Sigma$ is induced by $\Sigma^{\dd}$). Well-positioned subsets (resp. subschemes) satisfy good functoiral properties.
 \begin{proposition}[{Proposition \ref{lemma: pullback of well-positioned is well-positioned}, \ref{proposition: open-closed subschemes are well positioned}}]\leavevmode
     \begin{enumerate}
         \item Let $Y^{\dd} \subset \Shum{K^{\dd}, \bar{s}}$ be well-positioned, then its pullback $Y = \varphi^{-1}(Y^{\dd}) \subset \Shum{K, \bar{s}}(G, X)$ is well-positioned, and moreover $Y^{\tor}_{\Sigma} = \varphi^{\tor, -1}(Y^{\dd, \tor}_{\Sigma^{\dd}})$, $Y^{\min} = \varphi^{\min, -1}(Y^{\dd, \min})$.
         \item Let $Y = Y_1 \sqcup Y_2 \subset \Shum{K, \bar{s}}$ be a topologically disjoint union. If $Y$ is well-positioned, then $Y_1$, $Y_2$ are well-positioned, and $Y^{\tor}_{\Sigma} = Y^{\tor}_{1, \Sigma} \sqcup Y^{\tor}_{2, \Sigma}$, $Y^{\min} = Y^{\min}_1 \sqcup Y^{\min}_2$ are topologcally disjoint unions.
     \end{enumerate}
 \end{proposition}

   Given the integral canonical model $\Shum{K^{\dd}}(G^{\dd}, X^{\dd})$ of a Siegel variety with hyperspecial level structure, Boxer proved that EO strata are well-positioned in \cite{boxer2015torsion}, and Lan and Stroh proved that Newton strata and central leaves are well-positioned in \cite{lan2018compactifications}. Consequently, well-position of Newton strata, central leaves and EO strata in Theorems \ref{main theorem, KPZ} and \ref{main theorem, PZ} can be derived from the above proposition. 
   
 Moreover, since the partial minimal compactifications of EO strata and central leaves on the integral model of Siegel modular varieties are affine, see \cite[Theorem C]{boxer2015torsion} and \cite[Proposition 1.9]{caraiani2024generic} (\cite[Theorem 3.3.3]{santos2023generic}) respectively, we have

 \begin{corollary}\leavevmode
\begin{enumerate}
    \item Minimal compactifications of EO strata are affine, see Corollary \ref{cor: min EO are affine} (\cite[Proposition 6.3.1(c)]{goldring2019strata}).
    \item Minimal compactifications of central leaves are affine, see Corollary \ref{cor: minimal compactifications of central leaves are affine}.
\end{enumerate}
\end{corollary}

\subsection{}

In Kisin-Pappas integral models with stablizer quasi-parahoric subgroup, we approach Igusa varieties in a similar manner. Let $\CE^b \subset \Shum{K, \Bar{s}}(G, X)$ be a central leaf and $\CE^{b^{\dd}} \subset \Shum{K^{\dd}, \bar{s}}(G^{\dd}, X^{\dd})$ be the unique central leaf containing $\varphi(\CE^b)$. In \cite{hamacher2019adic}, Kim and Hamacher have proven that the Igusa variety $\Ig_{b}$ over $\CE^b$, defined by Hodge tensors related to $\GG$, is a closed union of connected components in the pullback of the Igusa variety $\Ig_{b^{\dd}}$ over $\CE^{b^{\dd}}$ under $\CE^b \to \CE^{b^{\dd}}$, we prove that this closed union and its complement do not intersect along the boundary. Since the Igusa variety over $\CE^{b^{\dd}}$ are well-positioned due to the work of \cite{caraiani2024generic} and \cite{santos2023generic}, we have
\begin{proposition}[{\ref{prop: Igusa varieties are well-positioned}, \ref{cor: min Igusa are affine}, \ref{cor: Ig tor closed embedding}, \ref{prop: pink to Igusa varieties}}]
    Igsua variety $\Ig_b$ is well-positioned in the sense of Definition \ref{def: well-positioned over well-positioned}. Moreover,
    \begin{enumerate}
        \item $\Ig_b^{\min}$ is affine.
        \item $\Ig^{\tor}_b \hookrightarrow \Ig^{\tor}_{b^{\dd}}$ is a closed embedding when $(\Sigma, \Sigma^{\dd})$ are chosen to be strictly compatible in the sense of \cite[Definition 4.6]{lan2022closed}.
        \item One can define Igusa variety $\Ig_{K(l^m), b}$ over $\CE^b \subset \Shum{K(l^m), \bar{s}}$ for $m\geq 0$ and a good prime $l\nmid p$, then $\varinjlim_m R\Gamma(\Ig_{K(l^m), b}^{\tor}, \F_l) \rightiso \varinjlim_m R\Gamma(\Ig_{K(l^m), b}, \F_l)$.
    \end{enumerate}
\end{proposition}

\subsection{}

 In Kisin-Pappas integral models, assume $\GG = \GGc$ is connected. Use the strategy in \cite{gortz2012supersingular}, with a crucial input being the quasi-affineness of KR strata which was proved in \cite[Theorem 3.5.9]{shen2021ekor}, we are able to prove the following corollary, which gives an affirmative answer to \cite[Conjecture 4.3.1]{van2024mod} under some mild conditions:
  \begin{corollary}[Corollary \ref{cor: minimal EKOR has zero dim}, \ref{corollary: Hoften 3.7.5 removed}]
     Assume $G^{\ad}$ is $\Q$-simple or has no $\PGL_2$-factor\footnote{We expect this condition can be easily removed once the theory of well-positioned subschemes in abelian-type Shimura varieties becomes available}. Let $Y$ be any irreducible component of the closure of an EKOR stratum at any parahoric level, then $Y$ intersects nontrivially with the minimal EKOR stratum. The minimal EKOR has $0$ dimension and has no boundary.

 \end{corollary}

\subsection{}
Let $(G, X)$ be a Hodge-type Shimura datum. Let us work with Pappas-Rapoport integral models. Let $K_{1, p} \subset K_{2, p}$ be two quasi-parahoric subgroups of $G(\Q_p)$, choose $K_1^p \subset K_2^p$, let $K_1 = K_{1, p}K_1^p$, $K_2 = K_{2, p}K_2^p$, then we have a morphism $\pi: \Shum{K_1} \to \Shum{K_2}$ which extends the finite morphism $\shu{K_1} \to \shu{K_2}$, due to \cite[Corollary 4.3.2]{pappas2024p} and \cite[Corollary 4.1.10]{daniels2024conjecture}. In \cite[Theorem 4.5]{mao2025boundary}, we show that $\pi$ extends uniquely to proper morphisms $\pi^{\min}: \Shumm{K_1} \to \Shumm{K_2}$ and $\pi^{\tor}: \Shumc{K_1, \Sigma_1} \to \Shumc{K_2, \Sigma_2}$ (let $\Sigma_1$ be induced by $\Sigma_2$), which in turns imply that $\pi$ is proper, and we describe the structures of $\pi^{\min}$ and $\pi^{\tor}$ in details in \cite[Proposition 4.11]{mao2025boundary}. This allows us to consider the change-of-parahoric morphism for well-positioned subschemes:
\begin{proposition}[{Proposition \ref{prop: extension of pullback, well-positioned}}]
    Let $Y_1 \subset \Shum{K_1, T}$ and $Y_2 \subset \Shum{K_2, T}$ be well-positioned subscheme such that $\pi: \Shum{K_1} \to \Shum{K_2}$ induces a morphism $\pi_Y: Y_1 \to Y_2$, then $\pi$ extends to morphisms 
    \begin{equation}
     \pi^{\min}_Y: Y^{\min}_1 \to Y^{\min}_2 \quad \pi^{\tor}_Y: Y^{\tor}_{1, \Sigma_1} \to Y^{\tor}_{2, \Sigma_2}.
    \end{equation}
    \begin{enumerate}
       \item If $\pi_Y$ is proper, then $\pi^{\min}_Y$ and $\pi^{\tor}_Y$ are proper.
        \item  If $\pi_Y$ is proper surjective, then $\pi^{\min}_Y$ and $\pi^{\tor}_Y$ are proper surjective.
        \item If $\pi_Y$ is smooth (resp. \'etale), then $\pi^{\tor}_Y$ is smooth (resp. \'etale).
        \item If $\pi_Y^{\min}$ is quasi-finite (resp. finite), then $\pi_Y^{\tor}$ is quasi-finite (resp. finite).
    \end{enumerate}
\end{proposition}

\begin{corollary}[{Corollary \ref{cor: central leaves extends, levels}, \ref{cor: EKOR extends, level}}]
    Let $K_{1, p} \subset K_{2, p}$ be parahoric subgroups. Assume the axioms in \cite{he2017stratifications} are varified.
    \begin{enumerate}
        \item Given central leaves $\CE^b_{K_1} \subset \Shum{K_1, \bar{s}}$, $\CE^b_{K_2} \subset \Shum{K_2, \bar{s}}$, then the finite surjection $\pi_b: \CE^b_{K_1} \to \CE^b_{K_2}$ extends to finite surjections
        \[ \pi_b^{\min}: \CE^{b, \min}_{K_1} \to \CE^{b, \min}_{K_2},\quad \pi_b^{\tor}: \CE^{b, \tor}_{K_1, \Sigma_1} \to \CE^{b, \tor}_{K_2, \Sigma_2}. \]
       If $\pi_b$ is \'etale, then $\pi_b^{\tor}$ is \'etale.
        \item Given $w \in \prescript{\KK_{2, p}}{}{\Adm(\lrbracket{\mu})} \subset \prescript{\KK_{1, p}}{}{ \Adm(\lrbracket{\mu})}$, the finite surjection $ \pi_{w}:  \EKOR_{K_1, w} \to \EKOR_{K_2, w}$ extends to finite surjections 
        \[ \pi_{w}^{\min}:  \EKOR_{K_1, w}^{\min} \to \EKOR_{K_2, w}^{\min},\quad \pi_{w}^{\tor}:  \EKOR_{K_1, w, \Sigma_1}^{\tor} \to \EKOR_{K_2, w, \Sigma_2}^{\tor}. \]
        If $\pi_{w}$ is \'etale, then $\pi_{w}^{\tor}$ is \'etale.
    \end{enumerate}
\end{corollary}

\subsection{}

   In PEL-type cases, Lan and Stroh showed that the boundary strata of central leaves, Newton strata, KR strata, and EO strata are again central leaves, Newton strata, KR strata, and EO strata, respectively, on the special fibers of integral models of smaller PEL-type Shimura varieties. In contrast to PEL-type cases, Hodge-type Shimura varieties lack moduli interpretations, and the descriptions of boundary strata are primarily described in a group-theoretical manner. 
   
       \begin{theorem}[{\cite[Theorem 3.58]{mao2025boundary}}]
           Minimal compactification of Pappas-Rapoport integral model (resp. Kisin-Pappas integral models) of Shimura variety of Hodge type with quasi-parahoric level structure is stratified by finite quotients of Pappas-Rapoport integral models (resp. Kisin-Pappas integral model) of smaller Shimura variety of Hodge type with quasi-parahoric level structure (resp. under some extra conditions).
       \end{theorem}

 Recall that the minimal compactification $\Shumm{K}(G, X)$ is stratified by:
 \begin{equation*}
     \bigsqcup_{[\Phi]  \in \Cusp_K(G, X)} \Lambda_{\Phi, K}\backslash\Shum{K_{\Phi, h}}(G_h, X_{\Phi, h}),
 \end{equation*}
 where $\Cusp_K(G, X)$ denotes the set of cusp labels, $\Shum{K_{\Phi, h}}(G_h, X_{\Phi, h})$ is the integral model of a Hodge-type Shimura variety obtained by taking relative normalization, and $\Lambda_{\Phi, K}$ is a group that acts on $\Shum{K_{\Phi, h}}(G_h, X_{\Phi, h})$ and factors through a finite quotient.

 We provide some partial answers regarding the description of the boundary of a specific stratum on $\Shum{K, \Bar{s}}(G, X)$. Let $Y \subset \Shum{K, \Bar{s}}(G, X)$ be a well-positioned subscheme,
\[ Y_{[\Phi]} = Y^{\min} \times_{\Shumm{K, \Bar{s}}(G, X)} (\Lambda_{\Phi, K}\backslash\Shum{K_{\Phi, h}, \Bar{s}}(G_h, X_{\Phi, h})) \] 
be its boundary stratum with respect to $[\Phi] \in \Cusp_K(G, X)$. Let $Y^{\natural}(\Phi) \subset \Shum{K_{\Phi, h}, \Bar{s}}(G_h, X_{\Phi, h})$ be the pullback of $Y_{[\Phi]}$. 
  \begin{proposition}[{Proposition \ref{prop: boundary of Newton strata}, \ref{prop: boundary of central leaves}, \ref{prop: boundary of KR is KR}, \ref{prop: boundary of EO strata}}]\leavevmode
     In Kisin-Pappas integral models,
        \begin{enumerate}
            \item If $Y$ is a Newton (resp. Kottwitz-Rapoport, Ekedahl-Oort) stratum, then $Y^{\natural}(\Phi)$ is also a Newton (resp. Kottwitz-Rapoport, Ekedahl-Oort) stratum. 
            \item If $Y$ is a central leaf, then $Y^{\natural}(\Phi)$ is a topologically disjoint finite union of central leaves.
        \end{enumerate}
\end{proposition}

\subsection*{Acknowledgement}

This article is based on the second half of the author's thesis, with many revisions and strengthened results. The author thanks his advisor, Kai-Wen Lan, for introducing this project and for his generous guidance and support. The author also thanks Alex Youcis, Peihang Wu, Pol van Hoften and Wansu Kim for many helpful discussions and expresses his gratitude to Georgios Pappas and Rong Zhou for their encouragement and assistance.

%% file: sections/well_position.tex
The set-up of well-positioned subschemes of Shimura varieties is axiomatic. We start with axiomatic descriptions of compactifications of integral models of Shimura varieties, and then work on general setting of well-position.

\subsection{Compactifications of Shimura varieties}

Let $\iota: (G_1, X_1) \to (G_2, X_2)$ be a morphism of Shimura data, $K_1 \subset G_1(\A_f)$, $K_2 \subset G_2(\A_f)$ be open compact subgroups.
\begin{assumption}\label{assumption: general set up}
    Suppose the morphism $(G_1, X_1, K_1) \to (G_2, X_2, K_2)$ is in one of the cases \cite[Assumption 2.1]{lan2022closed} (with less restrictive condition $K_1 \subset \iota^{-1}(K_2)$). In particular, this includes case $(4)$, where $(G_1, X_1) = (G, X)$ is a Hodge-type Shimura datum, $(G_2, X_2) = (G^{\dd}, X^{\dd})$ is a Siegel-type Shimura datum, $(G_1, X_1) \hookrightarrow (G_2, X_2)$ is a Siegel embedding, $K_{2, p} = K^{\dd}_p$ is hyperspecial, and $K_1 \subset K_2$.
\end{assumption}

Moreover, in later subsections, we consider the following case:
\begin{assumption}\label{assumption: adjusted Siegel embedding}
    Case $(4)'$: $(G_1, X_1, K_{p, 1}) \hookrightarrow (G_2, X_2, K_{p, 2})$ is an adjusted Siegel embedding in the sense of \cite[Definition 3.16, Remark 3.17]{mao2025boundary}. 
\end{assumption}

The Siegel embeddings appeared in \cite{kisin2018integral}, \cite{kisin2024independence}, \cite{kisin2024integral}, \cite{pappas2024p} all satisfied the above assumption after replacing the Siegel datum with a larger Siegel datum while combining chain of almost self-dual lattices to a self-dual lattice. We post this condition in order to use the constructions in \cite{pera2019toroidal} as well as the constructions for good integral models.

Under assumption \ref{assumption: general set up}, we have axiomatic descriptions of compactifications of Shimura varieties over a base scheme $S$ which was summarized as follows (which we recall is not a full list). In case ($4$), we take $S = \Spec \OO_{E(v)}$, where $E = E(G, X)$ is the reflex field, $v$ is a place over a fixed prime $p$.
\begin{theorem}{{\cite[Proposition 3.1]{lan2022closed}}}\label{theorem: axiomatic descriptions of compactifications}

    For each $i = 1, 2$, fix some projective smooth cone decompositions $\Sigma_i$. There are proper surjective morphisms between noetherian schemes $\oint_{K_i, \Sigma_i}: \Shumc{K_i, \Sigma_i}\to\Shumm{K_i}$ over $S$ that are functorial with respect to $\Sigma_i$. Denote by $\imin{K_i}: \Shum{K_i} \to \Shumm{K_i}$, $\itor{K_i, \Sigma_i}: \Shum{K_i} \to \Shumc{K_i, \Sigma_i}$, then $\oint_{K_i, \Sigma_i}\circ\imin{K_i} = \itor{K_i, \Sigma_i}$.

    The compactification $\Shumc{K_i, \Sigma_i}$ (resp. $\Shumm{K_i}$) is stratified by locally closed subschemes $\mathcal{Z}_i(\Upsilon_i)$ (resp. $\mathcal{Z}_i([\Phi_i])$) with respect to $\Upsilon_i \in \Cusp_{K_i}(\Sigma_i):=\Cusp_{K_i}(G_i, X_i, \Sigma_i)$ (resp. $[\Phi_i] \in \Cusp_{K_i}:=\Cusp_{K_i}(G_i, X_i)$), and has an open dense subscheme $\Shum{K_i}$. Fix a cusp label representative $\Phi_i$ in each equivalence class $[\Phi_i]$, then those $\Upsilon_i = [\Phi_i, \sigma_i] \in \Cusp_{K_i}(\Sigma_i)$ over $[\Phi_i]$ are indexed by $[\sigma_i] \in \Lambda_{K_i, \Phi_i}\backslash\Sigma_i(\Phi_i)^+$.
    
    The preimage $\oint_{K_i, \Sigma_i}^{-1}(\mathcal{Z}_i([\Phi_i]))=\bigsqcup_{[\sigma_i] \in \Lambda_{K_i, \Phi_i}\backslash\Sigma_i(\Phi_i)^+} \mathcal{Z}_i([\Phi_i, \sigma_i])$ (set theoretically) and the closure of $\mathcal{Z}_i([\Phi_i, \sigma_i])$ (resp. $\mathcal{Z}_i([\Phi_i])$) is the union of $\mathcal{Z}_i([\Phi'_i, \sigma'_i])$ (resp. $\mathcal{Z}_i([\Phi'_i])$) such that $[\Phi'_i, \sigma'_i] \preceq [\Phi_i, \sigma_i]$ (resp. $[\Phi'_i] \preceq [\Phi_i]$).
   
   For each $\Zb^{\bigsur}_i(\Phi_i)$, there is a projective surjective morphism $C_i(\Phi_i)\to\Zb_i^{\bigsur}(\Phi_i)$. Also, there is a $\mathbf{E}_{i, K_i}(\Phi_i)$-torsor $\Xi_i(\Phi_i)\to C_i(\Phi_i)$ where $\mathbf{E}_{i, K_i}(\Phi_i)$ is the pullback of a split torus over $\Spec \Z$. For each $\sigma_i\in\Sigma_i(\Phi_i)$, there is an affine toroidal embedding $\Xi_i(\Phi_i)\hookrightarrow\Xi_i(\Phi_i, \sigma_i)$ and a specified closed subscheme $\Xi_{i, \sigma_i}(\Phi_i)$ of $\Xi_i(\Phi_i,\sigma_i)$. $\Xi_i(\Phi_i, \sigma_i)$ has a stratification by $\{\Xi_{i, \tau_i}(\Phi_i)\}_{\tau_i}$ where $\tau_i$ runs over the set of all faces of $\sigma_i$. For each $\sigma_i\in\Sigma_i(\Phi_i)^+$, there is a canonical isomorphism $\Xi_{i, \sigma_i}(\Phi_i)\cong\mathcal{Z}_i([\Phi_i, \sigma_i])$. The following diagram is commutative:

\begin{equation}\label{diag: strata in two sides}
\begin{tikzcd}
	{\mathcal{Z}_i([\Phi_i, \sigma_i])} &&& {\mathcal{Z}_i([\Phi_i])} \\
	{\Xi_{i, \sigma_i}(\Phi_i)} & {C_i(\Phi_i)} & {\Zb^{\bigsur}_i(\Phi_i)} & {\Zb_i(\Phi_i)}
	\arrow["{\oint_{K_i, \Sigma_i}}", from=1-1, to=1-4]
	\arrow["\cong", from=2-1, to=1-1]
	\arrow[from=2-1, to=2-2]
	\arrow[from=2-2, to=2-3]
	\arrow[from=2-3, to=2-4]
	\arrow["\cong"', from=2-4, to=1-4]
\end{tikzcd}
\end{equation}
   
        The local structures of toroidal compactifications along boundary strata $\mathcal{Z}_i([\Phi_i, \sigma_i])$ could be described as follows:

     	Fix a $\Phi_i$ from now on. We omit $\Phi_i$ in the subscript to save notions. Consider the full toroidal embedding and its formal completion:
     	\[ \overline{\Xi}_i=\underset{\sigma_i\in\Sigma_i(\Phi_i)^+}{\bigcup} \Xi_i( \sigma_i), \quad  \mathfrak{X}_i=(\overline{\Xi}_i)^{\wedge}_{\underset{\tau_i\in\Sigma_i(\Phi_i)^{+}}{\cup}\Xi_{i, \tau_i}}  \]
     	Note that $\mathfrak{X}_i$ has open coverings $\mathfrak{X}_{i, \sigma_i}^{\circ}$ for $\sigma_i$ running over all $\sigma_i\in\Sigma_i(\Phi_i)^{+}$, where
     	\[ \mathfrak{X}_{i, \sigma_i}^{\circ}=(\Xi_i(\sigma_i))_{\underset{\tau_i\in\Sigma_i(\Phi_i)^{+},\ \bar{\tau}_i\subset\bar{\sigma}_i}{\cup}\Xi_{i, \tau_i}}^{\wedge}.  \]
         Let $\mathfrak{X}_{i, \sigma_i} = \Xi_i(\sigma_i)_{\Xi_{i, \sigma_i}}^{\wedge}$, then the composition $\mathfrak{X}_{i, \sigma_i} \to \mathfrak{X}_{i, \sigma_i}^{\circ}\to\mathfrak{X}_i\to\Shumc{K_i, \Sigma_i}$ induce isomorphisms
     	\[ \mathfrak{X}_{i, \sigma_i}^{\circ}\cong(\Shumc{K_i, \Sigma_i})_{\underset{\tau_i\in \Sigma_i(\Phi_i)^{+},\ \bar{\tau}_i\subset\bar{\sigma}_i}{\cup} \mathcal{Z}_i([\Phi_i, \tau_i])}^{\wedge}, \quad \mathfrak{X}_i \cong (\Shumc{K_i, \Sigma_i})_{\underset{[\sigma_i] \in \Lambda_{K_i, \Phi_i}\backslash\Sigma_i(\Phi_i)^+}{\cup} \mathcal{Z}_i([\Phi_i, \sigma_i])}^{\wedge},\]
        which extends $\mathfrak{X}_{i, \sigma_i}\cong(\Shumc{K_i, \Sigma_i})_{\mathcal{Z}_i([\Phi_i, \sigma_i])}^{\wedge}$, and extends the canonical boundary isomorphism $\Xi_{i, \sigma_i}\cong\mathcal{Z}_i([\Phi_i, \sigma_i])$.

        For any affine open formal subscheme $\mathfrak{W}_i=\Spf(R_i)\subset\mathfrak{X}_{i, \sigma_i}^{\circ}$, we have the canonical induced flat morphisms:
   \begin{equation}\label{def: W}
   	     W_i:=\Spec R_i \to \Xi_i(\sigma_i), \quad W_i \to \Shumc{K_i, \Sigma_i}. 
   \end{equation}
      The preimages of the boundary stratifications on the targets coincide in $W_i$. In particular, the preimages of the open strata $\Xi_i$ and $\Shum{K_i}$ coincide in $W_i$, we denote it by $W^{0}_i$. These two morphisms \ref{def: W} are the key ingredients in the definition of \emph{well position} of a locally closed subscheme.
     \end{theorem}

   \begin{proposition}{{\cite[Proposition 3.4]{lan2022closed}}}\label{proposition: functorialities on toroidal compactifications}
          Keep assumptions and notations as in Theorem \ref{theorem: axiomatic descriptions of compactifications}. 
           \begin{enumerate}
               \item The finite morphism $\iota: \Shum{K_1} \to \Shum{K_2}$ extends uniquely to a finite morphism $\iota^{\min}: \Shumm{K_1} \to \Shumm{K^{\ddagger}}$ over $S$. For each boundary stratum $\mathcal{Z}_1([\Phi_1])$ of $\Shumm{K_1}$, there exists a unique boundary stratum $\mathcal{Z}_2([\Phi_2])$ ($\Phi_2 = \iota_*\Phi_1$) of $\Shumm{K_2}$ such that $\iota^{\min}(\mathcal{Z}_1([\Phi_1]))\subset \mathcal{Z}_2([\Phi_2])$. Moreover, $\mathcal{Z}_1([\Phi_1])$ is open and closed in $\iota^{\min, -1}(\mathcal{Z}_2([\Phi_2]))$.
               \item There exists a finite morphism $\Zb_1(\Phi_1) \to \Zb_2(\Phi_2)$ compatible with $\mathcal{Z}_1([\Phi_1]))\to \mathcal{Z}_2([\Phi_2])$, and it extends to finite morphisms $\Zb^{\bigsur}_1(\Phi_1) \to \Zb^{\bigsur}_2(\Phi_2)$, $C_1(\Phi_1) \to C_2(\Phi_2)$ and $\Xi_1(\Phi_1) \to \Xi_2(\Phi_2)$, and the induced finite morphism $\Xi_1(\Phi_1) \to \Xi_2(\Phi_2) \times_{C_2(\Phi_2)} C_1(\Phi_1)$ is equivariant with the pullback of a group homomorphism $\mathbf{E}_{1, K_1}(\Phi_1) \to \mathbf{E}_{2, K_2}(\Phi_2)$. If $K_1 = \iota^{-1}(K_2)$, then $\mathbf{E}_{1, K_1}(\Phi_1) \to \mathbf{E}_{2, K_2}(\Phi_2)$ is a closed immersion. 
                \item Let $\Sigma_1$ be induced by $\Sigma_2$, which means for any such pair $(\Phi_1, \Phi_2)$, each $\sigma_1 \in \Sigma_1(\Phi_1)^+$ is exactly the preimage of a unique $\sigma_2 \in \Sigma_2(\Phi_2)^+$, let us denote $\sigma_2 = \iota_*\sigma_1$. Note that there always exists projective smooth $\Sigma_1$ induced by $\Sigma_2$ due to \cite[Proposition 4.10]{lan2022closed}. Then $\iota$ induce morphisms $\Xi_1(\Phi_1, \sigma_1) \to \Xi_2(\Phi_2, \sigma_2)$, $\ovl{\Xi}_1 \to \ovl{\Xi}_2$, and induce finite morphisms $\mathfrak{X}_1 \to \mathfrak{X}_2$, $\mathfrak{X}_{1, \sigma_1} \to \mathfrak{X}_{2, \sigma_2}$, $\mathfrak{X}_{1, \sigma_1}^{\circ} \to \mathfrak{X}_{2, \sigma_2}^{\circ}$, $\mathcal{Z}_1([\Phi_1, \sigma_1]) \to \mathcal{Z}_2([\Phi_2, \sigma_2])$.
                \item Let $\Sigma_1$ be induced by $\Sigma_2$, then $\iota$ extends uniquely to a finite morphism $\iota^{\tor}_{\Sigma_1, \Sigma_2}: \Shumc{K_1, \Sigma_1} \to \Shumc{K_2, \Sigma_2}$ over $S$ such that $\iota^{\min}\circ\imin{K_1} = \imin{K_2}\circ\iota$,  $\iota^{\min}\circ\oint_{K_1, \Sigma_1} = \oint_{K_2, \Sigma_2}\circ\iota^{\tor}_{\Sigma_1, \Sigma_2}$.
              \item Moreover, let $\mathfrak{W}_2=\Spf(R_2)\subset \mathfrak{X}^{\circ}_{2, \sigma_2}(\Phi_2)$ be any affine open formal subscheme, $\mathfrak{W}_1=\Spf(R_1)$ be the pullback of $\mathfrak{W}_2$ to $\mathfrak{X}^{\circ}_{1, \sigma_1}(\Phi_1)$. Let $W_2=\Spec(R_2)\to \Shumc{K_2, \Sigma_2}$, $W_2\to \Xi_2(\Phi_2, \sigma_2)$, $W_1=\Spec(R_1) \to \Shumc{K_1, \Sigma_1}$, $W_1\to \Xi(\Phi_1, \sigma_1)$ be the induced morphisms introduced in Theorem \ref{theorem: axiomatic descriptions of compactifications}. Then the preimages of $\Shum{K_2}$ and $\Xi_2(\Phi_2)$ (resp. $\Shum{K_1}$ and $\Xi_1(\Phi_1)$) coincide as an open subscheme $W^{0}_2$ (resp. $W^0_1$) of $W_2$ (resp. $W_1$), and the pullback of $W^{0}_2$ to $W_1$ coincides with $W^0_1$.  
  
           \end{enumerate}
       \end{proposition}

\subsection{Well-positioned subschemes}\label{sec: Well-positioned subscheme}

   In this subsection, we recall the definition of a well-positioned subscheme.

    \begin{definition-proposition}{{\cite[Definition 2.2.1, Lemma 2.2.2]{lan2018compactifications}}}\label{def: well-positioned}
    	Let $\Shum{K}$ be an integral model of a Shimura variety whose minimal and toroidal compactifications have the axiomatic descriptions listed in \cite[Proposition 2.1.2]{lan2018compactifications} (we recall a part of the list in Theorem \ref{theorem: axiomatic descriptions of compactifications}). Let $T$ be a locally noetherian scheme over $S$. A locally closed subset (resp. subscheme) $Y$ of $\Shum{K, T}:=(\Shum{K})_T$ is called \emph{well positioned}, if, for each $[\Phi] \in \Cusp_K(G, X)$, there exists a locally closed subset (resp. subscheme) $Y^{\natural}(\Phi) \subset \Zb(\Phi)_T \rightiso \mathcal{Z}([\Phi])_T$ such that for some (thus for all) cone decompositions $\Sigma$, and for each $\sigma\in\Sigma(\Phi)^+$, there are some (thus for all) open affine coverings $\mathfrak{W}$ of $\mathfrak{X}^{\circ}_{\sigma}$ satisfying the following property: the pullback of $Y^{\natural}(\Phi)$ along $W^0_T \to \Xi(\Phi)_T \to C(\Phi)_T \to \Zb^{\bigsur}(\Phi)_T \to \Zb(\Phi)_T$ (see (\ref{def: W})) coincides with the pullback of $Y$ along $W^0_T \to \Shum{K, T}$ as a subset (resp. subscheme). If this is the case, we say $Y$ is well positioned with respect to $Y^{\natural}:=\{Y^{\natural}(\Phi)\}_{\Phi}$, and $Y^{\natural}$ is associated with $Y$.
    \end{definition-proposition}
     In practice, in case $(4)$ when $S = \Spec \OO_{E(v)}$, we usually take $T = s = \Spec k_E$ or $T=\bar{s} = \Spec \Bar{k}_E$.
     \begin{remark}\label{remark: slightly weaker condition for being well-positioned}\leavevmode
      \begin{enumerate}
          \item If there exists a subset $Y^{\natural}(\Phi) \subset \Zb(\Phi)_T$ whose preimage in $W^0_T$ is the preimage of a locally closed subset $Y \subset \Shum{K, T}$ in $W^0_T$ for each $W$, then $Y^{\natural}(\Phi)$ is automatically locally closed, see \cite[Lemma 2.3.10]{lan2018compactifications}.
          \item When $Y$ is a well-positioned subset with respect to $Y^{\natural}:=\{Y^{\natural}(\Phi)\}_{\Phi}$, then $Y^{\natural}(\Phi)$ is uniquely determined by $Y$.
          \item If moreover $C(\Phi) \to \Zb(\Phi)$ is flat (thus faithfully flat) for each $[\Phi]$ and $Y$ is a well-positioned subscheme, then $Y^{\natural}(\Phi)$ (as a subscheme of $\Zb(\Phi)$) is uniquely determined by $Y$, see \cite[Theorem 2.3.2(4)]{lan2018compactifications}.
      \end{enumerate}
        
     \end{remark}
 
  \subsubsection{Well-positioned subsets and subschemes}
    Let us discuss the difference of being a well-positioned subset and being a well-positioned subscheme. First, it follows from definition that the underlying set of a well-positioned subscheme is a well-positioned subset. To show the other direction, we need following assumption:
     \begin{assumption}\label{assumption: C to Z}
         Assume $C(\Phi) \to \Zb(\Phi)$ is flat and has geometrically reduced fiber for each $[\Phi] \in \Cusp_K(G, X)$.
     \end{assumption}
     The assumption \ref{assumption: C to Z} holds if for each $[\Phi] \in \Cusp_K(G, X)$, $C(\Phi) \to \Zb^{\bigsur}(\Phi)$ is an abelian scheme torsor and $\Zb^{\bigsur}(\Phi) \to \Zb(\Phi)$ is finite \'etale. In particular, due to \cite[Theorem 3.58]{mao2025boundary}, \cite[Corollary 4.27]{wu2025arith} 
     \begin{proposition}\label{remark: when the assumption C to Z is true}
       The assumption \ref{assumption: C to Z} holds for both $\Shum{K_1}(G_1, X_1)$ and $\Shum{K_2}(G_2, X_2)$ in case $(4)'$, where $K_{1, p}$ is a quasi-parahoric level subgroup, $K_{2, p}$ is a hyperspecial level subgroup, and $(G_1, X_1, K_{1, p}) \to (G_2, X_2, K_{2, p})$ is an adjusted Siegel embedding, and $K^p, K^{\dd, p}$ properly chosen. In particular, when we work with Kisin-Pappas-Zhou integral models or Pappas-Rapoport integral models, after replacing the chain of almost self-dual lattices with a self-dual lattice, assumption \ref{assumption: C to Z} is automatic.
    \end{proposition}
    \begin{remark}\label{remark: subset well positioned implies subscheme well positioned}
        Under the assumption \ref{assumption: C to Z}, the pullbacks of reduced subschemes along the morphism $C \to \Zb$ are reduced. Let $Y \subset \Shum{K, T}$ be a locally closed subscheme endowed with the induced reduced subscheme structure, then $Y$ is a well-positioned subscheme if and only if the underlying set $Y$ is a well-positioned subset.
    \end{remark}

   \subsubsection{Partial compactifications}
    \begin{definition}\cite[Definition 2.3.1]{lan2018compactifications}
    	Let $Y$ be a locally closed subset (resp. subscheme) of $\Shum{K, T}$. Let $\overline{Y}$ be the closure (resp. scheme theoretical closure) of $Y$ in $\Shum{K, T}$, $Y_0=\overline{Y}\backslash Y$ the complement subset of $Y$ in $\overline{Y}$ (resp. the complement subscheme). Let $\overline{Y}^{\min}$ and $Y_0^{\min}$ be the set theoretical (resp. scheme theoretical) closure of $\overline{Y}$ and $Y_0$ in $\Shumm{K, T}$ respectively, $Y^{\min}$ be the complement subset of $Y_0^{\min}$ in $\overline{Y}^{\min}$ (resp. the complement subscheme). Similarly we define $\overline{Y}^{\tor}_{\Sigma}$, $Y_{0, \Sigma}^{\tor}$, $Y^{\tor}_{\Sigma}$. We call $Y^{\min}$ (resp. $Y^{\tor}_{\Sigma}$) the \emph{partial minimal compactification} (resp. the \emph{partial toroidal compactification} with respect to $\Sigma$) of $Y$. In the case of subsets, we always view these locally closed sets as schemes with induced reduced subscheme structures.
    \end{definition}
    \begin{remark}\label{remark: set/scheme theoretical closure does not matter}
        Since $\Shum{K, T}$, $\Shumc{K, \Sigma, T}$, $\Shumm{K, T}$ are locally noetheiran schemes (each of them is finitely presented over the locally noetherian base $T$), any locally closed immersion $Y \to \Shumm{K, T}$ or $Y \to \Shumc{K, \Sigma, T}$ is quasi-compact (\cite[\href{https://stacks.math.columbia.edu/tag/01OX}{Tag 01OX}]{stacks-project}), thus the underlying set of $Y$ is dense in its scheme theoretical closure (\cite[\href{https://stacks.math.columbia.edu/tag/01R8}{Lemma 01R8}]{stacks-project}), the underlying set of the scheme theoretical partial minimal and toroidal compactifications are exactly the set theoretical partial minimal and toroidal compactifications, respectively.
    \end{remark}
    
    Given a morphism $f: Y \to X$, let $Z$ be a locally closed subscheme of $X$, then the underlying topological space of scheme-theoretical pullback $f^{-1}(Z)$ is the set-theoretical pullback $f^{-1}(Z)$. We use the symbol $f^{-1}(Z)$ for both pullbacks if there is no confusion.
    
    \begin{proposition}\cite[Theorem 2.3.2]{lan2018compactifications}\label{prop: Y same as Shum}
    	Let $Y$ be a well-positioned subset (resp. subscheme) of $\Shum{K, T}$ with respect to $Y^{\natural}:=\{Y^{\natural}(\Phi)\}_{\Phi}$, consider the partial compactifications:
    	\[\begin{tikzcd}
    		Y && {Y^{\tor}_{\Sigma}} \\
    		&& {Y^{\min}}
    		\arrow["{J_{Y^{\tor}_{\Sigma}}}", from=1-1, to=1-3]
    		\arrow["{\pi_Y}", from=1-3, to=2-3]
    		\arrow["{J_{Y^{\min}}}"', from=1-1, to=2-3]
    	\end{tikzcd}\]
       where the morphism $\pi_Y$ is proper and surjective induced by the natural proper surjection:
       \[ \oint_{K, \Sigma}: \Shumc{K, \Sigma} \to \Shumm{K} \]
       The structure morphisms $Y^{\tor}_{\Sigma} \to T$, $Y^{\min} \to T$ are proper when $Y$ is closed in $\Shum{K, T}$.
       
        The stratification $\mathcal{Z}(\Upsilon)$ (resp. $\mathcal{Z}([\Phi])$) of $\Shumc{K, \Sigma}$ (resp. $\Shumm{K}$) induces a stratification $Y_{\mathcal{Z}(\Upsilon)}$ (resp. $Y_{[\Phi]}$) of $Y^{\tor}_{\Sigma}$ (resp. $Y^{\min}$), where 
        \[ Y_{\mathcal{Z}(\Upsilon)}=\mathcal{Z}(\Upsilon)\times_{\Shumc{K, \Sigma}}Y^{\tor}_{\Sigma},\quad  Y_{[\Phi]}=\mathcal{Z}([\Phi])\times_{\Shumm{K}}Y^{\min}.\]
        Each $\Zb(\Phi) \rightiso \mathcal{Z}([\Phi])$ induces a canonical morphism $Y^{\natural}(\Phi)\to Y_{[\Phi]}$ which is bijective on the underlying set.

        For each $\sigma \in \Sigma(\Phi)^+$, let $Y_C^{\natural}$ be the pullback of $Y^{\natural}(\Phi)$ along $C \to \Zb$. Let $(\ast)$ be any of the following: $\Xi$, $\Xi_{\sigma}$, $\Xi({\sigma})$, $\mathfrak{X}$, $\mathfrak{X}_{\sigma}$, $\mathfrak{X}_{\sigma}^{\circ}$, let $Y^{\natural}_{(\ast)}$ be the pullback of $Y^{\natural}_C$ along $(\ast) \to C$. In the case of subsets, in order to formulate results on schemes and formal schemes as in Theorem \ref{theorem: axiomatic descriptions of compactifications}, we endow $Y_C^{\natural}$ and $Y^{\natural}(\Phi)$ with the induced reduced subscheme structures. Then $Y \to Y_{\Sigma}^{\tor} \to Y^{\min}$ satisfy the descriptions listed in Theorem \ref{theorem: axiomatic descriptions of compactifications}.
    \end{proposition}
    \begin{remark}\label{remark: why need assumption?}
    The pullbacks in Proposition \ref{prop: Y same as Shum} start from $C$ instead of $\Zb$, since the pullback of $Y^{\natural}(\Phi)$ along $C \to \Zb$ might not be a reduced subscheme even if $Y^{\natural}(\Phi)$ itself is a reduced subscheme. The problem comes from the possibly non-reduceness of the morphism $C \to \Zb$. The assumption \ref{assumption: C to Z} assures that, if $Y$ is a well-positioned subset with respect to $Y^{\natural}:=\{Y^{\natural}(\Phi)\}_{\Phi}$ and if we endow both $Y^{\natural}(\Phi)$ and $Y^{\natural}_{C(\Phi)}$ the induced reduced subscheme structures, then $Y^{\natural}_{C(\Phi)}$ is indeed the pullback of $Y^{\natural}(\Phi)$, thus we can pull back everything from $\Zb$ instead of $C$. 
        
    \end{remark}

    \subsubsection{}
         Let $Y$ be a well-positioned subset (resp. subscheme) of $\Shum{K}$ with respect to $Y^{\natural} = \lrbracket{Y^{\natural}(\Phi)}_{\Phi}$, since $Y \to Y^{\tor}_{\Sigma} \to Y^{\min}$ satisfy the axiomatic descriptions listed in Theorem \ref{theorem: axiomatic descriptions of compactifications}, we could define well-positioned subschemes of $Y$.
    \begin{lemma}\label{lemma: well-position of well-position is well-positioned}
        When $Y$ is a well-positioned subset (resp. subscheme) of $\Shum{K}$ with respect to $Y^{\natural} = \lrbracket{Y^{\natural}(\Phi)}_{\Phi}$, let $Z \subset Y$ be a well-positioned subset (resp. subscheme) with respect to $Z^{\natural}=\lrbracket{Z^{\natural}(\Phi)}_{\Phi}$, $Z^{\natural}(\Phi) \subset Y^{\natural}(\Phi) \subset \Zb(\Phi)_T$, then $Z$ is a well-positioned subset (resp. subscheme) of $\Shum{K, T}$ with respect to $Z^{\natural}$.
    \end{lemma}
    \begin{proof}
        The proof is tautological. Let $\mathfrak{W}=\Spf(R) \to \mathfrak{X}_{\sigma}^{\circ}$ be an affine formal subscheme, its pullback along $Y^{\natural}_{\mathfrak{X}_{\sigma}^{\circ}} \to \mathfrak{X}_{\sigma}^{\circ}$ is a union of $\mathfrak{W}_i = \Spf(R_i) \to Y^{\natural}_{\mathfrak{X}_{\sigma}^{\circ}} $. The pullback of $W =\Spec(R) \to \Xi(\sigma)$ along $Y^{\natural}_{\Xi(\sigma)} \to \Xi(\sigma)$ is then a union of $W_i = \Spec (R_i) \to Y^{\natural}_{\Xi(\sigma)}$ due to the proof of \cite[Theorem 2.3.2]{lan2018compactifications}. Since $Y^{\natural}_{\Xi} \to Y^{\natural}_{\Xi(\sigma)}$ is the pullback of $\Xi \to \Xi(\sigma)$, then the pullback of $W^0 \to \Xi$ is a union of $W_i^0 \to Y^{\natural}_{\Xi}$. Since the pullback of $Z^{\natural} \to Y^{\natural} \to \Zb$ (resp. $Z \to Y \to \Shum{K, T}$) along $W^0 \to \Zb$ (resp. $W^0 \to \Shum{K, T}$) factors through $\cup W_i^0$, and these two pullbacks coincide over each $W_i^0$ since $Z \subset Y$ is well-positioned with respect to $Z^{\natural} \subset Y^{\natural}$, then the Lemma follows.
    \end{proof}

   \begin{lemma}\label{lemma: check equality over boundary strata}
       Given two locally closed subschemes $Y_1, Y_2 \subset \Shumc{K, \Sigma, T}$ and assume there is a morphism $f: Y_1 \to Y_2$ over $\Shumc{K, \Sigma, T}$. Assume for each boundary stratum $\mathcal{Z}(\Upsilon)$ of $\Shumc{K, \Sigma}$, $f$ induces an isomorphism $Y_1 \times_{\Shumc{K, \Sigma, T}} \mathcal{Z}(\Upsilon) \rightiso Y_2 \times_{\Shumc{K, \Sigma, T}} \mathcal{Z}(\Upsilon)$, then $f$ itself is an isomorphism.
   \end{lemma}
   \begin{proof}
       By assumption, $f$ is a bijection between the underlying topological spaces. It suffices to check that $f$ is an isomorphism \'etale locally. Due to \cite[Corollary 2.1.7]{lan2018compactifications}, given $x \in \mathcal{Z}([\Phi, \sigma])$, there is an \'etale neighbourhood $\ovl{U}$ of $\Shumc{K, \Sigma}$ at $x$ and an \'etale morphism $\ovl{U} \to E(\sigma) \times C$ such that the pullbacks of the stratifications of $\Shumc{K, \Sigma}$ and of $E(\sigma) \times C$ coincide, here $E \hookrightarrow E(\sigma)$ is an affine toroidal embedding. Then the lemma follows easily.
   \end{proof}

   \begin{lemma}\label{lemma: closure of well-pos is well-pos}
       Let $Y \subset \Shum{K, T}$ be a well-positioned subset (resp. subscheme) with repsect to $Y^{\natural} = \lrbracket{Y^{\natural}(\Phi)}_{\Phi}$. Assume \ref{assumption: C to Z} (in fact we only need the flatness of $C \to \Zb$), then $\ovl{Y}$ (resp. $Y_0$) is a well-positioned subset (resp. subscheme) with respect to $\lrbracket{ \ovl{Y^{\natural}(\Phi)}}_{\Phi}$ (resp. $\lrbracket{Y^{\natural}(\Phi)^0}_{\Phi}$). Moreover, $\ovl{Y}_{\Sigma}^{\tor} = \ovl{Y^{\tor}_{\Sigma}}$, $\ovl{Y}^{\min} = \ovl{Y^{\min}}$ (thus $Y_{0, \Sigma}^{\tor} = \ovl{Y^{\tor}_{\Sigma}} \setminus Y^{\tor}_{\Sigma}$, $Y_0^{\min} = \ovl{Y^{\min}} \setminus Y^{\min}$).
   \end{lemma}
   \begin{proof}
       The first part of the Lemma is \cite[Lemma 2.2.7]{lan2018compactifications}. We need to verify the last sentence. $\ovl{Y}^{\min} = \ovl{Y^{\min}}$ can be verified over the boundary $\Zb(\Phi)$, $\ovl{Y}^{\min} \cap \Zb(\Phi) = \ovl{Y}^{\natural}(\Phi) = \ovl{Y^{\natural}(\Phi)} = \ovl{Y^{\min}} \cap \Zb(\Phi)$. Now we verify $\ovl{Y}^{\tor} = \ovl{Y^{\tor}}$ with the help of Lemma \ref{lemma: check equality over boundary strata}: 
       
       Since $C(\Phi)\to \Zb(\Phi)$ is faithfully flat by assumption, $\Xi_{\sigma}(\Phi) \to C(\Phi) \to \Zb(\Phi)$ is faithfully flat, the closure of $Y^{\natural}_{\Xi_{\sigma}(\Phi)}$ is the pullback of $\ovl{Y^{\natural}(\Phi)}$. Use the canonical identification $Y_{\mathcal{Z}([\Phi, \sigma])} \cong Y_{\Xi_{\sigma}(\Phi)}$ and $\ovl{Y}_{\mathcal{Z}([\Phi, \sigma])} \cong \ovl{Y}_{\Xi_{\sigma}(\Phi)}$ for any $\mathcal{Z}([\Phi, \sigma])$, we have:
       \[ \ovl{Y}^{\tor}_{\Sigma} \cap \mathcal{Z}([\Phi, \sigma]) = \ovl{Y}_{\mathcal{Z}([\Phi, \sigma])} = \ovl{Y}_{\Xi_{\sigma}(\Phi)} = \ovl{Y_{\Xi_{\sigma}(\Phi)}} = \ovl{Y^{\tor}_{\Sigma}} \cap \mathcal{Z}([\Phi, \sigma]).\] 
   \end{proof}

   \begin{lemma}\label{lemma: intersection of closed well-positioned subschemes}
       Let $Y_1, Y_2 \subset \Shum{K, T}$ be \emph{closed} subsets (resp. subschemes), assume $Y_1$, $Y_2$ are well-positioned subsets (resp. subschemes) with respect to $Y^{\natural}_1 = \lrbracket{Y^{\natural}_1(\Phi)}_{\Phi}$ and $Y^{\natural}_2 = \lrbracket{Y^{\natural}_2(\Phi)}_{\Phi}$ respectively, then $Y:=Y_1 \cap Y_2$ is a well-positioned subset (resp. subscheme) with respect to $Y^{\natural} = \lrbracket{Y^{\natural}(\Phi)}_{\Phi}$, where $Y^{\natural}(\Phi) = Y^{\natural}_1(\Phi) \cap Y_2^{\natural}(\Phi)$. Moreover, $Y^{\tor}_{\Sigma} = Y^{\tor}_{1, \Sigma} \cap Y^{\tor}_{2, \Sigma}$ (resp. as subschemes), $Y^{\min} = Y_1^{\min} \cap Y_2^{\min}$.
   \end{lemma}
   \begin{proof}
      It follows immediately from Definition-Proposition \ref{def: well-positioned} that $Y=Y_1 \cap Y_2$ is a well-positioned subset (resp. subscheme) with respect to $Y^{\natural}$. Since \emph{closure of intersection is contained in intersection of closure}, then $Y^{\tor}_{\Sigma} \subset Y^{\tor}_{1, \Sigma} \cap Y^{\tor}_{2, \Sigma}$ (resp. $Y^{\min} \subset Y^{\min}_1 \cap Y^{\min}_2$). To check this is an isomorphism (resp. a bijection on underlying topological spaces), with the help of Lemma \ref{lemma: check equality over boundary strata}, it suffices to check on every boundary stratum. 
      
      $(Y^{\min})$: Consider underlying subsets, under the canonical identification $\Zb(\Phi) \rightiso \mathcal{Z}([\Phi])$, $Y_{[\Phi]} = Y^{\min} \cap \mathcal{Z}([\Phi])$ is $Y^{\natural}(\Phi) = Y^{\natural}_1(\Phi) \cap Y^{\natural}_2(\Phi)$, and 
      \[ (Y^{\min}_1 \cap Y^{\min}_2) \cap \mathcal{Z}([\Phi]) = Y_{1, [\Phi]} \cap Y_{2, [\Phi]} = Y^{\natural}_1(\Phi) \cap Y^{\natural}_2(\Phi) = Y^{\min} \cap \mathcal{Z}([\Phi]). \]

      $(Y^{\tor}_{\Sigma})$: As subsets (resp. subschemes) of $\mathcal{Z}([\Phi, \sigma])$, $Y^{\tor}_{\Sigma} \cap \mathcal{Z}([\Phi, \sigma]) = Y_{\mathcal{Z}([\Phi, \sigma])} \rightiso Y^{\natural}_{\Xi_{\sigma}(\Phi)}$, which is the pullback of $Y^{\natural}(\Phi)$ along $\Xi_{\sigma}(\Phi) \to \Zb(\Phi)$. Similarily, $(Y^{\tor}_{1, \Sigma} \cap Y^{\tor}_{2, \Sigma}) \cap \mathcal{Z}([\Phi, \sigma])$ is the pullback of $Y^{\natural}_1(\Phi) \cap Y^{\natural}_2(\Phi)$. The lemma follows from $Y^{\natural}(\Phi) = Y^{\natural}_1(\Phi) \cap Y^{\natural}_2(\Phi)$.
   \end{proof}

   \begin{lemma}\label{lemma: subset of well-positioned subschemes}
       Let $Y_1 \subset Y_2 \subset \Shum{K, T}$ be locally closed subsets (resp. subschemes), assume $Y_1$, $Y_2$ are well-positioned subsets (resp. subschemes) with respect to $Y^{\natural}_1 = \lrbracket{Y^{\natural}_1(\Phi)}_{\Phi}$ and $Y^{\natural}_2 = \lrbracket{Y^{\natural}_2(\Phi)}_{\Phi}$ respectively. Assume \ref{assumption: C to Z} for $\Shum{K, T}$, then $Y^{\natural}_1(\Phi) \subset Y^{\natural}_2(\Phi)$ for each $\Phi$, and $Y_{1, \Sigma}^{\tor} \subset Y_{2, \Sigma}^{\tor}$, $Y_1^{\min} \subset Y_2^{\min}$.
   \end{lemma}
   \begin{proof}
       We first show this is true when $Y_1 \subset Y_2$ are both \emph{closed} subsets (resp. subschemes) of $\Shum{K, T}$. In this case, we automatically have $Y_{1, \Sigma}^{\tor} \subset Y_{2, \Sigma}^{\tor}$, $Y_1^{\min} \subset Y_2^{\min}$ by taking closures. In particular, under the canonical isomorphism $\Zb(\Phi) \rightiso \mathcal{Z}([\Phi])$, $Y_{i, [\Phi]} = Y^{\min}_{i} \cap \mathcal{Z}([\Phi])$, $Y_{i, [\Phi]} \cong Y_i^{\natural}(\Phi)$, then $Y_{1, [\Phi]} \subset Y_{2, [\Phi]}$ implies that $Y_1^{\natural}(\Phi) \subset Y_2^{\natural}(\Phi)$.
       
       Now we consider the general case when $Y_1$, $Y_2$ are locally closed. Since $Y_1 \subset Y_2$, then $\ovl{Y}_1 \subset \ovl{Y}_2$, and $Y_{2, 0} \cap \ovl{Y}_1 \subset Y_{1, 0}$. Due to Lemma \ref{lemma: closure of well-pos is well-pos}, $Y_{2, 0}$, $\ovl{Y}_1$ and $Y_{1, 0}$ are well-positioned subsets (resp. subschemes). Due to Lemma \ref{lemma: intersection of closed well-positioned subschemes}, $Y_{2, 0, \Sigma}^{\tor} \cap \ovl{Y}_{1, \Sigma}^{\tor} = (Y_{2, 0} \cap \ovl{Y}_1)^{\tor}_{\Sigma}$, thus $Y_{2, 0, \Sigma}^{\tor} \cap \ovl{Y}_{1, \Sigma}^{\tor}  = (Y_{2, 0} \cap \ovl{Y}_1)^{\tor}_{\Sigma} \subset Y_{1, 0, \Sigma}^{\tor}$ (by the first paragraph), in particular, $Y_{1, \Sigma}^{\tor} \subset Y_{2, \Sigma}^{\tor}$. Similarly, $Y_1^{\min} \subset Y_2^{\min}$, and this implies $Y_1^{\natural}(\Phi) \subset Y_2^{\natural}(\Phi)$.
   \end{proof}
    
    \begin{lemma}\label{lemma: intersection of well-positioned subschemes}
           If a locally closed subset (resp. subscheme) $Y \subset \Shum{K, T}$ is an intersection of locally closed subsets (resp. subschemes) $Y_i \subset \Shum{K, T}$, $i \in I$, where each $Y_i$ is well-positioned with respect to the locally closed subsets (resp. subschemes) $Y^{\natural}_i = \lrbracket{Y^{\natural}_i(\Phi)}_{\Phi}$, then $Y$ is well-positioned with respect to the locally closed subsets (resp. subschemes) $Y^{\natural} = \lrbracket{\cap_{i\in I} Y^{\natural}_i(\Phi)}_{\Phi}$. Moreover, assume \ref{assumption: C to Z} for $\Shum{K, T}$, then $Y^{\tor}_{\Sigma} \cong \cap_{i \in I} Y^{\tor}_{i, \Sigma}$ and $Y^{\min} \cong \cap_{i \in I} Y^{\min}_i$ as subsets (in the case of subschemes, $Y^{\tor}_{\Sigma} \cong \cap_{i \in I} Y^{\tor}_{i, \Sigma}$ as subschemes).
       \end{lemma}
    \begin{proof}
         The first part follows trivially from the definition and Remark \ref{remark: slightly weaker condition for being well-positioned} (which says $\cap_{i\in I} Y^{\natural}_i(\Phi)$ is locally closed). We prove the second part: due to Lemma \ref{lemma: subset of well-positioned subschemes}, there are natural morphisms $Y^{\tor}_{\Sigma} \to \cap_{i \in I} Y^{\tor}_{i, \Sigma}$, $Y^{\min} \to \cap_{i \in I} Y^{\min}_i$. We can check the equalities as in the proof of Lemma \ref{lemma: intersection of closed well-positioned subschemes}.
    \end{proof}

     Similarly, we have
    \begin{lemma}\label{lem: union of well-positioned sets}
    	If a locally closed subset $Y \subset \Shum{K, T}$ is a union of well-positioned subsets $\{ Y_i \}_{i\in I}$ with respect to $Y_i^{\natural}=\{ Y_{i}^{\natural}(\Phi) \}_{\Phi}$, then $Y$ is a well-positioned subset with respect to $Y^{\natural}:=\{ \cup_{i\in I}Y_{i}^{\natural}(\Phi) \}_{\Phi}$. If assumption \ref{assumption: C to Z} holds, and $Y$, $Y^{\natural}$ are endowed with induced reduced subscheme structure, then $Y$ is a well-positioned subscheme, and moreover $Y^{\tor}_{\Sigma} \cong \cup_{i \in I} Y^{\tor}_{i, \Sigma}$ and $Y^{\min} = \cup_{i \in I} Y^{\min}_i$.
    \end{lemma}

    \subsubsection{}

    \begin{lemma}\label{lemma: pullback of well-positioned is well-positioned}
       Keep assumptions and notations as in Proposition  \ref{proposition: functorialities on toroidal compactifications}. Let $Y_2$ be a well-positioned subset (resp. subscheme) of $\Shum{K_2, T}$ with respect to subsets (resp. subschemes) $Y^{\natural}_2 = \{ Y^{\natural}_2(\Phi_2) \}_{\Phi_2}$, and let $Y_1$ be the pullback of $Y_2$ along $\varphi: \Shum{K_1, T}\to\Shum{K_2, T}$. Then $Y_1$ is a well-positioned subset (resp. subscheme) of $\Shum{K_1, T}$ with respect to the pullback subsets (resp. subschemes) $Y^{\natural}_1 = \{ Y^{\natural}_1(\Phi_1) \}_{\Phi_1}$ of $Y^{\natural}_2 = \{ Y^{\natural}_2(\Phi_2) \}_{\Phi_2}$ along those canonical morphisms $\Zb_1(\Phi_1)_T\to\Zb_2(\Phi_2)_T$. Moreover, if the assumption \ref{assumption: C to Z} holds for both $\Shum{K_1}$ and $\Shum{K_2}$ (which is always the case when we are in case $(4)'$), then $Y^{\tor}_{1, \Sigma_1}$ (resp. the underlying set of $Y^{\min}_1$) is the pullback (resp. preimage) of $Y^{\tor}_{2, \Sigma_2}$ (resp. the underlying set of $Y^{\min}_2$).
  
    \end{lemma}
    \begin{proof}
    	Due to Proposition \ref{proposition: functorialities on toroidal compactifications}, after refining the cone decompositions, we have the following commutative diagram (not necessarily cartesian):
\begin{equation}\label{graph:comm1}
\begin{tikzcd}
	{\Shum{K_1, T}} & {(W^0_1)_T} & {\Xi_1(\Phi_1)_T} & {C_1(\Phi_1)_T} & {\Zb_1(\Phi_1)_T} \\
	{\Shum{K_2, T}} & {(W^0_2)_T} & {\Xi_2(\Phi_2)_T} & {C_2(\Phi_2)_T} & {\Zb_2(\Phi_2)_T}
	\arrow[from=1-1, to=2-1]
	\arrow[from=1-2, to=1-1]
	\arrow[from=1-2, to=1-3]
	\arrow[from=1-2, to=2-2]
	\arrow[from=1-3, to=1-4]
	\arrow[from=1-3, to=2-3]
	\arrow[from=1-4, to=1-5]
	\arrow[from=1-4, to=2-4]
	\arrow[from=1-5, to=2-5]
	\arrow[from=2-2, to=2-1]
	\arrow[from=2-2, to=2-3]
	\arrow[from=2-3, to=2-4]
	\arrow[from=2-4, to=2-5]
\end{tikzcd}
\end{equation}
    Since the pullback of a covering $\mathfrak{W}_2$ of $\mathfrak{X}_{2, \sigma_2}^{\circ}$ induces a covering $\mathfrak{W}_1$ of $\mathfrak{X}_{1, \sigma_1}^{\circ}$, the first part of the lemma is a direct corollary of the Definition-Proposition \ref{def: well-positioned} and Proposition \ref{proposition: functorialities on toroidal compactifications}.

    We focus on the second statement. Fisrt, let us assume $Y_2 \subset \Shum{K_2, T}$ is closed, then $Y_1 \subset \Shum{K_1, T}$ is closed.  Let us denote by $Y^{\prime, \tor}_{1, \Sigma_1}$ (resp. $Y^{\prime, \min}_1$) the pullback of $Y^{\tor}_{2, \Sigma_2}$ (resp. the preimage of $Y^{\min}_2$). Since \emph{preimage of closure contains closure of preimage}, we have canonical inclusions $Y_{1, \Sigma_1}^{\tor} \to Y^{\prime, \tor}_{1, \Sigma_1}$, (resp. $Y_1^{\min} \subset Y_1^{\prime, \min}$). To show this is an isomorphism (resp. a homeomorphism on the underlying topological spaces), we apply Lemma \ref{lemma: check equality over boundary strata}, it suffices to check over the boundary stratification:

    In the case of toroidal compactifications, let $\mathcal{Z}_1([\Phi_1, \sigma_1]) \cong \Xi_{1, \sigma_1}(\Phi_1)$ be a stratum, the pullback of $Y^{\tor}_{1, \Sigma_1}$ over this stratum is $Y_{1, \mathcal{Z}_1([\Phi_1, \sigma_1])} \cong Y^{\natural}_{1, \Xi_{1, \sigma_1}(\Phi_1)}$, the pull back $Y_{1, \mathcal{Z}_1([\Phi_1, \sigma_1])}'$ of $Y^{\prime, \tor}_{1, \Sigma_1}$ is the pullback of $Y^{\tor}_{2, \Sigma_2}$ along $\mathcal{Z}_1([\Phi_1, \sigma_1]) \to \mathcal{Z}_2([\Phi_2, \sigma_2])$, which is canonically isomorphic to the pullback of $Y^{\natural}_{2, \Xi_{2, \sigma_2}(\Phi_2)}$ along $\Xi_{1, \sigma_1}(\Phi_1) \to \Xi_{2, \sigma_2}(\Phi_2)$. Note that $Y^{\natural}_{2, \Xi_{2, \sigma_2}(\Phi_2)}$ (resp. $Y^{\natural}_{1, \Xi_{1, \sigma_1}(\Phi_1)}$) is the pullback of $Y^{\natural}_2(\Phi_2)$ (resp. $Y^{\natural}_1(\Phi_1)$) along $\Xi_{2, \sigma_2}(\Phi_2)\to\Zb_2(\Phi_2)$ (resp. $\Xi_{1, \sigma_1}(\Phi_1) \to \Zb_1(\Phi_1)$), and $Y^{\natural}_1(\Phi_1)$ is the pullback of $Y^{\natural}_2(\Phi_2)$ along $\Zb_1(\Phi_1)\to\Zb_2(\Phi_2)$, thus $Y_{1, \mathcal{Z}_1([\Phi_1, \sigma_1])} \to Y_{1, \mathcal{Z}_1([\Phi_1, \sigma_1])}'$ is a canonical isomorphism. In the case of minimal compactifications, the calculation follows similarly. Note that $Y^{\natural}(\Phi) \to Y_{[\Phi]}$ is a only bijection on the underlying set, not necessarily an isomorphism as locally closed subschemes.

    Second, let $Y_1$, $Y_2$ be locally closed, this is the case where we need assumption \ref{assumption: C to Z}. Since $\varphi$ and $\varphi^{\min}$ are finite morphism (thus universally closed), then $\varphi(\ovl{Y}_1) \subset \ovl{Y}_2$ and $\varphi^{\min}(\ovl{Y}_1^{\min}) \subset \ovl{Y}_2^{\min}$. Since $\varphi^{-1}(Y_{2, 0}) \cap \ovl{Y}_1 = Y_{1, 0}$, and $Y_{2, 0}$ is well-positioned due to Lemma \ref{lemma: closure of well-pos is well-pos}, then $\varphi^{\min, -1}(Y_{2, 0}^{\min}) = \varphi^{-1}(Y_{2, 0})^{\min}$ by previous arguments ($Y_{2, 0}$ is closed), thus $\varphi^{\min, -1}(Y_{2, 0}^{\min}) \cap \ovl{Y}_1^{\min} = Y_{1, 0}^{\min}$ by Lemma \ref{lemma: intersection of closed well-positioned subschemes}. In particular, $\varphi^{\min}(Y_1^{\min}) \subset Y_2^{\min}$. Let $Y_1^{\prime, \min}$ be the preimage of $Y_2^{\min}$, then we have a canonical inclusion $Y_1^{\min} \to Y_1^{\prime, \min}$. By checking on boundary as before, we have $Y_1^{\min} = Y_1^{\prime, \min}$. Similarly, $Y_{1, \Sigma_1}^{\tor} = Y_{1, \Sigma_1}^{\prime, \tor}$.
    \end{proof}

    Well-position subsets and subschemes are also functorial with Hecke actions.
    \begin{lemma}\label{lemma: pullback of well positioned, different level}
	Let $Y$ be a locally closed subset (resp. subscheme) of $\Shum{K, T}$ with respect to $Y^{\natural}=\{Y^{\natural}(\Phi)\}_{\Phi}$. Let $g \in G(\A_f)$ such that $K' \subset gKg^{-1}$, $[g]: \Shum{K'}\to \Shum{K}$ be the Hecke action. Consider the pullback $Y'$ of $Y$ along $[g]$ and pullbacks $Y^{\prime, \natural}=\{Y^{\prime, \natural}(\Phi')\}_{\Phi'}$ of $Y^{\natural}=\{Y(\Phi)\}_{\Phi}$ along $\Zb^{\prime}(\Phi')\to \Zb(\Phi)$ induced by $[g]^{\min}: \Shumm{K^{\prime}}\to \Shumm{K}$ for each $\Zb^{\prime}(\Phi')$ above $\Zb(\Phi)$. In the case of subschemes, we assume moreover $[g]$ is faithfully flat. Then $Y$ is a well-positioned subset (resp. subscheme) with respect to $Y^{\natural}$ if and only if $Y^{\prime}$ is a well-positioned subset (resp. subscheme) with respect to $Y^{\prime, \natural}$. In this case, the preimage of $Y^{\min}$ is $Y^{\prime, \min}$, and the pullback of $Y_{\Sigma}^{\tor}$ is $Y_{\Sigma'}^{\prime, \tor}$ ($\Sigma'$ is a $g$-refinement of $\Sigma$, see \cite[Proposition 2.4.3]{lan2018compactifications}). If $Y$ is connected, then $Y_{\Sigma}^{\tor}$, $Y^{\min}$ are connected.
\end{lemma}
\begin{proof}
	The \emph{only-if} part of the lemma as well as the statement about the partial minimal and toroidal compactifications are exactly \cite[Proposition 2.4.3]{lan2018compactifications}. The \emph{if} part follows from the fact that $\cup W_i \to W$ (we use the notations from \emph{loc.cit.}) as the pullback of $\Shum{K'} \to \Shum{K}$ is faithfully flat, any two schemes over $W$ that are canonically isomorphic over each $W_i$, are also isomorphic over $W$. Since $Y$ is open and dense in both $Y_{\Sigma}^{\tor}$ and $Y^{\min}$, thus $Y_{\Sigma}^{\tor}$ and $Y^{\min}$ are connected when $Y$ is connected.
\end{proof}

    \begin{definition}\label{definition: insensitive to away from p part lift}
        Fix level $K_p \subset G(\Q_p)$, $S_{K_p}$ be a set, $\lrbracket{\phi_{K^p}: \Shum{K_pK^p, T}\to S_{K_p}}_{K^p}$ be a system of morphisms varying $K^p$ such that for all $g \in G(\A_f^p)$, $K^{p, \prime} \subset gK^pg^{-1}$, $\phi_{K^p} \circ [g] = \phi_{K^{p,\prime}}$. Let $Y_K \subset \Shum{K, T}$ be a locally closed subset (resp. subscheme), we say it is \emph{insensitive} to the away-from-$p$ level $K^p$, if $Y_K=\phi_{K^p}^{-1}(X)$ for some subset $X \subset S_{K_p}$. That is to say, $Y_{K_pK^{p, \prime}} = [g]^{-1}(Y_{K_pK^p})$.
    \end{definition}
    \begin{corollary}\label{corollary: insensitive to away from p part lift}
    	Let $Y_K$ be a locally closed subset (resp. subscheme) on $\Shum{K, T}$ that is insensitive to the away-from-$p$ level $K^p$. Then for any $K^p_1\subset K^p_2$, when $\Shum{K_pK_1^p} \to \Shum{K_pK_2^p}$ is finite \'etale (which is always true in case $(4)'$), $Y_{K_pK^p_1}$ is a well-positioned subset (resp. subscheme) of $\Shum{K_pK^p_1, T}$ if and only if $Y_{K_pK^p_2}$ is a well-positioned subset (resp. subscheme) of $\Shum{K_pK^p_2, T}$.
    \end{corollary}

    \subsubsection{}
    \begin{lemma}\label{lemma: pullback of stratifications}
        Let $Y$ be a well-positioned subset (resp. subscheme). Then the \'etale neighourhoods $\ovl{U}_Y \to Y^{\tor}_{\Sigma}$ and $\ovl{U}_Y \to Y^{\natural}_{\Xi(\sigma)}$ in \cite[Theorem 2.3.2(7)]{lan2018compactifications} can be taken as the pullbacks of \'etale neighourhoods $\ovl{U} \to \Shumc{K, \Sigma, T}$ and $\ovl{U} \to \Xi(\sigma)_T$ in \cite[Proposition 2.1.2(9)]{lan2018compactifications}.
        
        More precisely, let $x \in \mathcal{Z}([\Phi, \sigma])_T$, as in \cite[Proposition 2.1.2(9)]{lan2018compactifications}, there exists an \'etale neighourhoods $\ovl{U} \to \Shumc{K, \Sigma, T}$ of $x$ and an \'etale morphism $\ovl{U} \to \Xi(\sigma)_T$ respecting $x$ such that the stratification of $\ovl{U}$ induced by that of $\Shumc{K, \Sigma, T}$ coincides with the stratification of $\ovl{U}$ induced by that of $\Xi(\sigma)_T$ (see \emph{loc. cit.} for the precise meaning). 

        Assume $x \in Y_{\mathcal{Z}([\Phi, \sigma])}$. By refining $\ovl{U}$, let $\ovl{U}_Y \to Y^{\tor}_{\Sigma}$ be the pullback of the \'etale neighourhood $\ovl{U} \to \Shumc{K, \Sigma, T}$ along $Y^{\tor}_{\Sigma} \to \Shumc{K, \Sigma, T}$, then the induced morphism $\ovl{U}_Y \to \ovl{U} \to \Xi(\sigma)_T$ factors through $Y^{\natural}_{\Xi(\sigma)}$. In particular, the stratification of $\ovl{U}_Y$ induced by that of $Y^{\tor}_{\Sigma}$ (which is induced by $\Shumc{K, \Sigma}$) coincides with the stratification of $\ovl{U}_Y$ induced by that of $Y^{\natural}_{\Xi(\sigma)}$ (which is induced by $\Xi(\sigma)$), as stated in \cite[Theorem 2.3.2(7)]{lan2018compactifications}.
    \end{lemma}
    \begin{proof}

        By taking open subschemes of $\Shum{K, T}$, we can assume $Y$ is closed. Let $\ovl{U}_Y \to Y^{\tor}_{\Sigma}$, $\ovl{U}_Y \to Y^{\natural}_{\Xi(\sigma)}$ be the \'etale neighbourhoods of $x$ in \cite[Theorem 2.3.2(7)]{lan2018compactifications}. By shrinking $\ovl{U}_Y$, $\ovl{U}_Y$ is the restriction of an \'etale neighbourhood $\ovl{U}$ of $x \in \Shum{K, T}$. We furthur refine $\ovl{U}$ such that the neighbourhood $\ovl{U} \to \Shumc{K, \Sigma, T}$, $\ovl{U} \to \Xi(\sigma)_T$ satisfies \cite[Propsition 2.1.2(9)]{lan2018compactifications}. Then the restriction of $\ovl{U}$ on $Y_{\Sigma}^{\tor}$ furthur refines $\ovl{U}_Y$, and this $\ovl{U}$ is the \'etale neighourhood we need. 

   \end{proof}

    \begin{corollary}\label{corollay: cor 2.1.7 for well positioned Y}
        Let $x \in Y_{\mathcal{Z}([\Phi, \sigma])}$. Let $E := \mathbf{E}_K(\Phi)$, it is a split torus over $\Spec \Z$. Then there exists an \'etale neighborhood $\ovl{U}_Y \to Y^{\tor}_{\Sigma}$ of $x$ and an \'etale morphism $\ovl{U}_Y \to E(\sigma) \times_{\Spec \Z} Y^{\natural}_C$ such that the stratifications of $\ovl{U}$ induced by that of $Y^{\tor}_{\Sigma}$ and by that of $E(\sigma)$ coincide with each other, and the pullbacks of these \'etale morphisms to $\underset{\tau\in \Sigma(\Phi)^{+},\ \bar{\tau}\subset\bar{\sigma}}{\cup} Y_{\mathcal{Z}([\Phi, \tau])}$ and to $\underset{\tau\in \Sigma(\Phi)^{+},\ \bar{\tau}\subset\bar{\sigma}}{\cup} E_{\tau} \times_{\Spec \Z} Y^{\natural}_C$ are open immersions.
    \end{corollary}
    \begin{proof}
        This result is ture when $Y = \Shum{K}$, see \cite[Corollary 2.1.7]{lan2018compactifications}. For general $Y$, we apply Lemma \ref{lemma: pullback of stratifications}. Note that Zariski locally over $C$, $\Xi(\sigma) \to C$ is isomorphic to $E(\sigma) \times_{\Spec \Z} C \to C$ due to \cite[Lemma 2.3]{lan2018nearby}, thus Zariski locally over $Y^{\natural}_C$, $Y^{\natural}_{\Xi(\sigma)} \to Y^{\natural}_C$ is isomorphic to $E(\sigma) \times_{\Spec \Z} Y^{\natural}_C \to Y^{\natural}_C$.
    \end{proof}

    It is a formal consequence of \cite[Theorem 2.3.2]{lan2018compactifications}, Lemma \ref{lemma: pullback of well-positioned is well-positioned} and \ref{lemma: pullback of stratifications} that:
    \begin{lemma}\label{lemma: analog covering result for well-position}
        Keep the notations as in Lemma \ref{lemma: pullback of well-positioned is well-positioned} and assume $\Sigma_1$ is induced by $\Sigma_2$. Let $Y_2$ be a well-positioned subscheme of $\Shum{K_2, T}$, $Y_1$ be the pullback of $Y_2$ under $\Shum{K_1, T} \to \Shum{K_2, T}$, and assume \ref{assumption: C to Z} is true for both $\Shum{K_1}$ and $\Shum{K_2}$. Then the statements in \cite[Proposition 3.4, Corollary 3.7, 3.8]{lan2022closed} hold when we replace
        \begin{equation}\label{eq: Shum setting}
            \Shum{K_1},\ \Shumc{K_1, \Sigma_1},\ \Shumm{K_1},\ \Shum{K_2},\ \Shumc{K_2, \Sigma_2},\ \Shumm{K_2},
        \end{equation}
        with
        \begin{equation}\label{eq: well-position setting}
            Y_1,\ Y_{1, \Sigma_1}^{\tor},\ Y_1^{\min},\ Y_2,\ Y^{\tor}_{2, \Sigma_2},\ Y_2^{\min}.
        \end{equation}
        These statements are true when we furthur replace $Y_1$ with a clopen subscheme of $Y_1$ and replace $Y_{1, \Sigma_1}^{\tor}$ and $Y_1^{\min}$ respectively, using Proposition \ref{proposition: open-closed subschemes are well positioned}.
    \end{lemma}

\subsubsection{}
    The following simple observation is crucial in Proposition \ref{proposition: open-closed subschemes are well positioned}:
    \begin{lemma}\label{lemma: exclusiveness is redundant}
        Let $Y$ be a well-positioned subset (resp. subscheme) of $\Shum{K, T}$, then $\oint_{K, \Sigma}^{-1}(Y^{\min}) = Y^{\tor}_{\Sigma}$ as locally closed subsets. In particular, since $\oint_{K, \Sigma}$ is surjective and has geometrically connected fibers, then $\pi_Y := \oint_{K, \Sigma}|_{Y^{\tor}_{\Sigma}}: Y^{\tor}_{\Sigma} \to Y^{\min}$ is also surjective and has geometrically connected fibers.
    \end{lemma}
    \begin{proof}
        The proof bases on \cite[Theorem 2.3.2(3)]{lan2018compactifications}. Let $x \in \mathcal{Z}([\Phi, \sigma])$ with image in $Y_{[\Phi]} \subset \mathcal{Z}([\Phi])$, it suffices to show $x \in Y_{\mathcal{Z}([\Phi, \sigma])}$. Under the canonical identification $\Xi_{\sigma}(\Phi)\cong\mathcal{Z}([\Phi, \sigma])$, the pullback of $Y_{\mathcal{Z}([\Phi, \sigma])}$ identifies with the pullback of $Y_{C(\Phi)}^{\natural}$ along $\Xi_{\sigma}(\Phi) \to C(\Phi)$. On the other hand, due to the commutative diagram \ref{diag: strata in two sides}, the underlying topological space of $Y_{C(\Phi)}^{\natural}$ is the preimage of $Y_{[\Phi]}$ under $C(\Phi) \to \Zb(\Phi) \cong \mathcal{Z}([\Phi])$.
    \end{proof}

      \begin{definition}\label{def: topologically disjoint}
       Given locally closed subschemes (resp. subsets) $Y_1, Y_2 \subset X$ such that $Y_1 \cap Y_2 = \emptyset$, we say $Y = Y_1 \sqcup Y_2$ is a \emph{topologically disjoint union} (or $Y_1$ and $Y_2$ are \emph{topologically disjoint}), if $Y_1$, $Y_2$ are open and closed in $Y$.
   \end{definition}

    \begin{proposition}\label{proposition: open-closed subschemes are well positioned}
       Let $Y$ be a well-positioned locally closed subset (resp. subscheme) on $\Shum{K, T}$ with respect to $Y^{\natural}= \lrbracket{Y^{\natural}(\Phi)}_{\Phi}$, assume $Y$ is a topologically disjoint union of two locally closed subsets (resp. subschemes) $Y_1 \sqcup Y_2$, then $Y_1$ and $Y_2$ are well-positioned subsets (resp. subschemes). Moreover, $Y_{\Sigma}^{\tor} = Y_{1, \Sigma}^{\tor} \sqcup Y_{2, \Sigma}^{\tor}$, $Y^{\min} = Y^{\min}_1 \sqcup Y^{\min}_2$ as topologically disjoint unions.
    \end{proposition}
    \begin{proof}

    Let $Y_{i, \Sigma}^{\tor}$ (resp. $Y_i^{\min}$) be the closure of $Y_i$ in $Y^{\tor}_{\Sigma}$ (resp. $Y^{\min}$), then $\ovl{Y}^{\tor}_{\Sigma} = \ovl{Y}^{\tor}_{1, \Sigma} \cup \ovl{Y}^{\tor}_{2, \Sigma}$ (resp. $Y^{\min} =  Y^{\min}_{1} \cup Y^{\min}_{2}$). Moreover, when $Y = Y_1 \sqcup Y_2$ is a topological disjoint union, then $Y_{1, \Sigma}^{\tor} \cap Y_{2, \Sigma}^{\tor} = \emptyset$, see \cite[Proposition 2.3.11]{lan2018compactifications} (Although the statement of Proposition \cite[Proposition 2.3.11]{lan2018compactifications} requires that $Y^{\natural}_{C(\Phi)} \to Y^{\natural}(\Phi)$ has connected fibers for each $[\Phi]$, the first two paragraphs of the proof showing $Y_{1, \Sigma}^{\tor} \cap Y_{2, \Sigma}^{\tor} = \emptyset$ does not need this assumption.) In particular, $Y^{\tor}_{\Sigma} =  Y^{\tor}_{1, \Sigma} \sqcup Y^{\tor}_{2, \Sigma}$ is a topological disjoint union.

     The proper surjective morphism $\pi_Y$ is universally closed, thus $Y_i^{\min} \subset \pi_Y(Y_{i, \Sigma}^{\tor})$. On the other hand, $\pi_Y$ has geometrically connected fibers due to Lemma \ref{lemma: exclusiveness is redundant}, thus $\pi_Y(Y_{1, \Sigma}^{\tor}) \cap \pi_Y(Y_{2, \Sigma}^{\tor})$ is empty. This forces $Y^{\min}_i = \pi_Y(Y^{\tor}_i)$ and $Y^{\min}_1 \cap Y^{\min}_2 = \emptyset$, $Y^{\min} =  Y^{\min}_1 \sqcup Y^{\min}_2$ is a topological disjoint union.

    Let $Y_{i, \mathcal{Z}([\Phi, \sigma])} = Y_{i, \Sigma}^{\tor} \cap \mathcal{Z}([\Phi, \sigma])$, $Y_{i, \mathcal{Z}([\Phi])} = Y_i^{\min} \cap \mathcal{Z}([\Phi])$, then $Y_{\mathcal{Z}([\Phi, \sigma])} = Y_{1, \mathcal{Z}([\Phi, \sigma])} \sqcup Y_{2, \mathcal{Z}([\Phi, \sigma])}$, $Y_{\mathcal{Z}[\Phi]} = Y_{1, \mathcal{Z}([\Phi])} \sqcup Y_{2, \mathcal{Z}([\Phi])}$. Since given any point $y_i \in Y_{i, \mathcal{Z}([\Phi, \sigma])}$, there exists a point $x_i \in Y_i$ which specializes to $y_i$, and since $\pi_Y$ is a closed map (in topology), then $y_i$ maps to a point $z_i \in Y_{i, \mathcal{Z}([\Phi])}$ which is a specialization of $x_i$. In particular, since $\pi_Y|_{[\Phi, \sigma]}: Y_{\mathcal{Z}([\Phi, \sigma])} \to Y_{\mathcal{Z}[\Phi]}$ is surjective, then $\pi_Y|_{[\Phi, \sigma]}$ maps $Y_{i, \mathcal{Z}([\Phi, \sigma])}$ onto $Y_{i, \mathcal{Z}([\Phi])}$. The topologically disjoint union $Y_{\mathcal{Z}[\Phi]} = Y_{1, \mathcal{Z}([\Phi])} \sqcup Y_{2, \mathcal{Z}([\Phi])}$ induces canonically a disjoint union $Y^{\natural}(\Phi) = Y^{\natural}_1(\Phi) \sqcup Y^{\natural}_2(\Phi)$ under the canonical morphism $Y^{\natural}(\Phi) \to Y_{\mathcal{Z}[\Phi]}$. Under the canonical isomorphism $Y_{\mathcal{Z}([\Phi, \sigma])} \rightiso Y^{\natural}_{\Xi(\Phi)_{\sigma}}$, due to the above commutative diagram \ref{diag: strata in two sides}, $Y_{i, \mathcal{Z}([\Phi, \sigma])}$ is canonically isomorphic to the pullback $Y^{\natural}_{i, \Xi(\Phi)_{\sigma}}$ of $Y^{\natural}_i(\Phi)$ along $\Xi(\Phi)_{\sigma} \to \Zb(\Phi) = \mathcal{Z}([\Phi])$ (in the case of subsets, it is the pullback of $Y^{\natural}_{i, C(\Phi)}$ along $\Xi(\Phi)_{\sigma} \to C(\Phi)$, where $Y^{\natural}_{i, C(\Phi)}$ is the pullback of $Y^{\natural}_i(\Phi)$ along $C(\Phi) \to \Zb(\Phi)$ endowed with the induced reduced subscheme structure), thus the canonical isomorphism 
        \[Y^{\natural}_{\mathfrak{X}^{\circ}_{\sigma}(\Phi)} \rightiso (Y^{\tor}_{\Sigma})_{\underset{\tau\in \Sigma(\Phi)^{+},\ \bar{\tau}\subset\bar{\sigma}}{\cup} Y_{\mathcal{Z}([\Phi, \tau])}}^{\wedge}\] 
        induces an isomorphism 
        \[Y^{\natural}_{\mathfrak{X}^{\circ}_{\sigma}(\Phi)} \times_{Y^{\natural}(\Phi)} Y^{\natural}_i(\Phi) \rightiso (Y^{\tor}_{i, \Sigma})_{\underset{\tau\in \Sigma(\Phi)^{+},\ \bar{\tau}\subset\bar{\sigma}}{\cup} Y_{i, \mathcal{Z}([\Phi, \tau])}}^{\wedge}. \]
       The pullback of $Y^{\natural}_i(\Phi)$ over $W^0$ is equal to the pullback of $Y_i$, thus $Y_i$ is a well-positioned subset (resp. subscheme) with respect to $Y_i^{\natural} = \lrbracket{Y^{\natural}_{i}(\Phi)}_{\Phi}$.

       Let $\ovl{Y}$, $Y_0$, $\ovl{Y}_i$, $Y_{i, 0}$ be defined as before, and let $Y_i^{\min, \prime}$ be the partial minimal compactifications of $Y_i$ (we add a superscript $\prime$ since the usual notion $Y^{\min}_{i}$ has been occupied in the proof). It follows easily that $\ovl{Y} = \ovl{Y}_1 \cup \ovl{Y}_2$, $Y_0 = Y_{1, 0} \cup Y_{2, 0}$, $\ovl{Y}^{\min} = \ovl{Y}^{\min}_{1} \cup \ovl{Y}^{\min}_{2}$, $Y^{\min}_{0} = Y^{\min}_{1, 0} \cup Y^{\min}_{2, 0}$, $Y^{\min} = \ovl{Y}^{\min} \setminus Y^{\min}_{0} \subset Y^{\min, \prime}_1 \cup Y^{\min, \prime}_2$. Since the closure of $Y_i$ in $Y^{\min}$ is the intersection of $Y^{\min}$ with the closure of $Y_i$ in $\ovl{Y}^{\min}$, which is $\ovl{Y}_i^{\min} \cap Y^{\min} \subset Y_i^{\min, \prime}$. In particular, we have a natural inclusion $Y_i^{\min} \subset Y_i^{\min, \prime}$. This is an equality as subsets, since their intersections with each $\mathcal{Z}([\Phi])$ are both $Y_{i, \mathcal{Z}([\Phi])}$: for $Y_i^{\min}$, this comes from definition, and for $Y_i^{\min, \prime}$, this comes from the fact that $Y_i$ is well-positioned  with respect to $Y_i^{\natural} = \lrbracket{Y^{\natural}_{i}(\Phi)}_{\Phi}$. In particular, $Y_i^{\min}$ is the partial minimal compactification of $Y_i$. Similarly, $Y_{i, \Sigma}^{\tor}$, defined as the closure of $Y_i$ in $Y_{\Sigma}^{\tor}$, is the partial toroidal compactification of $Y_i$.
       
    \end{proof}
    \begin{remark}
        This Proposition \ref{proposition: open-closed subschemes are well positioned} shows that the assumption $Y^{\natural}_{C(\Phi)} \to Y^{\natural}(\Phi)$ has connected fibers for each $[\Phi]$ in \cite[Proposition 2.3.11]{lan2018compactifications} is redundant.
    \end{remark}
    \begin{corollary}
        Let $Y$ be a well-positioned subscheme, then $\pi_0(Y) \to \pi_0(Y_{\Sigma}^{\tor})$, $\pi_0(Y) \to \pi_0(Y^{\min})$ are bijection. In particular, if $\Shum{K}$ has geometrically reduced fibers, then all geometric fibers of $\Shum{K}$ have the same number of connected components.
    \end{corollary}
    \begin{proof}
        $\pi_0(Y) \to \pi_0(Y^{\tor}_{\Sigma})$ is surjective since $Y \hookrightarrow Y^{\tor}_{\Sigma}$ is open and dense, it is injective due to Proposition \ref{proposition: open-closed subschemes are well positioned}. Similarly for $\pi_0(Y) \to \pi_0(Y^{\min})$. We show the last statement: Since $\Shumc{K, \Sigma} \to \Spec \OO_{E(v)}$ is flat, proper, and has geometrically reduced fibers, then due to \cite[\href{https://stacks.math.columbia.edu/tag/0E0N}{Tag 0E0N}]{stacks-project}, all geometric fibers of $\Shumc{K, \Sigma}$ have the same number of connected components. Note that $\Shum{K}$ itself is a well-positioned subscheme ($T = \Spec \OO_{E(v)}$), its geometric fibers are all well-positioned. Due to the first part of the statement, all geometric fibers of $\Shum{K}$ have the same number of connected components.
    \end{proof}

\subsection{Generalized definition}\label{sec: generalized definition of well position}

   In \cite{caraiani2024generic}, the authors proved that Igusa varieties over certain PEL type Shimura varieties are well-positioned. In order to adapt the situations there, and prove that Igusa varieties over Hodge-type Shimura varieties are well-positioned, we should slightly generalize the definition of being a well-positioned scheme.

   Let $Y \subset \Shum{K, T}$ be a well-positioned subscheme with partial toroidal compactification $Y^{\tor}_{\Sigma} \subset \Shum{K, \Sigma, T}$. Let $\wdt{Y}^{\tor}_{\Sigma}$ be a scheme over $Y^{\tor}_{\Sigma}$, and $\wdt{Y} \subset \wdt{Y}^{\tor}_{\Sigma}$ be the pullback of $Y \subset Y^{\tor}_{\Sigma}$. Given each boundary stratum $Y_{\mathcal{Z}([\Phi, \sigma])} \subset Y^{\tor}_{\Sigma}$, let $\wdt{Y}_{\mathcal{Z}([\Phi, \sigma])} \subset \wdt{Y}^{\tor}_{\Sigma}$ be the its pullback.

   Let $\wdt{Y}^{\natural}_C$ be a scheme over $Y^{\natural}_C$. Given each $\sigma \in \Sigma(\Phi)^+$, let $(\ast)$ be any of the following: $\Xi$, $\Xi_{\sigma}$, $\Xi({\sigma})$, $\mathfrak{X}$, $\mathfrak{X}_{\sigma}$, $\mathfrak{X}_{\sigma}^{\circ}$, let $\wdt{Y}^{\natural}_{(\ast)}$ be the pullback of $\wdt{Y}^{\natural}_C$ along $\wdt{Y}^{\natural}_C \to Y^{\natural}_C$ (since $\wdt{Y}^{\natural}_{C} \to Y^{\natural}_C$ might not be locally of finitely presented, we take the pullback in the category of formal schemes, i.e., using completed tensor products instead of the usual tensor products).

  Proposition \ref{prop: Y same as Shum} implies that
   \begin{equation}\label{eq: Y well-position}
       Y_{\Xi_{\sigma}}^{\natural} \cong Y_{\mathcal{Z}([\Phi, \sigma])},\quad Y^{\natural}_{\mathfrak{X}_{\sigma}} \cong (Y^{\tor}_{\Sigma})^{\wedge}_{Y_{\mathcal{Z}([\Phi, \tau])}},
   \end{equation}
   \begin{equation*}
       Y^{\natural}_{\mathfrak{X}_{\sigma}^{\circ}} \cong (Y^{\tor}_{\Sigma})_{\underset{\tau\in \Sigma(\Phi)^{+},\ \bar{\tau}\subset\bar{\sigma}}{\cup} Y_{\mathcal{Z}([\Phi, \tau])}}^{\wedge}, \quad Y^{\natural}_{\mathfrak{X}} \cong (Y^{\tor}_{\Sigma})_{\underset{[\sigma] \in \Lambda_{K, \Phi}\backslash\Sigma(\Phi)^+}{\cup} Y_{\mathcal{Z}([\Phi, \sigma])}}^{\wedge}, 
   \end{equation*}

   Let $\mathfrak{W} = \Spf R \subset Y^{\natural}_{\mathfrak{X}_{\sigma}^{\circ}}$ be an open formal scheme, consider the induced maps $W = \Spec R \to Y_{\Xi(\sigma)}^{\natural}$ and $W \to Y^{\tor}_{\Sigma}$, the preimages of the boundary stratifications on the targets coincide in $W$.

   \begin{lemma}\label{lemma: fpqc W}
       Such coverings $\bigcup W \to Y_C^{\natural}$ (resp. $\bigcup W^0 \to Y_C^{\natural}$) are fpqc coverings up to refining $\bigcup \mathfrak{W} \subset Y^{\natural}_{\mathfrak{X}_{\sigma}^{\circ}}$. If moreover $C \to \Zb^{\bigsur}$ is flat, then same statements are true when we replace $Y_C^{\natural}$ with $Y^{\natural}_{\Zb^{\bigsur}}$.
   \end{lemma}
   \begin{proof}
      Each $W \to Y_{\Xi(\sigma)}^{\natural}$ (resp. $W^0 \to Y_{\Xi}^{\natural}$) is flat, and $Y_{\Xi(\sigma)}^{\natural} \to Y_{C}^{\natural}$ (resp. $Y_{\Xi}^{\natural} \to Y_{C}^{\natural}$) is flat since it is the pullback of the smooth morphism $\Xi(\sigma) \to C$ (resp. $\Xi \to C$). Here Zariski locally on $C$, $\Xi(\sigma) \to C$ (resp. $\Xi \to C$) is the projection $E(\sigma) \times C \to C$ (resp. $E \times C \to C$) with $E \to E(\sigma)$ a toroidal embedding of a split torus $E$ over $\Spec \Z$.

      Given each point $y \in Y^{\natural}_C$, since $\Xi \to C$ and $\mathfrak{X}_{\sigma}^{\circ} \to C$ are surjections, we can always find a point $y_1 \in Y_{\Xi}^{\natural}$ which specializes to a point $y_2 \in Y_{\mathfrak{X}_{\sigma}^{\circ}}^{\natural}$ with image $y$, and such that there is an open formal subscheme $\mathfrak{W} \subset Y_{\mathfrak{X}_{\sigma}^{\circ}}^{\natural}$, $y_1$ and $y_2$ have images in the image of $W \to Y^{\natural}_{C}$, in particular, $y_1$ is in the image of $W^0 \to Y^{\natural}_{C}$, $W$ and $W^0$ have same image in $Y^{\natural}_{C}$.

      Finally, since $\bigcup \mathfrak{W} \subset Y_{\mathfrak{X}_{\sigma}^{\circ}}^{\natural}$ is an open affine covering, and $Y_{\mathfrak{X}_{\sigma}^{\circ}}^{\natural} \to Y_C^{\natural}$ is quasi-compact, the second condition of being quasi-compact in \cite[\href{https://stacks.math.columbia.edu/tag/03NW}{Tag 03NW}]{stacks-project} is also verified for $\bigcup W \to Y_C^{\natural}$. Thus $\bigcup W \to Y_C^{\natural}$ is a fpqc covering. Since each $W$ and $W^0$ have same image in $Y^{\natural}_{C}$, then $\bigcup W^0 \to Y_C^{\natural}$ is a fpqc covering. 
   \end{proof}

   \begin{definition}\label{def: well-positioned over well-positioned}
       We say $\wdt{Y} \to Y$ is well-positioned with respect to morphisms $\lrbracket{\wdt{Y}^{\natural}_{C(\Phi)} \to Y^{\natural}_{C(\Phi)}}_{[\Phi]}$, if $\wdt{Y} \to Y$ is affine, flat, and for some (thus for all) coverings $\mathfrak{W}$, the pullbacks of $\wdt{Y}^{\tor}_{\Sigma} \to Y^{\tor}_{\Sigma}$ and $\wdt{Y}^{\natural}_C \to Y^{\natural}_C$ coincide.
   \end{definition}
   \begin{lemma}\label{lemma: fpqc descent for Ytor and Yc}
       If $\wdt{Y} \to Y$ is well-positioned with respect to $\lrbracket{\wdt{Y}^{\natural}_{C(\Phi)} \to Y^{\natural}_{C(\Phi)}}_{[\Phi]}$, then each $\wdt{Y}^{\natural}_{C(\Phi)} \to Y^{\natural}_{C(\Phi)}$ is also affine and flat, and $\wdt{Y}^{\tor}_{\Sigma} \to Y^{\tor}_{\Sigma}$ is affine and flat.
   \end{lemma}
   \begin{proof}
       Since $\wdt{Y}^{\tor}_{\Sigma} \to Y^{\tor}_{\Sigma}$ and $\wdt{Y}^{\natural}_C \to Y^{\natural}_C$ coincide over $W$, then $\wdt{Y} \to Y$ and $\wdt{Y}^{\natural}_C \to Y^{\natural}_C$ coincide over $W^0$. Since $\bigcup W^0 \to Y^{\natural}_C$ is a fpqc covering due to Lemma \ref{lemma: fpqc W}, then $\wdt{Y}^{\natural}_C \to Y^{\natural}_C$ is affine and flat by fpqc descent. Since $Y \cup (\bigcup W) \to Y^{\tor}_{\Sigma}$ is a fqpc covering, then by fpqc descent, $\wdt{Y}^{\natural}_C \to Y^{\natural}_C$ is affine and flat 
 implies that $\wdt{Y}^{\tor}_{\Sigma} \to Y^{\tor}_{\Sigma}$ is affine and flat.
   \end{proof}

   \begin{remark}
       Flatness ensures that $\wdt{Y}_{\mathcal{Z}([\Phi, \sigma])}$ give a stratification of the boundary of $\wdt{Y}^{\tor}_{\Sigma}$, and $\wdt{Y}$ is dense in $\wdt{Y}^{\tor}_{\Sigma}$, see \cite[\href{https://stacks.math.columbia.edu/tag/081H}{Tag 081H}]{stacks-project} (note that all the base schemes related to $Y$ are locally noetherian, thus the quasi-compactness comes for free, see \cite[\href{https://stacks.math.columbia.edu/tag/01OX}{Tag 01OX}]{stacks-project}).
   \end{remark}
   \begin{remark}
       Affineness of $\tilde{Y} \to Y$ and $\lrbracket{\wdt{Y}^{\natural}_{C(\Phi)} \to Y^{\natural}_{C(\Phi)}}_{[\Phi]}$ is to ensure the effectiveness of fpqc descent of morphisms. Here we regard these schemes $\tilde{Y}$, $\wdt{Y}^{\natural}_{C(\Phi)}$ as quasi-coherent sheaves over corresponding bases. Since we apply this definition to Igusa varieties which are group torsors over the base, the affineness is guaranteed.
   \end{remark}
 
   \begin{proposition}\label{prop: wdt Y satisfy same prop}
        Let $\wdt{Y} \to Y$ be well-positioned with respect to $\lrbracket{\wdt{Y}^{\natural}_{C(\Phi)} \to Y^{\natural}_{C(\Phi)}}_{[\Phi]}$, then \cite[Theorem 2.3.2 (5)(6)(7)]{lan2018compactifications} holds for $(\wdt{Y}, \lrbracket{\wdt{Y}^{\natural}_{C(\Phi)}}_{[\Phi]})$. That is to say
        \begin{enumerate}
            \item We have canonical isomorphisms
            \begin{equation}
                \wdt{Y}_{\Xi_{\sigma}}^{\natural} \cong \wdt{Y}_{\mathcal{Z}([\Phi, \sigma])},\quad \wdt{Y}^{\natural}_{\mathfrak{X}_{\sigma}} \cong (\wdt{Y}^{\tor}_{\Sigma})^{\wedge}_{\wdt{Y}_{\mathcal{Z}([\Phi, \tau])}},
            \end{equation}
            \begin{equation*}
           \wdt{Y}^{\natural}_{\mathfrak{X}_{\sigma}^{\circ}} \cong (\wdt{Y}^{\tor}_{\Sigma})_{\underset{\tau\in \Sigma(\Phi)^{+},\ \bar{\tau}\subset\bar{\sigma}}{\cup} \wdt{Y}_{\mathcal{Z}([\Phi, \tau])}}^{\wedge} \quad \wdt{Y}^{\natural}_{\mathfrak{X}} \cong (\wdt{Y}^{\tor}_{\Sigma})_{\underset{[\sigma] \in \Lambda_{K, \Phi}\backslash\Sigma(\Phi)^+}{\cup} \wdt{Y}_{\mathcal{Z}([\Phi, \sigma])}}^{\wedge} 
            \end{equation*}

            that are compactible with \ref{eq: Y well-position}.
            \item For each $\sigma \in \Sigma(\Phi)^+$, and for each open formal subscheme $\mathfrak{W} = \Spf R \subset \wdt{Y}^{\natural}_{\mathfrak{X}_{\sigma}^{\circ}}$ coming from an open formal subscheme in $Y^{\natural}_{\mathfrak{X}^{\circ}_{\sigma}}$ along the affine morphism $\tilde{Y}^{\natural}_{\mathfrak{X}^{\circ}_{\sigma}} \to Y^{\natural}_{\mathfrak{X}^{\circ}_{\sigma}}$. Consider the induced maps $W = \Spec R \to \wdt{Y}_{\Xi(\sigma)}^{\natural}$ and $W \to \wdt{Y}^{\tor}_{\Sigma}$, the preimages of the boundary stratifications (induced by the boundary stratification of $Y_{\Xi(\sigma)}^{\natural}$ and $Y^{\tor}_{\Sigma}$ respectively) on the targets coincide in $W$.
            \item Let $x$ be a point of $\wdt{Y}^{\natural}_{\Xi_{\sigma}(\Phi)} \cong \wdt{Y}_{\mathcal{Z}([\Phi, \sigma])}$, there exists an \'etale neighborhood $\ovl{U} \to \wdt{Y}^{\tor}_{\Sigma}$ of $x$ and an \'etale morphism $\ovl{U} \to \wdt{Y}^{\natural}_{\Xi(\sigma)}$ respecting $x$ such that the stratifications of $\ovl{U}$ induced by $\wdt{Y}^{\tor}_{\Sigma}$ and by $\wdt{Y}^{\natural}_{\Xi(\sigma)}$ coincides. Moreover, the pullbacks of these \'etale morphisms to $\underset{\tau\in \Sigma(\Phi)^{+},\ \bar{\tau}\subset\bar{\sigma}}{\cup} \wdt{Y}_{\mathcal{Z}([\Phi, \tau])}$ and to $\underset{\tau\in \Sigma(\Phi)^{+},\ \bar{\tau}\subset\bar{\sigma}}{\cup} \wdt{Y}^{\natural}_{\Xi_{\tau}}$ are open immersions.
            \item (This is slightly stronger than $(3)$) Let $x \in \wdt{Y}_{\mathcal{Z}([\Phi, \sigma])}$. Let $E := \mathbf{E}_K(\Phi)$, it is a split torus over $\Spec \Z$. Then there exists an \'etale neighborhood $\ovl{U}_{\wdt{Y}} \to \wdt{Y}^{\tor}_{\Sigma}$ of $x$ and an \'etale morphism $\ovl{U}_{\wdt{Y}} \to E(\sigma) \times_{\Spec \Z} \wdt{Y}^{\natural}_C$ such that the stratifications of $\ovl{U}$ induced by that of $\wdt{Y}^{\tor}_{\Sigma}$ and by that of $E(\sigma)$ coincide with each other, and the pullbacks of these \'etale morphisms to $\underset{\tau\in \Sigma(\Phi)^{+},\ \bar{\tau}\subset\bar{\sigma}}{\cup} \wdt{Y}_{\mathcal{Z}([\Phi, \tau])}$ and to $\underset{\tau\in \Sigma(\Phi)^{+},\ \bar{\tau}\subset\bar{\sigma}}{\cup} E_{\tau} \times_{\Spec \Z} \wdt{Y}^{\natural}_C$ are open immersions.
        \end{enumerate}
   \end{proposition}
   \begin{proof}
       \begin{enumerate}
           \item First of all, the pullback of $W \to Y^{\tor}_{\Sigma}$ along $Y_{\mathcal{Z}([\Phi, \sigma])} \to Y^{\tor}_{\Sigma}$ is canonically isomorphic to the pullback of $W \to Y^{\natural}_{\Xi(\sigma)}$ along $Y^{\natural}_{\Xi_{\sigma}} \to Y^{\natural}_{\Xi(\sigma)}$, under the canonical isomorphism $Y_{\mathcal{Z}([\Phi, \sigma])} \cong Y^{\natural}_{\Xi_{\sigma}}$. Such pullbacks give a fpqc covering over the stratum $Y_{\mathcal{Z}([\Phi, \sigma])} \cong Y^{\natural}_{\Xi_{\sigma}}$. By fpqc descent, we have $\wdt{Y}_{\Xi_{\sigma}}^{\natural} \cong \wdt{Y}_{\mathcal{Z}([\Phi, \sigma])}$. The second isomorphism comes from the third, and the third isomorphism can also be checked over $W$: since the stratifications of $\wdt{Y}_{\Sigma}^{\tor}$ and $\wdt{Y}_{\Xi(\sigma)}^{\natural}$ are induced by the stratifications of $Y_{\Sigma}^{\tor}$ and $Y_{\Xi(\sigma)}^{\natural}$, over $W$, every single pair of strata $\wdt{Y}_{\Xi_{\sigma}}^{\natural} \cong \wdt{Y}_{\mathcal{Z}([\Phi, \sigma])}$ is matched, the base $\wdt{Y}^{\natural}_{\Xi(\sigma)} \cong \wdt{Y}^{\tor}_{\Sigma}$ (over $W$) is matched, and more importantly, these stratifications are the pullbacks of the stratification of $W$ (induced by $Y^{\natural}_{\Xi(\sigma)}$ and $Y^{\tor}_{\Sigma}$). We do formal completions along the boundary of $W$ which gives the isomorphism $\wdt{Y}^{\natural}_{\mathfrak{X}_{\sigma}^{\circ}} \cong (\wdt{Y}^{\tor}_{\Sigma})_{\underset{\tau\in \Sigma(\Phi)^{+},\ \bar{\tau}\subset\bar{\sigma}}{\cup} \wdt{Y}_{\mathcal{Z}([\Phi, \tau])}}^{\wedge} $ on $\mathfrak{W}$. Run $\mathfrak{W}$ as an open covering of $Y^{\natural}_{\mathfrak{X}_{\sigma}^{\circ}}$, and patch the isomorphisms, we have the wanted isomorphism over $Y^{\natural}_{\mathfrak{X}_{\sigma}^{\circ}} \cong (Y^{\tor}_{\Sigma})_{\underset{\tau\in \Sigma(\Phi)^{+},\ \bar{\tau}\subset\bar{\sigma}}{\cup} Y_{\mathcal{Z}([\Phi, \tau])}}^{\wedge}$. The other isomorphism follows similarly.
           \item We pullback the corresponding results on $Y^{\natural}_{\mathfrak{X}^{\circ}_{\sigma}}$ from last subsection.
           \item Same as $(2)$.
           \item This is a standard result follows from the third part of the Proposition, see the proof of Corollary \ref{corollay: cor 2.1.7 for well positioned Y}.
       \end{enumerate}
   \end{proof}

   We follow the notations in Lemma \ref{lemma: pullback of well-positioned is well-positioned}:
   \begin{lemma}\label{lemma: pullback of well-positioned is well-positioned, second}
     Assume \ref{assumption: C to Z} hold for both $\Shum{K_1}$ and $\Shum{K_2}$. Let $Y_2 \subset \Shum{K_2, T}$ be well-positioned with respect to $\lrbracket{Y^{\natural}_{2}(\Phi_2)}_{[\Phi_2]}$, and $\wdt{Y}_2 \to Y_2$ be well-positioned with respect to $\lrbracket{\wdt{Y}^{\natural}_{2, C_2(\Phi_2)} \to Y^{\natural}_{2, C_2(\Phi_2)}}_{[\Phi_2]}$. Let $Y_1$ (resp. $Y^{\natural}_{1}(\Phi_1)$) be the pullback of $Y_2$ (resp. $Y^{\natural}_{2}(\Phi_2)$) along $\Shum{K_1, T} \to \Shum{K_2, T}$ (resp. $\Zb_1(\Phi_1)_T \to \Zb_2(\Phi_2)_T$), and let $\wdt{Y}^{\tor}_{1, \Sigma_1} \to Y^{\tor}_{1, \Sigma_1}$ (resp. $\wdt{Y}^{\natural}_{1, C_1(\Phi_1)} \to Y^{\natural}_{1, C_1(\Phi_1)}$) be the pullback of $\wdt{Y}^{\tor}_{2, \Sigma_2} \to Y^{\tor}_{2, \Sigma_2}$ (resp. $\wdt{Y}^{\natural}_{2, C_2(\Phi_2)} \to Y^{\natural}_{2, C_2(\Phi_2)}$) along $Y^{\tor}_{1, \Sigma_1} \to Y^{\tor}_{2, \Sigma_2}$ (resp. $Y^{\natural}_{1, C_1(\Phi_1)} \to Y^{\natural}_{2, C_2(\Phi_2)}$). Then $\wdt{Y}_1 \to Y_1$ is well-positioned with respect to $\lrbracket{\wdt{Y}^{\natural}_{1, C_1(\Phi_1)} \to Y^{\natural}_{1, C_1(\Phi_1)}}_{[\Phi_1]}$.
   \end{lemma}
   \begin{proof}
       Due to Lemma \ref{lemma: pullback of well-positioned is well-positioned}, $Y_1$ is well-positioned with respect to $\lrbracket{Y^{\natural}_{1}(\Phi_1)}_{[\Phi_1]}$, and $Y^{\tor}_{1, \Sigma_1}$ is the pullback of $Y^{\tor}_{2, \sigma_2}$ along the finite morphism $\Shumc{K_1, \Sigma_1, T} \to \Shumc{K_2, \Sigma_2, T}$. Let $\mathfrak{W}_2 \to Y^{\natural}_{2, \mathfrak{X}^{\circ}_{2, \Sigma_2}}$ be an open formal subscheme, and $\mathfrak{W}_1 \to Y^{\natural}_{1, \mathfrak{X}^{\circ}_{1, \sigma_1}}$ be its pullback. Due to Lemma \ref{lemma: analog covering result for well-position}, we have the following commutative diagram:
\[\begin{tikzcd}
	{Y_{1, \Sigma_1}^{\tor}} & {W_1} & {Y^{\natural}_{1, \Xi_1(\Phi_1)(\sigma_1)}} & {Y^{\natural}_{1, C_1(\Phi_1)}} \\
	{Y_{2, \Sigma_2}^{\tor}} & {W_2} & {Y^{\natural}_{2, \Xi_2(\Phi_2)(\sigma_2)}} & {Y^{\natural}_{2, C_2(\Phi_2)}}
	\arrow[from=1-1, to=2-1]
	\arrow[from=1-2, to=1-1]
	\arrow[from=1-2, to=1-3]
	\arrow[from=1-2, to=2-2]
	\arrow[from=1-3, to=1-4]
	\arrow[from=1-3, to=2-3]
	\arrow[from=1-4, to=2-4]
	\arrow[from=2-2, to=2-1]
	\arrow[from=2-2, to=2-3]
	\arrow[from=2-3, to=2-4]
\end{tikzcd}\]
       Since the pullback of a covering $\mathfrak{W}_2$ of $\mathfrak{X}_{2, \sigma_2}^{\circ}$ induces a covering $\mathfrak{W}_1$ of $\mathfrak{X}_{1, \sigma_1}^{\circ}$, we are done.
   \end{proof}

   Similar as Lemma \ref{lemma: pullback of well positioned, different level}, we have
   \begin{lemma}\label{emma: pullback of well positioned over well position, different level}
       Let $Y \subset \Shum{K, T}$ be a well-positioned subscheme with respect to $\lrbracket{Y^{\natural}(\Phi)}_{[\Phi]}$, and $\wdt{Y} \to Y$ be well-positioned with respect to $\lrbracket{\wdt{Y}^{\natural}_{C(\Phi)} \to Y^{\natural}_{C(\Phi)}}_{[\Phi]}$. Let $g \in G(\A_f)$ such that $K' \subset gKg^{-1}$, $[g]: \Shum{K', T}\to \Shum{K, T}$ be the Hecke action. Consider the pullback $Y'$ of $Y$ along $[g]$ and pullbacks $Y^{\prime, \natural}=\{Y^{\prime, \natural}(\Phi')\}_{\Phi'}$ of $Y^{\natural}=\{Y(\Phi)\}_{\Phi}$ along $\Zb^{\prime}(\Phi')\to \Zb(\Phi)$ induced by $[g]^{\min}: \Shumm{K^{\prime}, T}\to \Shumm{K, T}$ for each $\Zb^{\prime}(\Phi')$ above $\Zb(\Phi)$, and the pullback $\wdt{Y}^{\prime, \tor}_{\Sigma}$ of $\wdt{Y}^{\tor}_{\Sigma}$ along $Y^{\prime, \tor}_{\Sigma'} \to Y^{\tor}_{\Sigma}$ and pullbacks $\wdt{Y}^{\prime, \natural}_{C(\Phi')}$ of $\wdt{Y}^{\natural}_{C(\Phi)}$ along $Y^{\prime, \natural}_{C(\Phi')} \to Y^{\natural}_{C(\Phi)}$ for each $[\Phi']$ above $[\Phi]$. Then $Y'$ is a well-positioned subscheme with respect to $\lrbracket{Y^{\prime, \natural}(\Phi')}_{[\Phi']}$, and $\wdt{Y}' \to Y'$ is well-positioned with respect to $\lrbracket{\wdt{Y}^{\prime, \natural}_{C(\Phi')} \to Y^{\prime, \natural}_{C(\Phi')}}_{[\Phi']}$.
   \end{lemma}

   It follows from Proposition \ref{proposition: open-closed subschemes are well positioned} that
   \begin{lemma}\label{lemma: open closed subschemes are well positioned, 2}
       Let $Y \subset \Shum{K, T}$ be a well-positioned subscheme with respect to $\lrbracket{Y^{\natural}(\Phi)}_{[\Phi]}$, and $\wdt{Y} \to Y$ be well-positioned with respect to $\lrbracket{\wdt{Y}^{\natural}_{C(\Phi)} \to Y^{\natural}_{C(\Phi)}}_{[\Phi]}$. Let $Y_1 \subset Y$ be an open and closed subscheme, $\wdt{Y}_{1, \Sigma}^{\tor} \to Y_{1, \Sigma}^{\tor}$ (resp. $\wdt{Y}^{\natural}_{1, C(\Phi)} \to Y^{\natural}_{1, C(\Phi)}$) be the pullback of $\wdt{Y}_{\Sigma}^{\tor} \to Y_{\Sigma}^{\tor}$ (resp. $\wdt{Y}^{\natural}_{C(\Phi)} \to Y^{\natural}_{C(\Phi)}$). Then $Y_1 \subset \Shum{K, T}$ is well-positioned subscheme with respect to $\lrbracket{Y^{\natural}_{1}(\Phi)}_{[\Phi]}$, and $\wdt{Y}_1 \to Y_1$ is well-positioned with respect to $\lrbracket{\wdt{Y}^{\natural}_{1, C(\Phi)} \to Y^{\natural}_{1, C(\Phi)}}_{[\Phi]}$.
   \end{lemma}

   \begin{lemma}\label{lemma: open closed subschemes are well positioned 3}
       Let $Y \subset \Shum{K, T}$ be a well-positioned subscheme with respect to $\lrbracket{Y^{\natural}(\Phi)}_{[\Phi]}$,  and $\wdt{Y} \to Y$ be well-positioned with respect to $\lrbracket{\wdt{Y}^{\natural}_{C(\Phi)} \to Y^{\natural}_{C(\Phi)}}_{[\Phi]}$. Let $\wdt{Y}_1 \subset \wdt{Y}$ be an open and closed subscheme with complement $\wdt{Y}_2$. Let $i = 1, 2$, $\wdt{Y}^{\tor}_{i, \Sigma}$ be the closure of $\wdt{Y}_i$ in $\wdt{Y}^{\tor}_{\Sigma}$, then there exists open and closed subscheme $\wdt{Y}^{\natural}_{i, C(\Phi)} \subset \wdt{Y}^{\natural}_{i, C(\Phi)}$ for each $[\Phi]$ such that $\wdt{Y}_i$ is well-positioned with respect to $\lrbracket{\wdt{Y}^{\natural}_{i, C(\Phi)} \to Y^{\natural}_{C(\Phi)}}_{[\Phi]}$. Moreover, $\wdt{Y}^{\natural}_{C(\Phi)} = \wdt{Y}^{\natural}_{1, C(\Phi)} \sqcup \wdt{Y}^{\natural}_{2, C(\Phi)}$ is a topologically disjoint union.
   \end{lemma}
\begin{proof}
    We first show that $\wdt{Y}^{\tor}_{\Sigma} = \wdt{Y}^{\tor}_{1, \Sigma} \sqcup \wdt{Y}^{\tor}_{2, \Sigma}$ (if it is disjoint, then it is automatically topologically disjoint by definition). We follow the proof of \cite[Proposition 2.3.11]{lan2018compactifications}: If $\wdt{Y}^{\tor}_{1, \Sigma} \cap \wdt{Y}^{\tor}_{2, \Sigma}$ is nonempty, we could find a boundary stratum $\mathcal{Z}([\Phi, \sigma])$ such that $\wdt{Y}^{\tor}_{1, \Sigma} \cap \wdt{Y}^{\tor}_{2, \Sigma} \cap \wdt{Y}_{\mathcal{Z}([\Phi, \sigma])}$ is nonempty and let $x$ be any point in it. Using the identification $\wdt{Y}_{\mathcal{Z}([\Phi, \sigma])}^{\natural} \cong \wdt{Y}_{\Xi_{\sigma}(\Phi)}$, we also denoty by $x \in \wdt{Y}_{\Xi_{\sigma}(\Phi)}^{\natural}$. Due to Proposition \ref{prop: wdt Y satisfy same prop}, there exist \'etale morphisms $\ovl{U} \to \wdt{Y}^{\tor}_{\Sigma}$ and $\ovl{U} \to \wdt{Y}^{\natural}_{\Xi(\sigma)}$ respecting $x$ such that the stratification on $\ovl{U}$ induced by $\wdt{Y}^{\tor}_{\Sigma}$ (which is induced by $Y^{\tor}_{\Sigma}$) coincides with the stratification induced by $\wdt{Y}^{\natural}_{\Xi(\sigma)}$ (which is induced by $Y^{\natural}_{\Xi(\sigma)}$), and their pullbacks to $\wdt{Y}_{\mathcal{Z}([\Phi, \sigma])}$ and to $\wdt{Y}_{\Xi_{\sigma}(\Phi)}$ are both open immersions. 
    
    The remaining part is similar as the second paragraph of \cite[Proposition 2.3.11]{lan2018compactifications}. We write it out in more details for the reader's convinience.

    Replace $\ovl{U}$ by an open embedding if necessary, we assume $\ovl{U} \to \wdt{Y}^{\tor}_{\Sigma}$ and $\ovl{U} \to \wdt{Y}^{\natural}_{\Xi(\sigma)}$ have connected fibers. We identify the open subscheme $U = \ovl{U} \times_{\wdt{Y}^{\tor}_{\Sigma}} \wdt{Y}$ with $\ovl{U} \times_{\wdt{Y}^{\natural}_{\Xi(\sigma)}} \wdt{Y}^{\natural}_{\Xi}$, and $\wdt{Y} = \wdt{Y}_1 \sqcup \wdt{Y}_2$ induces $U = U_1 \sqcup U_2$. Let $\ovl{U}_i$ be the closure of $U_i$ in $\ovl{U}$. Since replacing $\wdt{Y}^{\tor}_{\Sigma}$ by an open subscheme $\ovl{V}$ (the image of $\ovl{U}$) does not change the closure relations in $\wdt{Y}^{\tor}_{\Sigma}$, due to Lemma \ref{lemma: simple topo lemma, etale covering}, $x$ is in the image of $\ovl{U}_i$. Since $\wdt{Y}$ is open dense in $\wdt{Y}^{\tor}_{\Sigma}$, $U_1$ and $U_2$ are non-empty.

    If there exists $x_1 \in U_1$, $x_2 \in U_2$ such that $x_i$ specializes to $x$ and have same image $y$ in $\wdt{Y}^{\natural}_{C}$, since $ \wdt{Y}^{\natural}_{\Xi} \to \wdt{Y}^{\natural}_{C}$ is a torsor under a torus, $U \to \wdt{Y}^{\natural}_{\Xi} \to \wdt{Y}^{\natural}_{C}$ has connected fibers, this forces the fiber of $y$ is connected, which leads to a contradiction, thus such $x$ does not exist.  On the other hand, we show that we could take such $x_i$, which leads to a contradiction: Note that $\wdt{Y}^{\natural}_{\Xi} \to \wdt{Y}^{\natural}_{\Xi(\sigma)}$ is an qcqs toroidal embedding, and is Zariski locally $E \times_{\Spec \Z} \wdt{Y}^{\natural}_{C} \hookrightarrow E(\sigma) \times_{\Spec \Z} \wdt{Y}^{\natural}_{C}$, see Proposition \ref{prop: wdt Y satisfy same prop}. Therefore, one could take the specialiazations $x_i \mapsto x$ coming from $E(\Phi) \hookrightarrow E(\Phi)(\sigma)$ part.
  
  Since $\mathcal{Z}([\Phi, \sigma]) \cong \Xi_{\sigma}(\Phi) \to C(\Phi)$ (thus $Y_{\mathcal{Z}([\Phi, \sigma])} \to Y_{C(\Phi)}^{\natural}$) is faithfully flat, thus there exists open and closed subscheme $\wdt{Y}^{\natural}_{i, C(\Phi)} \subset \wdt{Y}^{\natural}_{i, C(\Phi)}$ whose pullbacks are $\wdt{Y}_{i, \mathcal{Z}([\Phi, \sigma])} \cong \wdt{Y}^{\natural}_{i, \Xi_{\sigma}}$ (under the isomorphism $\wdt{Y}_{\mathcal{Z}([\Phi, \sigma])} \cong \wdt{Y}^{\natural}_{\Xi_{\sigma}}$) along $Y_{\mathcal{Z}([\Phi, \sigma])} \cong Y^{\natural}_{\Xi_{\sigma}} \to Y^{\natural}_{C}$. $\wdt{Y}^{\natural}_{i, C(\Phi)}$ is independent of choice of $\sigma \in \Sigma(\Phi)^+$, this can be checked over $W$, note that points in $\wdt{Y}_{i, \mathcal{Z}([\Phi, \sigma])}$ come from specializations of points in $\wdt{Y}_{i}$. Since $\wdt{Y}_i \subset \wdt{Y}$ is open and closed, thus $\wdt{Y}_i \to Y$ is affine and flat, we are done.
\end{proof}
\begin{lemma}\label{lemma: simple topo lemma, etale covering}
    Let $\varphi: \ovl{U} \to \ovl{V}$ be an open surjection, $V \subset \ovl{V}$ be an open dense subscheme such that $V = V_1 \sqcup V_2$ is a topologically disjoint union. Let $U_i = \varphi^{-1}(V_i)$, $U = \varphi^{-1}(V)$. Let $\ovl{U}_i$ be the closure of $U_i$ in $\ovl{U}$, $\ovl{V}_i$ be the closure of $V_i$ in $\ovl{V}$, then $\varphi(\ovl{U}_i) = \ovl{V}_i$.
\end{lemma}
\begin{proof}
    On one hand, since $\varphi(\ovl{U}\setminus\ovl{U}_i)$ is open and does not intersect with $V_i$, thus $\varphi(\ovl{U}\setminus\ovl{U}_i)^c \supset \ovl{V}_i$. Note that $U$ is open dense in $\ovl{U}$, otherwise the image of the complementary open subscheme of the closure of $U$ is open in $\ovl{V}$ and does not intersect with the open dense subscheme $V$, which is absurd. Use the tautological inclusion $\varphi({\ovl{U}_i}) \supset \varphi(\ovl{U}\setminus\ovl{U}_i)^c$, we have $\varphi(\ovl{U}_i) \supset \ovl{V}_i$. On the other hand, $\varphi^{-1}(\ovl{V}_i) \supset \ovl{U}_i$, $\ovl{V}_i = \varphi(\varphi^{-1}(\ovl{V}_i)) \supset \varphi(\ovl{U}_i)$. Therefore, $\varphi(\ovl{U}_i) = \ovl{V}_i$.
\end{proof}

\begin{definition}\label{def: minimal compactification of well-positioned, second}
      Let $Y \subset \Shum{K, T}$ be a well-positioned subscheme with respect to $\lrbracket{Y^{\natural}(\Phi)}_{[\Phi]}$,  and $\wdt{Y} \to Y$ be well-positioned with respect to $\lrbracket{\wdt{Y}^{\natural}_{C(\Phi)} \to Y^{\natural}_{C(\Phi)}}_{[\Phi]}$. We define the so-called minimal compactification $\wdt{Y}^{\min}$ of $\wdt{Y}$ as the relative normalization of $Y^{\min}$ in $\wdt{Y}^{\tor}_{\Sigma}$ along $\wdt{Y}^{\tor}_{\Sigma} \to Y^{\tor}_{\Sigma} \to Y^{\min}$. 
\end{definition}
\begin{lemma}\label{lemma: minimal compactifications disjoint}
     Keep notations from \ref{def: minimal compactification of well-positioned, second}. Let $\wdt{Y} = \wdt{Y}_1 \sqcup \wdt{Y}_2$ be a topologically disjoint union, then $\wdt{Y}^{\tor}_{\Sigma} = \wdt{Y}^{\tor}_{1, \Sigma} \sqcup \wdt{Y}^{\tor}_{2, \Sigma}$ and $\wdt{Y}^{\min} = \wdt{Y}^{\min}_1 \sqcup \wdt{Y}^{\min}_2$ are topologically disjoint unions.
\end{lemma}
\begin{proof}
    The result for toroidal compactifications comes from Lemma \ref{lemma: open closed subschemes are well positioned 3}, and the result for minimal compactifications comes from \cite[\href{https://stacks.math.columbia.edu/tag/03GO}{Tag 03GO}]{stacks-project} ($\wdt{Y}^{\tor}_{\Sigma} \to Y^{\tor}_{\Sigma}$ is qcqs, $Y^{\tor}_{\Sigma} \to Y^{\min}$ is proper, thus $\wdt{Y}^{\tor}_{\Sigma} \to Y^{\min}$ is qcqs).
\end{proof}

\begin{lemma}\label{lemma: limit of well-positioned}
     Let $Y \subset \Shum{K, T}$ be a well-positioned subscheme with respect to $\lrbracket{Y^{\natural}(\Phi)}_{[\Phi]}$,  and $\wdt{Y}_i \to Y$ be well-positioned with respect to $\lrbracket{\wdt{Y}^{\natural}_{i, C(\Phi)} \to Y^{\natural}_{C(\Phi)}}_{[\Phi]}$ for $i \in I$. Assume for each $i \to j \in I$ and each $\Phi$, $\wdt{Y}_i^{\tor} \to \wdt{Y}_j^{\tor}$, $\wdt{Y}^{\natural}_{i, C(\Phi)} \to Y^{\natural}_{j, C(\Phi)}$ is affine, then $\wdt{Y} \to Y$ be well-positioned with respect to $\lrbracket{\wdt{Y}^{\natural}_{C(\Phi)} \to Y^{\natural}_{C(\Phi)}}_{[\Phi]}$, where $\wdt{Y}  = \prolim_{i \in I} \wdt{Y}_i $, $\wdt{Y}^{\natural}_{C(\Phi)} = \prolim_{i \in I} \wdt{Y}^{\natural}_{i, C(\Phi)}$.
\end{lemma}
\begin{proof}
    The affiness of transition maps ensure the existence of $\wdt{Y}^{\tor}$ and $\wdt{Y}^{\natural}_{C(\Phi)}$ that are affine over $Y^{\tor}$ and $Y^{\natural}_{C(\Phi)}$ respectively. Since all $\wdt{Y}_i^{\tor} \to Y^{\tor}$, $\wdt{Y}^{\natural}_{i, C(\Phi)} \to Y^{\natural}_{C(\Phi)}$ are flat, then $\wdt{Y}^{\tor}$ and $\wdt{Y}^{\natural}_{C(\Phi)}$ are flat over the corresponding base. Then this lemma follows from the fact that taking inverse limits commutes with taking products.
\end{proof}

   \begin{lemma}\label{lemma: reverse, well-position}
       Let $Y \subset \Shum{K, T}$ be a well-positioned subscheme with respect to $\lrbracket{Y^{\natural}(\Phi)}_{[\Phi]}$, and let $\wdt{Y}^{\tor}_{\Sigma} \to Y^{\tor}_{\Sigma}$ be an affine flat morphism and $\lrbracket{\wdt{Y}^{\natural}_{C(\Phi)} \to Y^{\natural}_{C(\Phi)}}_{[\Phi]}$ be a collection of affine flat morphisms. If for each $[\Phi, \sigma]$, there is an isomorphim $\wdt{Y}_{\Xi_{\sigma}}^{\natural} \cong \wdt{Y}_{\mathcal{Z}([\Phi, \sigma])}$ (compatible with  $Y_{\Xi_{\sigma}}^{\natural} \cong Y_{\mathcal{Z}([\Phi, \sigma])}$) which extends to an isomorphism $\wdt{Y}^{\natural}_{\mathfrak{X}_{\sigma}} \cong (\wdt{Y}^{\tor}_{\Sigma})^{\wedge}_{\wdt{Y}_{\mathcal{Z}([\Phi, \sigma])}}$ (compatible with $Y^{\natural}_{\mathfrak{X}_{\sigma}} \cong (Y^{\tor}_{\Sigma})^{\wedge}_{Y_{\mathcal{Z}([\Phi, \sigma])}}$), then $\wdt{Y} \to Y$ be well-positioned with respect to $\lrbracket{\wdt{Y}^{\natural}_{C(\Phi)} \to Y^{\natural}_{C(\Phi)}}_{[\Phi]}$.
   \end{lemma}
   \begin{proof}

       By standard limit argument, we might assume both $\wdt{Y}^{\tor}_{\Sigma} \to Y^{\tor}_{\Sigma}$ and $\wdt{Y}^{\natural}_{C} \to Y^{\natural}_{C}$ are locally of finitely presented. It follows from the condition that for affine open coverings $\mathfrak{W}$ on $Y_{\mathfrak{X}_{\sigma}}^{\natural} \cong (Y^{\tor}_{\Sigma})^{\wedge}_{\mathcal{Z}([\Phi, \sigma])}$, the pullbacks of $\wdt{Y}^{\tor}_{\Sigma}$ and $\wdt{Y}^{\natural}_{C}$ coincides over $W$.
       
       Note that there are only finitely many strata in the union $\Xi(\sigma)^+ = \underset{\tau\in \Sigma(\Phi)^{+},\ \bar{\tau}\subset\bar{\sigma}}{\cup} \Xi_{\tau}$. Consider the faithfully flat morphism:
       \[  \underset{\tau\in \Sigma(\Phi)^{+},\ \bar{\tau}\subset\bar{\sigma}}{\bigsqcup} Y^{\natural}_{\mathfrak{X}_{\tau}} := \underset{\tau\in \Sigma(\Phi)^{+},\ \bar{\tau}\subset\bar{\sigma}}{\bigsqcup} (Y^{\natural}_{\Xi(\sigma)})^{\wedge}_{Y^{\natural}_{\Xi_{\tau}}} \to (Y^{\natural}_{\Xi(\sigma)})^{\wedge}_{Y^{\natural}_{\Xi(\sigma)^+}} =: 
       Y^{\natural}_{\mathfrak{X}_{\sigma}^{\circ}} \]
       Let $\mathfrak{W}$ be a affine open covering over $Y^{\natural}_{\mathfrak{X}_{\sigma}^{\circ}}$, its pullback can be refined to an affine open covering $\bigsqcup \mathfrak{W}_{\tau}$ on $\bigsqcup Y^{\natural}_{\mathfrak{X}_{\tau}}$. Then the result follows from fpqc descent along $\bigsqcup W_{\tau} \to W$.
   \end{proof}
\subsubsection{}
\begin{definition-proposition}{\cite[Definition 4.6, Proposition 4.10]{lan2022closed}}
    Let 
    \[ \Sigma_1 = \lrbracket{\Sigma_{1}(\Phi_1)}_{[\Phi_1] \in \Cusp_{K_1}(G_1, X_1)}, \quad \Sigma_2 = \lrbracket{\Sigma_{2}(\Phi_2)}_{[\Phi_2] \in \Cusp_{K_2}(G_2, X_2)} \] be cone decompositiones. We say $\Sigma_1$ and $\Sigma_2$ are strictly compatible if the image of each $\sigma_1 \in \Sigma_1(\Phi_1)^+$ under $\mathbf{B}_{K_1}(\Phi_1)_{\R} \hookrightarrow \mathbf{B}_{K_2}(\Phi_2)_{\R}$ is exactly some $\Sigma_2 \in \Sigma(\Phi_2)^+$. There always exists such strictly compatible $(\Sigma_1, \Sigma_2)$ such that $\Sigma_1, \Sigma_2$ are smooth and projective.
\end{definition-proposition}
\begin{remark}
    The proof in \cite[Proposition 4.10]{lan2022closed} does not require the full strength of \cite[Assumption 2.1]{lan2022closed}, the weaker assumption in \cite[Proposition 3.4]{lan2022closed} is enough.
\end{remark}
\begin{proposition}\label{prop: tor closed embedding}
    Let $(\Sigma_1, \Sigma_2)$ be strictly compatible. Assume \ref{assumption: C to Z} hold for both $\Shum{K_1}$ and $\Shum{K_2}$. Let $Y_2 \subset \Shum{K_2, T}$ be a well-positioned subscheme with respect to $\lrbracket{Y^{\natural}_{2}(\Phi_2)}_{[\Phi_2]}$, and $\wdt{Y}_2 \to Y_2$ be well-positioned with respect to $\lrbracket{\wdt{Y}^{\natural}_{2, C_2(\Phi_2)} \to Y^{\natural}_{2, C_2(\Phi_2)}}_{[\Phi_2]}$. Let $Y_1 \subset \Shum{K_1, T}$ be an open and closed subscheme of the pullback of $Y_2$ along $\Shum{K_1, T} \to \Shum{K_2, T}$, and $\wdt{Y}_1$ be an open and closed subscheme of the pullback of $\wdt{Y}_2 \to Y_2$ along $Y_1 \to Y_2$.

    If the finite morphism $\wdt{Y}_1 \to \wdt{Y}_2$ is indeed a closed embedding, then $\wdt{Y}^{\tor}_{1, \Sigma_1} \to \wdt{Y}^{\tor}_{2, \Sigma_2}$ and $\wdt{Y}_{1, C_1(\Phi_1)}^{\natural} \to \wdt{Y}_{2, C_2(\Phi_2)}^{\natural}$ are closed embedding. In particular, if the finite morphism $Y_1 \to Y_2$ is indeed a closed embedding, then $Y^{\tor}_{1, \Sigma_1} \to Y^{\tor}_{2, \Sigma_2}$, $\wdt{Y}_{1, C_1(\Phi_1)}^{\natural} \to \wdt{Y}_{2, C_2(\Phi_2)}^{\natural}$ are closed embeddings.
\end{proposition}
\begin{proof}
    The proof is similar to the one of \cite[Proposition 4.9]{lan2022closed}. Due to Lemma \ref{lemma: analog covering result for well-position}, \ref{lemma: pullback of well-positioned is well-positioned, second}, \ref{lemma: open closed subschemes are well positioned, 2}, \ref{lemma: open closed subschemes are well positioned 3}, and Proposition \ref{prop: wdt Y satisfy same prop}, given any point $x \in \wdt{Y}^{\natural}_{1, \mathcal{Z}_1([\Phi_1, \sigma_1])} \cong \wdt{Y}^{\natural}_{1, \Xi_{1, \sigma_1}}$, there exist \'etale morphisms $\ovl{U} \to \wdt{Y}^{\tor}_{2, \Sigma_2}$, $\ovl{U} \to E_2(\sigma_2) \times_{\Spec \Z} \wdt{Y}^{\natural}_{2, C_2} \subset \wdt{Y}^{\natural}_{\Xi_2(\sigma_2)}$ with respect to $x$ such that the left diagram is identified with the right diagram over $\ovl{U}$:
\[\begin{tikzcd}
	{\wdt{Y}_1} & {\wdt{Y}_{1, \Sigma_1}^{\tor}} & {E_1\times_{\Spec \Z}\wdt{Y}^{\natural}_{1, C_1}} & {E_1(\sigma_1)\times_{\Spec \Z}\wdt{Y}^{\natural}_{1, C_1}} \\
	{\wdt{Y}_2} & {\wdt{Y}_{2, \Sigma_2}^{\tor}} & {E_2\times_{\Spec \Z}\wdt{Y}^{\natural}_{2, C_2}} & {E_2(\sigma_2)\times_{\Spec \Z}\wdt{Y}^{\natural}_{2, C_2}}
	\arrow[hook, from=1-1, to=1-2]
	\arrow[hook, from=1-1, to=2-1]
	\arrow[from=1-2, to=2-2]
	\arrow[hook, from=1-3, to=1-4]
	\arrow[from=1-3, to=2-3]
	\arrow[from=1-4, to=2-4]
	\arrow[hook, from=2-1, to=2-2]
	\arrow[hook, from=2-3, to=2-4]
\end{tikzcd}\]

    Note that $E_1\times_{\Spec \Z}\wdt{Y}^{\natural}_{1, C_1} \to E_2\times_{\Spec \Z}\wdt{Y}^{\natural}_{2, C_2}$ is a closed embedding over $U \subset \ovl{U}$ by assumption. Due to \cite[Lemma 4.3]{lan2022closed}, when $(\Sigma_1, \Sigma_2)$ are strictly compatible, then $E_1(\sigma_1) \to E_2(\sigma_2)$ is a closed embedding, in particular, $E_1(\sigma_1)\times_{\Spec \Z}\wdt{Y}^{\natural}_{1, C_1} \to E_2(\sigma_2)\times_{\Spec \Z}\wdt{Y}^{\natural}_{2, C_2}$ is also a closed embedding over $\ovl{U}$. By varing $x$ and by \'etale descent, $\wdt{Y}^{\tor}_{1, \Sigma_1} \to \wdt{Y}^{\tor}_{2, \Sigma_2}$ is a closed embedding.

    Notice that for each $[\Phi]$, and $\sigma$ of top dimensional in $\Sigma(\Phi)$, we have $Y_{\mathcal{Z}([\Phi, \sigma])} \cong Y^{\natural}_C$ canonically, thus $\wdt{Y}_{\mathcal{Z}([\Phi, \sigma])} \cong \wdt{Y}^{\natural}_C$ canonically. We can take $\sigma_2 \in \Sigma_2(\Phi_2)$ being of top dimensional, then $\sigma_1 \in \Sigma_1(\Phi_1)$ will be of top-dimensional, thus $\wdt{Y}^{\tor}_{1, \Sigma_1} \to \wdt{Y}^{\tor}_{2, \Sigma_2}$ is a closed embedding implies that $\wdt{Y}_{\mathcal{Z}_1([\Phi_1, \sigma_1])} \to \wdt{Y}_{\mathcal{Z}_2([\Phi_2, \sigma_2])}$ is a locally closed embedding, thus the finite morphism $\wdt{Y}_{1, C_1(\Phi_1)}^{\natural} \to \wdt{Y}_{2, C_2(\Phi_2)}^{\natural}$ is a closed embedding.
\end{proof}

\subsubsection{}
Also, let us recall a result from \cite[\S 2.7]{pink1992}. Consider an inverse system of schemes:
\[\begin{tikzcd}
	{X_m} & {X_n} & X \\
	{\ovl{X}_m} & {\ovl{X}_n} & {\ovl{X}} \\
	{S_m} & {S_n} & S
	\arrow["{\varphi_{nm}}", from=1-1, to=1-2]
	\arrow["{\varphi_n}", from=1-2, to=1-3]
	\arrow["{\bar{\varphi}_{nm}}", from=2-1, to=2-2]
	\arrow["{\bar{\varphi}_n}", from=2-2, to=2-3]
	\arrow["{\omega_{nm}}", from=3-1, to=3-2]
	\arrow["{\omega_n}", from=3-2, to=3-3]
	\arrow["j"', hook, from=1-3, to=2-3]
	\arrow["{\bar{f}}"', from=2-3, to=3-3]
	\arrow["{j_n}"', hook, from=1-2, to=2-2]
	\arrow["{\bar{f}_n}"', from=2-2, to=3-2]
	\arrow["{j_m}"', hook, from=1-1, to=2-1]
	\arrow["{\bar{f}_m}"', from=2-1, to=3-1]
\end{tikzcd}\]
and assume
\begin{enumerate}
    \item All horizontal morphisms are finite and flat.
    \item A prime $l$ invertible over $S$.
    \item Every $\bar{f}_n$ is smooth and $j_n$ is an open embedding where the complement is a divisor with normal crossings over $S_n$.
    \item For every $n$ and every irreducible component $D_n$ of $\ovl{X}_n\setminus X_n$, there exists $m > n$ such that the ramification degree of $\bar{\varphi}_{nm}$ in $D_n$ is divisible by $l$,
\end{enumerate}
then
\begin{equation}\label{eq: pink's isomorphism}
    \varinjlim \varphi_{n*} \F_l \to \varinjlim \varphi_{n*} Rj_{n*} \F_l
\end{equation}
is an isomorphism.

We apply this to the setting of well-positioned subschemes. Fix a prime $l$ different from $p$. Fix a level group $K = K(l^0) = K_pK(l^0)K^{p, l}$, $K(l^n) = K_pK(l^n)K^{p, l}$ where $K(l^n) \subset G(\Q_l)$ be the principal level groups when $G$ is hyperspecial at level $l$, or more general we can take Moy-Prasad level group of level $n$ at prime $l$, $K^{p, l}$ be the level group away from primes $p, l$. Let $Y_K \subset \Shum{K, T}$ be a well-positioned subscheme with respect to $\lrbracket{Y^{\natural}_{K}(\Phi)}_{[\Phi]}$, and $\wdt{Y}_K \to Y_K$ be well-positioned with respect to $\lrbracket{\wdt{Y}^{\natural}_{C_{K}(\Phi)} \to Y^{\natural}_{C_{K}(\Phi)}}_{[\Phi]}$. Let $Y_{K(l^n)}$ (resp. $\wdt{Y}^{\tor}_{K(l^n), \Sigma}$) be the pullback of $Y_K$ (resp. $\wdt{Y}^{\tor}_{K, \Sigma}$) Consider the following two diagrams:
\begin{equation}\label{eq: pink formula diagram}
\begin{tikzcd}
	{Y_{K(l^m)}} & {Y_{K(l^n)}} & {Y_{K}} & {\wdt{Y}_{K(l^m)}} & {\wdt{Y}_{K(l^n)}} & {\wdt{Y}_{K}} \\
	{Y_{K(l^m), \Sigma}^{\tor}} & {Y_{K(l^n), \Sigma}^{\tor}} & {Y_{K, \Sigma}^{\tor}} & {\wdt{Y}_{K(l^m), \Sigma}^{\tor}} & {\wdt{Y}_{K(l^n), \Sigma}^{\tor}} & {\wdt{Y}_{K, \Sigma}^{\tor}}
	\arrow["{\varphi_{nm}}", from=1-1, to=1-2]
	\arrow["{j_m}"', hook, from=1-1, to=2-1]
	\arrow["{\varphi_n}", from=1-2, to=1-3]
	\arrow["{j_n}"', hook, from=1-2, to=2-2]
	\arrow["j"', hook, from=1-3, to=2-3]
	\arrow["{\tilde{\varphi}_{nm}}", from=1-4, to=1-5]
	\arrow["{\tilde{j}_{m}}"', hook, from=1-4, to=2-4]
	\arrow["{\tilde{\varphi}_n}", from=1-5, to=1-6]
	\arrow["{\tilde{j}_{n}}"', hook, from=1-5, to=2-5]
	\arrow["{\tilde{j}}"', from=1-6, to=2-6]
	\arrow["{\bar{\varphi}_{nm}}", from=2-1, to=2-2]
	\arrow["{\bar{\varphi}_{n}}", from=2-2, to=2-3]
	\arrow["{\tilde{\bar{\varphi}}_{nm}}", from=2-4, to=2-5]
	\arrow["{\tilde{\bar{\varphi}}_{n}}", from=2-5, to=2-6]
\end{tikzcd}
\end{equation}

\begin{proposition}\label{prop: pink's formula}\leavevmode
     In case $(4)'$ (Assumption \ref{assumption: adjusted Siegel embedding}), we have:
    \begin{enumerate}
        \item $\varinjlim \varphi_{n*} \F_l \rightiso \varinjlim \varphi_{n*} Rj_{n*} \F_l$
        \item $\varinjlim \tilde{\varphi}_{n*} \F_l \rightiso \varinjlim \tilde{\varphi}_{n*} R\tilde{j}_{n*} \F_l$
    \end{enumerate}
\end{proposition}
\begin{proof}
    It suffices to work \'etale locally at boundary. Let $x\in Y_{\mathcal{Z}_K([\Phi, \sigma])}$, and we lift $x$ to a system of points $x \in Y_{\mathcal{Z}_{K(l^n)}([\Phi, \sigma])}$ with the a system of cusp labels (which we still denote by $[\Phi, \sigma]$). \'Etale locally, the above two diagrams \ref{eq: pink formula diagram} (here we only draw out the one associated with $\wdt{Y}$) can be upgraded to
\[\begin{tikzcd}
	{E_{K(l^m)} \times \wdt{Y}^{\natural}_{C_{K(l^m)}}} & {E_{K(l^n)} \times \wdt{Y}^{\natural}_{C_{K(l^n)}}} & {E_{K} \times \wdt{Y}^{\natural}_{C_{K}}} \\
	{E_{K(l^m)}(\sigma) \times \wdt{Y}^{\natural}_{C_{K(l^m)}}} & {E_{K(l^n)}(\sigma) \times \wdt{Y}^{\natural}_{C_{K(l^n)}}} & {E_{K}(\sigma) \times \wdt{Y}^{\natural}_{C_{K}}} & {} \\
	{\wdt{Y}^{\natural}_{C_{K(l^m)}}} & {\wdt{Y}^{\natural}_{C_{K(l^n)}}} & {\wdt{Y}^{\natural}_{C_{K}}}
	\arrow["{\tilde{\varphi}_{nm}}", from=1-1, to=1-2]
	\arrow["{\tilde{j}_m}"', from=1-1, to=2-1]
	\arrow["{\tilde{\varphi}_n}", from=1-2, to=1-3]
	\arrow["{\tilde{j}_n}"', from=1-2, to=2-2]
	\arrow["{\tilde{j}}"', from=1-3, to=2-3]
	\arrow["{\tilde{\bar{\varphi}}_{nm}}", from=2-1, to=2-2]
	\arrow["{\tilde{\bar{f}}_m}"', from=2-1, to=3-1]
	\arrow["{\tilde{\bar{\varphi}}_{n}}", from=2-2, to=2-3]
	\arrow["{\tilde{\bar{f}}_n}"', from=2-2, to=3-2]
	\arrow["{\tilde{\bar{f}}}"', from=2-3, to=3-3]
	\arrow["{\tilde{\omega}_{nm}}", from=3-1, to=3-2]
	\arrow["{\tilde{\omega}_{n}}", from=3-2, to=3-3]
\end{tikzcd}\]

    Note that in case $(4)$, by the construction in \cite[\S 2, \S 4]{pera2019toroidal}, $E_{K(l^m)} \to E_{K(l^n)} \to E_K$ (resp. $C_{K(l^m)} \to C_{K(l^n)} \to C_K$) are relative normalizations of $E_{K_2}$ (resp. $C_{K_2}$, here we omit the subscript $(\ast)_1$) in $E_{K(l^m), \eta} \to E_{K(l^n), \eta} \to E_{K, \eta}$ (resp. $C_{K(l^m), \eta} \to C_{K(l^n), \eta} \to C_{K, \eta}$). In particular, $E_{K(l^m)} \to E_{K(l^n)} \to E_K$ (resp. $C_{K(l^m)} \to C_{K(l^n)} \to C_K$) are finite and surjective due to \cite[Lemma 3.44]{mao2025boundary}. By miracle flatness, these finite surjections between normal noetherian schemes are flat morphisms. Also,
    \[ E_{K(l^m)} \times C_{K(l^m)} \to (E_{K(l^n)} \times C_{K(l^n)}) \times_{C_{K(l^n)}} C_{K(l^m)} = E_{K(l^n)} \times C_{K(l^m)} \]
    is the multiplication of $l^{m-n}$ over $C_{K(l^m)}$.  We then apply the isomorphism \ref{eq: pink's isomorphism} to the above diagram.
\end{proof}

\subsection{General cases}\label{subsec: general cases}

Let $\iota: (G_1, X_1, K_1) \to (G_2, X_2, K_2)$ be a map of Shimura data in Assumption \ref{assumption: general set up}. Let $\varphi: \Shum{K_1} \to \Shum{K_2}$ be the induced morphism (we omit the $\otimes E_2$ part on $\Shum{K_2}$ side).

\begin{definition}
    Let $S$ be a scheme, we say locally closed subschemes $\lrbracket{Y_i}_{i \in I}$ in $S$ form a (weak) stratification on $S$ if $S = \sqcup_{i \in I} Y_i$. In this case, we call $\lrbracket{Y_i}_{i \in I}$ strata in $S$.
\end{definition}

\subsubsection{Basic assumption}

We first deal with the axiomatic setting for central leaves, Newton strata and Ekedahl-Oort strata.

\begin{assumptionaxiom}\label{assumption: fiber discreteness}
Let $\lrbracket{Y_i}_{i \in I}$, $\lrbracket{Y_j}_{j \in J}$ be weak stratifications of $\Shum{K_1, T}$, $\Shum{K_2, T}$ respectively, such that:
    \begin{itemize}
        \item Assume \ref{assumption: C to Z} hold for both $\Shum{K_1}$ and $\Shum{K_2}$.
        \item Given any type $Y_i \subset \Shum{K_1, T}$, there exists a unique $Y_j \subset \Shum{K_2, T}$ such that $\varphi(Y_i) \subset Y_j$, and moreover $Y_i \subset \varphi^{-1}(Y_j)$ is open and closed in $\varphi^{-1}(Y_j)$. We call this \emph{discrete fiber property} for strata $\lrbracket{Y_i}_{i \in I}$, $\lrbracket{Y_j}_{j \in J}$.
        \item $Y_j \subset \Shum{K_2, T}$ are well-positioned for all $j \in J$.
    \end{itemize}
\end{assumptionaxiom}

\begin{proposition}\label{prop: general well-position}
    Assume \ref{assumption: fiber discreteness} holds, then $Y_i \subset \Shum{K, T}$ are well-positioned subsets (resp. subschemes) for all $i \in I$. Moreover, let $Y \subset \Shum{K_1, T}$ be a locally closed subsets that is a union of some connected components of some collection $\lrbracket{Y_i}_{i \in I_Y \subset I}$, then $Y$ is a well-positioned subset, and is a well-positioned subscheme with respect to induced reduced subscheme structure. 
\end{proposition}
\begin{proof}\leavevmode
    \begin{itemize}
     \item Pullbacks of well-positioned subsets (resp. subschemes) under $\varphi$ are again well-positioned subsets (resp. subschemes), see Lemma \ref{lemma: pullback of well-positioned is well-positioned}.
    \item Open and closed subsets (resp. subschemes) of well-positioned subsets (resp. subschemes) are again well-positioned subsets (resp. subschemes), see Proposition \ref{proposition: open-closed subschemes are well positioned}.
     \item Locally closed subsets which are unions of well-positioned subsets are again well-positioned subsets, see Lemma \ref{lem: union of well-positioned sets}.
     \item Well-positioned subsets are well-positioned subschemes with induced reduced subscheme structures, see Remark \ref{remark: subset well positioned implies subscheme well positioned}.
 \end{itemize}
\end{proof}

\begin{proposition}\label{prop: affineness in general}
     Assume \ref{assumption: fiber discreteness} holds. Assume the partial minimal compactifications $Y_j^{\min} \subset \Shumm{K_2, T}$ are affine for all $j \in J$, then the partial minimal compactifications $Y_i^{\min} \subset \Shumm{K_1, T}$ are affine for all $i \in I$.
\end{proposition}
\begin{proof}
    This follows from Lemma \ref{lemma: pullback of well-positioned is well-positioned} and Proposition \ref{proposition: open-closed subschemes are well positioned}, note that the finite morphism $\varphi^{\min}$ pulls back an affine scheme to an affine scheme, and an open and closed subscheme of an affine scheme is again affine.
\end{proof}

\begin{proposition}\label{prop: general well-position, disjoint}
    Assume \ref{assumption: fiber discreteness} holds. We furthur assume that any two $Y_{j_1}, Y_{j_2} \subset \Shum{K_2, T}$ has no intersection at boundary, i.e. $Y_{j_1}^{\min}\cap Y_{j_2}^{\min} = \emptyset$, $Y_{j_1, \Sigma_2}^{\tor}\cap Y_{j_2, \Sigma_2}^{\tor} = \emptyset$
    
    Let $\lrbracket{Z_p}_{p \in P}$ be a weak stratification of $\Shum{K_1, T}$ such that each $Z_p$ is a union (not necessarily a finite union) of some collection $\lrbracket{Y_i}_{i \in I_p \subset I}$ (for example, take $\lrbracket{Z_p}_{p \in P} = \lrbracket{Y_i}_{i \in I}$). Then 
    \[ \Shumm{K_1, T} = \bigsqcup_{p \in P} Z_p^{\min},\quad  \Shumc{K_1, \Sigma_1, T} = \bigsqcup_{p \in P} Z^{\tor}_{p, \Sigma_1}. \]
    In particular, if $p \neq q$, then $Z^{\tor}_{p, \Sigma_1}$ and $Z^{\tor}_{q, \Sigma_1}$ do not intersect, $Z^{\min}_p$ and $Z^{\min}_q$ do not intersect.
    
    If moreover the closure of any $Z_p$ is a union of some collection $\lrbracket{Y_i}_{i \in I_{\bar{p}} \subset I}$, then for any $p, q \in P$, we have 
    \begin{equation}
              \ovl{Z_p} \cap Z_q \neq \emptyset \Longleftrightarrow \ovl{Z_p^{\min}} \cap Z_q^{\min} \neq \emptyset \Longleftrightarrow \ovl{Z^{\tor}_{p, \Sigma_1}} \cap Z^{\tor}_{q, \Sigma_1} \neq \emptyset.
    \end{equation}
\end{proposition}
\begin{proof}
    
    The assumption here together with Lemma \ref{lemma: pullback of well-positioned is well-positioned}, Proposition \ref{proposition: open-closed subschemes are well positioned} imply that each $Y_i$ ($i\in I$) is well-positioned, and 
    \[ \Shumm{K_1, T} = \bigsqcup_{i \in I} Y_i^{\min}, \quad \Shumc{K_1, \Sigma_1, T} = \bigsqcup_{i \in I} Y_{i, \Sigma_1}^{\tor}. \]
    Due to Lemma \ref{lem: union of well-positioned sets}, we have
      \[ \Shumm{K_1, T} = \bigsqcup_{p \in P} Z_p^{\min},\quad  \Shumc{K_1, \Sigma_1, T} = \bigsqcup_{p \in P} Z^{\tor}_{p, \Sigma_1}. \]
    
  We consider the second part of the statement. It suffices to show $\ovl{Z_p^{\min}} \cap Z_q^{\min}  \neq \emptyset \Longrightarrow \ovl{Z_p} \cap Z_q \neq \emptyset$. The case of toroidal compactifications follows similarly. Due to Lemma \ref{lemma: closure of well-pos is well-pos} and \ref{lem: union of well-positioned sets}, $\ovl{Z_p}$ is well-positioned and $\ovl{Z_p^{\min}} = \ovl{Z_p}^{\min}$. By assumption (resp. by definition), $\ovl{Z_p} = \sqcup_{i \in I_{\bar{p}}} Y_i$ (resp. $Z_q = \sqcup_{j \in I_q} Y_j$). Use Lemma \ref{lem: union of well-positioned sets} $\ovl{Z_p^{\min}} = \ovl{Z_p}^{\min} = \sqcup_{i \in I_{\bar{p}}} Y_i^{\min}$ (resp. $Z_q^{\min} = \sqcup_{j \in I_q} Y_j^{\min}$). If $\ovl{Z_p^{\min}} \cap Z_q^{\min} \neq \emptyset$, then there exists $Y_i \subset \ovl{Z_p}$ and $Y_j \subset Z_q$ such that $Y_i^{\min} \cap Y_j^{\min}$ is nontrivial. Then $Y_i \cap Y_j$ is nontrivial due to the first part of this lemma. Therefore, $\ovl{Z_p} \cap Z_q \neq \emptyset$.
\end{proof}

\subsubsection{Enhanced assumption}

We moreover assume we are in case $(4)'$ (Assumption \ref{assumption: adjusted Siegel embedding}). Let 
\[\iota_{[\Phi]}: (G_{1, h}, X_{1, \Phi_1, h}, K_{1, \Phi_1, h}) \to (G_{2, h}, X_{2, \Phi_2, h}, K_{2, \Phi_2, h}),\]
\[\varphi_{[\Phi]}: \Shum{K_{1, \Phi_1, h}}(G_{1, h}, X_{1, \Phi_1, h}) \to \Shum{K_{2, \Phi_2, h}}(G_{2, h}, X_{2, \Phi_2, h})\] 
be the morphisms at boundary introduced in Proposition \ref{proposition: functorialities on toroidal compactifications}. Then due to \cite[Theorem 3.58]{mao2025boundary}, $\iota_{[\Phi]}$ are again in case $(4)'$.

\begin{assumptionaxiom}\label{assumption: fiber discreteness, enhanced}
    Let $\lrbracket{Y_i}_{i \in I}$, $\lrbracket{Y_j}_{j \in J}$ be weak stratifications of $\Shum{K_1, T}$, $\Shum{K_2, T}$ respectively that satisfy the assumption \ref{assumption: fiber discreteness} under $\varphi$. For each $[\Phi_1] \in \Cusp_{K_1}(G_1, X_1)$, $\Phi_2 = \iota_*\Phi_1$, let $\lrbracket{Z_i}_{i \in I_{\Phi_1}}$, $\lrbracket{Z_j}_{j \in J_{\Phi_2}}$ be weak stratifications of $\Shum{K_{1, \Phi_1, h}, T}$, $\Shum{K_{2, \Phi_2, h}, T}$ respectively that satisfy the assumption \ref{assumption: fiber discreteness} under $\varphi_{[\Phi]}$, such that
    \begin{enumerate}
        \item For each $Y_j \subset \Shum{K_2, T}$ ($j \in J$), its boundary $Y^{\natural}_{j, \Zb^{\bigsur}_2(\Phi_2)} \subset \Zb^{\bigsur}_2(\Phi_2)_T = \Shum{K_{2, \Phi_2, h}, T}$ is $Z_{l_j}$ for some $l_j \in J_{\Phi_2}$.
        \item Any two $Y_{j_1}, Y_{j_2} \subset \Shum{K_2, T}$ has no intersection at boundary. In particular, $Y_{j_1}$, $Y_{j_2}$ correspond to different $Z_{l_{j_1}}$, $Z_{l_{j_2}}$.
    \end{enumerate}
\end{assumptionaxiom}

   \begin{proposition}\label{prop: boundaries of fiber discrete strata}
         Assume \ref{assumption: fiber discreteness, enhanced} holds. Let $Y \subset \Shum{K_1, T}$ be open and closed in some $Y_i \subset \Shum{K_1, T}$ ($i \in I$), then its boundary $Y_{1, \Zb^{\bigsur}_1}^{\natural} \subset \Zb^{\bigsur}_{1, T}$ is a union of open closed subschemes of some $\lrbracket{Z_i}_{i \in I_{\Phi_1, l_i}}$, where $Z_i \subset \Zb^{\bigsur}_{1, T} = \Shum{K_{1, \Phi_1, h}, T}$. 
     \end{proposition}
    \begin{proof}
        Due to Proposition \ref{proposition: functorialities on toroidal compactifications}, both $\varphi: \Shum{K_1} \to \Shum{K_2}$ and $\varphi^{\min}: \Shumm{K_1} \to \Shumm{K_2}$ are finite. Let $Y_j \subset \Shum{K_2, T}$ be the unique stratum containing $\varphi(Y_i)$, then by assumption, $Y$ is open and closed in $\varphi^{-1}(Y_j)$.

        By assumption, the pullback $Z_{l, j} = Y^{\natural}_{j, \Zb^{\bigsur}_2(\Phi_2)} \subset \Zb^{\bigsur}_2(\Phi_2)_T$ along $\varphi^{\natural}_{\Phi}: \Zb^{\bigsur}_1(\Phi_1) \to \Zb^{\bigsur}_2(\Phi_2)$ is a topologically disjoint union of some $\lrbracket{Z_i}_{i \in I_{\Phi_1, l_i}}$ on $\Zb_1^{\bigsur}(\Phi_1)_{T}$. The proof of Proposition \ref{proposition: open-closed subschemes are well positioned} implies that, $Y^{\natural}_{\Zb^{\bigsur}_1(\Phi_1)}$ is an open and closed subscheme of $\varphi^{\natural, -1}_{[\Phi]}(Z_{l_j})$, thus is a union of some open and closed subschemes of $\lrbracket{Z_i}_{i \in I_{\Phi_1, l_i}}$.
     \end{proof}

\subsubsection{Igusa varieties}
     Let us consider the axiomatic setting for Igusa varieties
      \begin{assumptionaxiom}\label{assumption: well position over well position}
         Assume \ref{assumption: fiber discreteness} holds. Let $\lrbracket{\wdt{Y}_i \to Y_i}_{i \in I}$ and $\lrbracket{\wdt{Y}_j \to Y_j}_{j \in J}$ be some affine and flat morphisms such that
         \begin{enumerate}
             \item For each $i \in I$, let $j \in J$ be the unique element such that $\varphi(Y_i) \subset Y_j$, then $\wdt{Y}_i$ is an open and closed subscheme in the pullback of $\wdt{Y}_j \to Y_j$ along $Y_i \to Y_j$.
             \item $\lrbracket{\wdt{Y}_j \to Y_j}_{j \in J}$ are well-positioned in the sense of \ref{def: well-positioned over well-positioned} for all $j \in J$.
         \end{enumerate}
    
     \end{assumptionaxiom}

     \begin{proposition}\label{prop: well position over well position, and affine}
         Assume \ref{assumption: well position over well position} holds. Then $\wdt{Y}_i \to Y_i$ are well-positioned in the sense of \ref{def: well-positioned over well-positioned} for all $i \in I$. Moreover, if the partial minimal compactifications $Y_j^{\min} \subset \Shumm{K_2, T}$ are affine for all $j \in J$, then the partial minimal compactifications $\wdt{Y}_i^{\min}$ and $\wdt{Y}_j^{\min}$ are affine for all $i \in I$, $j \in J$.
     \end{proposition}
     \begin{proof}
         This follows from Proposition \ref{prop: general well-position} and Lemma \ref{lemma: pullback of well-positioned is well-positioned, second}, \ref{lemma: open closed subschemes are well positioned, 2} and \ref{lemma: open closed subschemes are well positioned 3} that $\wdt{Y}_i \to Y_i$ are well-positioned in the sense of \ref{def: well-positioned over well-positioned} for all $i \in I$.

        By definition, the partial minimal compactification $\wdt{Y}_i^{\min}$ of $\wdt{Y}_i$ is the relative normalization of $Y^{\min}_i$ in $\wdt{Y}^{\tor}_{i, \Sigma_1}$. Due to Proposition \ref{prop: affineness in general}, $Y^{\min}_i$ is affine, thus $\wdt{Y}^{\min}$ is affine, see \cite[\href{https://stacks.math.columbia.edu/tag/035H}{Tag 035H}]{stacks-project}.
     \end{proof}

\subsubsection{Hecke actions}
     Now we consider the Hecke action. To save notations, we omit the subscript $(\ast)_1$. Fix a $\Phi$. Consider $\Shum{K_{\Phi, h, p}} = \prolim_{K^{p}_{\Phi, h}} \Shum{K_{\Phi, h, p}K^{p}_{\Phi, h}}$, let $\pi_h: \Shum{K_{\Phi, h, p}} \to \Shum{K_{\Phi, h}} =: \Zb^{\bigsur}(\Phi)$ be the projection.

     \begin{proposition}\label{prop: Hecke disc, well positioned}
       Assume \ref{assumption: fiber discreteness, enhanced}. Let $Y_1, Y_2 \subset \Shum{K, T}$ be two such strata that are insensitive to the away-from-$p$ level $K^p$ in the sense of Definition \ref{definition: insensitive to away from p part lift}, assume there exists $y_i \in Y_{i, \Zb^{\bigsur}}^{\natural} \subset \Zb^{\bigsur}_{T}$ ($i = 1, 2$) for some boundary stratum $\Zb^{\bigsur}$ such that $y_1$ and $y_2$ are in $Z_i \subset \Zb^{\bigsur}_{T}$ for some $i \in I_{\Phi}$. If $G_h(\A_f^p)$ acts transitively on $\pi_0(\pi_h^{-1}(Z_i))$, then $Y_1 = Y_2$.
     \end{proposition}
     \begin{proof}
          Let $y_1', y_2' \in \pi_h^{-1}(Z_i)$ be liftings of $y_1$, $y_2$ respectively. Since $G_h(\A_f^p)$ acts transitively on $\pi_0(\pi_h^{-1}(Z_i))$, we could find $g_0 \in G_h(\A_f^p)$ such that $g_0y_2'$ and $y_1'$ are in the same connected component of $\pi_h^{-1}(Z_i)$. Let $\Tilde{y}_2 \in Z_i$ be the image of $g_0y_2'$, then $y_1$ and $\Tilde{y}_2$ are in the same same connected component of $Z_i$.

         Recall $\Phi = ((P, X_{P, +}), g)$. Let $g' = g^{-1}g_0g \in G(\A_f^p)$, and let $K' \subset G(\A_f)$ be an open compact subgroup such that $K^{\prime, p} \subset g^{\prime, -1}K^{p}g'$, $K'_p = K_p$, then by the functoriality of relative normalizations, we have a Hecke action $[\cdot g']: \Shum{K'} \to \Shum{K}$ which induces 
         \begin{equation}
             [\cdot g']^{\tor}: \Shumc{K', \Sigma'} \to \Shumc{K, \Sigma}, \quad  [\cdot g']^{\min}: \Shumm{K'} \to \Shumm{K},
         \end{equation}
         such that the resitriction of $[\cdot g']^{\min}$ on $\mathcal{Z}_{K'}(\Phi) \to \mathcal{Z}_{K}(\Phi)$ is compatible with the Hecke action $[\cdot g_0]: \Zb_{K'}^{\bigsur}(\Phi) \to \Zb_{K}^{\bigsur}(\Phi)$.

         Let $Y_2' \subset \Shum{K', T}$ be the pullback of $Y_2$ along the projection $\Shum{K'} \to \Shum{K}$, since $Y_2$ is insensitive to the away-from-$p$ Hecke action, $Y_2'$ coincides with the pullback of $Y_2$ along the Hecke action $[\cdot g']: \Shum{K'} \to \Shum{K}$. Due to \cite[Proposition 2.4.3]{lan2018compactifications} and Lemma \ref{lemma: pullback of well positioned, different level}, the canonical morphisms
         \begin{equation}
           Y_2^{\min} \times_{\Shumm{K}, [\cdot \identity]^{\min}} \Shumm{K'} \longleftarrow  Y_2^{\prime, \min} \longrightarrow Y_2^{\min} \times_{\Shumm{K}, [\cdot g']^{\min}} \Shumm{K'}
         \end{equation}
         induce isomorphisms between the reduced subschemes. Therefore, $y_2$ and $\Tilde{y}_2$ are in $Y_{2, \Zb^{\bigsur}}^{\natural}$. Since $\Tilde{y}_2$ and $y_1$ are in the same same connected of $Z_i$, then $Y_1 = Y_2$ due to the proof of Proposition \ref{prop: general well-position, disjoint}.
     \end{proof}
     \begin{corollary}\label{cor: Hecke disc, well positioned}
          Keep assumptions from Proposition \ref{prop: Hecke disc, well positioned}. We pullback these strata to $\Shum{K_p, T} = \prolim_{K^p} \Shum{K_pK^p, T}$. If the away-from-$p$ Hecke action acts transitively on the set of connected components of each $Y_i \subset \Shum{K_p, T}$ ($i \in I$), then each $Y^{\natural}_{i, \Zb^{\bigsur}(\Phi)}$ is a full union of some $\lrbracket{Z_i}_{i \in I_{\Phi_1, l_i}}$, where $Z_i \subset \Zb^{\bigsur}_{T} = \Shum{K_{\Phi, h}, T}$ (compared with Proposition \ref{prop: boundaries of fiber discrete strata}).
     \end{corollary}

\subsubsection{Nearby cycles}

   In this part, we only need to consider a single Shimura datum instead of a morphism between Shimura data. Let $\Shum{K}(G, X)$ be an integral model of a Shimura varitey $\shu{K}(G, X)$ which has a good compactification theory (i.e. satisfies Theorem \ref{theorem: axiomatic descriptions of compactifications}) and such that $C(\Phi) \to \Zb^{\bigsur}(\Phi)$ are smooth with geometrically connected fibers for all $[\Phi] \in \Cusp_K(G, X)$. For example, in case $(4)'$ (Assumption \ref{assumption: adjusted Siegel embedding}), this is true for $(G, X, K_p) := (G_1, X_1, K_{1, p})$. More generally, it is true when we replace the quasi-parahoric subgroup $K_p$ with any deep level $K_p(n) := K_p \cap \ker(\GSP(\Z_p) \to \GSP(\Z_p/p^n\Z_p))$, see \cite[Corollary 3.56]{mao2025boundary}.

    Let $\Lambda = \Q_l$ or $\ovl{\Q}_l$ ($l \neq p$), we regard it as a constant sheaf on $\shu{K}$ and on $\shu{K_{\Phi, h}}$, and $d = \dim \shu{K}$, $d_{\Phi} = \dim \shu{K_{\Phi, h}}$. We fix $T = \Bar{s}$ in this section. Let us recall
      \begin{definition-proposition}
   	Given any perverse sheaf $\mathcal{F}$ on a finite type scheme $X$, it has a Jordan-Holder series. Each constitute has the form $j_{!*}i_*\mathcal{L}([\dim Y^{\circ}])$, where $Y$ is an integral closed subscheme of $X$, $Y^{\circ}$ is a smooth open dense subscheme of $Y$, $i: Y\to X$ is the closed embedding, $j: Y^{\circ}\to Y$ is the open dense embedding, $\mathcal{L}$ is a irreducible lisse sheaf on $Y^{\circ}$. We define the \emph{Support} of $\mathcal{F}$: $\Supp(\mathcal{F})$ to be the collection of all these $Y$.
  \end{definition-proposition}

    \begin{proposition}\label{prop: boundary of nby is nby}
    Let $Y \subset \Shum{K, T}$ be a support of $\nby(\Lambda[d])$ and moreover assume $Y$ is well-positioned with respect to $\lrbracket{Y^{\natural}(\Phi)}_{\Phi}$. Assume $C(\Phi) \to \Zb^{\bigsur}(\Phi)$ are smooth with geometrically connected fibers for all $[\Phi] \in \Cusp_K(G, X)$. Then for any boundary stratum $\Zb$, $Y^{\natural}_{\Zb^{\bigsur}}$ is either empty or a support of $\nby(\Lambda[d_{\Phi}])$ on $\Zb^{\bigsur}_{T}$, and $Y^{\tor}_{\Sigma} \subset \Shumc{K, \Sigma, T}$ is a support of $\nby^{\Shumc{K, \Sigma}}(\Lambda[d])$.
\end{proposition}
\begin{remark}
    \cite[Lemma 3.7.9, Proposition 3.7.13]{lan2018compactifications} proved that to study the nearby cycle near the boundary is essentially the same as studying the nearby cycles on the boundary strata. In particular, when $C \to \Zb$ is smooth and has geometrically connected fibers (\cite[Assumption 3.7.7]{lan2018compactifications}), then any support of $\nby(\Lambda[d])$ is well-positioned. Since $C \to \Zb$ is not smooth in general, we can not use the statements there directly.
\end{remark}
\begin{remark}
    The assumption that $Y$ is well-positioned is easy to verify: due to Lemma \ref{lem: union of well-positioned sets}, since supports of nearby cycles are locally closed, it suffices to show they are unions of well-positioned subsets (for examples, in quasi-parahoric levels, they are unions of central leaves, and central leaves are well-positioned).
\end{remark}
\begin{proof}
    
    Due to \cite[Corollary 2.1.7, Lemma 4.1.2]{lan2018compactifications}, there exists a collection 
    \begin{equation}
        \mathcal{U}_{\Zb}:= \lrbracket{(\ovl{U}_i, a_{\ovl{U}_i}, a_{\ovl{U}_i, C}^{\natural})}_{i \in I_{\Zb}},
    \end{equation}
    where
    \begin{equation}
        a_{\ovl{U}_i}: \ovl{U}_i \to \Shumc{K, \Sigma}, \quad a_{\ovl{U}_i, C}^{\natural}: \ovl{U}_i \to E(\sigma_i) \times_{\Z} C\ (E := \mathbf{E}_K(\Phi))
    \end{equation}
    are \'etale morphisms such that the pullback of stratifications of $\Shumc{K, \Sigma}$ along $a_{\ovl{U}_i}$ match the pullback of stratifications of $E(\sigma_i)$ along $a_{\ovl{U}_i, C}^{\natural}$, and moreover, let $\wdt{\mathcal{Z}}([\Phi]) \subset \Shumc{K, \Sigma}$ be the preimage of $\mathcal{Z}([\Phi]) \subset \Shumm{K}$ along $\oint_{K, \Sigma}: \Shumc{K, \Sigma} \to \Shumm{K}$, then 
    \begin{equation}\label{eq: covering U bar}
        \lrbracket{\ovl{U}_i^+:= \ovl{U}_i \times_{\Shumc{K, \Sigma}} \mathcal{Z}([\Phi]) \to \mathcal{Z}([\Phi])}_{i \in I_Z}
    \end{equation}
    is an open covering for each $\mathcal{Z}([\Phi])$. Also, for each $i \in I_{\Zb}$ and each $\Zb$, the preimage of $Y^{\tor}_{\Sigma}$ under $a_{\ovl{U}_i}$ coincides with the preimage of $Y^{\natural}(\Phi)$ under the composition
    \begin{equation*}
        a_{\ovl{U}_i, \Zb}^{\natural}: \ovl{U}_i \stackrel{a_{\ovl{U}_i, C}^{\natural}}{\to} E(\sigma_i) \times_{\Z} C  \to C \to \Zb^{\bigsur} \to \Zb.
    \end{equation*}

    Let $a_{\ovl{U}_i, \Zb^{\bigsur}}^{\natural}: \ovl{U}_i \stackrel{a_{\ovl{U}_i, C}^{\natural}}{\to} E(\sigma_i) \times_{\Z} C  \to C \to \Zb^{\bigsur}$, it is smooth since $C \to \Zb^{\bigsur}$ is smooth by assumption. Let $U_i \to \ovl{U}_i$ be the common preimages of $\Shum{K} \to \Shumc{K, \Sigma}$ and of $E \times_{\Z} C \to E(\sigma)\times_{\Z} C$. Let $a_{U_i}$ (resp. $a_{U_i, \Zb^{\bigsur}}^{\natural}$) be the composition of $U_i \to \ovl{U}_i$ and $a_{\ovl{U}_i}$ (resp. $a_{\ovl{U}_i, \Zb^{\bigsur}}^{\natural}$). Then
    \begin{equation}\label{eq: pullback of Y, Y diamond, KR}
        a_{U_i}^{-1}(Y) = a_{U_i, \Zb^{\bigsur}}^{\natural, -1}(Y^{\natural}_{\Zb^{\bigsur}}).
    \end{equation}

    On the other hand, since both $a_{U_i}$ and $a_{U_i, \Zb^{\bigsur}}^{\natural}$ are smooth, up to shifting degrees, we have
    \begin{equation*}
      a_{U_i}^{-1}(\nby^{\Shum{K}}(\Lambda)) \cong \nby^{U_i}(\Lambda) \cong a_{U_i, \Zb^{\bigsur}}^{\natural, -1} (\nby^{\Zb^{\bigsur}}(\Lambda)).
    \end{equation*}
    In particular, we could refine $\ovl{U}_i$ such that $\ovl{U}_i \to \Shumc{K, \Sigma}$ and $\ovl{U}_i \to \Zb^{\bigsur}$ have geometrically connected fibers, then
    \begin{equation}\label{eq: pullback of nby}
         a_{U_i}^{-1}(\Supp \nby^{\Shum{K}}(\Lambda)) \cong \Supp \nby^{U_i}(\Lambda) \cong a_{U_i, \Zb^{\bigsur}}^{\natural, -1}(\Supp \nby^{\Zb^{\bigsur}}(\Lambda)).
    \end{equation}

    Since \ref{eq: covering U bar} forms a full covering on the boundary strata, by comparing \ref{eq: pullback of Y, Y diamond, KR} and \ref{eq: pullback of nby}, we see that $Y^{\natural}_{\Zb^{\bigsur}}$ is a support of $\nby(\Lambda[d_{\Phi}])$ on $\Zb^{\bigsur}_{T}$ if nonempty.

   Since both $a_{\ovl{U}_i}$ and $a_{\ovl{U}_i, \Zb^{\bigsur}}^{\natural}$ are smooth, up to shifting degrees, we have
    \begin{equation*}
      a_{\ovl{U}_i}^{-1}(\nby^{\Shumc{K, \Sigma}}(\Lambda)) \cong \nby^{\ovl{U}_i}(\Lambda) \cong a_{\ovl{U}_i, \Zb^{\bigsur}}^{\natural, -1} (\nby^{\Zb^{\bigsur}}(\Lambda)),
    \end{equation*}
    in particular, up to shifting degrees,
     \begin{equation*}
         a_{\ovl{U}_i}^{-1}(\Supp \nby^{\Shumc{K, \Sigma}}(\Lambda)) \cong \Supp \nby^{\ovl{U}_i}(\Lambda) \cong a_{U_i, \Zb^{\bigsur}}^{\natural, -1}(\Supp \nby^{\Zb^{\bigsur}}(\Lambda)).
    \end{equation*}
    Then $Y^{\tor}_{\Sigma}$ is a support of $\nby^{\Shumc{K, \Sigma}}(\Lambda)$ over each $\ovl{U}_i$. Since $\bigcup \ovl{U}_i$ together with $\Shum{K, T}$ form an \'etale covering of $\Shumc{K, \Sigma}$, we are done.
\end{proof}

%% file: sections/various_stratification.tex
\subsection{Index sets}

We recall the well-known index sets of various strata.

\subsubsection{$B(G)$}\label{subsubsec: intro of B(G)}

Let $F$ be a finite extension of $\Q_p$, $\bF$ be the completion of the maximal unramified extension of $F$, $\Sigma = \Gal(\ovl{F}|F)$, $I = \Gal(\ovl{F}|\bF)$, $\Sigma_0 = \Gal(\bF|F)$. Let $G$ be a $F$-reductive group, $T$ be a maximal torus of $G$ defined over $F$. Let $B(G)$ be the set of $G(\bF)$-$\sigma$-conjugate classes of $G(\bF)$. $X_*(T)_{I, \Q} = X_*(T)_{\Q}^I$, let $X_*(T)_{I, \Q}^+$ be the set of dominant elements, the action of $\sigma$ on $X_*(T)_{I, \Q}$ induces an action on $X_*(T)_{I, \Q}^+$, let $\NE(G)=(X_*(T)_{I, \Q}^+)^{\sigma}$. There are a Newton map $\nu$ as well as a Kottwitz map $\kappa$ which induce an injective Kottwitz morphism:
    \begin{equation}\label{eq: Kottwtiz map}
        B(G)\stackrel{(\nu, \kappa)}{\longrightarrow} \NE(G) \times \pi_1(G)_{\Sigma} 
    \end{equation}
  Fix a Weyl chamber, we define an order in $X_*(T)_{\Q}$ as follows: $\nu_1\leq\nu_2$ if and only if $\nu_2-\nu_1$ is a non-negative $\Q$-linear combination of positive roots. This order induces an order on $\NE(G)$. Finally, we define an order on $B(G)$ using Kottwitz morphism: $[b_1]\leq[b_2]$ if and only if $\nu([b_1])\preceq\nu([b_2])$, $\kappa([b_1])=\kappa([b_2])$. 
  
  Let $\mu \in X_*(T)$ be a cocharacter, $\bar{\mu} \in \NE(G)$ be the Galois average of the image of $\mu$ in $X_*(T)_{I, \Q}^+$ under the action of $\lrangle{\sigma}$, and
  \[ B(G, \{\mu\}) =\{[b]\in B(G)\mid \nu([b])\leq\bar{\mu}, \kappa([b])=\mu^{\sharp} \} \]

  Given a quasi-parahoric group scheme $\GG$ of $G$ defined over $\OO_F$, let $C(\GG)$ be the set of $\GG(\OO_{\bF})$-$\sigma$-conjugate classes of $G(\bF)$, then we have a natural surjection $C(\GG) \to B(G)$. We denote $C(\GG, \lrbracket{\mu})$ be the preimage of $B(G, \lrbracket{\mu})$.
\subsection{Jacobson schemes}\label{section: Jacobson schemes}

 Let $X$ be a Jacobson scheme over $s = \Spec L$, where $L$ is a field of $\Char p$, let $X_0$ be the collection of all closed points of $X$. By \cite[\href{https://stacks.math.columbia.edu/tag/005Z}{Tag 005Z}]{stacks-project}, the collection of locally closed subsets $E$ of $X$ is bijective to the collection of locally closed subsets $E_0$ of $X_0$ under $E \mapsto E_0 = E \cap X_0$. For any locally closed subset $E$ of $X$, there is a unique induced reduced subscheme structure on $E$. Let $E \to X$ be a morphism, $L'|L$ be a field extension, $\bar{s}=\Spec (L')$, $E' \to X'$ denotes the base change of $E \subset X$ along $\bar{s} \to s$. Since $\bar{s} \to s$ is a fppf covering, $E\to X$ is an immersion if and only if $E'\to X'$ is an immersion, see \cite[\href{https://stacks.math.columbia.edu/tag/02YM}{Tag 02YM}]{stacks-project}.

 Since $\Shum{K, \ovl{s}}$ is of finite type over the algebraically closed field $k = \ovl{\F}_q$, and the strata we are considering in this section are defined using geometric points, it is sufficient to focus on those strata defined in the set of $k$-points of $\Shum{K, \ovl{s}}$. This set is the set of closed points of $\Shum{K, \ovl{s}}$.

 When $Y_0$ is a locally closed subset of $\Shum{K, \ovl{s}}(k)$, it can be regarded as the $k$-points of a unique locally closed subset $Y \subset \Shum{K, \ovl{s}}$. We can equip $Y$ with the induced reduced subscheme structure.

\subsection{Stratifications in Kisin-Pappas models}

Let $(G, X)$ be a Hodge-type Shimura datum, $(G, X) \hookrightarrow (G^{\dd}, X^{\dd}) = (\GSp(V, \psi), S^{\pm})$ be a Siegel embedding, $\mu$ be a Hodge cocharacter of $(G, X)$, $\mu^{\dd}$ be the composition of $\mu$ and $G \to \GSp(V, \psi)$. Let $K_p$ be a stablizer quasi-parahoric subgroup, i.e., $K_p = \GG(\Z_p):= \GG_x(\Z_p)$ for some $x \in \Bui_{\ext}(G, \Q_p)$. Fix an embedding $\ovl{\Q} \to \ovl{\Q}_p$, we also regard $\mu$ as a cocharacter of $G_{\Q_p}$.

\subsubsection{Set up}\label{subsubsec: set up KP}

\begin{assumption}\label{general condition}
    We assume $(G, \mu, \GG)$ satisfies the following conditions: $G$ splits over a tamely ramified extension of $\Q_p$, $p>2$, and $p\nmid |\pi_1(G^{\der})|$. Moreover, assume $G_{\Q_p} \to \GL(V_{\Q_p})$ extends to a very good embedding $(\GG, \mu) \to (\GL(\Lambda), \mu_d)$, in the sense of \cite[Definition 5.2.5]{kisin2024integral}, where $\Lambda \subset V_{\Q_p}$ is a $\Z_p$-lattice that is contained in its dual under $\psi$.
\end{assumption}

\begin{definition}
    Given a morphism between Shimura data $(G_1, X_1) \to (G_2, X_2)$, we say $(G_1, X_1, K_{1, p}) \hookrightarrow (G_2, X_2, K_{2, p})$ is \emph{good} (resp. \emph{very good}) if $(G_{1, \Q_p}, \mu_{1, \ovl{\Q}_p}, K_{1, p}) \hookrightarrow (G_{2, \Q_p}, \mu_{2, \ovl{\Q}_p}, K_{2, p})$ is \emph{good} in the sense of \cite[Definition 3.1.6]{kisin2024independence} (resp. \emph{very good} in the sense of \cite[Definition 5.2.5]{kisin2024integral}), here $\mu_1$ is a Hodge-cocharacter of $(G_1, X_1)$ and $\mu_2$ is the composition of $\mu_1$ and $G_1 \to G_2$.
\end{definition}

Under the assumption \ref{general condition}, \cite{kisin2018integral} constructed a Siegel embedding $(G, X) \hookrightarrow (G^{\dd}, X^{\dd})$, together with a level group $K_p^{\dd} = \GG^{\dd}(\Z_p) \subset G^{\dd}(\Q_p)$ such that $(G, X, K_p) \hookrightarrow (G^{\dd}, X^{\dd}, K_p^{\dd})$ is a good embedding. As explained in \cite{kisin2024integral}, choosing a good embedding is not sufficient for the construction in \cite{kisin2018integral}. In order to consider central leaves, Newton strata, EKOR strata, KR strata at boundary on Kisin-Pappas integral models, we need to furthur assume the good embeddings $\iota_{\Phi}$ are very good. 

Under the assumption \ref{general condition}, we can always replace the Siegel embedding by an adjusted Siegel embedding $(G, X, K_p) \hookrightarrow (G^{\dd}, X^{\dd}, K_p^{\dd})$ (case $(4)'$, assumption \ref{assumption: adjusted Siegel embedding}) that is a very good, see \cite[\S 3.3]{mao2025boundary}, note that after applying Zarhin's trick we still get a very good embedding from a very good embedding, due to \cite[Lemma 5.3.7, 5.3.8]{kisin2024integral}. We fix such an adjusted Siegel embedding that is very good from now on.

Let $\Shum{K}(G, X)$ be the Kisin-Pappas integral model defined via $(G, X, K_p) \hookrightarrow (G^{\dd}, X^{\dd}, K_p^{\dd})$, that is to say, by choosing $K^p \subset K^{\dd, p}$ properly, let $\Shum{K}(G, X)$ be the relative normalization of $\Shum{K^{\dd}}(G^{\dd}, X^{\dd})$ in $\shu{K}(G, X)$.

\begin{proposition}{\cite[Theorem 3.58]{mao2025boundary}}\label{prop: boundary of KP models}
    Let $\iota: (G, X, K_p) \to (G^{\dd}, X^{\dd}, K_p^{\dd})$, for each $[\Phi] \in \Cusp_K(G, X)$, let 
    \[ \iota_{\Phi}: (G_{\Phi, h}, X_{\Phi, h}, K_{\Phi, h}) \to (G_{\Phi^{\dd}, h}^{\dd}, X_{\Phi^{\dd}, h}^{\dd}, K^{\dd}_{\Phi^{\dd}, h}) \]
     be the induced morphism at boundary, we have the following: 
    \begin{enumerate}
        \item If $\iota$ is an adjusted Siegel embedding, then $\iota_{\Phi}$ is an adjusted Siegel embedding.
        \item If $(G, \mu, \GG)$ satisfies the assumption \ref{general condition}, then $(G_{\Phi, h}, \mu_{\Phi, h}, \GG_{\Phi, h})$ satisfies the assumption \ref{general condition}.
        \item If $(G, \mu, \GG)$ satisfies assumption \ref{general condition}, if $\iota$ is a good embedding, then $\iota_{\Phi}$ is a good embedding.
        \item Let $\Shum{K}(G, X)$ is the Kisin-Pappas integral model defined via the good adjusted Siegel embedding $\iota$. Then $\Shum{K_{\Phi, h}}(G_h, X_{\Phi, h})$ is the Kisin-Pappas integral model defined via the good adjusted Siegel embedding $\iota_{\Phi}$.
    \end{enumerate}
  
\end{proposition}

  \subsubsection{Hodge tensors}

\begin{remark}\label{rmk: KPZ models vs KP models}
     The only reason we work with Kisin-Pappas integral models constructed in \cite{kisin2018integral} rather than the more recent Kisin-Pappas-Zhou integral models constructed in \cite{kisin2024independence} and \cite{kisin2024integral} is that, a good theory of stratifications on the special fiber of the integral model has only been written down for Kisin-Pappas integral models. Once the globally defined crystalline tensors in \cite{hamacher2019adic} has been written down in Kisin-Pappas-Zhou integral models, then statements and arguments in this article works verbatim in that case. The constructions of the globally defined crystalline tensors in \cite{hamacher2019adic} should work almost verbatim in Kisin-Pappas-Zhou case, with the help of arguments in \cite[\S 7]{kisin2024integral}.
 \end{remark}

  We briefly recall the notations in \cite{kisin2018integral} that we need. Let $\ab$ be the universal abelian scheme over $\Shum{K^{\ddagger}}$, $h: \ab_{\GG}\to\Shum{K}$ be the pullback of $\ab\to\Shum{K^{\ddagger}}$.
   
   Since $(G, X, K_p) \hookrightarrow (G^{\dd}, X^{\dd}, K_p^{\dd})$ is chosen to be a good embedding, then $\GG \to \GG^{\dd}$ is a closed embedding of group schemes. Let $\Lambda:= V_{\Z_{(p)}}$ be the self-dual lattice in $V_{\Q}$ (we also denote $\Lambda \otimes \Z_p$ by $\Lambda$ if there is no confusion) whose stablizer subgroup in $G^{\dd}(\Q_p)$ is $K_p^{\dd}$. Then $\GG$ is defined by a family of tensors $(s_{\alpha}) \in \Lambda^{\otimes}$. 
   
  Let $x\in\Shum{\Kk}(\cF_p)$, lift $x$ to some $\tilde{x}\in\Shum{\Kk}(\OO_F)$ for a finite field extension $F|\bE_v$, where $E_v$ is the completion of $E$ at a prime $v$ over $p$, $\bE_v$ is the completion of composition of $E_v$ and $\unQ_p$. We also regard $\tilde{x}$ as a $F$-point. We denote by $((t_{\alpha, \et, \tilde{x}}), (t_{\alpha, \dR, \tilde{x}}))$ the fiber of the absolute Hodge cycles induced by $(s_{\alpha})$ at $\tilde{x}$. By crystalline \'etale comparison, $(t_{\alpha, \et, \tilde{x}})$ give rise to a family of tensors on Diedonn\'e display $(t_{\alpha, \tilde{x}})\in\mM(\ab_{\GG}[p^{\infty}]_x)_{\Q}^{\otimes}$ and they are actually integral, i.e. $(t_{\alpha, \tilde{x}})\in \mM(\ab_{\GG}[p^{\infty}]_x)^{\otimes}$. Moreover, $(t_{\alpha, \tilde{x}})$ are independent of the lifting $\tilde{x}$, therefore we denote them as $(t_{\alpha, x})$. Using the isomorphism between Dieudonn\'e displays and Dieudonn\'e modules, $(t_{\alpha, x})$ give rise to $(t_{\alpha, x})\in\DD(\ab_{\GG}[p^{\infty}]_x)^{\vee \otimes}_{\bZ_p}$. We call $(t_{\alpha, x})$ the family of \emph{crystalline-Tate tensors}.

  Moreover, in \cite{hamacher2019adic}, there is a family of globally defined Tate-crystalline tensors. This allows us to construct crystalline tensors over any geometric point over $\Shum{K, \bar{s}}$. Let $\tilde{x} \to \Shum{K, \bar{s}}$ be a geometric point, since $\Shum{K, \OO_{\Breve{E}_v}}$ is flat over $\Spec \OO_{\Breve{E}_v}$, we could lift $\tilde{x}$ to a discrete valuation ring $\Spec \OO_F \to \Shum{K, \OO_{\Breve{E}_v}}$, let $C$ be the completion of the algebraic closure of $F$, and $\bar{\eta}$ be the associated geometric point over $\shu{K}$, use the absolute Hodge tensors $(t_{\alpha, \et}) \in T_p(\XX_{\bar{\eta}})$ and comparisons in \cite[\S 4.2]{hamacher2019adic}, one has a family of tensors $(t_{\alpha, \bar{x}}) \in (\DD(\XX)^{\vee}_{W(\OO_F)})^{\otimes}$, this coincides with the pullback of the global tensors $(t_{\alpha})$ due to \cite[Corollary 4.5]{hamacher2019adic}, and induces a family of crystalline $(t_{\alpha, \bar{x}}) \in (\DD(\XX)^{\vee}_{W(k(\bar{x}))}[1/p])^{\otimes}$ by the arguments above \cite[Lemma 4.4]{hamacher2019adic}. Over a perfect scheme, if the $p$-divisible group is geometrically constant, then globally defined family of tensors $(t_{\alpha})$ is also geometrically constant due to \cite[Lemma 2]{hamacher2017almost}. Let $x$ be a closed point that specializes $\bar{x}$ in the same central leaf, the Tata-crystalline tensors $(t_{\alpha, x})$ is defined integrally, thus $(t_{\alpha, \bar{x}}) \in (\DD(\XX)^{\vee}_{W(k(\bar{x}))})^{\otimes}$.

 \subsubsection{Local construction}\label{subsubsec: local construction}
 
 We follow \cite[\S 3]{shen2021ekor}. Given any geometric point $x\to\Shum{K, \Bar{s}}$, consider the following functor $I_x$:
 \begin{equation}\label{eq: I_x def}
     I_x := \Isom_{W(k(x))}((\Lambda^{\vee}\otimes W(k(x)), (s_{\alpha} \otimes \identity)), (\DD(\ab_{\GG}[p^{\infty}]_{x})_{W(k(x))}^{\vee}, (t_{\alpha, x})  ) ),
 \end{equation}
  it is a $\GG$-torsor due to \cite[Corollary 2.3.14]{shen2021ekor} where we could check the flatness over maximal points, and is trivial over the strictly henselian local ring $W(k(x))$. Let $t \in I_x(W(k(x)))$, this gives an isomorphism $\DD(\ab_{\GG}[p^{\infty}]_{x})_{W(k(x))}^{\vee}\cong \Lambda^{\vee}\otimes_{\Z_p} W(k(x))$ which maps the family of crystalline tensors $(t_{\alpha, x})$ to $(s_{\alpha}\otimes\mathrm{1})$, the pullback of the Frobenoius morphism $\varphi_x$ on the Dieudonn\'e module is of the form $(\identity\otimes \sigma) g_{x, t}$, where $g_{x, t}\in G(W(k(x))_{\Q})$. 
  
  Let $[g_{x, t}]$ be the image of $g_{x, t}$ in $B(G)(W(k(x))_{\Q})$, where $B(G)(W(k(x))_{\Q})$ is the set of $G(W(k(x))_{\Q})$-$\sigma$-conjugate class in $G(W(k(x))_{\Q})$. $[g_{x, t}]$ is independent of the choice of $t$, thus we write $[g_x] = [g_{x, t}]$. $[g_x]$ is well defined in the sense that if there is a geometric point $y$ over $x$, then $[g_x]$ maps to $[g_y]$ under the bijection $B(G)(W(k(x))_{\Q}) \to B(G)(W(k(y))_{\Q})$, see \cite[Lemma 1.3]{rapoport1996classification}. We simply write $B(G)$ for $B(G)(\bQ)$. 
  
  Let $\llbracket g_x\rrbracket \in C(G, \GG(W(k(x))))$ be the image of $g_{x, t}$ in the set of $\GG(W(k(x)))$-$\sigma$-conjugacy classes of $G(W(k(x)))_{\Q})$, it is independent of the choice of $t$. Since isomorphism classes of Dieudonne modules are independent of the base field (i.e. given $k \to k'$ an extension between algebraically closed fields, then isomorphism classes of Dieudonne modules over $W(k)$ map bijective to isomorphism classes of Dieudonne modules over $W(k')$ by field extension), and since the Hodge tensors are essentially defined over $\Z_p$, thus $C(G, \GG(W(k(x)))) \rightiso  C(G, \GG(W(k(y))))$, we simply denote it by $C(\GG)$. 
  
  The isogeny class of $(\ab_{\GG}[p^{\infty}]_{x}, (t_{\alpha, x}))$ defines a class $[g_{x}]\in B(G)$, and the isomorphism class of $(\ab_{\GG}[p^{\infty}]_{x}, (t_{\alpha, x}))$ defines a class $\llbracket g_x\rrbracket \in C(\GG)$. In particular, we have a morphism:
   \begin{equation}\label{eq: newton, K}
       \Shum{K, \bar{s}}\longrightarrow C(\GG, \{ \mu\}) \longrightarrow B(G, \{\mu\})
   \end{equation}
  We call the fibers of $\Shum{K, \bar{s}} \to C(\GG, \{ \mu\})$ and $B(G, \{\mu\})$ central leaves and Newton strata respectively. Due to \cite[Corollary 4.12]{hamacher2019adic}, Newton strata are locally closed and central leaves are closed in Newton strata.

\subsection{Stratifications in Pappas-Rapoport models}\label{subsubsec: set up PR}

In this subsection, we do not need assume \ref{general condition}. Let $(G, X, K_p) \hookrightarrow (G^{\dd}, X^{\dd}, K_p^{\dd})$ be an adjusted Siegel embedding (not necessarily a good embedding) in case $(4)'$ (assumption \ref{assumption: adjusted Siegel embedding}), where $K_p = \GG_x(\Z_p)$ is a stablizer quasi-parahoric subgroup, let $\Kf_p = \GGf_x(\Z_p)$ be a general quasi-parahoric subgroup, $\Kc_p = \GGc_x(\Z_p)$ be the parahoric subgroup.  Due to \cite[Theorem 4.3.1, 4.5.2]{pappas2024p} (for $K_p$) and \cite[Theorem 4.1.8, 4.1.12]{daniels2024conjecture} (for $\Kf_p$), the relative normalization $\Shum{K}(G, X)$ of $\Shum{K^{\dd}}(G^{\dd}, X^{\dd})$ in $\shu{K}(G, X)$ (varing $K^p$) verifies \cite[Conjecture 4.2.2]{pappas2024p}. In particular, $\Shum{K}(G, X)$ is independent of choice of Siegel embeddings due to \cite[Theorem 4.2.4]{pappas2024p} and \cite[Corollary 4.1.9]{daniels2024conjecture}.

\begin{proposition}{\cite[Theorem 3.58]{mao2025boundary}}\label{prop: boundary of PR models}
  We follow the notations in Proposition \ref{prop: boundary of KP models}. Let $\Shum{\Kf}(G, X)$ be the Pappas-Rapoport integral model with quasi-parahoric level structure defined via the adjusted Siegel embedding $\iota$. Then $\Shum{K_{\Phi, h}}(G_h, X_{\Phi, h})$ is the Pappas-Rapoport integral model with quasi-parahoric level structure defined via the adjusted Siegel embedding $\iota_{\Phi}$.
\end{proposition}

Due to \cite[Theorem 4.1.8, 4.1.12]{daniels2024conjecture}, we have commutative diagrams (see \emph{loc. cit.} for notations, we do not recall here)

\begin{equation}
\begin{tikzcd}
	{\Shum{\Kc}(G, X)^{\diamondsuit/}} & {\Shum{\Kf}(G, X)^{\diamondsuit/}} & {\Shum{K}(G, X)^{\diamondsuit/}} & {\Shum{K^{\dd}}(G^{\dd}, X^{\dd})^{\diamondsuit/}} \\
	{\Sht_{\GGc, \mu}} & {\Sht_{\GGf, \mu}} & {\Sht_{\GG, \mu}} & {\Sht_{\GG^{\dd}, \mu^{\dd}}}
	\arrow[from=1-1, to=1-2]
	\arrow[from=1-1, to=2-1]
	\arrow[from=1-2, to=1-3]
	\arrow[from=1-2, to=2-2]
	\arrow[from=1-3, to=1-4]
	\arrow[from=1-3, to=2-3]
	\arrow[from=1-4, to=2-4]
	\arrow[from=2-1, to=2-2]
	\arrow[from=2-2, to=2-3]
	\arrow[from=2-3, to=2-4]
\end{tikzcd}
\end{equation}

\subsubsection{Central leaves and Newton strata}
By applying the reduction functors introduced in \cite{gleason2024specialization}, we have
\begin{equation}\label{eq: G-shtuka, comm diagram}
\begin{tikzcd}
	{\Shum{\Kc,\bar{s}}(G, X)^{\perf}} & {\Shum{\Kf,\bar{s}}(G, X)^{\perf}} & {\Shum{K,\bar{s}}(G, X)^{\perf}} & {\Shum{K^{\dd},\bar{s}}(G^{\dd}, X^{\dd})^{\perf}} \\
	{\Sht_{\GGc, \mu}^W} & {\Sht_{\GGf, \mu}^W} & {\Sht_{\GG, \mu}^W} & {\Sht_{\GG^{\dd}, \mu^{\dd}}^W}
	\arrow[from=1-1, to=1-2]
	\arrow[from=1-1, to=2-1]
	\arrow[from=1-2, to=1-3]
	\arrow[from=1-2, to=2-2]
	\arrow[from=1-3, to=1-4]
	\arrow[from=1-3, to=2-3]
	\arrow[from=1-4, to=2-4]
	\arrow[from=2-1, to=2-2]
	\arrow[from=2-2, to=2-3]
	\arrow[from=2-3, to=2-4]
\end{tikzcd}
\end{equation}
where $\Sht_{\GG, \mu}^W$ is the stack of Witt vector $(\GG, \mu)$-shtukas, see \cite[\S 3.1]{daniels2024igusa}. Moreover, $\Sht_{\GGc, \mu}^W$ is the $\Sht_{\mu^{-1}, \KKc}^{\loc}$ (defined similarly as $\Sht_{\mu^{-1}, \KK}^{\loc}$ using $\KKc$ instead of $\KK$) used in \cite{shen2021ekor}, see \cite[Remark 3.1.8]{daniels2024igusa}.

By the main results of \cite{daniels2024conjecture}, both $\Shum{\Kf}(G, X)^{\diamondsuit/} \to \Sht_{\GGf, \mu}$ and $\Sht_{\GGc, \mu} \to \Sht_{\GGf, \mu}$ factors through an open and closed substack $\Sht_{\GGf, \mu, \delta = 1} \subset \Sht_{\GGf, \mu}$, and $\Sht_{\GGc, \mu} \to \Sht_{\GGf, \mu, \delta = 1}$ is an \'etale torsor for the abelian group $\pi_0(\GGf)^{\phi} = \Kf_p/\Kc_p$.

Let $\Bun_G$ be the $v$-stack of $G$-torsors on Fargues-Fontaine curves, let $B(G)$ be endowed with the topology from the \emph{opposite} of the partial order introduced in subsection \ref{subsubsec: intro of B(G)}. Then \cite[Theorem 1]{viehmann2024newton} gives an homeomorphism $|\Bun_G| \cong B(G)$. Let $\Bun_{G, \mu^{-1}}$ be the open substack associated with $B(G, \lrbracket{\mu^{-1}}) \subset B(G)$.

Let $G\textit{-}\Isoc$ be the $v$-stack of $G$-isocrystals (see \cite[\S 3.2.1]{daniels2024igusa}), then there is a homeomorphism $|G\textit{-}\Isoc| \cong B(G)$ ($B(G)$ endowed with the topology recalled in subsection \ref{subsubsec: intro of B(G)}). Moreover, due to \cite[Theorem 7.14]{gleason2025meromorphic}, $G\textit{-}\Isoc \cong \Bun_G^{\red}$, the reduction of $\Bun_G$. Let $G\textit{-}\Isoc_{\leq\mu^{-1}}$ be the closed substack associated with $B(G, \lrbracket{\mu^{-1}}) \subset B(G)$.

Recall that we have integral Beauville–Laszlo morphism $\BL^{\circ}: \Sht_{\GGf, \mu} \to \Sht_{\GGf} \to \Bun_G$. By \cite[Lemma 3.1.9]{daniels2024conjecture}, $\BL^{\circ}: \Sht_{\GGc, \mu} \to \Bun_{G}$ factors through $\Bun_{G, \mu^{-1}}$, thus $\BL^{\circ}: \Sht_{\GGf, \mu, \delta = 1} \to \Bun_{G}$ factors through $\Bun_{G, \mu^{-1}}$. Moreover, taking reductions, $\Sht^W_{\GGf, \mu, \delta = 1} \to G\textit{-}\Isoc_{\leq\mu^{-1}}$ is a $v$-surjection due to \cite[Proposition 3.2.3]{daniels2024igusa}.

Given $[b] \in B(G)$, let $\Sht_{\GGf, \mu, \delta = 1}^{[b]} \subset \Sht_{\GGf, \mu, \delta = 1}$, $\Sht_{\GGc, \mu}^{[b]} \subset \Sht_{\GGc, \mu}$, $G\textit{-}\Isoc^{[b]} \subset G\textit{-}\Isoc$ be the locally closed Newton strata. Since the composition of reduction maps
\begin{equation}\label{eq: isocrystal on special fiber}
  \Shum{\Kf, \bar{s}}^{\perf} \to \Sht_{\GGf, \mu}^W \to \Bun_G^{\red} \cong G\textit{-}\Isoc    
\end{equation}
gives $F$-Isocrystal structures with $G$-structures on $\Shum{\Kf, \bar{s}}^{\perf}$, then use \cite[Theorem 3.6]{rapoport1996classification}, and note that taking perfections does not change the underlying topological space, we have a semi-lower continuous map $\Shum{\Kf, \bar{s}} \to B(G, \lrbracket{\mu^{-1}})$, given $[b] \in B(G)$, we define the Newton straum $\NE^{[b]} \subset \Shum{\Kf, \bar{s}}$ as the locally closed subscheme whose perfection is $\Shum{\Kf, \bar{s}}^{\perf} \times_{G\textit{-}\Isoc} G\textit{-}\Isoc^{[b]}$.

Now, we consider the central leaves. Let $G \to \GSp(V, \psi) \to \GL(V)$ be the composition map, $\Lambda \subset V_{\Q_p}$ be the self-dual lattice associated with $K_p^{\dd}$, let $\ovl{\GG}$ be the closure of $G$ in $\GLL = \GL(\Lambda)$, then $\ovl{\GG}$ is cut out by a family of tensors $(s_{\alpha}) \in \Lambda^{\otimes}$. Then $\GG \to \ovl{\GG}$ is the smoothening, and $\GG(W(R)) = \ovl{\GG}(W(R))$ for all perfect $k$-algebra $R$, see \cite[Lemma 4.6.1]{pappas2024p}. 

We summarize some arguments in \cite[\S 4.6.3]{pappas2024p}. Let $\ab_{\GG} \to \Shum{K}$ be the pullback of the universal abelian scheme. The $\GG$-shtuka $\mathcal{P}_E$ over $\shu{K}^{\diamondsuit}$ gives a vector bundle shtuka $(\VS, \phi_{\VS})$ over $\shu{K}$. Then $(s_{\alpha, E}) \in V_{\Q_p}^{\otimes}$ give rises to a family of tensors $(t_{\alpha, E}) \in (\VS, \phi_{\VS})^{\otimes}$. Given $S = \Spa(R, R^+)$ an affinoid perfectoid space over $k$, and a map $S \to (\wdh{\Shum{K}})^{\diamondsuit}$ given by the untilt $\Spa(R^{\sharp}, R^{\sharp+}) \to (\wdh{\Shum{K}})^{\ad}$. Let $M_{\inf}(R^{\sharp+})$ be the BKF-module (with leg along $\phi(\xi) = 0$) of the pullback to $\Spec R^{\sharp+}$ of $\ab_{\GG}[p^{\infty}]$. This is a finite locally free module over $A_{\inf}(R^{\sharp+}) = W(R^+)$. Denote $(\VS_S, \phi_{\VS_S})$ the corresponding shtuka over $S$ with a leg at $S^{\sharp}$, given by the restriction of $M(W(R^+)) = (\phi^{-1})^*M_{\inf}(R^{\sharp+})$ to $\Spa(W(R^+))\setminus\lrbracket{[\mathfrak{\pi}]=0}$. The tensors $(t_{\alpha, E}) \in (\VS, \phi_{\VS})^{\otimes}$ extend uniquely to integral tensors $(t_{\alpha}) \in (\VS_S, \phi_{\VS_S})^{\otimes}$.

If $R$ is a perfect algebra, and let $Z = \Spec R \to \Shum{K}$ inducing $\Spd(R) \to (\wdh{\Shum{K}})^{\diamondsuit}$. Note that $M(W(R))$ can be identified with $\DD^{\sharp}(W(R)) = \phi^{-1}(\DD(W(R))^*)$, where $\DD(W(R))^*$ is the linear dual of the contravariant Dieudonn\'e module of $\ab_{\GG}|_{\Spec R}$, and from $(t_{\alpha})  \in (\VS_S, \phi_{\VS_S})^{\otimes}$ one obtain tensors $(t_{\alpha, \crys}) \in \DD^{\sharp}(W(R))$. Using the result that \cite[Equation 4.6.5]{pappas2024p} is a $\GG$-torsor, we have a $\GG$-torsor over $W(R)$:
\begin{equation}\label{eq: crystalline, PR}
    \mathscr{T}_{\crys}(R):= \Isom_{(t_{\alpha, \crys}), (s_{\alpha}\otimes \identity)}(\DD^{\sharp}(W(R)), \Lambda \otimes_{\Z_p} W(R)).
\end{equation}
\begin{equation}
    T_{\crys}(R):= \mathscr{T}_{\crys}(R)[1/p] = \Isom_{(t_{\alpha, \crys}), (s_{\alpha}\otimes \identity)}(\DD^{\sharp}(W(R))[1/p], \Lambda \otimes_{\Z_p} W(R)[1/p])
\end{equation}
is a $G$-torsor over $\Spec(W(R)[1/p])$, by Tannakian formalism, $T_{\crys}(R)$ gives a Frobenius $G$-isocrystal over the perfect scheme $Z = \Spec R$. These produce the morphism $\Shum{K, \bar{s}}^{\perf} \to \Sht^W_{\GG, \mu}$ and $\Shum{K, \bar{s}}^{\perf} \to G\textit{-}\Isoc$ in \ref{eq: isocrystal on special fiber}.

Since $\Sht^W_{\GGc, \mu} = \Sht^{\loc}_{\mu^{-1}, \KKc}$, we have 
\[\Sht^W_{\GGc, \mu}(k) = \Sht^{\loc}_{\mu^{-1}, \KKc}(k) = C(\GGc, \lrbracket{\mu^{-1}}) = \KKc \Adm(\lrbracket{\mu^{-1}})\KKc/\KKc_{\sigma}. \]
Let $b \in G(\bQ)$, also denote $b: \Spd k \to \Sht_{\GGc, \mu}^W$. As in \cite[\S 2.14]{hamacher2025point}, there is a locally closed subscheme $\CE^b \subset \Shum{\Kc, \bar{s}}$ whose set of $k$-points consist of those $x$ such that $x: \Spd k \to \Shum{\Kc, \bar{s}}^{\perf} \to \Sht^W_{\GGc, \mu}$ is isomorphic to $b$. Then $\CE^b \subset \NE^{[b]}$ is closed, due to \cite[Proposition 2.15(3)]{hamacher2025point}.

Let $b \in \Sht_{\GG, \mu, \delta = 1}(k)$, also regard it as an element in $G(\bQ)$, we can define $\CE^b \subset \Shum{K, \bar{s}}$ as the subscheme whose $k$-points are those  $x$ such that $x: \Spd k \to \Shum{K, \bar{s}}^{\perf} \to \Sht^W_{\GG, \mu}$ is isomorphic to $b$, later in Lemma \ref{lemma: central leaves disjoint union, different levels} we show $\CE^b(k)$ is closed in $\NE^{[b]}(k)$, which implies that $\CE^b(k)$ is locally closed in $\Shum{K, \bar{s}}(k)$, thus $\CE^b$ is the unique locally closed subscheme whose set of $k$-points is $\CE^b(k)$, see subsection \ref{section: Jacobson schemes}. We use $\CE^b_{K} \subset \Shum{K, \bar{s}}$ and $\CE^b_{\Kc} \subset \Shum{\Kc, \bar{s}}$ to distinguish the central leaves defined in different levels.

With the help of the diagram \ref{eq: G-shtuka, comm diagram}, given $b \in G(\bQ)$, we see that $\CE^b_{\Kc}$ has image contained in $\CE^b_K$ under $\varphi: \Shum{\Kc} \to \Shum{K}$.
\begin{lemma}\label{lemma: central leaves disjoint union, different levels}
  $\CE^b_{\Kc}$ is mapped onto $\CE^b_K$ under $\varphi$, $\CE^b_K$ is closed in $\NE^{[b]}$ thus locally closed in $\Shum{K, \bar{s}}$. Moreover, $\CE^b_{\Kc}$ is open and closed in $\varphi^{-1}(\CE^b_{K})$.

\end{lemma}
\begin{proof}
    Due to \cite[Proposition 2.3.1]{daniels2024conjecture}, $\varphi$ is an \'etale torsor under a finite abelian group $\pi_0(\GG)^{\phi}$, and is in particular a finite \'etale surjection, $\varphi(\CE^b_{\Kc}) \subset \NE^{[b]}$ is closed. Due to the proof of \cite[Theorem 4.1.12]{daniels2024conjecture}, we have following cartesian diagram
\[\begin{tikzcd}
	{\Shum{\Kc, \bar{s}}^{\perf}} & {\Shum{K, \bar{s}}^{\perf}} \\
	{\Sht_{\GGc, \mu}^W} & {\Sht_{\GG, \mu, \delta = 1}^W}
	\arrow[from=1-1, to=1-2]
	\arrow[from=1-1, to=2-1]
	\arrow[from=1-2, to=2-2]
	\arrow[from=2-1, to=2-2]
\end{tikzcd}\]
    where both horizontal lines are  \'etale torsor under $\pi_0(\GG)^{\phi}$. In particular, given different $b_1, b_2 \in \Sht_{\GGc, \mu}(k)$ with same image in $\Sht_{\GG, \mu}(k)$, there exists $\delta \in \pi_0(\GG)^{\phi}$ which maps $b_1$ to $b_2$, then the action $\delta$ on $\Shum{\Kc, \bar{s}}(k)$ maps $\CE^{b_1}_{\Kc}$ isomorphically onto $\CE^{b_2}_{\Kc}$, thus $\CE^{b_1}_{\Kc}$ and $\CE^{b_2}_{\Kc}$ have same image in $\Shum{K, \bar{s}}(k)$. In particular, $\varphi(\CE^b_{\Kc}) = \CE^b_K$. Since $\CE^b_{\Kc}$ is closed in $\varphi^{-1}(\NE^{[b]})$, thus it is closed in $\varphi^{-1}(\CE^b_K)$. Note that $\pi_0(\GG)^{\phi}$ is a finite abelian group which swap different central leaves $\CE^b \subset \varphi^{-1}(\CE^b_K)$, thus $\CE^b_{\Kc}$ is open and closed in $\varphi^{-1}(\CE^b_K)$.
\end{proof}
\begin{remark}\label{rmk: central leaves disjoint union, different levels}
   One can similarly define $\CE^b_{\Kf} \subset \Shum{\Kf, \bar{s}}$ as the subscheme whose $k$-points are those $x$ such that $x: \Spd k \to \Shum{\Kf, \bar{s}}^{\perf} \to \Sht^W_{\GGf, \mu}$ is isomorphic to $b$, then it is closed in $\NE^b$. Note that \cite[Proposition 2.3.1, Theorem 4.1.12]{daniels2024conjecture} work for general quasi-parahoric subgroups, thus by the arguments in Lemma \ref{lemma: central leaves disjoint union, different levels}, the finite-\'etale map $\varphi_{\Kf}: \Shum{\Kc} \to \Shum{\Kf}$ maps $\CE^b_{\Kc}$ onto $\CE^b_{\Kf}$, and $\CE^b_{\Kf}$ is open and closed in $\varphi^{\Kf, -1}(\CE^b_K)$, where $\varphi^{\Kf}: \Shum{\Kf} \to \Shum{K}$.
\end{remark}

\subsubsection{Compatible constructions}

We denote by $(t_{\alpha, \mathrm{HK}, \crys}) \in \DD^{\otimes}$ the crystalline tensors defined in \cite{hamacher2019adic}, due to \cite[Proposition 5.3.3]{daniels2024igusa}, given any perfect algebra $R$, over $\DD(W(R))^* = \phi(\DD^{\sharp}(W(R)))$, for each $\alpha$, we have an equality
\begin{equation}\label{eq: two crystalline tensors equal}
    \phi^*(t_{\alpha, \crys}) = t_{\alpha, \mathrm{HK}, \crys}
\end{equation}

Note that the morphism $\Shum{K, \bar{s}}^{\perf} \to \Sht_{\GG, \mu}^W$ is independent of the choice of the Siegel embedding, thus central leaves and Newton strata defined in Pappas-Rapoport integral models are independent of the choice of Siegel embedding. This is also true for Kisin-Pappas integral models.
\begin{proposition}\label{prop: strata independent of Siegel}
    Let $\Shum{K}(G, X)$ be a Kisin-Pappas integral model with stablier quasi-parahoric level structure (resp. with connected parahoric level structure), then the central leaves and Newton strata (resp. together with KR strata, EKOR strata) are independent of the choice of Siegel embeddings.
\end{proposition}
\begin{proof}
    For central leaves and Newton strata, we simply compare the definitions \ref{eq: I_x def} and \ref{eq: crystalline, PR} using above comparison on tensors \ref{eq: two crystalline tensors equal}. For KR strata and EKOR strata, note that in group theoretical aspect, these strata are determined by the index sets under the natural projection
    \[ \Shum{K, \bar{s}}(G, X) \to \CE(\GG, \lrbracket{\mu}) \to \prescript{\KK_p}{}{\Adm(\lrbracket{\mu})} \to \Adm(\lrbracket{\mu})_{\KK_p}. \]
\end{proof}

\subsubsection{Local model diagram}

In \cite[\S 4.9.1]{pappas2024p}, the authors constructed a functorial $v$-sheaf theoretical local model diagram
\[  \Sht_{\GG, \mu} \to [\mathbb{M}^{v}_{\GG, \mu}/\GG^{\diamondsuit}], \]
where $\mathbb{M}^{v}_{\GG, \mu}$ is the $v$-sheaf local model, $\GG^{\diamondsuit}$ is associated with $\GG$: given an affinoid perfectoid space $S = \Spa(R, R^+)$ of charactersitic $p$, $\GG^{\diamondsuit}(S) = (S^{\sharp}, g)$, where $S^{\sharp} = \Spa(R^{\sharp}, R^{\sharp+})$ is an untilt of $S$, $g \in \GG(R^{\sharp})$.

Due to \cite[Theorem 1.14]{anschutz2022p} and \cite[Theorem 1.3]{gleason2024tubular}, $\mathbb{M}^{v}_{\GG, \mu}$ is representable by a unique normal scheme $\mathbb{M}_{\GG, \mu}$ flat and proper over $\OO_{E_v}$ with reduced special fibers. \cite{pappas2024p} conjectured that the canonical models $\lrbracket{\Shum{K}(G, X)}_{K^p}$ admit schematic local model diagram for all sufficiently small $K^p \subset G(\A_f^p)$, that is to say, for each such $K^p$, there exists a smooth morphism of algebraic stack
\[ \pi_{\dR, \GG}: \Shum{K}(G, X) \to [\mathbb{M}_{\GG, \mu}/\GG] \]
whose generic fiber comes from standard principal bundle, and together with a $2$-commutative diagram ($\mathbb{M}^{\diamondsuit}_{\GG, \mu} = \mathbb{M}^{\diamondsuit/}_{\GG, \mu}$ since $\mathbb{M}_{\GG, \mu}$ is proper):
\begin{equation}\label{eq: conj, local model diagram}
    \begin{tikzcd}
	{\Shum{K}(G, X)^{\diamondsuit/}} & {\Sht_{\GG, \mu}} \\
	{[\mathbb{M}^{\diamondsuit}_{\GG, \mu}/\GG^{\diamondsuit/}]} & {[\mathbb{M}^{\diamondsuit}_{\GG, \mu}/\GG^{\diamondsuit}]}
	\arrow["{\pi_{\crys}}", from=1-1, to=1-2]
	\arrow["{\pi_{\dR, \GG}^{\diamondsuit/}}"', from=1-1, to=2-1]
	\arrow[from=1-2, to=2-2]
	\arrow[from=2-1, to=2-2]
\end{tikzcd}
\end{equation}

\begin{theorem}{\cite[Proposition 4.3.3, Theorem 4.3.6]{daniels2024conjecture}}\label{thm: local model diagram exists}
    If $p>2$, $p\nmid \pi_1(G^{\der})$, and $G$ is acceptable and $R$-smooth, then for all small enough $K^p$, $\lrbracket{\Shum{K}(G, X)}_{K^p}$, $\lrbracket{\Shum{\Kf}(G, X)}_{K^p}$ admit schematic local model diagram, and we have the following commutative diagram $(\mathbb{M}_{\GGf, \mu} = \mathbb{M}_{\GG, \mu})$:
    \begin{equation}\label{eq: local model diagram, comm}
\begin{tikzcd}
	{\Shum{\Kf}(G, X)} & {[\mathbb{M}_{\GG, \mu}/\GGf] } \\
	{\Shum{K}(G, X)} & {[\mathbb{M}_{\GG, \mu}/\GG] }
	\arrow["{\pi_{\dR, \GGf}}", from=1-1, to=1-2]
	\arrow[from=1-1, to=2-1]
	\arrow[from=1-2, to=2-2]
	\arrow["{\pi_{\dR, \GG}}", from=2-1, to=2-2]
\end{tikzcd}
    \end{equation}
\end{theorem}

\subsubsection{Kottwitz-Rapoport strata}\label{subsec: def KR strata}

In this subsection, we assume schematic local model diagram \ref{eq: conj, local model diagram} exists, and assume $p\nmid \pi_1(G^{\der})$. The existence of $\pi_{\dR, \GGc}$ is a consequence of the existence of $\pi_{\dR, \GG}$ due to \cite[Proposition 4.3.3]{daniels2024conjecture}. Note that the generic fiber of $\pi_{\dR}$ is compatible with standard principal bundle, one needs to take de Rham homology instead of de Rham cohomology in the construction of local models in \cite{kisin2018integral}, see \cite[\S A 3.2]{daniels2024conjecture} for details.

Recall that, let $\fl_{\GGc}$ be the Witt vector affine flag variety, it is representable by an ind-(perfectly projective) and ind-(perfect scheme). For each $w \in \wdt{W}_{\KKc}\setminus \wdt{W}/\wdt{W}_{\KKc}$, one can associate a Schubert perfect $k$-scheme $\fl_{\GGc, w} \subset \fl_{\GGc}$ as the closure of the Schubert perfect orbit $\fl_{\GGc, w}^{\circ} \subset \fl_{\GGc}$, here $\fl_{\GGc, w}$ is normal and $\fl_{\GGc, w}^{\circ}$ is perfectly smooth and open dense in $\fl_{\GGc, w}$, see \cite[Proposition 3.7]{anschutz2022p}. Let $\ab_{\GGc, \mu^{-1}} \subset \fl_{\GGc}$ be the $\mu^{-1}$-admissable locus, under the assumption $p\nmid \pi_1(G^{\der})$, we have 
\begin{equation}
    \ab_{\GGc, \mu^{-1}} = \bigcup_{w\in\Adm(\lrbracket{\mu^{-1}})_{\KKc}} \fl_{\GGc, w}.
\end{equation}

By \cite[Theorem 6.16]{anschutz2022p} and the fully faithfulness of the diamond functor on the category of perfect schemes, we have
\[  \mathbb{M}_{\GGc, \mu, k}^{\perf} = \ab_{\GGc, \mu^{-1}} = \bigcup_{w\in\Adm(\lrbracket{\mu^{-1}})_{\KKc}} \fl_{\GGc, w}, \]
In particular, $|[\mathbb{M}_{\GGc, \mu, k}/\GGc_k]| = \Adm(\lrbracket{\mu^{-1}})_{\KKc}$. Note that
 \[ \fl_{\GGc, w}=\textrm{the closure of}\ \fl_{\GGc, w}^{\circ}=\bigsqcup_{w'\leq w} \fl_{\GGc, w'}^{\circ},  \]
 then with the help of smoothness of $\pi_{\dR, \GGc}$, given $w \in \Adm(\lrbracket{\mu^{-1}})_{\KKc}$, let  $\KR_{\Kc, w}$ (resp. $\KR_{\Kc, \leq w}$) be the locally closed subscheme (resp. closed subscheme) whose underlying topological space is $\pi_{\dR, \GGc}^{-1}(w)$ (resp. $\bigcup_{w'\leq w} \pi_{\dR, \GGc}^{-1}(w')$), endowed with induced reduced subscheme structure, then we have
     \begin{equation}\label{equation: closure relations of KR}
        \overline{\KR_{\Kc, w}}= \KR_{\Kc, \leq w} = \bigsqcup_{w' \leq w}\KR_{\Kc, w'}
    \end{equation}

   \begin{remark}\label{remark: nearby cycles and KR strata}
       Note that $\lrbracket{\fl_{\GGc, w}}_{w \in \Adm(\lrbracket{\mu^{-1}})_{\KKc}} \subset \fl_{\GGc}$ are the supports of the nearby cycle $\nby(\Lambda[d])$ of the constant sheaf, where $d$ is the dimension of flag variety $\Fl_{G, \mu}$. Since nearby cycles commutes with smooth base change, we can also define $\KR_{\Kc, w}$ using nearby cycles to local model diagram coming from $\pi_{\dR, \GGc}$:
    \begin{equation}
        \Shum{\Kc} \stackrel{\pi}{\longleftarrow} \wdt{\Shum{\Kc}} \stackrel{q}{\longrightarrow} \mathbb{M}_{\GGc, \mu}.
    \end{equation}
   \end{remark}

%% file: sections/compactifications_of_hodge_type.tex
\subsection{Mixed Hodge structures}\label{sec: mixed Hodge structures}

We follow \cite[\S 2]{pera2019toroidal}. Fix a Shimura datum $(G, X)$ and a level group $K \subset G(\A_f)$. Fix a cusp label representative $[\Phi] \in \Cusp_K(G, X)$. Recall that, fix a $\Q$-admissible parabolic subgroup $Q \subset G$, for any $x\in X$, there exists a unique morphism $h_{\infty, x}: \DS_{\CC} \to Q_{\CC}$ ($h_{\infty, x} := \omega_x\circ h_{\infty}$ in Pink's notations) which attaches a weight filtration associated with $Q$ to the Hodge filtration induced by $h_x: \DS \to G_{\R}$, see \cite[Proposition 4.6, Proposition 4.12]{pink1989arithmetical}. Let $P$ be the smallest normal $\Q$-subgroup of $Q$ such that $h_{\infty, x}$ factors through $P_{\CC}$, due to \cite[Claim 4.7]{pink1989arithmetical}, $P$ is independent of the choice of $x \in X$, and there is a $Q(\R)$-equivariant morphism $\varphi: X \to \Hom(\DS_{\CC}, P_{\CC})$ maps $x$ to $h_{\infty, x}$. 

Consider an algebraic representation $M$ of $P$ over $\Q$. It has a weight filtration $W_{\bullet}M$ given by the weight cocharacter $w\circ h_{\infty, x}$ attached to $\Phi$ (independent of the choice of $x$). See \cite[\S 4]{pink1989arithmetical} for the relation between this weight cocharacter and the weight cocharacter attached to $(G, X)$.

Consider the quotient
\begin{equation}
    W_{\bullet}\mathbf{M}_{\betti}(\Phi):= P(\Q) \backslash (X_{\Phi} \times W_{\bullet}M(\Q)) \times (P(\A_f)/K_{\Phi, P}),
\end{equation}
it is a $\Q$-local system over $\shu{K_{\Phi, P}}(P, X_{\Phi})(\CC)$, equipped with a variation of mixed Hodge structure
\begin{equation}
    \mathbf{M}_{\MH}(\Phi) = (\mathbf{M}_{\betti}(\Phi), W_{\bullet}\mathbf{M}_{\betti}(\Phi), F^{\bullet}(\OO^{\an}_{\shu{K_{\Phi, P}}(P, X_{\Phi})(\CC)} \otimes \mathbf{M}_{\betti}(\Phi))).
\end{equation}

Let $\shu{}(P, X_{\Phi})(\CC) = \prolim_{K_{\Phi, P}}\shu{K_{\Phi, P}}(P, X_{\Phi})(\CC)$, $\shu{K_{\Phi, P}}(P, X_{\Phi})(\CC) = \shu{}(P, X_{\Phi})(\CC)/K_{\Phi, P}$. Set $\mathbf{M}_{\A_f}(\Phi) = \A_f \otimes_{\Q} \mathbf{M}_{\betti}(\Phi)$, as in \cite[\S 2.1.9.4]{pera2019toroidal}, we have an isomorphism of $\A_f$-sheaves over $\shu{}(P, X_{\Phi})(\CC)$:
\begin{equation}\label{eq: local system and sheaf}
    \underline{M}(\A_f) \cong \mathbf{M}_{\A_f}(\Phi)|_{\shu{}(P, X_{\Phi})(\CC)}.
\end{equation}

Let $M(\hat{\Z}) \subset M(\A_f)$ be a $\hat{\Z}$-lattice fixed by $K_{\Phi, P}$, then it induces a $K_{\Phi, P}$-invariant  $\hat{\Z}$-lattice in $\mathbf{M}_{\A_f}(\Phi)|_{\shu{}(P, X_{\Phi})(\CC)}$, which descends to a $\hat{\Z}$-lattice 
 $\mathbf{M}_{\hat{\Z}}(\Phi) \subset \mathbf{M}_{\A_f}(\Phi)$ over $\shu{K_{\Phi, P}}(P, X_{\Phi})(\CC) $, and induces a $\Z$-lattice
 \begin{equation}
     \mathbf{M}_{\betti}(\Phi)_{\Z} = \mathbf{M}_{\betti}(\Phi) \cap \mathbf{M}_{\hat{\Z}}(\Phi) \subset \mathbf{M}_{\betti}(\Phi).
 \end{equation}
 This gives the mixed Hodge structure $\mathbf{M}_{\MH}(\Phi)$ a $\Z$-lattice structure.

\subsection{Siegel side}

We recall some calculation from \cite[\S 2.2]{pera2019toroidal} (the following subsections follow tightly the presentations there):

Let $G^{\dd} = \GSp(V, \psi)$, $(G^{\dd}, X^{\dd})$ be a Siegel Shimura datum, an admissable parabolic subgroup $Q^{\dd} \subset G^{\dd}$ is the stablizer of an isotropic subspace $I \subset V$, let $W_{<-2}V = 0$, $W_{-2}V = I$, $W_{-1}V = I^{\bot}$, $W_{\geq 0}V = V$. Fixing a connected component $X^{\dd, +}$ of $X^{\dd} = S^{\pm}$ corresponds to the choice of $i = \sqrt{-1}$, or equivalently a point in $\pi_0(S^{\pm})$. Since $P^{\dd}(\R)$ acts transitively on $\pi_0(X^{\dd})$, the $P^{\dd}(\R)U^{\dd}(\CC)$-orbit $X_{\Phi^{\dd}}^{\dd} := X_{P^{\dd}, +}^{\dd}$ does not depend on the choice of $X^{\dd, +}$.

Let $x \in X^{\dd}$, then $h_x$ induces a Hodge filtration $F^{\bullet}_xV(\CC)$ on $V(\CC)$, let $\lambda: \Q \rightiso \Q(1)$ be the unique isomorphism such that $F^{\bullet}_xV(\CC)$ is polarized by $\lambda \circ \psi$, then $(F_x^{\bullet}V(\CC), W_{\bullet}V(\Q))$ is the polarized mixed Hodge structure on $V$ with weights $(-1, -1)$, $(0, -1)$, $(-1, 0)$, $(0, 0)$ associated with the image of $x$ in $X_{\Phi^{\dd}}^{\dd}$, since \cite[Proposition 4.12]{pink1989arithmetical} implies that the Hodge filtration on $V(\CC)$ induced by $h_{\infty, x}$ is the same as the one induced by $h_x$. That is to say, the open immersion $X^{\dd} \to X_{\Phi^{\dd}}^{\dd}$ induces a weight filtration on $V$ while keeping the Hodge filtration on $V(\CC)$.

Now we take $g \in G^{\dd}(\A_f)$ into account. We fix the self-dual integral lattice $V_{\Z}$. Let $V_{\hat{\Z}}:= V_{\Z} \otimes \hat{\Z}$, $V_{\hat{\Z}}^g = gV_{\hat{\Z}} \subset V_{\A_f}$, $V_{\Z}^g = V_{\Q} \cap V_{\hat{\Z}}^g$. Note that by approximation, $g = g'k$, $g' \in G^{\dd}(\Q)$, $k$ stablizes $V_{\hat{\Z}}$, thus the action of $g$ is indeed rationally. In particular, $V^g_{\Z} \otimes \hat{\Z} = V^g_{\hat{\Z}}$.
 
Let $\psi^g = N(g)\psi$, $N(g) \in \Q^{\times}_{> 0}$ is the unique element with image $\nu(g)^{-1}$ under $\Q^{\times} \to \A_f^{\times}/\hat{\Z}^{\times}$. Let $V^{\vee}(\nu)$ be the twist of the dual representation $V^{\vee}$ by the similitude character $\nu: G^{\dd} \to \Gm$, then $\psi^g$ induces an isomorphism of $G^{\dd}$-representations $f^g: V \rightiso V^{\vee}(\nu)$. Let $W_iV^{\vee}(\nu) \subset V^{\vee}(\nu)$ be the annihilator of $W_{-3-i}V$, then $f^g$ keeps the weight filtration.

Let $V_{\Z}^{g, \vee} \subset V^{\vee}$ be the lattice dual to $V_{\Z}^g$, and let $V_{\Z}^{g, \vee}(\nu) = N(g)V_{\Z}^{g, \vee} \subset V^{\vee}(\nu)$. The weight filtration on $V$ induces a weight filtration on $V_{\Z}^g$, and we define a dual filtration on $V_{\Z}^{g, \vee}$. Since $V_{\Z}$ is self-dual under $\psi$, then $V_{\Z}^{g}$ is self-dual under $\psi^g$.

\subsubsection{}

Apply the construction in subsection \ref{sec: mixed Hodge structures} to $V_{\Z}^g \subset V$ gives a variation of mixed $\Z$-Hodge structure $\mhg{\Phi^{\dd}}_{\Z}$ (we use $\mhg{\Phi^{\dd}}$ instead of $\mathbf{V}_{\MH}(\Phi^{\dd})$ to avoid the confusion between the representations of the vector space $V$ and the unipotent group $V = W/U$) of weights $(-1, -1), (-1, 0), (0, -1), (0, 0)$ over $\shu{K_{\Phi^{\dd}, P^{\dd}}^{\dd}}(P^{\dd}, X_{\Phi^{\dd}}^{\dd})(\CC)$. The pairing $\psi^g$ gives a polarization on $\mhg{\Phi^{\dd}}_{\Z}$: $\mhg{\Phi^{\dd}}_{\Z} \times \mhg{\Phi^{\dd}}_{\Z} \to \mathbf{1}_{\Z}(1)$. Using the equivalence between $1$-motive and variation of mixed Hodge-structure of weights $(-1, -1), (-1, 0), (0, -1), (0, 0)$, we have polarized $1$-motif $(\mQ, \lambda)$ over $\shu{K_{\Phi^{\dd}, P^{\dd}}^{\dd}}(P^{\dd}, X_{\Phi^{\dd}}^{\dd})(\CC)$. We have canonical identifications $\mdr{\Phi^{\dd}} = H^1_{\dR}(\mQ)$, $\hat{T}(\mQ) = \mathbf{H}_{\hat{\Z}}(\Phi^{\dd})$.

\subsubsection{}\label{subsubsec: 1 motive, Siegel side, comparison}

Moreover, let $K_p^{\dd}$ be the stablizer of the self-dual lattice $V_{\Z_p}$. Consider the integral model $\Shum{K_{\Phi^{\dd}, P^{\dd}}^{\dd}}(P^{\dd}, X_{\Phi^{\dd}}^{\dd})$, it is the stack of groupoids over $\Z_{(p)}$ such that its $S$-point parametrizes the following data $(\mQ, \lambda, \alpha, \alpha^{\vee}, \epsilon)$, where:
\begin{compactitem}
    \item $(\mQ, \lambda)$ is a polarized $1$-motif over $S$,
    \item $\alpha: \gr^0_W \underline{V}^g(\Z) \rightiso \mQ^{\et}$, $\alpha^{\vee}: \gr^0_W \underline{V}^{g, \vee}(\nu) \rightiso \mQ^{\vee, \et} = \mQ^{\mult}$,
    \item $\epsilon \in H^0(S, \Isom(\underline{V}^g(\hat{\Z}^p), \hat{T}^p(\mQ))/K_{\Phi^{\dd}}^{\dd, p})$ which is compatible with $\lambda, \alpha, \alpha^{\vee}$.
\end{compactitem}
Note that $\Shum{K_{\Phi^{\dd}, P^{\dd}}^{\dd}}(P^{\dd}, X_{\Phi^{\dd}}^{\dd})$ is a smooth quasi-projective scheme over $\Z_{(p)}$.

The integral model $\Shum{K_{\Phi^{\dd}, h}^{\dd}}(G_h^{\dd}, X_{\Phi^{\dd}, h}^{\dd})$ parametrizes $(\mB, \lambda^{\abe}, \epsilon^{\abe})$ over a $\Z_{(p)}$-scheme $S$, where $(\mB, \lambda^{\abe})$ is a principally polarized abelian scheme over $S$, and 
\[\epsilon^{\abe}\in H^0(S, \Isom(\gr^{-1}_W\underline{V}^g(\hat{\Z}^p), \hat{T}(\mB))/K_{\Phi^{\dd}, h}^{\dd, p}).\] 
The natural projection $\Shum{K_{\Phi^{\dd}, P^{\dd}}^{\dd}}(P^{\dd}, X_{\Phi^{\dd}}^{\dd}) \to \Shum{K_{\Phi^{\dd}, h}^{\dd}}(G_h^{\dd}, X_{\Phi^{\dd}, h}^{\dd})$ maps $(\mQ, \lambda, \alpha, \alpha^{\vee}, \epsilon)$ to $(\mQ^{\abe}, \lambda^{\abe}, \epsilon^{\abe})$ ($\epsilon^{\abe} = \gr^{-1}_W\epsilon$). Now we use $\Xi^{\dd}(\Phi^{\dd}), \Zb^{\bigsur, \dd}(\Phi^{\dd})$ instead of $\Shum{K_{\Phi^{\dd}, P^{\dd}}^{\dd}}(P^{\dd}, X_{\Phi^{\dd}}^{\dd})$ and $\Shum{K_{\Phi^{\dd}, h}^{\dd}}(G_h^{\dd}, X_{\Phi^{\dd}, h}^{\dd})$ to simplify the notations.

Given any point $x \in \Xi^{\dd}_{\sigma^{\dd}}(\Phi^{\dd}) \cong \mathcal{Z}([\Phi^{\dd}, \sigma^{\dd}])$, let $R_x$ be the complete local ring of $\Xi^{\dd}(\sigma^{\dd})(\Phi^{\dd})$ at $x$, let $S = \Spec R_x$, $U = S \times_{\Xi^{\dd}(\sigma^{\dd})(\Phi^{\dd})}\Xi^{\dd}(\Phi^{\dd})$, then the polarized $1$-motif $(\mQ, \lambda)|_U$ is associated with a canonical polarized abelian scheme $(\ab_U', \psi_U')$ over $U$, and we have following properties:
\begin{proposition}\label{prop: relations between 1 motif and abe}\leavevmode
    \begin{enumerate}
    \item $\ab_U^{\prime, \vee}$ is associated with $\mQ^{\vee}.$
    \item $\ab_U'[n]\rightiso\mQ|_U[n]$ as finite flat group schemes over $U$, and is compatible with polarizations.
    \item $H^1_{\dR}(\ab_U')\rightiso H^1_{\dR}(\mQ|_U)$ as filtered vector bundles with integrable connections.
    \item $\hat{T}(\ab_U') \rightiso \hat{T}(\mQ|_U)$ as $\hat{\Z}$-local systems compatible with the Weil pairings induced by polarizations.
\end{enumerate}
\end{proposition}

The tautological $K^p_{\Phi^{\dd}}$-level structure $\epsilon \in H^0(S, \Isom(\underline{V}^g(\hat{\Z}^p), \hat{T}^p(\mQ))/K_{\Phi^{\dd}}^{\dd, p})$ induces a $K^{\dd, p}$-level structure $\eta'\in H^0(S, \Isom(\underline{V}(\hat{\Z}^p), \hat{T}^p(\ab_U))/K^{\dd, p})$ via $\underline{V}(\hat{\Z}^p) \stackrel{g}{\to} \underline{V}^g(\hat{\Z}^p) \stackrel{\epsilon}{\to} \hat{T}^p(\mQ) \rightiso \hat{T}^p(\ab_U')$. $(\ab_U', \psi_U', \eta')$ induces a canonical morphism $U \to \Shum{K^{\dd}}$, and is isomorphic to the pullback $(\ab_U, \psi_U, \eta)$ of tautological triple over $\Shum{K^{\dd}}$. 

\subsection{Hodge tensors}\label{subsec: Hodge tensors on 1-motive}

\subsubsection{Rational tensors}

Fix a pair of cusp label representatives $(\Phi, \Phi^{\dd})$ as usual, $\Phi = ((P, X_{P, +}, g))$, $\Phi^{\dd} = ((P^{\dd}, X^{\dd}_{P^{\dd}, +}, g))$. One can pull back the universal $1$-motive $(\mQ, \lambda)$ over $\shu{K_{\Phi^{\dd}, P^{\dd}}^{\dd}}$ to $\shu{K_{\Phi, P}}$. Let $G^{\dd, g} = \GSp(V^g_{\Q}, \psi^g)$, then under the canonical identification $V^g_{\Q} = V_{\Q}$, $G^{\dd, g} \cong G^{\dd}$. Let $G^g \hookrightarrow G^{\dd, g}$ be associated with $G \hookrightarrow G^{\dd}$ under such identification $V^g_{\Q} \cong V_{\Q}$.

Let $(s_{\alpha, \Phi}) \in (V_{\Q}^g)^{\otimes}$ and $(s_{\alpha}) \in V_{\Q}^{\otimes}$ be the collection of tensors defining the closed subgroup $G^{g} \hookrightarrow G^{\dd, g}$ and $G \hookrightarrow G^{\dd}$ respectively. Then under the canonical identification $V^g_{\Q} = V_{\Q}$, $(s_{\alpha, \Phi}) = (s_{\alpha})$. Due to \cite{brylinski19831} (see \cite[Proposition 3.1.2]{pera2019toroidal}), over $\shu{K_{\Phi, P}}$, we have:
\[
    (s_{\alpha, \Phi, \dR}) \in \mdr{\Phi}^{\otimes} := \mdr{\Phi^{\dd}}^{\otimes},\quad (s_{\alpha, \Phi, \et}) \in \mathbf{H}_{\A_f}(\Phi)^{\otimes} := \mathbf{H}_{\A_f}(\Phi^{\dd})^{\otimes}.
\]
Since $\mdr{\Phi^{\dd}} = H^1_{\dR}(\mQ)$, $\mathbf{H}_{\hat{\Z}}(\Phi^{\dd}) = \hat{T}(\mQ)$, we also write (over $\shu{K_{\Phi, P}}$)
\begin{equation}
    (s_{\alpha, \Phi, \dR}) \in H^1_{\dR}(\mQ)^{\otimes}, \quad  (s_{\alpha, \Phi, \et})\in \hat{T}(\mQ)^{\otimes}_{\A_f}.
\end{equation}
Keep notations from subsection \ref{subsubsec: 1 motive, Siegel side, comparison}. Let $U_0$ be the generic fiber of $U$. The pullback $\ab_{U_0}$ of the universal abelian scheme over $U_0$ is equipped with tensors
\begin{equation}
    (t_{\alpha, \dR}) \in H^1_{\dR}(\ab_{U_0})^{\otimes}, \quad  (t_{\alpha, \et}) \in \hat{T}(\ab_{U_0})^{\otimes}.
\end{equation}
Using the identification (due to \cite{lan2012comparison}) between the analytic construction (by Pinks) and the algebraic construction (by Chai-Faltings), the canonical isomorphisms 
\begin{equation}\label{eq: comparison over U_0, without tensors}
    \hat{T}(\mQ_{U_0}) \rightiso \hat{T}(\ab'_{U_0}) \rightiso \hat{T}(\ab_{U_0}), \quad H^1_{\dR}(\mQ_{U_0}) \rightiso H^1_{\dR}(\ab'_{U_0}) \rightiso H^1_{\dR}(\ab_{U_0})
\end{equation}
induce
\begin{proposition}{\cite[Proposition 3.1.6]{pera2019toroidal}}\label{prop: comparison, 1-motif and ab}
   Over $U_0$, $(s_{\alpha, \Phi, \et}) \mapsto (t_{\alpha, \et})$, $(s_{\alpha, \Phi, \dR}) \mapsto (t_{\alpha, \dR})$.
\end{proposition}

\subsubsection{Integral tensors}

From now on, we work under the setting of Kisin-Pappas integral models, assumption \ref{general condition}, let $(G, X, K_p) \hookrightarrow (G^{\dd}, X^{\dd}, K_p^{\dd})$ be an adjusted Siegel embedding that is a very good embedding.
      
Pull back the tautological polarized $1$-motif $(\mQ, \lambda)$ to $\Shum{K_{\Phi, P}}$, and the universal abelian scheme $(\ab, \lambda)$ to $\Shum{K}$. Let $x \in \Xi_{\sigma}(\Phi) \cong \mathcal{Z}([\Phi, \sigma])$, let $R_{G, x}$ be the complete local ring of $\Xi(\sigma)(\Phi)$ at $x$, let $S = \Spec R_{G, x}$, $U = S \times_{\Xi(\sigma)(\Phi)}\Xi(\Phi)$.

Recall that we have a family of Hodge tensors $(s_{\alpha}) \in V_{\Z_{(p)}}^{\otimes}$, which cuts out the closed subgroup scheme $\GG \to \GSP$ over $\Z_{(p)}$. The tensors at boundary $\Phi$ are twisted by $g$. We denote by $(s_{\alpha, \Phi}) \in V_{\Z_{(p)}}^{g, \otimes}$ (where $V^g_{\Z_{(p)}} = gV_{\Z_p} \cap V_{\Q}$) be the collection of Hodge tensors which cuts out the closed subgroup scheme $\GG^g \to \GSP^g$, here $\GG^g$ and $\GSP^g$ are the $g$-conjugation of $\GG$ and $\GSP$ respectively. There is a canonical identification $V^g_{\Z_{(p)}} \otimes \Q = V_{\Z_{(p)}} \otimes \Q$, under which we have $(s_{\alpha, \Phi}) = (s_{\alpha})$ rationally. $\GG^{g} \to \GSP^{g}$ has generic fiber $G^{g} \to \GSp^{g}$ that can be identified with $G \to \GSp$ via identify $V^g_{\Z_{(p)}} \otimes \Q = V_{\Z_{(p)}} \otimes \Q$. 

In the $\A_f$-module
\[ \underline{V}(\A_f) \stackrel{g}{\rightiso} \underline{V}^g(\A_f) \rightiso \hat{T}(\mQ_{U_0})_{\A_f} \rightiso \hat{T}(\ab'_{U_0})_{\A_f} \rightiso \hat{T}(\ab_{U_0})_{\A_f}, \]
by identifying $\underline{V}^g(\A_f) = \underline{V}(\A_f)$, the $\hat{\Z}$-lattice $\hat{T}(\mQ|_{U_0})$ corresponds to the twisted lattice $\underline{V}^g(\hat{\Z})$ stablized by $K_{\Phi}$ and $\hat{\Z}$-lattice $\hat{T}(\ab_{U})$ corresponds to the standard lattice $\underline{V}(\hat{\Z})$ stablized by $K$.

Recall that the \'etale tensors over $\shu{K_{\Phi, P}}$ comes from the Betti tensors on $\mathbf{H}_{\betti}(\Phi) \cong \underline{V}^g(\Q)$. The integral tensors $(s_{\alpha, \Phi})$ supports on the integral structure on $\mathbf{H}_{\betti}(\Phi)_{\Z_{(p)}} \cong \underline{V}^g(\Z_{(p)})$, thus we have integral tensors on $(s_{\alpha, \et}) \in T_p(\mQ)^{\otimes}$ which maps to integral Hodge tensors $(t_{\alpha, \et}) \in T_p(\ab_{U_0})^{\otimes}$ under the identification over $U_0$.

\subsubsection{}

By construction of the Hodge tensors $(s_{\alpha, \Phi}) \in \mb{\Phi}^{\otimes}$, at each point $\bar{t} \in \shu{K_{\Phi, P}}(\CC)$, we have an isomorphism
\begin{equation}\label{eq: interior isomorphism on MH}
(\mb{\Phi}^{\otimes}_{\Z_{(p)}, \bar{t}}, (s_{\alpha, \Phi, \bar{t}})) \cong (V_{\Z_{(p)}}^g, (s_{\alpha}))    
\end{equation}
which keeps the weight filtrations and the pairings on both sides. Note that the weight filtrations on both sides determine and are determined by the fixed admissable parabolic subgroup $Q_{\Phi}$.

Such an isomorphism induces an isomorphism on the subquotient $\mb{\mQ^{\abe}}_{\Z_{(p)}} \cong V_{1, \Z_{(p)}}^g$ which keeps the polarization on the left and the perfect paring on the right, here $\mb{\mQ^{\abe}}_{\Z_{(p)}}$ is the Betti realization of the abelian scheme $\mQ^{\abe} = \mB$ together with the integral structure determined by $K_{\Phi, h}$, and $V_{1, \Z_{(p)}}^g = \gr^{-1}_W V^g_{\Z_{(p)}}$. 

Recall that in Proposition \ref{prop: boundary of KP models}, the induced morphism between boundary Shimura data $(G_{\Phi, h}, X_{\Phi, h}, K_{\Phi, h, p}) \hookrightarrow (G_{\Phi, h}^{\dd}, X_{\Phi, h}^{\dd}, K_{\Phi^{\dd}, h, p}^{\dd})$ is again an adjusted Siegel embedding, if $K_p = \GG(\Z_p)$ is a stablizer quasi-parahoric subgroup, then $K_{\Phi, h, p} = \GG_h^g(\Z_p)$ is also a stablizer quasi-parahoric subgroup, here $\GG_h^g$ is a smooth model of $G^g_h \cong G_h$. We can apply the constructions of the absolute Hodge tensors on this pair. Since $(G_{\Phi, h}, X_{\Phi, h}, K_{\Phi, h, p}) \hookrightarrow (G_{\Phi, h}^{\dd}, X_{\Phi, h}^{\dd}, K_{\Phi^{\dd}, h, p}^{\dd})$ is a good embedding, let $(s_{\alpha, \Phi}^{\star}) \in (V_{1, \Z_{(p)}}^g)^{\otimes}$ be the Hodge tensors defining $\GG^g_h(\Phi) \hookrightarrow \GG_h^{\dd, g}(\Phi^{\dd})$ (with $\GG_h^{\dd, g}(\Phi^{\dd}) = \GSP(V_{1, \Z}^g)$), we similarly have 
\[ (t_{\alpha, \Phi, B}^{\star}) \in \mb{\mQ^{\abe}}^{\otimes}_{\Z_{(p)}} = H^1_B(\mB(\CC))_{\Z_{(p)}},\quad (t_{\alpha, \Phi, \dR}^{\star}) \in H^1_{\dR}(\mB)^{\otimes},\quad (t_{\alpha, \Phi, \et}^{\star}) \in \hat{T}(\mB)^{\otimes}. \]
\begin{remark}
    Note that the collection $(s_{\alpha, \Phi}^{\star})$ can be strictly larger then the restriction of $(s_{\alpha, \Phi})$ on the graded piece $\gr^{-1}_W$, since $P \subset G \cap P^{\dd}$ can be strict.
\end{remark}

Note that the subquotient $V_{1, \Z_{(p)}}^g$ of the local system $W^{\bullet}V_{\Z_{(p)}}^g$ (constructed using $W^{\bullet}V_{\Z_{(p)}}^g$) on $\shu{K_{\Phi, P}}(\CC)$ descends over $\shu{K_{\Phi, h}}(\CC)$ that coincides with the local system constructed by $V_{1, \Z_{(p)}}^g$. By regarding $\bar{t} \in \shu{K_{\Phi, h}}(\CC)$, by construction, the isomorphism \ref{eq: interior isomorphism on MH} induces an isomorphism
\begin{equation}
   (H^1_B(\mB(\CC))_{\Z_{(p)}, \bar{t}}, (t_{\alpha, \Phi, B, \bar{t}}^{\star})) \cong (V_{1, \Z_{(p)}}^g, (s_{\alpha, \Phi}^{\star})).
\end{equation}

Due to the theory of absolute Hodge cycles, we have the following proposition:

\begin{proposition}\label{prop: dr tensors, 1-motive}
    For every geometric point $\bar{t} \to \shu{K_{\Phi, P}}$, there exists an isomorphism 
    \begin{equation}\label{eq: dr Hodge trivialization}
        (W^{\bullet}\mdr{\Phi}_{\bar{t}}, (s_{\alpha, \Phi, \dR, \bar{t}})) \cong (W^{\bullet}V^g, (s_{\alpha}))
    \end{equation}
    which induces an isomorphism on the subquotient graded piece
    \begin{equation}\label{eq: dr Hodge trivialization, over B}
        (H^1_{\dR}(\mB)_{\bar{t}}, (t_{\alpha, \Phi, \dR, \bar{t}}^{\star})) \cong (V_{1}^g, (s_{\alpha, \Phi}^{\star})).
    \end{equation}
    Same statement is true with the de-Rham realization $H_{\dR}$ replaced by the \'etale realization $T_p$.
\end{proposition}

Since $W^0$ is flat over $S$, Given a geometric point $\bar{t} \to W^0_{\bar{s}}$, there exists a discrete valuation ring $\Spec \OO_F \to W^0$ with residue field $k(\bar{t})$ and generic point $\eta$ of characteristic $0$. Since closed points of $W$ concentrate on $W \setminus W^0$, then the image of $\bar{t}$ in $W^0$ specializes to some $x \in W \backslash W^0$. Assume $x$ has image in $\mathcal{Z}([\Phi, \sigma]) \cong \Xi_{\sigma}(\Phi)$. Recall that we let $S = \Spec R_{G, x}$ be the complete local ring of $\Xi(\sigma)(\Phi)$ at $x$, and $U = S \times_{\Xi(\Phi)(\sigma)} \Xi(\Phi)(\sigma)$. Since $U \to W^0$ is flat, we can lift $\eta \to W^0$ to $\eta \to U$ and $\bar{t} \to W^0$ to $\bar{t} \to U$ up to field extensions, thus we assume $\Spec \OO_F \to W^0$ factors through $U$. By the comparisons \ref{prop: comparison, 1-motif and ab}, comparisons \cite[\S 4]{hamacher2019adic}, uniqueness of the extensions of Hodge tensors, and Proposition \ref{prop: dr tensors, 1-motive}, we have
\begin{proposition}
    There exists an isomorphism 
    \begin{equation}\label{eq: de rham Hodge tri, abelian sch}
        (W^{\bullet}H^1_{\dR}(\ab)_{\OO_F}, (t_{\alpha, \dR, \OO_F})) \cong (W^{\bullet}V \otimes \OO_F, (s_{\alpha}))
    \end{equation}
    which induces an isomorphism on the (conjugated) subquotient graded piece
    \begin{equation}\label{eq: de rham Hodge tri, over B, abelian sch}
        (H^1_{\dR}(\mB)_{\OO_F}, (t_{\alpha, \Phi, \dR, \OO_F}^{\star})) \cong (V_{1}^g \otimes \OO_F, (s_{\alpha, \Phi}^{\star})).
    \end{equation}
\end{proposition}
\begin{proposition}
    There exists an isomorphism 
    \begin{equation}\label{eq: crys Hodge tri, abelian sch}
        (W^{\bullet}\DD(\ab_{\bar{t}})_{W(k(\bar{t}))}, (t_{\alpha, \bar{t}})) \cong (W^{\bullet}V \otimes W(k(\bar{t})), (s_{\alpha}))
    \end{equation}
    which induces an isomorphism on the (conjugated) subquotient graded piece
    \begin{equation}\label{eq: crys Hodge tri, over B, abelian sch}
        (\DD(\mB_{\bar{t}})_{W(k(\bar{t}))}, (t_{\alpha, \Phi, \bar{t}}^{\star})) \cong (V_{1}^g \otimes W(k(\bar{t})), (s_{\alpha, \Phi}^{\star})).
    \end{equation}
\end{proposition}

%% file: sections/newton_strata.tex
From this section, we prove that various strata on the special fiber of the Kisin-Pappas's integral model of a Hodge-type Shimura variety with stablizer quasi-parahoric level group are well-positioned, and prove some useful propositions along the way. We also give some remarks regarding Pappas-Rapoport integral models.

\subsection{Kisin-Pappas integral models}\label{subsec: well-position, KP model}

From now on, unless specified, we use the set-up introduced in subsection \ref{subsubsec: set up KP}, in particular, we assume \ref{general condition}. Fix an adjusted Siegel embedding $(G, X, K_p) \hookrightarrow (G^{\dd}, X^{\dd}, K_p^{\dd})$ which is very good, where $K_p = \GG_x(\Z_p)$ be a stablizer quasi-parahoric subgroup. Let $K^p \subset G(\A_f^p)$ and $K^{\dd, p} \subset G^{\dd}(\A_f^p)$ be properly chosen, and let $\Shum{K}(G, X)$ is the Kisin-Pappas integral model of the Hodge-type Shimura datum $\shu{K}(G, X)$. Under this setting, both $\Shum{K}(G, X)$ and $\Shum{K^{\dd}}(G^{\dd}, X^{\dd})$ satisfy assumption \ref{assumption: C to Z} due to Proposition \ref{remark: when the assumption C to Z is true}.

  \begin{proposition}\label{proposition: newton strata are well-positioned}
   	Newton strata satisfy the assumption \ref{assumption: fiber discreteness}. In particular, let $b \in B(G)$ $\NE^{[b]}$, $\NE^{\leq [b]}$, and their connected components are well-positioned subsets of $\Shum{K, \Bar{s}}$, and they are well-positioned subschemes with respect to the induced reduced subscheme structures.
   \end{proposition}
   \begin{proof}
 Recall that in the construction \ref{subsubsec: local construction}, we have a commutative diagram given by forgetting the structure of Hodge tensors:
  \begin{equation}\label{eq: Newton strata, comm diag, KP}
\begin{tikzcd}
	{\Shum{K, \bar{s}}(G, X)} & {\Shum{K^{\dd}, \bar{s}}(G^{\dd}, X^{\dd})} \\
	{B(G, \lrbracket{\mu})} & {B(G^{\dd}, \lrbracket{\mu^{\dd}})}
	\arrow["\varphi", from=1-1, to=1-2]
	\arrow[from=1-1, to=2-1]
	\arrow[from=1-2, to=2-2]
	\arrow[from=2-1, to=2-2]
\end{tikzcd}
  \end{equation}
        Since $\Shum{K, \Bar{s}} \to B(G, \lrbracket{\mu} )$ and $\Shum{K^{\dd}, \Bar{s}} \to B(\GSp,  \lrbracket{\mu^{\dd}} )$ are lower semi-continuous due to \cite[Theorem 3.6]{rapoport1996classification} and \cite[Corollary 4.12]{hamacher2019adic}, given $[b] \in B(G, \lrbracket{\mu})$ with image $[b^{\dd}] \in B(G^{\dd}, \lrbracket{\mu^{\dd}})$, we need to show $\NE^{[b]}$ is open and closed in $\varphi^{-1}(\NE^{[b^{\dd}]})$: let $J \subset B(G, \lrbracket{\mu})$ be the collection of elements with image $[b^{\dd}] \in B(G^{\dd}, \lrbracket{\mu^{\dd}})$, this is a finite set, due to Lemma \ref{lemma: B(G) quasi-finite}. Note that $\varphi^{-1}(\NE^{[b^{\dd}]}) = \bigsqcup_{[b] \in J} \NE^{[b]}$, it suffices to show given different $[b_1], [b_2] \in J$, then $\ovl{\NE^{[b^1]}}$ and $\ovl{\NE^{[b^2]}}$ has no intersection in $\varphi^{-1}(\NE^{[b^{\dd}]})$. Recall that, $\ovl{\NE^{[b_i]}} \subset \NE^{\leq [b_i]}$, if $\ovl{\NE^{[b_1]}} \cap \ovl{\NE^{[b_2]}} \cap \varphi^{-1}(\NE^{[b^{\dd}]})$ is non-empty, then there exists $b \in J$ such that $[b] < [b_i]$ for $i = 1, 2$. Due to Proposition \ref{prop: BG main prop}, this is impossible. Therefore, the preimage of Newton strata on $\Shum{K^{\dd}, \bar{s}}$ are topologically disjoint union of Newton strata on $\Shum{K, \bar{s}}$. 

        We follow the notations 
        The Newton strata on $\Shum{K^{\ddagger}, \bar{s}}$ here are the same as the one appeared in \cite[\S 3.3]{lan2018compactifications}. Due to \cite[Proposition 3.3.9]{lan2018compactifications}, Newton strata on $\Shum{K^{\ddagger}, \bar{s}}$ are well-positioned. These verify the the assumption \ref{assumption: fiber discreteness}. Then $\NE^{b}$ and its connected components are well-positioned subsets of $\Shum{K, \bar{s}}$, and are well-positioned subscheme endowed with induced reduced subscheme structures, due to Proposition \ref{prop: general well-position}.

         It might happen that for some $[b_1], [b_2]\in \iota^{-1}(\iota([b]))$, $\NE^{\leq [b_1]}\cap \NE^{\leq [b_2]}$ is non-empty, and $\overline{\NE^{[b_i]}}\subsetneq \NE^{\leq [b_i]}$. Nevertheless, the closed subset $\NE^{\leq [b]} = \sqcup_{[b_i] \in B(G), [b_i] \leq [b]} \NE^{[b_i]}$ is a well-positioned subset due to Lemma \ref{lem: union of well-positioned sets}, is a well-positioned subscheme with induced reduced subscheme structure due to Remark \ref{remark: subset well positioned implies subscheme well positioned}, its connected components are well-positioned due to Proposition \ref{proposition: open-closed subschemes are well positioned}.
\end{proof}

\begin{proposition}\label{proposition: central leaves are well--positioned}
   	Central leaves satisfy the assumption \ref{assumption: fiber discreteness}. In particular, central leaves and their connected components are well-positioned subsets of $\Shum{K, \Bar{s}}$, and they are well-positioned subschemes with respect to the induced reduced subscheme structures.
   \end{proposition}
   \begin{proof}
       Given a central leaf $\CE^b \subset \Shum{K, \bar{s}}$, let $b^{\dd}$ be the image of $b$ in $G^{\dd}(\bQ)$, then $\CE^b \subset \varphi^{-1}(\CE^{b^{\dd}})$. The central leaves on $\Shum{K^{\ddagger}, \bar{s}}$ are the same as the one used in \cite[\S 3.4]{lan2018compactifications}. Due to \cite[Proposition 3.4.2]{lan2018compactifications}, central leaves on $\Shum{K^{\ddagger}, \bar{s}}$ are well-positioned. To apply the proof of Proposition \ref{proposition: newton strata are well-positioned}, it suffices to show $\CE^b \subset \varphi^{-1}(\CE^{b^{\dd}})$ is open and closed. Due to the proof of Proposition \ref{proposition: newton strata are well-positioned}, $\NE^{[b]}$ is closed in $\varphi^{-1}(\NE^{[b^{\dd}]})$, and $\CE^b$ is closed in $\NE^{[b]}$, thus $\CE^b$ is closed in $\varphi^{-1}(\CE^{b^{\dd}})$.
       
       Let $S = \Shum{K, \Bar{s}}^{\perf}$ ($S^{\dd} = \Shum{K^{\dd}, \Bar{s}}^{\perf}$) be the perfection of $\Shum{K, \Bar{s}}$ (resp. $\Shum{K^{\dd}, \Bar{s}}$), let $\XX^{\perf}$ be the pullback of the universal $p$-divisible group along $S \to S^{\dd} \to \Shum{K^{\dd}, \Bar{s}}$. If $\XX^{\perf}$ is geometrically constant on a scheme $T$ over $k$ (i.e. given any two points $x, y \in T$, $\XX^{\perf}_{\Bar{x}} \cong \XX^{\perf}_{\Bar{y}}$ for some geometric points $\Bar{x}$, $\Bar{y}$ over $x, y$ respectively), then the globally defined family of tensors $(t_{\alpha})$ is also geometrically constant on $T$ due to \cite[Lemma 2]{hamacher2017almost}. The proof there shows that one can even take a connected and perfect scheme $T'$ such that $T' \to T^+$ ($T^+$ is any fixed connected component of $T$) is surjective and $\XX^{\perf}$ and $(t_{\alpha})$ are constant over $T'$.
    
     The perfection of $\varphi^{-1}(\CE^{b^{\dd}})$ is exactly the locus of $S$ where $\XX^{\perf}$ (without the tensors structures) is geometrically constant and isomorphic to $\XX^{\perf}_{\Bar{x}}$ for some geometric point $\Bar{x}$ over $x\in \CE^b$. By \cite[Corollary 4.11]{hamacher2019adic}, $\XX^{\perf}$ has a globally defined Tata-crystalline tensors $(t_{\alpha})$ over $S$. By definition, the perfection of $\CE^b$ is the locus of $S$ where $(\XX^{\perf}, (t_{\alpha}))$ is geometrically isomorphic to $(\XX^{\perf}_{\bar{x}}, (t_{\alpha, \bar{x}}))$ for some geometric point $\bar{x}$ over $x \in \CE^b$. In particular, $\CE^b \subset \varphi^{-1}(\CE^{b^{\dd}})$ is open.

     Note that $\varphi^{-1}(\CE^{b^{\dd}})$ is locally closed in the Noetherian space $\Shum{K, \bar{s}}$, thus itself is Noetherian, there are only finitely many connected components. In particular, $\varphi^{-1}(\CE^{b^{\dd}})$ is a finite union of central leaves on $\Shum{K, \bar{s}}$.
       \end{proof}
    We find that such fiber discreteness of central leaves has already been worked out in Hesse's thesis \cite[Theorem 2.11]{hesse2020central}. We choose to keep the above proof here, since it is also used in the more general case, see the proof of Proposition \ref{prop: Newton central are well-posiitoned, PR}.
       
 Since the partial minimal compactifications of central leaves in $\Shum{K^{\dd}, \Bar{s}}$ are affine (\cite[Proposition 1.9]{caraiani2024generic}, \cite[Theorem 3.3.3]{santos2023generic}), due to Proposition \ref{prop: affineness in general}, we have:
   \begin{corollary}\label{cor: minimal compactifications of central leaves are affine}
       Partial minimal compactifications of central leaves in $\Shum{K, \Bar{s}}$ are affine.
   \end{corollary}

   \subsubsection{Boundary has no intersections}

  \begin{lemma}\label{lemma: closure of central leaves}
        Let $\CE^b \subset \Shum{K, \bar{s}}$ be a central leaf, $\NE^{[b]} \subset \Shum{K, \bar{s}}$ be a Newton stratum, then the closures $\ovl{\CE^b} \subset \Shum{K, \bar{s}}$ and $\ovl{\NE^{[b]}} \subset \Shum{K, \bar{s}}$ are unions of central leaves in $\Shum{K, \bar{s}}$.
    \end{lemma}
    \begin{proof}
       It suffices to work with closed points. Given a closed point $x \in \Shum{K, \bar{s}}$, we can consider its formal local ring $\OO_{\Shum{K}, x}^{\wedge}$. Let $x \in \ovl{\CE^b}(k)$ (resp. $x \in \ovl{\NE^{[b]}}(k)$), then $x \in \CE^{b_1}(k)$ for some $b_1 \in G(\bQ)$. Let $y \in \CE^{b_1}(k)$ be another point in the central leaf. By definition, we have $(\ab_{\GG}[p^{\infty}]_x, (t_{\alpha, x})) \cong (\ab_{\GG}[p^{\infty}]_y, (t_{\alpha, y}))$. We have $\OO_{\Shum{K}, x}^{\wedge} \cong \OO_{\Shum{K}, y}^{\wedge}$ as versal deformation space $\Def_{\GG}(\ab_{\GG}[p^{\infty}]_x; (t_{\alpha, x})) \cong \Def_{\GG}(\ab_{\GG}[p^{\infty}]_y; (t_{\alpha, y}))$. Note that taking formal neighbourhoods preserves central leaves stratification (resp. Newton stratification), $\Spec \OO_{\Shum{K}, x}^{\wedge} \cap \CE^b$ (resp. $\Spec \OO_{\Shum{K}, x}^{\wedge} \cap \NE^{[b]}$) is nonempty implies that $\Spec \OO_{\Shum{K}, y}^{\wedge} \cap \CE^b$ (resp. $\Spec \OO_{\Shum{K}, y}^{\wedge} \cap \NE^{[b]}$) is nontrivial, thus $\CE^{b_1}$ is contained in $\ovl{\CE^b}$ (resp. $\ovl{\NE^{[b]}}$).
    \end{proof}

    \begin{lemma}\label{lemma: boundary of Newton no intersection, siegel}
    Let $\NE^{[b^{\dd}_1]}, \NE^{[b^{\dd}_2]} \subset \Shum{K^{\dd}, \bar{s}}$ be Newton strata, assume $[b^{\dd}_1] \neq [b^{\dd}_2]$, then 
    \[  \NE^{[b^{\dd}_1], \min} \cap \NE^{[b^{\dd}_2], \min} = \emptyset,\quad \NE^{[b^{\dd}_1], \tor}_{\Sigma^{\dd}} \cap \NE^{[b^{\dd}_2], \tor}_{\Sigma^{\dd}} = \emptyset.  \]
\end{lemma}
\begin{proof}
   We omit the $(\ast)^{\dd}$ to save notations. It suffices to show $\NE^{[b_1], \min} \cap \NE^{[b_2], \min} = \emptyset$. The case of toroidal compactifications follows from the case of minimal compactifications since $\oint_{K, \Sigma}(\NE^{[b_i], \tor}_{\Sigma}) = \NE^{[b_i], \min}$.

  Since Newton strata are insensitive to the away-from-$p$ Hecke action, due to Lemma \ref{lemma: pullback of well positioned, different level}, it suffices to prove the case when $K^{p}$ is properly chosen such that $\Zb^{\bigsur}(\Phi) \cong \Zb(\Phi)$ for every $[\Phi] \in \Cusp_K(G, X)$ (this is always possible in Siegel case). This can simplify the arguments. Assume $\NE^{[b_1], \min} \cap \NE^{[b_2], \min}$ is nontrivial, then there exists a boundary stratum $\Zb(\Phi) \cong \mathcal{Z}([\Phi])$ such that $\NE^{[b_1], \min} \cap \NE^{[b_2], \min} \cap \Zb(\Phi)_{\bar{s}}$ is nontrivial. Due to \cite[Proposition 3.3.9]{lan2018compactifications}, $\NE^{[b_i]}$ is well-positioned, $\NE^{[b_i], \min} \cap \Zb(\Phi)_{\bar{s}} = \NE^{[b_i], \natural}_{\Zb(\Phi)}$, where $\NE^{[b_i], \natural}_{\Zb(\Phi)} = \NE^{[b_i^{\natural}]}$ is itself a Newton stratum of $\Zb^{\bigsur}(\Phi)_{\Bar{s}} = \Shum{K_{\Phi, h}, \Bar{s}}(G_h, X_{\Phi, h})$ for some $[b_i^{\natural}] \in B(G_h)$. In the proof of \cite[Proposition 3.3.9]{lan2018compactifications}, we see that $[b_1^{\natural}] = [b_2^{\natural}]$ if and only if $[b_1] = [b_2]$, $\NE^{[b_1], \min} \cap \NE^{[b_2], \min} \cap \Zb(\Phi)_{\bar{s}}$ is nontrivial implies that $[b_1^{\natural}] = [b_2^{\natural}]$, thus $\NE^{[b_1]}  = \NE^{[b_2]}$, contradiction.
\end{proof}

With the help of \cite[Proposition 3.4.2]{lan2018compactifications}, we have an analogue of Lemma \ref{lemma: boundary of Newton no intersection, siegel}:
 \begin{lemma}\label{lemma: boundary of central leaves no intersection, siegel}
    Let $\CE^{b^{\dd}_1}, \CE^{b^{\dd}_2} \subset \Shum{K^{\dd}, \bar{s}}$ be central leaves, assume $[[b^{\dd}_1]] \neq [[b^{\dd}_2]]$, then 
    \[  \CE^{b^{\dd}_1, \min} \cap \CE^{b^{\dd}_2, \min} = \emptyset,\quad \CE^{b^{\dd}_1, \tor}_{\Sigma^{\dd}} \cap \CE^{b^{\dd}_2, \tor}_{\Sigma^{\dd}} = \emptyset.  \]
\end{lemma}

     \begin{corollary}\label{corollary: central leaves, exclusive}\leavevmode
       \begin{enumerate}
          \item $\Shum{K, \bar{s}}$ (resp. $\Shumm{K, \bar{s}}$, $\Shumc{K, \Sigma, \bar{s}}$) is a disjoint union of connected components of $\CE^b$ (resp. $\CE^{b, \min}$, $\CE^{b, \tor}_{\Sigma}$) where $b$ runs over the index set $C(\GG, \lrbracket{\mu})$.
          \item Given any central leaves $\CE^{b_1}$, $\CE^{b_2}$,
          \begin{equation}
              \ovl{\CE^{b_1}} \cap \CE^{b_2} \neq \emptyset \Longleftrightarrow \ovl{\CE^{b_1, \min}} \cap \CE^{b_2, \min} \neq \emptyset \Longleftrightarrow \ovl{\CE^{b_1, \tor}_{\Sigma}} \cap \CE^{b_2, \tor}_{\Sigma} \neq \emptyset.
          \end{equation}
      \end{enumerate}
   \end{corollary}
   \begin{proof}
         Due to Remark \ref{remark: set/scheme theoretical closure does not matter}, we focus on the underlying topological spaces. We apply Proposition \ref{prop: general well-position, disjoint}, with the help of Lemma \ref{lemma: boundary of central leaves no intersection, siegel} and \ref{lemma: closure of central leaves}.
   \end{proof}

    Similarly, use Lemma \ref{lemma: boundary of Newton no intersection, siegel} in place of Lemma \ref{lemma: boundary of central leaves no intersection, siegel}, we have
     \begin{corollary}\label{corollary: Newton strata, exclusive}\leavevmode
      \begin{enumerate}
          \item $\Shum{K, \bar{s}}$ (resp. $\Shumm{K, \bar{s}}$, $\Shumc{K, \Sigma, \bar{s}}$) is a disjoint union of connected components of $\NE^{[b]}$ (resp. $\NE^{b, \min}$, $\NE^{b, \tor}_{\Sigma}$) where $b$ runs over the index set $B(G, \lrbracket{\mu})$.
          \item Given any Newton strata $\NE^{b_1}$, $\NE^{b_2}$,
          \begin{equation}
              \ovl{\NE^{b_1}} \cap \NE^{b_2} \neq \emptyset \Longleftrightarrow \ovl{\NE^{b_1, \min}} \cap \NE^{b_2, \min} \neq \emptyset \Longleftrightarrow \ovl{\NE^{b_1, \tor}_{\Sigma}} \cap \NE^{b_2, \tor}_{\Sigma} \neq \emptyset.
          \end{equation}
      \end{enumerate}
   \end{corollary}

   \subsubsection{Boundary descriptions}\label{subsec: boundary descriptions, newton}

 In this part, we furthur determine the boundary strata of given Newton strata and central leaves. In order to define Newton strata and central leaves at boundary, we need to furthur assume all the good embeddings $\iota_{\Phi}$ are very good in Proposition \ref{prop: boundary of KP models}. We keep this assumption from now on.

Apply Lemma \ref{lemma: boundary of Newton no intersection, siegel} and Proposition \ref{prop: boundaries of fiber discrete strata}:
\begin{corollary}
    Newton strata satisfy assumption \ref{assumption: fiber discreteness, enhanced}, where we take both $\lrbracket{Z_i}_{i \in I_{\Phi}}$, $\lrbracket{Z_j}_{j \in J_{\Phi^{\dd}}}$ be Newton strata of $\Shum{K_{\Phi, h}, \bar{s}}$, $\Shum{K_{\Phi^{\dd}, h}^{\dd}, \bar{s}}$ respectively. Given a Newton stratum $\NE^{[b]} \subset \Shum{K, \bar{s}}$, its boundary stratum $\NE^{[b], \natural}_{\Zb^{\bigsur}(\Phi)} \subset \Zb^{\bigsur}(\Phi)_{\bar{s}}$ (for each $[\Phi] \in \Cusp_K(G, X)$) is a union of connected components of some Newton strata in $\Zb^{\bigsur}(\Phi)_{\bar{s}} = \Shum{K_{\Phi, h}, \bar{s}}$.
\end{corollary}

Apply Lemma \ref{lemma: boundary of central leaves no intersection, siegel} and Proposition \ref{prop: boundaries of fiber discrete strata}:
\begin{corollary}
    Central leaves satisfy assumption \ref{assumption: fiber discreteness, enhanced}, where we take both $\lrbracket{Z_i}_{i \in I_{\Phi}}$, $\lrbracket{Z_j}_{j \in J_{\Phi^{\dd}}}$ be central leaves of $\Shum{K_{\Phi, h}, \bar{s}}$, $\Shum{K_{\Phi^{\dd}, h}^{\dd}, \bar{s}}$ respectively. Then given a central leaf $\CE^{b} \subset \Shum{K, \bar{s}}$, its boundary stratum $\CE^{b, \natural}_{\Zb^{\bigsur}(\Phi)} \subset \Zb^{\bigsur}(\Phi)_{\bar{s}}$ (for each $[\Phi] \in \Cusp_K(G, X)$) is a union of connected components of some central leaves in $\Zb^{\bigsur}(\Phi)_{\bar{s}} = \Shum{K_{\Phi, h}, \bar{s}}$.
\end{corollary}

 We can say something more about the boundary strata of given Newton strata and central leaves.
 \begin{remark}\label{rmk: fix mu and mu_h}
     Given a point $x \in X_{\Phi, P}$, we can associated it with a Hodge cocharacter $\mu_P = h_x \circ w: (\Gm)_{\CC} \to P_{\CC}$. Different $x$ give conjugated $\mu_P$. We can find a $\mu_P$ defined over $\ovl{\Q}$. On the other hand, we take $\mu$ be a Hodge cocharacter associated with a point $x \in X$ and is defined over $\ovl{\Q}$ after conjugation. Due to \cite[Proposition 12.1]{pink1989arithmetical}, $\mu_P$ is $G(\CC)$-conjugated to $\mu$, thus is $G(\ovl{\Q})$-conjugated to $\mu$. In particular $\lrbracket{\mu} = \lrbracket{\mu_P}$ in $G(\CC)$ and in $G(\ovl{\Q})$. Fix an embedding $\ovl{\Q} \to \ovl{\Q}_p$, $\mu$ and $\mu_P$ define the same conjugate class in $G(\ovl{\Q}_p)$. Also, consider the projeciton $\mu_h: (\Gm)_{\ovl{\Q}} \to P_{\ovl{\Q}} \to G_{h, \ovl{\Q}}$ to its Levi part, since the image of $\mu_P$ is semisimple, we can take a section $G_h \to P$ and assume $\mu_h = \mu_P$ if necessary.
 \end{remark}
 
 Given a geometric point $\bar{t} \to W^0_{\bar{s}}$ as usual, due to \ref{eq: crys Hodge tri, abelian sch} and \ref{eq: crys Hodge tri, over B, abelian sch}, there exists an isomorphism 
    \begin{equation}\label{eq: Frob use 1}
        (W^{\bullet}\DD(\ab_{\bar{t}})_{W(k(\bar{t}))}, (t_{\alpha, \bar{t}})) \cong (W^{\bullet}V_{W(k(\bar{t}))}, (s_{\alpha}))
    \end{equation}
    which induces an isomorphism on the (conjugated) subquotient graded piece
    \begin{equation}\label{eq: Frob use 2}
        (\DD(\mB_{\bar{t}})_{W(k(\bar{t}))}, (t_{\alpha, \Phi, \bar{t}}^{\star})) \cong (V_{1, W(k(\bar{t}))}^g, (s_{\alpha, \Phi}^{\star})).
    \end{equation}
   Under such isomorphism, the $\sigma$-linear Frobenious morphism on the Dieudonne module $\DD(\ab_{\bar{t}})_{W(k(\bar{t}))}$ (resp. $\DD(\mB_{\bar{t}})_{W(k(\bar{t}))}$) is sent to an element $b\sigma$ (resp. $b_h\sigma$) where $b \in G(W(k(\bar{t}))_{\Q})$ (resp. $b \in G_h(W(k(\bar{t}))_{\Q})$). In particular, the image of $\bar{t} \to \Shum{K, \bar{s}}$ (resp. $\bar{t} \to \Shum{K_{h}, \bar{s}}$) falls into the Newton stratum $\NE^{[b]}$ (resp. $\NE^{[b_h]}$). This implies that $[b] \in B(G, \lrbracket{\mu})$ (resp. $ [b_h] \in B(G_h, \lrbracket{\mu_h})$)
   
   Now we determine the relations between $b$ and $b_h$. Since Frobenious preserves the weight filtration, $b \in Q(W(k(\bar{t}))_{\Q})$ ($Q = Q^{\dd} \cap G$). Consider the graded pieces, 
   \[ \gr^W_{-2}\DD(\ab_{\bar{t}}) = \DD(\mT_{\bar{t}}), \quad \gr^W_{-1}\DD(\ab_{\bar{t}}) = \DD(\mB_{\bar{t}}), \quad \gr^W_{0}\DD(\ab_{\bar{t}}) = \DD(Y_{\bar{t}} \otimes \Q_p/\Z_p).\]
   Since the Frobenious morphism determines and is determined by its action on the isoclinic subquotients (induced by the slope filtration), and Frobenious on slope $0$ (the \'etale part) and slope $1$ (the multiplicative part) are constant, therefore the Frobenious morphism determines and is determined by its action on the graded pieces induced by the weight filtration. First, we need to replace $b$ by $b^g:=gb\sigma(g)^{-1} \in G(W(k(\bar{t}))_{\Q})$ under the canonical idetification $G^g \cong G$. Let $b_L$ be the projection of $b^g$ under $Q \to L$ (here we use the identification $Q^g \cong Q, L^g \cong L$, $G^g_h \cong G_h$), we get $b_L \in L(W(k(\bar{t}))_{\Q})$. 
   
   Let us forget about the Hodge tensors for now. Recall that we have a canonical isomorphism $L^{\dd} = G_h^{\dd} \times G_l^{\dd}$, where $G_h^{\dd}$ is associated with $\gr^W_{-1}$, and $G_l^{\dd}$ is associated with $\gr^W_{-0}$ (here $\gr^W_{-2}$-part is determined by the $\gr^W_{-0}$-part due to the existence of the principal polarization). The Frobenious action on the $\gr^W_{-0}$-part (the \'etale part) is simply the identity morphism, thus the projection of $b^g$ under $L \to L^{\dd} \to G_l^{\dd}$ is trivial ($G_l^{\dd} = \GL(\gr^W_{-2}V) \cong \GL(\gr^W_{0}V)$). In particular, let $b_h'$ be the image of $b_L$ under $L \to L^{\dd} \to G_h^{\dd}$, we have $b_L = b_h' \in L \cap G_h^{\dd}$. Note that the isomorphism \ref{eq: Frob use 2} keeps the tensors, we moreover have $b_L \in G_h(W(k(\bar{t}))_{\Q})$, and this element is $b_h$. In other words, by choosing the isomorphism \ref{eq: Frob use 1}, we force $b^g \in P(W(k(\bar{t}))_{\Q})$, and $b_h \in G_h(W(k(\bar{t}))_{\Q})$ is simply the projection of $b^g$ under $P \to G_h$.

      Given another point $\bar{t}' \in W^0$, without lose of generality, assume $k(\bar{t}') = k(\bar{t})$, repeat this process, we have corresponding Newton strata $\NE^{[b']}$, $\NE^{[b_h']}$. Note that $[b]=[b'] \in B(G)$ if and only if $[b^g] = [b^{\prime, g}] \in B(G^g) = B(G)$. Assume $[b_h] = [b_h'] \in B(G_h)$. Then $b_h = hb_h'\sigma(h)^{-1}$ for some $h \in G_h(W(k(\bar{t}))_{\Q})$. Lift $h$ to an element $\tilde{h} \in P(W(k(\bar{t}))_{\Q})$. $b_L = hb_L'\sigma(h)^{-1}$ implies that $b^g = \tilde{h}b^{\prime, g}\sigma(\tilde{h})^{-1}$ modulo $W(W(k(\bar{t})_{\Q}))$. Since $B(Q) \to B(L)$ is a canonical bijection, we have $[b] = [b'] \in B(G)$.
      On the other hand, assume $[b] = [b'] \in B(G)$. Then $b^g = hb^{\prime,g}\sigma(h)^{-1}$ for some $h \in G(W(k(\bar{t}))_{\Q})$. Due to the isomorphism \ref{eq: Frob use 1} we can furthur request $h \in Q(W(k(\bar{t}))_{\Q})$ since it keeps the weight filtration. Let $h_L$ be the image of $h$ under $Q \to L$, then $b_L = h_L b_L' \sigma(h_L)^{-1} \in G_h(W(k(\bar{t}))_{\Q})$. In other words, $b_h = h_L b_h' \sigma(h_L)^{-1} \in G_h(W(k(\bar{t}))_{\Q})$.
    With the help of Lemma \ref{lem: B(G) bijection}, we arrive at the following proposition:

     \begin{proposition}\label{prop: boundary of Newton strata}
        Let $\NE^{[b]} \subset \Shum{K, \Bar{s}}$ be a Newton stratum, then its boundary stratum $\NE^{[b], \natural}_{\Zb^{\bigsur}(\Phi)} \subset \Zb^{\bigsur}(\Phi)_{\Bar{s}}$ (for each $[\Phi]\in \Cusp_K(G, X)$) is either empty or a Newton stratum $\NE^{[b_h]} \subset \Zb^{\bigsur}(\Phi)_{\Bar{s}}$, where $[b_h] \in B(G_h)$ is associated with $[b] \in B(G)$ in the way presented above.
     \end{proposition}

   Let $\bar{t}, \bar{t}' \in W^0$ be two geometric points that lying over $\CE^{b}$, $\CE^{b'}$ respectively (along $W^0 \to \Shum{K, \bar{s}}$), 
and over $\CE^{b'}$, $\CE^{b_h'}$ respectively (along $W^0 \to \Zb^{\bigsur}(\Phi)_{\bar{s}}$). Take $k(\bar{t}) = k(\bar{t}')$ to simplify notations. Note that in the arguments above Proposition \ref{prop: boundary of Newton strata}, we showed that that the Frobenious morphisms which keep the weight filtration on the Dieudonne modules only factor through the graded pieces of the weight filtration, thus one can pick representatives $b, b' \in P(W(k(\bar{t}))_{\Q})$, $b_h, b_h' \in G_h(W(k(\bar{t}))_{\Q})$ in the equivalent classes $[[b]], [[b']] \in \CE(\GG, \lrbracket{\mu})$, $[[b_h]], [[b'_h]] \in \CE(\GG_h^g, \lrbracket{\mu_h})$ such that the projections of $b^g = b_h$, $b^{\prime, g} = b_h'$. Note that $[[b]] = [[b']] \in \CE(\GG)$ if and only if $[[b^g]] = [[b^{\prime, g}]] \in \CE(\GG^g)$.

Assume $[[b_h]] = [[b_h']] \in C(\GG_h^g, \lrbracket{\mu_h})$. Then $b_h = hb_h'\sigma(h)^{-1}$ for some $h \in \GG_h^g(W(k(\bar{t})))$, since $\KK_{\Phi, h, p} = \GG_h^g(W(k(\bar{t}))) \subset \KK_{\Phi, p} = \GG^g(W(k(\bar{t})))$, thus $b^g = hb^{\prime, g}\sigma(h)^{-1}$ with $h \in \GG^g(W(k(\bar{t})))$, we have $[[b]] = [[b']] \in C(\GG, \lrbracket{\mu})$.
      On the other hand, assume $[[b]] = [[b']] \in C(\GG, \lrbracket{\mu})$. Then $b^g = hb^{\prime, g}\sigma(h)^{-1}$ for some $h \in \GG^g(W(k(\bar{t})))$. Due to the isomorphism \ref{eq: Frob use 1} we can furthur request $h \in Q(W(k(\bar{t}))_{\Q}) \cong Q^g(W(k(\bar{t}))_{\Q})$ since it keeps the weight filtration. Let $h_L$ be the image of $h$ under $Q \to L$, then $b_L = h_L b_L' \sigma(h_L)^{-1} \in G_h(W(k(\bar{t}))_{\Q})$. In other words, $b_h = h_L b_h' \sigma(h_L)^{-1} \in G_h(W(k(\bar{t}))_{\Q})$ and $h_L \in \KK_{\Phi, L, p} = \pi(\KK_{\Phi} \cap Q(W(k(\bar{t}))_{\Q}))$. With the help of Corollary \ref{cor: GGc to LLc injective when Z connected}, we have:
  \begin{proposition}\label{prop: boundary of central leaves}
        Let $\CE^b \subset \Shum{K, \Bar{s}}$ be a central leaf, then its boundary stratum $\CE^{b, \natural}_{\Zb^{\bigsur}(\Phi)} \subset \Zb^{\bigsur}(\Phi)_{\Bar{s}}$ (for each $[\Phi]\in \Cusp_K(G, X)$) is either empty or a topologically disjoint finite union of some $\CE^{b_h} \subset \Zb^{\bigsur}(\Phi)_{\Bar{s}}$, where the collection $\lrbracket{b_h} \subset G_h(\bQ)$ is associated with $b \in G(\bQ)$ in the way presented above. Moreover, $[[b_h]] \in C(\GG^g_h)$ is uniquely determined by $[[b]] \in C(\GG)$ when $1 \to G_h \to L \to G_l \to 1$ splits as $L \cong G_h \times G_l$, or the closure of $Z_h$ (the center of $G_h$) in $\GG_h$ (assume $\GGc_h = \GG_h$) has connected fibers.
     \end{proposition}

 \subsection{Pappas-Rapoport integral models}\label{subsec: well-posiiton, PR models}

From now on, unless specified, we use the set-up introduced in subsection \ref{subsubsec: set up PR}. Under this setting, both $\Shum{\Kf}(G, X)$ and $\Shum{K^{\dd}}(G^{\dd}, X^{\dd})$ satisfy assumption \ref{assumption: C to Z} due to Proposition \ref{remark: when the assumption C to Z is true}. Here $\Kf$ is a general quasi-parahoric level subgroup.

 In subsection \ref{subsubsec: set up PR}, we recall the definitions of central leaves and Newton strata on Pappas-Rapoport integral models. We show that many statements in subsection \ref{subsec: well-position, KP model} involving Kisin-Pappas integral models work in this setting. To be more precise:
 \begin{proposition}\label{prop: Newton central are well-posiitoned, PR}\leavevmode
     \begin{enumerate}
         \item Newton strata and central leaves on $\Shum{\Kf, \bar{s}}$ satisfy the assumption \ref{assumption: fiber discreteness}. In particular, let $b \in G(\bQ)$, then $\CE^b$, $\NE^{[b]}$, $\NE^{\leq [b]}$, and their connected components are well-positioned subsets of $\Shum{\Kf, \Bar{s}}$, and they are well-positioned subschemes with respect to the induced reduced subscheme structures.
         \item Newton strata and central leaves on $\Shum{\Kf, \bar{s}}$ satisfy the assumption \ref{assumption: fiber discreteness, enhanced}. In particular, let $b \in G(\bQ)$, let $\NE^{[b]} \subset \Shum{\Kf, \Bar{s}}$ (resp. $\CE^{b} \subset \Shum{\Kf, \Bar{s}}$) be a Newton stratum (resp. a central leaf), then its boundary stratum $\NE^{[b], \natural}_{\Zb^{\bigsur}(\Phi)} \subset \Zb^{\bigsur}(\Phi)_{\Bar{s}}$ (resp. $\CE^{b, \natural}_{\Zb^{\bigsur}(\Phi)}$) (for each $[\Phi]\in \Cusp_K(G, X)$) is a topologically disjoint union of some connected components of some finite collections of $\NE^{[b_h]} \subset \Zb^{\bigsur}(\Phi)_{\Bar{s}}$ (resp. $\CE^{[b_h]} \subset \Zb^{\bigsur}(\Phi)_{\Bar{s}}$).
         \item $\Shum{K, \bar{s}}$ (resp. $\Shumm{K, \bar{s}}$, $\Shumc{K, \Sigma, \bar{s}}$) is a disjoint union of connected components of $\CE^b$ (resp. $\CE^{b, \min}$, $\CE^{b, \tor}_{\Sigma}$) where $b$ runs over the index set $C(\GG, \lrbracket{\mu^{-1}})$, and is a disjoint union of connected components of $\NE^{[b]}$ (resp. $\NE^{b, \min}$, $\NE^{b, \tor}_{\Sigma}$) where $b$ runs over the index set $B(G, \lrbracket{\mu^{-1}})$.
         \item Partial minimal compactifications $\CE^{b, \min}$ of central leaves $\CE^{b}$ are affine.
     \end{enumerate}
 \end{proposition}
\begin{proof}\leavevmode
    \begin{enumerate}
        \item We apply the proofs in Proposition \ref{proposition: newton strata are well-positioned} and \ref{proposition: central leaves are well--positioned} with the help of Lemma \ref{lemma: central leaves are open and closed, PR}.
        \item Same as the arguments below Lemma \ref{lemma: boundary of Newton no intersection, siegel} and \ref{lemma: boundary of central leaves no intersection, siegel}. Note that due to Proposition \ref{prop: boundary of PR models}, $\iota_{[\Phi]}$ are again adjusted Siegel embeddings and $\Shum{K_{\Phi, h}}$ are again Pappas-Rapoport integral models.
        \item Same as Corollary \ref{corollary: central leaves, exclusive} and \ref{corollary: Newton strata, exclusive}.
        \item Same as Corollary \ref{cor: minimal compactifications of central leaves are affine}.
    \end{enumerate}
\end{proof}

\begin{remark}
    Part $(2)$ of the Proposition \ref{prop: Newton central are well-posiitoned, PR} can be made stronger as Proposition \ref{prop: boundary of Newton strata} and \ref{prop: boundary of central leaves}, but this is enough for our purpose in later sections.
\end{remark}

\begin{lemma}\label{lemma: central leaves are open and closed, PR}
    Let $\varphi: \Shum{\Kf}(G, X) \to \Shum{K^{\dd}}(G^{\dd}, X^{\dd})$ be the induced finite morphism. Let $b \in G(\bQ)$ and $b^{\dd} \in G^{\dd}(\bQ)$ be the image of $b$. Then $\NE^{[b]}$ (resp. $\CE^b$) is open and closed in $\varphi^{-1}(\NE^{[b^{\dd}]})$ (resp. $\varphi^{-1}(\CE^{b^{\dd}})$).
\end{lemma}
\begin{proof}
    For Newton strata, we apply the first paragraph of the proof of Proposition \ref{proposition: newton strata are well-positioned} using the diagram \ref{eq: G-shtuka, comm diagram}.

    For central leaves, we only need to deal with the case when $\Kf_p = K_p$ is a stablizer quasi-parahoric subgroup, since the case of general quasi-parahoric subgroup follows from this and Lemma \ref{lemma: central leaves disjoint union, different levels} and Remark \ref{rmk: central leaves disjoint union, different levels}.

    Now, assume $\Kf_p = K_p$. We use the notations from subsection \ref{subsubsec: set up PR}. Apply arguments in the proof of Proposition \ref{proposition: central leaves are well--positioned}, we see that $\CE^b \subset \varphi^{-1}(\CE^{b^{\dd}})$ is closed. Moreover, note that given any test object $Z = \Spec R \to (\wdh{\Shum{K}})^{\diamondsuit}$ that factors through the perfection of $\varphi^{-1}(\CE^{b^{\dd}})$, the crystalline tensors $(t_{\alpha, \crys}) \in \DD^{\sharp}(W(R))$ are geometrically constant, again due to \cite[Lemma 2]{hamacher2017almost}, thus $\CE^b \subset \varphi^{-1}(\CE^{b^{\dd}})$ is a union of connected components thus is  open. 
\end{proof}

 \subsection{Index set $B(G)$}

 Given any embedding of $F$-reductive groups $G\to G'$, we have a well-defined morphism $\rho: B(G)\to B(G')$. The Kottwitz map $B(\ast)\to \NE(\ast) \times \pi_1(\ast)_{\Sigma}$ is functorial on $F$-reductive groups $(\ast)$ (\cite[\S 1.9, 1.15]{rapoport1996classification}).

\subsubsection{Fiber discreteness}
 The main proposition we prove in this subsection is the following:
 
   \begin{proposition}\label{prop: BG main prop}
   	 If $b_1\leq b_2$ in $B(G)$, then $\rho(b_1)\leq\rho(b_2)$ in $B(G')$. Moreover, if $b_1\lneqq b_2$, then $\rho(b_1)\lneqq\rho(b_2)$. In particular, the fiber of $B(G) \to B(G')$ is discrete.
   \end{proposition}

   Note that the Bruhat order $b_1 \leq b_2$ in $B(\ast) \stackrel{(\nu_G, \kappa_G)}{\longrightarrow} \NE(\ast) \times \pi_1(\ast)_{\Sigma}$ forces $\kappa_G(b_1) = \kappa_G(b_2)$, it suffices to consider the map $\nu_G$. 
   
   The first part of the Proposition \ref{prop: BG main prop} comes from the functorialty of the Kottwitz map and the following proposition:
   \begin{lemma}{{\cite[Lemma 2.2]{rapoport1996classification}}}
   	  Given $[\nu_1]$, $[\nu_2]$ in $\NE(G)$, $[\nu_1]\preceq[\nu_2]$ if and only if for all representations $\rho: G\to\GL(V)$, $\rho_*([\nu_1])\preceq\rho_*([\nu_2])$.
   \end{lemma}
   To prove the second part of Proposition \ref{prop: BG main prop}, it suffices to consider the case over $\ovl{F}$. Thus we base change everything to $\ovl{F}$, and assume $G$ is split. Take a maximal torus $T$ in $G$ (resp. $T'$ in $G'$) with $T\subset T'$, then $\NE(G) = (X_*(T)_{\Q}/W)$ (resp. $\NE(G') = (X_*(T')_{\Q}/W')$), where $W$ (resp. $W'$) is the absolute Weyl group of $G$ (resp. $G'$). Also, we will only consider the algebraic representations in this section, i.e. an algebraic representation is an algebraic group homomorphism $\rho: G\to\GL(V)$ over $\ovl{F}$.
   \begin{lemma}\label{lem: reduce repr to GL with inequality}
   	 Given $[\nu_1]\precneqq[\nu_2]$ in $\NE(G)$, then there exists a finite dimensional representation $\rho: G\to\GL(V)$ such that $\rho_*([\nu_1])\precneqq\rho_*([\nu_2])$ in $\NE(\GL(V))$, i.e. there does not exsit an $\omega$ in the Weyl group of $\GL(V)$ such that $\omega\rho_*([\nu_1])=\rho_*([\nu_2])$.
   \end{lemma}
   \begin{proof}
   	 Let $\nu_1$, $\nu_2$ the dominant representatives of $[\nu_1]$ and $[\nu_2]$ in $X_*(T)$ respectively, $\nu_1\lneqq\nu_2$ implies that $\exists\lambda$ an integral dominant weight with $\langle\nu_2-\nu_1, \lambda\rangle>0$. By highest weight theory, there exsits an irreducible representation $V=V(\lambda)$ of $G$ with the given highest weights $\lambda$. Assume $\rho: G\to\GL(V)$ has weight $\lambda=\lambda_1, \lambda_2, \dots, \lambda_n$, then $\lambda_1>\lambda_j$, $2\leq j\leq n$. Up to conjugation, assume $T\subset G$ maps into the diagonal torus in $\GL(V)$. For $t\in T$, $t_i=\langle t, \lambda_i\rangle$, then $\rho(t)=\mathrm{diag}\langle t_1, t_2, \dots, t_n\rangle$. Let $x\in\Gm$, for $i=1, 2$, $\rho(\nu_i(x))=\mathrm{diag}( x^{\langle \nu_i, \lambda_1\rangle}, x^{\langle \nu_i, \lambda_2\rangle}, \dots, x^{\langle \nu_i, \lambda_n\rangle} )$, $\langle \nu_2, \lambda_1\rangle$ is the unique largest element in all $\langle \nu_i, \lambda_j\rangle$, $i\in\{1, 2\}$, $j\in\{1, 2, \dots, n\}$. The lemma follows from the fact that the Weyl group of $\GL(V)$ acts on the diagonal entries as permutations.
   \end{proof}
   Proposition \ref{prop: BG main prop} follows from the following Proposition:
   \begin{proposition}
   	Let $\rho: G \to G'$ be a group embedding. Given $[\nu_1]$, $[\nu_2]$ in $\NE(G)$, if If $[\nu_1]\precneqq[\nu_2]$, then $\rho_*([\nu_1])\precneqq\rho_*([\nu_2])$.
   \end{proposition}
   \begin{proof}
   	 Let $V$ be the representation defined in Lemma \ref{lem: reduce repr to GL with inequality}, and let $V'=\Ind^{G'}_G(V)$ be the induced $G'$-representation (see \cite[\S 3.3]{jantzen2003representations} for definition) which satisfies the Frobenious reciprocity:
   	 \[ \Hom_G(V, \Res^{G'}_G(V'))\cong\Hom_{G'}(V', V') \]
   	 As a representation of $G$, $V\subset \Res^{G'}_G(V')$. Consider the minimal $G'$-invariant subrepresentation $W$ of $V'$ which includes $V$. $W$ is a finite dimensional subspace by \cite[(Local finiteness) 2.13]{jantzen2003representations}. Since any finite dimensional algebraic representation of a reductive group over an algebraically closed field is semisimple, $W\cong V_0\oplus V$ is a semisimple $G$-module, $G \to G' \to \GL(W)$ factors through $G \to \GL(V_0) \times \GL(V) \to \GL(W)$. Fix the diagonal tori $T_{V_0}$, $T_V$, $T_W$ and the standard uppertriangular Borel groups in diagonal embedding $\GL(V_0)\times\GL(V)\to\GL(W)$. Assume $\dim W=\dim V_0+\dim V=n+m$. Any cocharactor $\nu$ in $X_*(T_W)$ can be regarded as a vector in $W$. In this way, pick the dominant representatives $\rho_1$ and $\rho_2$ in $\rho_*([\nu_1])$ and $\rho_*([\nu_2])$ respectively, write $\rho_1=(a_1, \dots, a_n, b_1, \dots, b_m)$, $\rho_2=(a_1', \dots, a_n', b_1', \dots, b_m')$, since $\rho_1\leq\rho_2$, these two tuples of numbers satisfies the following conditions:
   	 	\begin{compactitem}\label{numbercondition}
   	 		\item $a_1\geq a_2\geq\dots\geq a_n$, $a_1'\geq a_2'\geq\dots\geq a_n'$, $b_1\geq b_2\geq\dots\geq b_m$, $b_1'\geq b_2'\geq\dots\geq b_m'$
   	 		\item $\Sigma_{i=1}^s a_i\leq \Sigma_{i=1}^s a_i'$ for all $s\in\{1, 2, \dots, n\}$, $\Sigma_{j=1}^r b_j\leq \Sigma_{j=1}^r b_j'$ for all $r\in\{1, 2, \dots, m\}$
   	 	\end{compactitem}
   	   By the choice of $V$ in Lemma \ref{lem: reduce repr to GL with inequality}, $b_1<b_1'$. To prove $\rho_1\neq\rho_2$ in $X_*(T_W)_{\Q, \dom}$, it's equivalent to prove there is no permutation in $S_{n+m}$ acting on the index sending $(a_1, \dots, a_n, b_1, \dots, b_m)$ to $(a_1', \dots, a_n', b_1', \dots, b_m')$. This is a pure elementary combinatoric proposition. To prove this, we do induction on the size of $n+m$. When $n = 0$, then $W = V$, $\rho_1 \lneqq \rho_2$ by construction, thus we assume $n \geq 1$.
   	   
   	If $n+m=2$ ($n = m = 1$), then $a_1 \leq a_1'$, $b_1 < b_1'$. If $(a_1, b_1) = (a_1', b_1')$ up to permuation, then $a_1 = b_1'$, $a_1' = b_1$, but then $a_1 \leq a_1' = b_1 < b_1'$, absurd.
   	   
       Now assume the lemma is true for $n+m=k$, $k\geq 2$, consider the case $n+m=k+1$:
       
       If $a_1=a_1'$, we can assume $a_1$ is fixed under the permutation. Remove $a_1$ and $a_1'$, the remaining tuples of numbers still satisfy the condition \ref{numbercondition}, then use the induction step on $n+m=k$. Therefore we assume $a_1<a_1'$. $b_1'>b_1\geq b_2\geq\dots\geq b_m$, so we must have $a_1\geq b_1'$, otherwise there is no element in $(a_1, \dots, a_n, b_1, \dots, b_m)$ can match $b_1'$. Now $a_1'>a_1\geq a_2\geq\dots\geq a_n$, $a_1'>a_1\geq b_1'>b_1\geq b_2\geq\dots\geq b_m$, no element in $(a_1, \dots, a_n, b_1, \dots, b_m)$ can match $a_1'$.
   \end{proof}

\subsubsection{Quasi-finiteness}
Let us prove another useful lemma:
\begin{lemma}\label{lemma: B(G) quasi-finite}
    $B(G) \to B(G')$ is quasi-finite, i.e., given any $[g'] \in B(G')$, $\rho^{-1}([g']) \subset B(G)$ is either empty or a finite set.
\end{lemma}
\begin{proof}
    Due to \cite[Theorem 1.15]{rapoport1996classification}, we have following functorial diagram with exact rows (in the sense of pointed sets)
\[\begin{tikzcd}
	{H^1(F, G)} & {B(G)} & {\NE(G)} \\
	{(\pi_1(G)_{\Gamma})_{\mathrm{tors}}} & {\pi_1(G)_{\Gamma}} & {\pi_1(G)^{\Gamma}\otimes \Q}
	\arrow[from=1-1, to=1-2]
	\arrow[from=1-1, to=2-1]
	\arrow[from=1-2, to=1-3]
	\arrow[from=1-2, to=2-2]
	\arrow[from=1-3, to=2-3]
	\arrow[from=2-1, to=2-2]
	\arrow[from=2-2, to=2-3]
\end{tikzcd}\]
    To show $B(G) \to B(G')$ is quasi-finite, it suffices to show $\NE(G) \to \NE(G')$ is quasi-finite and $H^1(F, G)$ is a finite set. Recall that $\NE(G) = (X_*(T)_{I, \Q}^+)^{\sigma}$, since $X_*(T)_{I, \Q} = X_*(T)^I_{\Q}$, $X_*(T)_{I, \Q} \to X_*(T')_{I, \Q}$ is injective, then $X_*(T)_{I, \Q}^+ \to X_*(T')_{I, \Q}^+$ is quasi-finite, thus $\NE(G) \to \NE(G')$ is quasi-finite. On the other side, $H^1(F, G)$ is a finite set, which is moreover bijective to the torsion part of the finitely generated algebraic fundamental group $(\pi_1(G)_{\Gamma})_{\mathrm{tors}}$.
\end{proof}

\begin{lemma}\label{lem: B(G) bijection}
     Let $\mu_h$ be a cocharacter of $G_h$, thus of $L$, then $B(G_h, \lrbracket{\mu_h}) \to B(L, \lrbracket{\mu_h})$ is a bijection.
\end{lemma}
\begin{proof}
  Since $L^{\ad} \cong G_h^{\ad} \times G_l^{\ad}$, the projection of $\mu_h$ to $G_l$ is trivial, thus
      \[ B(L^{\ad}, \lrbracket{\mu_h^{\ad}}) = B(G_h^{\ad}, \lrbracket{\mu_h}) \times B(G_l^{\ad}, \lrbracket{0}) = B(G_h^{\ad}, \lrbracket{\mu_h}) \]
      Recall that \cite[\S 6.5.1]{kottwitz1997isocrystals} says $B(G, \lrbracket{\mu}) \cong B(G^{\ad}, \lrbracket{\mu^{\ad}})$ for general reductive groups and cocharacters, the natural inclusion $B(G_h, \lrbracket{\mu_h}) \to B(L, \lrbracket{\mu_h})$ induces an isomorphism:
\[\begin{tikzcd}
	{B(G_h, \lrbracket{\mu_h})} & {B(L, \lrbracket{\mu_h})} \\
	{B(G_h^{\ad}, \lrbracket{\mu_h^{\ad}})} & {B(L^{\ad}, \lrbracket{\mu_h^{\ad}})}
	\arrow["\cong", from=1-1, to=1-2]
	\arrow["\cong"', from=1-1, to=2-1]
	\arrow["\cong"', from=1-2, to=2-2]
	\arrow["\cong", from=2-1, to=2-2]
\end{tikzcd}\]
     Since $[b_h], [b_h'] \in B(G_h, \lrbracket{\mu_h})$ thus $\in B(L, \lrbracket{\mu_h})$, and since $[b_h] = [b_h'] \in B(L)$, thus $[b_h] = [b_h'] \in B(G_h, \lrbracket{\mu_h})$. 
\end{proof}

\subsection{Index set $C(\GG)$}

\subsubsection{Adjoint action}

We want to consider Lemma \ref{lem: B(G) bijection} with $B(G, \lrbracket{\mu})$ replaced with $C(\GGc, \lrbracket{\mu})$. Here we use the following notations: let $\mu$ be a cocharacter of $G$, recall that we have defined $C(\GGc, \lrbracket{\mu}) = \KKc\Adm(\lrbracket{\mu})\KKc/\KKc_{\sigma}$, due to the compatibility of Newton maps (for example, see \cite[\S 1.3.1]{shen2021ekor}), $C(\GGc, \lrbracket{\mu})$ maps surjectively to $B(G, \lrbracket{\mu})$ under the natural projection $C(\GGc) \to B(G)$.

The following lemma is indeed \cite[Lemma 5.5.2]{shen2021ekor}, but we might need to replace $C(\GGc) \to C(\GG^{\ad, \circ})$ in that statement with $C(\GGc, \lrbracket{\mu}) \to C(\GG^{\ad, \circ}, \lrbracket{\mu^{\ad}})$ for some $\mu$. Since if we do not restrict the index sets using cocharacter, when we take $G = T$ being a torus, then $C(\mathcal{T}^{\circ})$ is definitely non-trivial.
\begin{lemma}{{\cite[Lemma 5.5.2]{shen2021ekor}}}\label{lem: ekor central leaves}
    Assume the center $Z$ of $G$ is connected. Let $\mu$ be a cocharacter of $G$. When the closure of $Z$ in $\GGc$ has connected fibers, the natural projection $C(\GGc, \lrbracket{\mu}) \to C(\GG^{\ad, \circ}, \lrbracket{\mu^{\ad}})$ is bijective.
\end{lemma}
\begin{proof}
    First of all, by assumption, the center $Z$ of $G$ is connected, then $G(\bF) \to G^{\ad}(\bF)$ is surjective since $H^1(\bF, Z)$ is trivial, thus $\GGc \to \GG^{\ad, \circ}$ is surjective, see \cite[Proposition 2.59]{mao2025boundary}. Given $\llbracket g \rrbracket, \llbracket g' \rrbracket \in C(\GGc, \lrbracket{\mu})$ with same image in $C(\GG^{\ad, \circ}, \lrbracket{\mu^{\ad}})$, we can adjust $\KKc_{\sigma}$ such that $g$ and $g'$ have same image in $G^{\ad}(\bF)$. We write $g' = z g$, $z \in Z(\bF)$. We claim $z \in \KKc$: Since $C(\GGc, \lrbracket{\mu}) = \KKc\Adm(\lrbracket{\mu})\KKc/\KKc_{\sigma}$, since the image of $\Adm(\lrbracket{\mu})$ under $\kappa_G$ is a single element, thus $\kappa_G(g) = \kappa_G(gz)$, $z \in \ker \kappa_G$. Due to \cite[Proposition 2.84]{mao2025boundary}, $z \in \KKc$. Apply Lemma \ref{lemma: Z connected then action stablize central leaf}, we have $z \in Z(\bF) \cap \KKc = Z(\bF)_0$. Since $H^1(\lrangle{\sigma}, Z(\bF)_0) = 1$, thus we can write $z = z_1 \sigma(z_1)^{-1}$ for some $z_1 \in Z(\bF)_0$, thus $\llbracket gz \rrbracket = \llbracket z_1g\sigma(z_1)^{-1} \rrbracket = \llbracket g \rrbracket \in C(\GGc, \lrbracket{\mu})$, $C(\GGc, \lrbracket{\mu}) \to C(\GG^{\ad, \circ}, \lrbracket{\mu^{\ad}})$ is injective. The surjectivity comes from the fact that $Z$ is connected and $H^1(\bF, Z)$ is trivial.
\end{proof}
\begin{lemma}\label{lemma: Z connected then action stablize central leaf}
   Assume the center $Z$ is connected. The closure of $Z$ in $\GGc$ has connected fibers if and only if $Z(\bF)_0 = Z(\bF) \cap \KKc$.
\end{lemma}
\begin{proof}
     Recall that $Z(\bF) \cap \KKc = Z(\bF)^1 \cap \KKc$, for example, see \cite[Proposition 2.84]{mao2025boundary}. Let $\ZZ'$ be the integral closure of $Z$ in $\GGc$, then $\ZZ'(\OO_{\bF}) = Z(\bF)^1 \cap \KKc$. Take the smoothening $\wdt{\ZZ}\to \ZZ'$, $\wdt{\ZZ}$ is the unique smooth model of $Z$ such that $\wdt{\ZZ}(\OO_{\bF}) = \ZZ'(\OO_{\bF})$. Note that the unique parahoric group scheme $\ZZ^{\circ} $ of $Z$ is also the neutral connected component of $\wdt{\ZZ}$. $\ZZ'$ has connected fibers if and only if $\wdt{\ZZ}$ has connected fiber, if and only if $\wdt{\ZZ} = \ZZ^{\circ}$. Since both $\ZZ^{\circ} \subset \wdt{\ZZ}$ are smooth, $\ZZ^{\circ} = \wdt{\ZZ}$ if and only if $\ZZ^{\circ}(\OO_{\bF}) = \wdt{\ZZ}(\OO_{\bF})$, if and only if $Z(\bF)_0 = Z(\bF)^1 \cap \KKc$.
\end{proof}

\begin{corollary}\label{cor: GGc to LLc injective when Z connected}
   Let $\mu_h$ be a cocharacter of $G_h$. In either one of the following case, $C(\GGc_h, \lrbracket{\mu_h}) \to C(\LLc, \lrbracket{\mu_h})$ is injective.
    \begin{enumerate}
        \item The closure of $Z_h$ (the center of $G_h$) in $\GGc_h$ has connected fibers.
        \item $1 \to G_h \to L \to G_l \to 1$ splits as $L \cong G_h \times G_l$.
    \end{enumerate}
\end{corollary}
\begin{proof}
  \begin{enumerate}
      \item $L^{\ad} \cong G_h^{\ad} \times G_l^{\ad}$ induces $\LL^{\ad, \circ} \cong \GG_h^{\ad, \circ} \times \GG_l^{\ad, \circ}$. In particular, $C(\GG_h^{\ad, \circ}) \to C(\LL^{\ad, \circ})$ is injective. Given $\llbracket g_1 \rrbracket, \llbracket g_2 \rrbracket \in C(\GGc_h)$ such that $\llbracket g_1 \rrbracket = \llbracket g_2 \rrbracket \in C(\LLc)$, then $\llbracket g_1 \rrbracket = \llbracket g_2 \rrbracket \in C(\GG_h^{\ad, \circ})$. Due to Lemma \ref{lem: ekor central leaves} and \ref{lemma: Z connected then action stablize central leaf}, $\llbracket g_1 \rrbracket = \llbracket g_2 \rrbracket \in C(\GGc_h)$.
      \item We simply use the fact that $\LLc = \GGc_h \times \GGc_l$.
  \end{enumerate}
    
\end{proof}

%% file: sections/igusa_varieties.tex
In this section we prove the well-position of Igusa varieties. Since Igusa varieties are not subschemes of $\Shum{K, \Bar{s}}$, the definition of well-positioned subschemes will be modified, see Subsection \ref{sec: generalized definition of well position}. We mainly follow the notations from \cite{caraiani2017generic}, \cite{caraiani2024generic}, \cite{santos2023generic} and \cite{hamacher2019adic}. In this section, we work with Kisin-Pappas integral models, and we work with stablizer quasi-parahoric group $K_p$.

\subsection{Siegel-type Igusa varieties}

We recall constructions and theorems in \cite{caraiani2017generic}.

Let $b' \in \GSp(\bQ)$, $\CE^{b'}$ be the central leaf on $\Shum{K^{\dd}, \Bar{s}}$, we assume it is nonempty. Given any point $x \in \CE^{b'}$, $(\ab[p^{\infty}]_{\bar{x}}, \lambda_{\ab, \Bar{x}})$ is isomorphic to $(\XXp, \lambda_{\XXp})$ up to some field extension for some geometric point $\Bar{x}$ over $x$, where $\XXp$ is a $p$-divisible group over $k$ with principal polarization $\lambda_{\XXp}: \XXp\rightiso (\XXp)^{\vee}$, $(\XXp, \lambda_{\XXp})$ is determined by $b'$. Let us denote by $\XXp_{\ab}:= \ab[p^{\infty}]|_{\CE^{b'}}$. We assume $b'$ is completely divisible, then $\XXp$ is completely slope divisible and the filtration $0 \subset \XXp_{\leq \lambda_1} \subset \cdots \XXp_{\leq \lambda_r} = \XXp$ satisfies $\XXp_{i} \cong \XXp_{j}^{\vee}$ for those $\lambda_i + \lambda_j = 1$. Then $\XXp_{\ab}$ is completely slope divisible over the regular scheme $\CE^{b'}$, see \cite[Lemma 2.3.4]{caraiani2024generic}. In particular, there exists a slope filtration $0 \subset \XXp_{\ab, \leq \lambda_1} \subset \cdots \XXp_{\ab, \leq \lambda_r} = \XXp_{\ab}$ satisfies $\XXp_{\ab, i} \cong \XXp_{\ab, j}^{\vee}$ for those $\lambda_i + \lambda_j = 1$.
Let $\AUT(\XXp)$ be the group scheme parametrizing automorphisms of $\XXp$ that are compatible with the given polarization. $\AUT(\XXp)$ is a highly non-reduced scheme. The non-reduceness measures the failure of $\XXp$ being isoclinic. Its group of components is a profinite group $\Gamma_{\XXp} = \AUT(\XXp)(k)$.

\begin{definition-proposition}{{\cite[\S 4]{caraiani2017generic}}}\label{def: Igusa Siegel}\leavevmode
    \begin{enumerate}
        \item Let $\IGUSA_{b'} \to \CE^{b'}$ be the functor parametrizing isomorphisms $(\XXp, \lambda_{\XXp}) \cong (\ab[p^{\infty}], \lambda_{\ab})$, it is represented by a $\AUT(\XXp)$-torsor over $\CE^{b'}$.
        \item Let $\Ig_{b'} \to \CE^{b'}$ be the functor parametrizing isomorphisms $(\rho_i)_{i=1}^r$ $\rho_i: \XXp_i \cong \XXp_{\ab, i}$ such that $(\rho_i)_{i=1}^r$ commutes with polarizations up to a constant in $\Z_p^{\times}$ independent of $i$. $\Ig_{b'} \to \CE^{b'}$ is a pro-\'etale cover with Galois group $\Gamma_{\XXp}:=\Aut(\XXp)$.
        \item Let $\Ig_{m, b'} \to \CE^{b'}$ be the functor parametrizing isomorphisms $(\rho_{m, i})_{i=1}^r$: $\rho_{m, i}: \XXp_i[p^m] \cong \XXp_{\ab, i}[p^m]$ such that $(\rho_{m, i})_{i=1}^r$ commute with polarizations up to a constant in $(\Z_p/p^m\Z_p)^{\times}$ independent of $i$ and such that $(\rho_{m, i})_{i=1}^r$ lift to arbitary $m' \geq m$ keeping the polarizations. $\Ig_{m, b'} \to \CE^{b'}$ is a finite \'etale cover with Galois group $\Gamma_{m, \XXp}$, where $\Gamma_{m, \XXp}$ is the finite group of automorphisms of $\XXp[p^m]/k$ that lifts to $\XXp/k$ keeping the polarization. $\Gamma_{\XXp} = \prolim \Gamma_{m, \XXp}$.
        \item The natural morphism $\IGUSA_{b'} \to \Ig_{b'}$ is an isomorphism after perfection since the slope filtration splits canonically over a perfect base.
    \end{enumerate}
\end{definition-proposition}

\begin{remark}
    Another way to view $\IGUSA_{b'}$ is that it parametrizes quasi-isogenies of $(\XXp, \lambda_{\XXp}) \cong (\ab[p^{\infty}], \lambda_{\ab})$ up to $p$-power isogenies of $(\ab[p^{\infty}], \lambda_{\ab})$, see \cite[Lemma 4.3.4]{caraiani2017generic}, and in this way we can easily see that $\IGUSA_{b'}$ is a perfect scheme (\cite[Corollary 4.3.5]{caraiani2017generic}).
\end{remark}

    Fix a cusp label representative $\Phi^{\dd} \in \Cusp_{K^{\dd}}(G^{\dd}, X^{\dd})$. For simplicity, we take $K^{\dd, p}$ as principal congruence subgroups such that $\Zb^{\bigsur, \dd} \to \Zb^{\dd}$ is an isomorphism for each $\Phi^{\dd}$. From now on, we fix a $\Phi^{\dd}$ (and omit it in the index) and consider the corresponding boundary stratum.

    Note that central leaves $\CE^{b'}$ are well-positioned with respect to $\lrbracket{\CE^{b', \natural}(\Phi^{\dd})}_{\Phi^{\dd}}$, each $\CE^{b', \natural}_{\Zb^{\dd}}:=\CE^{b', \natural}_{\Zb^{\dd}(\Phi^{\dd})} :=\CE^{b', \natural}(\Phi^{\dd})$ is a single central leaf $\CE^{b_h'}$ over the Siegel Shimura variety $\Zb^{\dd}(\Phi^{\dd})_{\Bar{s}} = \Shum{K_{\Phi, h}^{\dd}, \Bar{s}}$, see \cite[\S 3.4]{lan2018compactifications}.

Over $C^{\dd}$, there is a semiabelian scheme $G$: $1 \to \mT \to \mG \to \mB \to 1$. Here $T$ is a constant torus with cocharacter group $Y$, and $\mB$ is the (pullback of) universal abelian scheme over $\Zb^{\dd} = \Shum{K_{\Phi, h}^{\dd}}$. Let us denote $\mG[p^{\infty}]^{\circ}$ by $\mH$. Since $\CE^{b'_h}$ and $C^{\dd} \to \Zb^{\dd}$ are smooth, and $\CE^{b', \natural}_{C^{\dd}}$ is the pullback of $\CE^{b'_h}$ along $C^{\dd} \to \Zb^{\dd}$, then $\CE^{b', \natural}_{C^{\dd}}$ is smooth, and $\mB[p^{\infty}]$ and $\mH$ are geometrically fiberwise constant over $\CE^{b', \natural}_{C^{\dd}}$. As a result, $\mB[p^{\infty}]$ and $\mH$ are completely slope divisible over $\CE^{b', \natural}_{C^{\dd}}$. 

As in \cite[\S 3.2]{caraiani2024generic}, we have a (connected) $p$-divisible group $\ab[p^{\infty}]^{\circ} \to \Shumc{K^{\dd}, \Sigma^{\dd}}$ coming from the universal semiabelian scheme, and $\ab[p^{\infty}]^{\circ} \cong \XXp^{\circ}$ over geometric points of $\CE^{b', \tor}$. In fact, over $W^{\dd}$, by considering the formal completions along identities, we have $\hat{\ab} \cong \hat{\mG}$ (by construction), which implies that $\ab[p^{\infty}]^{\circ} \cong \mH$. Over such $W^{\dd}$, the biconnected part $\ab[p^{\infty}]^{\circ}_{(0, 1)} = \ab[p^{\infty}]^{\circ}/\ab[p^{\infty}]^{\circ}_{\mu}$ together with the principal polarization (which extended the one on interior), it is isomorphic to the pullback of the restriction of the bi-connected part $\mB[p^{\infty}]_{(0, 1)} $. When restricted to the pullback of the central leaf $\CE^{b', \tor}$, they are geometrically fiberwise isomorphic to $\XXp_{(0, 1)}= \oplus_{0<\lambda_i<1} \XXp_{i}$. Since $\CE^{b'}$ is a regular scheme, then $\CE^{b', \tor}$ is a regular scheme, $\ab[p^{\infty}]^{\circ}$ is completely slope divisible over $\CE^{b', \tor}$, it admits a slope filtration with graded pieces $\ab_i[p^{\infty}]$ which are iso-clinic of slope $\lambda_i$.

\begin{definition-proposition}\leavevmode
\begin{enumerate}
    \item Let $\Ig_{m, b'}^{\tor} \to \CE^{b', \tor}$ be the functor that classifies the data $(\rho_{m, i})$ (for those $\lambda_i > 0$), where $\rho_{m, i}: \XXp_i[p^m] \cong \ab_i[p^m]$ such that $(\rho_{m, i})$ commute (explicitly, for those $\rho_{m, i}$, $\rho_{m, j}$ with $\lambda_i + \lambda_j = 1$) with polarizations up to a constant in $(\Z_p/p^m\Z_p)^{\times}$ independent of $i$ and such that $(\rho_{m, i})$ lift to arbitary $m' \geq m$ keeping the polarizations. Let $\Ig_{b'}^{\tor} = \prolim_{m} \Ig_{m, b'}^{\tor}$. Then \cite[Theorem 3.2.4]{caraiani2024generic} says that $\Ig_{b'}^{\tor} \to \CE^{b', \tor}$ (resp. $\Ig_{m, b'}^{\tor} \to \CE^{b', \tor}$) is a $\Gamma_{\XXp}$-(resp. $\Gamma_{m, \XXp}$-)torsor.
    \item Let $\IGUSA_{b'}^{\tor} \to \CE^{b', \tor}$ be the functor that classifies the isomorphism $\rho: \XXp[p^{\infty}]^{\circ} \cong \ab[p^{\infty}]^{\circ}$ such that the induced isomorphism $\rho: \XXp[p^{\infty}]^{\circ}_{(0, 1)} \cong \ab[p^{\infty}]^{\circ}_{(0, 1)}$ commutes with the polarization up to a constant in $(\Z_p/p^m\Z_p)^{\times}$. Then \cite[Theorem 3.2.8]{caraiani2024generic} says that  $\IGUSA_{b'}^{\tor} \to \CE^{b', \tor}_{\perf}$ is a $\Gamma_{\XXp}$-torsor (the \'etale part is uniquely determined by the multiplicative part).
\end{enumerate} 
\end{definition-proposition}

 Due to the modular interpretations in \cite[Definition 3.2.5, 3.2.9]{caraiani2024generic}, one has a natural morphism $\IGUSA_{b'}^{\tor} \to \Ig_{b'}^{\tor}$ which extends $\IGUSA_{b'} \to \Ig_{b'}$. It induces an isomorphism $\IGUSA_{b'}^{\tor} \rightiso \Ig_{b'}^{\tor, \perf}$, see \cite[Theorem 3.2.8]{caraiani2024generic} and arguments below it.

We are able to determine the boundary structures of Igusa varieties. They are described in \cite[Theorem 3.2.6]{caraiani2024generic}, we summarize it as follows (the language here is slightly different from the cited article). 
\begin{definition-proposition}\label{def: Igusa, 1-motive, Siegel}
 One can similarly define a functor (see the lines above \cite[Theorem 3.2.6]{caraiani2024generic}, where $\CE^{b', \natural}_{C^{\dd}}$ is denoted by $C_{Z} \times_{\Shum{Z}} \CE_{Z}^{\XXp}$) which is representated by a $\Gamma_{m, \XXp^{\circ}}$-torsor $\Ig^{\natural}_{m, b', C^{\dd}} \to \CE^{b', \natural}_{C^{\dd}}$: it classifies $(\rho_{m, i})_i$ (for those $\lambda_i > 0$): $\rho_{m, i}: \XXp_i[p^m] \cong \mH_{i}[p^m]$ such that $(\rho_{m, i})$ commute with polarizations up to a constant in $(\Z_p/p^m\Z_p)^{\times}$ independent of $i$ and such that $(\rho_{m, i})$ lift to arbitary $m' \geq m$ keeping the polarizations.  Let $\Ig^{\natural}_{b', C^{\dd}} = \prolim_{m} \Ig^{\natural}_{m, b', C^{\dd}}$.   
\end{definition-proposition}
\begin{remark}\label{rmk: Igusa, 1-motive, Siegel}
    Recall that $\XXp$ has a principal polarization on it, fix a constant in $\Z_p^{\times}$, we have $\Gamma_{(m), \XXp} = \Gamma_{(m), \XXp^{\circ}}$, that is to say, the isomorphism between the \'etale part determines and is determined by the isomorphism between the multiplicative part. 
\end{remark}

   Fix a $\Phi^{\dd}$. Consider the $\Gamma_{(m), \XXp}$-torsor $\Ig^{\natural}_{(m),  b', C^{\dd}} \to \CE^{b', \natural}_{C^{\dd}}$ defined above. Let $\Ig^{\natural}_{(m), b', \Xi^{\dd}} \to \CE^{b', \natural}_{\Xi^{\dd}}$ be the pullback of $\Ig^{\natural}_{(m), b', C^{\dd}} \to \CE^{b', \natural}_{C^{\dd}}$ along $\Xi^{\dd} \to C^{\dd}$. We similarly deonte the pullbacks of $\Ig^{\natural}_{(m), b', C^{\dd}} \to \CE^{b', \natural}_{C^{\dd}}$ along $\Xi^{\dd}(\sigma^{\dd}) \to C^{\dd}$, $\Xi^{\dd}_{\sigma^{\dd}} \to C^{\dd}$, $\ovl{\Xi}^{\dd}_{\Sigma^{\dd}(\Phi^{\dd})} \to C^{\dd}$ by $\Ig^{\natural}_{(m), b', \Xi^{\dd}(\sigma^{\dd})}$, $\Ig^{\natural}_{(m), b', \Xi^{\dd}_{\sigma^{\dd}}}$, $\Ig^{\natural}_{(m), b', \ovl{\Xi}^{\dd}_{\Sigma^{\dd}(\Phi^{\dd})}}$ respectively.

    For any $\sigma^{\dd} \in \Sigma^{\dd}(\Phi^{\dd})^+$, we denote by $\Ig_{(m), b', \mathcal{Z}([\Phi^{\dd}, \sigma^{\dd}])}$ the boundary stratum of $\Ig_{(m), b'}^{\tor}$ over $\CE^{b', \tor}_{\mathcal{Z}([\Phi^{\dd}, \sigma^{\dd}])}$.

    Let $(\Ig_{(m), b'}^{\tor})^{\wedge}_{\mathcal{Z}([\Phi^{\dd}, \sigma^{\dd}])}$ (resp. $\Ig^{\natural}_{(m), b', \mathfrak{X}^{\dd}_{\sigma^{\dd}}}$) be the completion of $\Ig_{(m), b'}^{\tor}$ (resp. $\Ig^{\natural}_{(m), b', \Xi^{\dd}(\sigma^{\dd})}$) along $\Ig_{(m), b', \mathcal{Z}([\Phi^{\dd}, \sigma^{\dd}])}$ (resp. $\Ig^{\natural}_{(m), b', \Xi^{\dd}_{\sigma^{\dd}}}$). Since all the schemes involved are finite type over locally noetherian base, and $\Xi^{\dd}(\sigma^{\dd}) \to C^{\dd}$ is a finite type morphism, one can also regard $\Ig^{\natural}_{(m), b', \mathfrak{X}^{\dd}_{\sigma^{\dd}}}$ as the pullback of $\Ig^{\natural}_{(m), b', C^{\dd}}$ along $\mathfrak{X}^{\dd}_{\sigma^{\dd}} \to C^{\dd}$.
    
    Let $(\Ig_{(m), b'}^{\tor})^{\wedge}_{\Phi^{\dd}, \sigma^{\dd}}$ (resp. $(\Ig_{(m), b'}^{\tor})^{\wedge}_{\Phi^{\dd}}$) be the completion of $\Ig_{(m), b'}^{\tor}$ along
    \[ \bigcup_{\tau^{\dd} \in \Sigma^{\dd}(\Phi^{\dd})^+,\ \Bar{\tau}^{\dd} \subset \Bar{\sigma}^{\dd}}  \Ig_{(m), b', \mathcal{Z}([\Phi^{\dd}, \tau^{\dd}])},\quad  (\textit{resp.}\ \bigcup_{[\tau^{\dd}] \in \Sigma^{\dd}(\Phi^{\dd})^+/\Gamma_{\Phi^{\dd}}}  \Ig_{(m), b', \mathcal{Z}([\Phi^{\dd}, \tau^{\dd}])}), \]
    and let $\Ig^{\natural}_{(m), b', \mathfrak{X}^{\dd, \circ}_{\sigma^{\dd}}}$ (resp.  $\Ig^{\natural}_{(m), b', \mathfrak{X}^{\dd}_{\Sigma^{\dd}(\Phi^{\dd})}}$) be the completion of $\Ig^{\natural}_{(m), b', \Xi^{\dd}(\sigma^{\dd})}$ (resp. $\Ig^{\natural}_{(m), b', \ovl{\Xi}^{\dd}_{\Sigma^{\dd}(\Phi^{\dd})}}$) along 
    \[ \bigcup_{\tau^{\dd} \in \Sigma^{\dd}(\Phi^{\dd})^+,\ \Bar{\tau}^{\dd} \subset \Bar{\sigma}^{\dd}}  \Ig^{\natural}_{(m), b', \Xi^{\dd}_{\tau^{\dd}}},\quad  (\textit{resp.}\ \bigcup_{\tau^{\dd} \in \Sigma^{\dd}(\Phi^{\dd})^+} \Ig^{\natural}_{(m), b', \Xi^{\dd}_{\tau^{\dd}}}).  \]

        Since $\mathfrak{X}_{\sigma}^{\circ} \hookrightarrow \mathfrak{X}$ is an open covering when $[\sigma]$ runs over $\Lambda_{\Phi, K}\backslash \Sigma(\Phi)^+$, it follows from \cite[Theorem 3.2.6]{caraiani2024generic} and Lemma \ref{lemma: limit of well-positioned} that
  
    \begin{proposition}\label{prop: well-position, Igusa, Siegel case}

        $\Ig_{(m), b'} \to \CE^{b'}$ (resp. $\Ig_{b'} \to \CE^{b'}$) is well-positioned with respect to $\lrbracket{\Ig^{\natural}_{(m), b', C^{\dd}(\Phi^{\dd})} \to \CE^{b', \natural}_{C^{\dd}(\Phi^{\dd})}}_{[\Phi^{\dd}]}$ (resp. $\lrbracket{\Ig^{\natural}_{b', C^{\dd}(\Phi^{\dd})} \to \CE^{b', \natural}_{C^{\dd}(\Phi^{\dd})}}_{[\Phi^{\dd}]}$) in the sense of Definition \ref{def: well-positioned over well-positioned}. Note that $\Ig_{b'} = \varprojlim_m \Ig_{m, b'}$, $\Ig^{\natural}_{b', C^{\dd}} = \varprojlim_m \Ig^{\natural}_{m, b', C^{\dd}}$.
    \end{proposition}

    Due to Proposition \ref{prop: wdt Y satisfy same prop},
    \begin{corollary}
        The canonical isomorphism $\Ig_{(m), b', \mathcal{Z}([\Phi^{\dd}, \sigma^{\dd}])} \cong \Ig^{\natural}_{(m), b', \Xi^{\dd}_{\sigma^{\dd}}}$ extends to isomorphisms of formal schemes:
        \[(\Ig_{(m), b'}^{\tor})^{\wedge}_{\mathcal{Z}([\Phi^{\dd}, \sigma^{\dd}])} \cong \Ig^{\natural}_{(m), b', \mathfrak{X}^{\dd}_{\sigma^{\dd}}},\quad (\Ig_{(m), b'}^{\tor})^{\wedge}_{\Phi^{\dd}, \sigma^{\dd}} \cong  \Ig^{\natural}_{(m), b', \mathfrak{X}^{\dd, \circ}_{\sigma^{\dd}}},\quad (\Ig_{(m), b'}^{\tor})^{\wedge}_{\Phi^{\dd}} \cong \Ig^{\natural}_{(m), b', \mathfrak{X}^{\dd}_{\Sigma^{\dd}(\Phi^{\dd})}}/\Gamma_{\Phi^{\dd}}. \]
    \end{corollary}

\subsection{Hodge-type Igusa varieties}

Let $b \in G(\bQ)$, and $b' \in G^{\dd}(\bQ)$ be its image. Let us denote by $\XXp_{\ab_{\GG}} = \ab[p^{\infty}]|_{\CE^{b}}$. Due to \cite[Lemma 2.4.3]{kim2019central}, each $[b] \in B(G, \lrbracket{\mu})$ is completely slope divisible in the sense of \cite[Definition 2.4.1]{kim2019central}, we pick a representative that is completely slope divisible. 

\subsubsection{Relative Normalizations}

        \begin{lemma}\label{lemma: rela norm of rela norm is rela norm}
         Let $Y \to Z \to X$ be qcqs morphisms, let $Z_X \to X$ (resp. $Y_Z$, $Y_X$) be the relative normalization of $X$ in $Z$ (resp. $Z_X$ in $Y$, $X$ in $Y$), then there is a canonical isomorphism $Y_X = Y_Z$. 

     \end{lemma}
     \begin{proof}
         Due to \cite[\href{https://stacks.math.columbia.edu/tag/035I}{Tag 035I}]{stacks-project}, there is a unique morphism $Y_X \to Y_Z$ such that $Y \to Y_Z$ factors through $Y_X$, and $Y_X \to Y_Z$ is the relative normalization of $Y_Z$ in $Y$. On the other hand, due to \cite[\href{https://stacks.math.columbia.edu/tag/0BXA}{Tag 0BXA}]{stacks-project}, the relative normalization of $Y_Z$ in $Y$ is $Y_Z$ itself, thus $Y_X = Y_Z$.
     \end{proof}
     
 \begin{lemma}\label{lemma: rela norm of limits}
         Let $\lrbracket{f_i: X_i \to Y}_{i \in I}$ be a direct inverse system of qcqs morphisms, and assume the transition morphisms $X_i \to X_j$ are affine. Let $Z_i$ be the relative normalization of $Y$ in $X_i$, then $Z = \prolim Z_i$ is the relative normalization of $Y$ in $X = \prolim X_i$.
     \end{lemma}
     \begin{proof}
         Let $f: X \to Y$ be the morphism induced by $\lrbracket{f_i}_{i \in I}$, $p_i: X \to X_i$ are affine (\cite[\href{https://stacks.math.columbia.edu/tag/01YX}{Tag 01YX}]{stacks-project}), $f: X \to Y$ is qcqs. Let $Z'$ be the relative normalization of $Y$ in $X$, then $Z' = \SPEC_Y \OO'$, where $\OO'$ is the integral closure of $\OO_Y$ in $f_*\OO_X$. By construction, $f_*\OO_X = \injlim f_{i, *}\OO_{X_i}$ as quasi-coherent sheaves over $\OO_Y$. The property follows from the fact that taking integral closure of modules commutes with taking colimit.
     \end{proof}

\subsubsection{Hodge-type case, step 1}

\begin{definition}\label{definition: pullback of Igusa variety}\leavevmode
\begin{enumerate}
    \item Let $\Ig^{\square}_{m, b} = \Ig_{m, b'} \times_{\CE^{b'}} \CE^b$, $\Ig^{\square, \tor}_{m, b} = \Ig_{m, b'}^{\tor} \times_{\CE^{b', \tor}} \CE^{b, \tor}$.
    \item Let $\Ig^{\square, \min}_{m, b} \to \CE^{b, \min}$ be the relative normalization of $\CE^{b, \min}$ in $\Ig^{\square}_{m, b}$. 
\end{enumerate}
    We define $\Ig^{\square}_b$, $\Ig^{\square, \tor}_b$, $\Ig^{\square, \min}_b$ in the same way. This is compatible with \cite[Definition 3.3.7]{caraiani2024generic}.
\end{definition}

\begin{definition}\leavevmode
\begin{enumerate}
    \item Let $\IGUSA^{\square}_b = \IGUSA_{b'} \times_{\CE^{b', \perf}} \CE^{b, \perf}$, $\IGUSA^{\square, \tor}_b = \IGUSA_{b'}^{\tor} \times_{\CE^{b', \tor, \perf}} \CE^{b, \tor, \perf}$, these are perfection of $\Ig^{\square}_b$ and $\Ig^{\square, \tor}_b$ respectively.
    \item Let $\IGUSA^{\square, \min}_{b} \to \CE^{b, \min}$ be the relative normalization of $\CE^{b, \min}$ in $\IGUSA^{\square}_b$.
\end{enumerate}
\end{definition}

\begin{remark}\label{remark: limit of pullback Igusa}
    Due to \cite[\href{https://stacks.math.columbia.edu/tag/01YZ}{Tag 01YZ}]{stacks-project} and Lemma \ref{lemma: rela norm of limits}, $\Ig^{\square}_b = \prolim_{m \geq 0} \Ig^{\square}_{m, b}$, $\Ig^{\square, \min}_b=\prolim_{m \geq 0} \Ig^{\square, \min}_{m, b}$, $\Ig^{\square, \tor}_b = \prolim_{m \geq 0} \Ig^{\square, \tor}_{m, b}$. Since $\Ig^{\square}_{m, b}$ (resp. $\Ig^{\square}_b$) is finite \'etale (resp. pro-\'etale) over the smooth scheme $\CE^b$, thus $\Ig^{\square, \min}_{m, b}$ (resp. $\Ig^{\square, \min}_{b}$) is a normal scheme.
\end{remark}
\begin{remark}\label{remark: open dense Igusa variety}
    Since both $\Ig^{\square}_{(m), b} \to \CE^b$ and $\Ig^{\square, \tor}_{(m), b} \to \CE^{b, \tor}$ are $\Gamma_{(m), \XXp}$-torsors, $\Ig^{\square}_{(m), b}$ is open dense in $\Ig^{\square, \tor}_{(m), b}$.
\end{remark}

\begin{lemma}
    $\Ig^{\square, \min}_{m, b}$ (resp. $\Ig^{\square, \min}_b$) is the relative normalization of $\CE^{b, \min}$ in $\Ig^{\square, \tor}_{m, b}$ (resp. $\Ig^{\square, \tor}_b$). In particular, $\Ig^{\square, \tor}_{m, b} \to \Ig^{\square, \min}_{m, b}$ is the Stein factorization of $\Ig^{\square, \tor}_{m, b} \to \CE^{b, \min}$, thus it has geometrically connected fibers.
\end{lemma}
\begin{proof}
    Due to Lemma \ref{lemma: rela norm of limits} and Remark \ref{remark: limit of pullback Igusa}, it suffices to consider the finite level. Consider the commutative diagram, where the morphism $\Ig^{\square, \min}_{m, b} \to \Ig^{\min}_{m, b'}$ is induced by functoriality, see \cite[\href{https://stacks.math.columbia.edu/tag/035J}{Tag 035J}]{stacks-project}:
\[\begin{tikzcd}[sep=tiny]
	&& {\Ig^{\square, \min}_{m, b}} \\
	& {\Ig^{\square, \tor}_{m, b}} & {\CE^{b, \min}} & {\Ig^{\min}_{m, b'}} \\
	{\Ig^{\square}_{m, b}} & {\CE^{b, \tor}} & {\Ig^{\tor}_{m, b'}} & {\CE^{b', \min}} \\
	{\CE^b} & {\Ig_{m, b'}} & {\CE^{b', \tor}} \\
	& {\CE^{b'}}
	\arrow[from=3-1, to=4-1]
	\arrow[from=2-2, to=3-2]
	\arrow[from=1-3, to=2-3]
	\arrow[from=3-1, to=2-2]
	\arrow[from=2-2, to=1-3]
	\arrow[from=4-1, to=3-2]
	\arrow[from=3-2, to=2-3]
	\arrow[from=1-3, to=2-4]
	\arrow[from=2-4, to=3-4]
	\arrow[from=3-3, to=4-3]
	\arrow[from=5-2, to=4-3]
	\arrow[from=4-3, to=3-4]
	\arrow[from=4-2, to=3-3]
	\arrow[from=3-3, to=2-4]
	\arrow[from=2-3, to=3-4]
	\arrow[from=3-2, to=4-3]
	\arrow[from=4-1, to=5-2]
	\arrow[from=2-2, to=3-3]
	\arrow[from=3-1, to=4-2]
	\arrow[from=4-2, to=5-2]
\end{tikzcd}\]
Since $\Ig_{m, b}^{\square, \tor} \to \CE^{b, \tor}$ is finite and $\CE^{b, \tor} \to \CE^{b, \min}$ is proper, then $\Ig_{m, b}^{\square, \tor} \to \CE^{b, \min}$ is proper, then the relative normalization of $\CE^{b, \min}$ in $\Ig_{m, b}^{\square, \tor}$ is exactly the Stein factorization of $\Ig_{m, b}^{\square, \tor} \to \CE^{b, \min}$ (\cite[\href{https://stacks.math.columbia.edu/tag/03GQ}{Tag 03GQ}]{stacks-project}), we denote it by $\Ig_{m, b}^{\square, \min, \star}$. Due to Lemma \ref{lemma: rela norm of rela norm is rela norm} and the above discussion, $\Ig_{m, b}^{\square, \min}$ is the relative normalization of $\Ig_{m, b}^{\square, \min, \star}$ in $\Ig_{m, b}^{\square}$. Since $\Ig_{m, b}^{\square, \tor}$ is smooth, $\Ig_{m, b}^{\square, \min, \star}$ is normal and contains $\Ig_{m, b}^{\square}$ as an open subscheme, thus $\Ig_{m, b}^{\square, \min} \rightiso \Ig_{m, b}^{\square, \min, \star}$ canonically.
\end{proof}

\begin{corollary}\label{corollary: Ig is rela norm of Cmin in Igtor}
     $\IGUSA^{\square, \min}_b$ is the relative normalization of $\CE^{b, \min}$ in $\IGUSA^{\square, \tor}_b$. The natural morphism $\IGUSA_{b}^{\square, \min} \to \Ig_{b}^{\square, \min}$ induces an isomorphism $\IGUSA_{b}^{\square, \min} \to \Ig_{b}^{\square, \min, \perf}$. $\IGUSA^{\square, \tor}_b \to \IGUSA^{\square, \min}_b$ has geometrically connected fibers.
\end{corollary}

\subsubsection{Hodge-type case, step 2}

We follow the construction in \cite{hamacher2019adic}, there is a sub-functor $\IGUSA$ of $\IGUSA^{\square}$ which serves as the real Igusa variety on $\CE^b$. 
\begin{definition-proposition}{\cite[Definition-Lemma 6.1]{hamacher2019adic}}
    Let $\IGUSA_b \subset \IGUSA_b^{\square}$ be the sub-functors such that for all $\CE^{b}$-scheme $S$, $\IGUSA_b(S)$ parametrizes isomorphisms $\rho: \XXp|_S \to \XXp_{\ab_{\GG}}|_S$ which commute with polarizations and preserve the Hodge tensors point-wisely. $\IGUSA_b \subset \IGUSA_b^{\square}$ is a closed union of connected components.
\end{definition-proposition}

 Similarly, the statement above \cite[Definition 2.3.3]{caraiani2024generic} can be applied here. Note that one can always take a quasi-isogeny $\XXp_{b_1} \to \XXp_{b_2}$ which keeps the \'etale part, multiplicative part and biconnected part.
\begin{lemma}\label{lemma: Igusa isomorphism interior}
    Assume $[b_1] = [b_2] \in B(G, \lrbracket{\mu})$, then we have an isogeny $\XXp_{b_1} \to \XXp_{b_2}$ keeping the extra structures that induces an isomorphism $\IGUSA_{b_1} \rightiso \IGUSA_{b_2}$.
\end{lemma}

Note that a completely slope divisible $p$-divisible group over a perfect scheme has a canonical splitting given by slope isoclinic subquotients.
\begin{definition}\label{def: Igusa, Hodge}
\leavevmode
\begin{enumerate}
    \item  We define $\Ig_{m, b}(S) \subset \Ig_{m, b}^{\square}(S)$ be the subfunctor such that for schemes $S$ over $\CE^b$, and over the perfection $S^{\perf}$ of $S$, the data $(\rho_{m, i})_{i = 1}^r$ satisfies the following property: there exists a lifting $(\rho_{m', i})_{m' \geq m}$ such that
    \[
  \prolim_{m'} \oplus_{i = 1}^r \rho_{m', i, S^{\perf}}: \XXp[p^{\infty}]_{S^{\perf}} \cong \bigoplus \XXp_i[p^{\infty}]_{S^{\perf}} \cong \bigoplus \XXp_{\ab_{\GG}, i}[p^{\infty}]_{S^{\perf}} \cong \XXp_{\ab_{\GG}}[p^{\infty}]_{S^{\perf}} \]
  keeps the Hodge tensors.
  \item  We define $\Ig_{b}(S) \subset \Ig_{b}^{\square}(S)$ be the subfunctor such that for schemes $S$ over $\CE^b$, and over the perfection $S^{\perf}$ of $S$,  the data $(\rho_{i})_{i = 1}^r$ satisfies the following property: 
    \[
  \oplus_{i = 1}^r \rho_{i, S^{\perf}}: \XXp[p^{\infty}]_{S^{\perf}} \cong \bigoplus \XXp_i[p^{\infty}]_{S^{\perf}} \cong \bigoplus \XXp_{\ab_{\GG}, i}[p^{\infty}]_{S^{\perf}} \cong \XXp_{\ab_{\GG}}[p^{\infty}]_{S^{\perf}} \]
  keeps the Hodge tensors.
\end{enumerate}
\end{definition}
It follows easily that the natural projection $\Ig_{(m'), b}^{\square} \to \Ig_{(m), b}^{\square}$ induces surjections $\Ig_{(m'), b} \to \Ig_{(m), b}$ for all $m' > m$, $\Ig_{b} = \prolim_m \Ig_{m, b}$, and $\IGUSA_b \to \Ig_b$ induces an isomorphism $\IGUSA_b \to \Ig_b^{\perf}$. Taking perfection does not change the underlying topological space, $\Ig_b \subset \Ig_b^{\square}$ has the underlying topological space with $\IGUSA_b \subset \IGUSA^{\square}_b$.

Moreover, let $\Gamma_{m, \XXp, (s_{\alpha})}$ be the finite group of automorphisms of $\XXp[p^m]/k$ that lifts to $\XXp/k$ keeping the polarization and Hodge tensors, let $\Gamma_{\XXp, (s_{\alpha})} = \prolim \Gamma_{m, \XXp, (s_{\alpha})}$. Then $\Ig_{m, b}$ (resp. $\Ig_b$) is a $\Gamma_{m, \XXp, (s_{\alpha})}$- (resp. $\Gamma_{\XXp, (s_{\alpha})}$-) torsor over $\CE^{b}$. In particular, $\Ig_{m, b}$ (resp. $\Ig_b$) is finite \'etale (resp. pro-\'etale) over $\CE^{b}$.

\begin{proposition}\label{prop: open and closed, finite Igusa variety}
    $\Ig_{(m), b} \subset \Ig_{(m), b}^{\square}$ is a closed union of connected components. In particular, since $\Ig_{m, b}$ is a noetherian scheme, all its connected components are open, any union of connected components are clopen.
\end{proposition}
\begin{proof}
    This mainly follows from the same reason as in \cite[Definition-Lemma 6.1]{hamacher2019adic}. Another way to see this is as follows: $\Ig_b \subset \Ig_b^{\square}$ has the underlying topological space with $\IGUSA_b \subset \IGUSA^{\square}_b$, thus $\Ig_{b} \subset \Ig_{b}^{\square}$ is a closed union of connected components. At finite level, we pass to the perfection of the bases. Given a perfect connected scheme $S$, if the abelian scheme over that is geometrically fiberwise constant, then the crystalline Hodge tensors are also geometrically fiberwise constant, see the proof of Lemma \ref{lemma: central leaves are open and closed, PR}. In particular, if a point in $S$ is in the image of the projection $\Ig_b \to \Ig_{m, b}^{\square}$, then whole $S$ is in the projection $\Ig_b \to \Ig_{m, b}^{\square}$, thus $\Ig_{m, b} \subset \Ig_{m, b}^{\square}$ is a union of connected components. 
\end{proof}

\begin{definition-proposition}\label{def: tor Igusa variety, level m}
    Let $\Ig_{m, b}^{\tor} \subset \Ig_{m, b}^{\square, \tor}$ be the partial toroidal compactification of $\Ig_{m, b}$ inside $\Ig_{m, b}^{\square, \tor}$, which is simply the closure of $\Ig_{m, b}$ in $\Ig_{m, b}^{\square, \tor}$ since $\Ig_{m, b} \subset \Ig_{m, b}^{\square}$ is closed. It follows from Lemma \ref{lemma: open closed subschemes are well positioned 3} that $\Ig_{m, b}^{\tor} \subset \Ig_{m, b}^{\square, \tor}$ is open and closed.
\end{definition-proposition}

We can use full strength of Proposition \ref{prop: well position over well position, and affine}, with the help of Proposition \ref{prop: well-position, Igusa, Siegel case},
\begin{proposition}\label{prop: Igusa varieties are well-positioned}
        Igusa varieties $\Ig_{m, b}$ over central leaves at stablizer quasi parahoric levels satisfy the assumption \ref{assumption: well position over well position}. In particular, for each $[\Phi] \in \Cusp_K(G, X)$, there exist finite \'etale $\Ig_{m, b, C(\Phi)}^{\natural} \to \CE^{b, \natural}_{C(\Phi)}$ such that $\Ig_{m, b} \to \CE^b$ is well-positioned with respect to $\lrbracket{\Ig_{m, b, C(\Phi)}^{\natural} \to \CE^{b, \natural}_{C(\Phi)}}_{[\Phi]}$. 
\end{proposition}
\begin{remark}
    $\Ig_{m, b, C}^{\natural} \subset \Ig_{m, b, C}^{\square, \natural}$ is open and closed in $\Ig_{m, b, C}^{\square, \natural}$ due to Lemma \ref{lemma: open closed subschemes are well positioned 3}, where $\Ig_{m, b, C}^{\square, \natural}$ is the pullback of $\Ig_{m, b', C^{\dd}}^{\natural} \to \CE^{b', \natural}_{C}$ along $\CE^{b, \natural}_{C} \to \CE^{b', \natural}_{C^{\dd}}$.
\end{remark}

\begin{definition-proposition}\label{def: tor Igusa variety}
    Let $\Ig_{b}^{\tor} \subset \Ig_{b}^{\square, \tor}$ be the partial toroidal compactification of $\Ig_{b}$ inside $\Ig_{b}^{\square, \tor}$, which is simply the closure of $\Ig_{b}$ in $\Ig_{b}^{\square, \tor}$ since $\Ig_{b} \subset \Ig_{b}^{\square}$ is closed. We similarly define $\IGUSA_b^{\tor} \subset \IGUSA_b^{\square, \tor}$ as the subscheme with underlying topological space of $\Ig_b^{\tor} \subset \Ig_b^{\square, \tor}$.
\end{definition-proposition}

For all $m' \geq m$, $\Ig_{m', b}^{\square, \tor} \to \Ig_{m, b}^{\square, \tor}$ and $\Ig_{m', b} \to \Ig_{m, b}$ are affine, flat and surjective, thus $\Ig_{m', b}^{\tor} \to \Ig_{m, b}^{\tor}$ is affine, flat and surjective. By fpqc descent, $\Ig_{m', b, C(\Phi)}^{\natural} \to \Ig_{m, b, C(\Phi)}^{\natural}$ are also affine, flat and surjective. Due to Lemma \ref{lemma: limit of well-positioned},
\begin{proposition}\label{prop: limit Igusa varieties are well-positioned}
        For each $[\Phi] \in \Cusp_K(G, X)$, there exist pro-\'etale $\Ig_{b, C(\Phi)}^{\natural} \to \CE^{b, \natural}_{C(\Phi)}$ such that $\Ig_b \to \CE^b$ is well-positioned with respect to $\lrbracket{\Ig_{b, C(\Phi)}^{\natural} \to \CE^{b, \natural}_{C(\Phi)}}_{[\Phi]}$, where $\Ig_{b, C(\Phi)}^{\natural} = \prolim \Ig_{m, b, C(\Phi)}^{\natural}$.
\end{proposition}

\begin{definition}
    Let $\Ig_{(m), b}^{\min}$ be the relative normalization of $\CE^{b, \min}$ in $\Ig_{(m), b}$.
\end{definition}
Lemma \ref{lemma: minimal compactifications disjoint} shows that $\Ig_{m, b}^{\min} \subset \Ig_{m, b}^{\square, \min}$ is open and closed, and Corollary \ref{cor: minimal compactifications of central leaves are affine} together with Proposition \ref{prop: well position over well position, and affine} imply that
\begin{corollary}\label{cor: min Igusa are affine}
    Minimal compactifications $\Ig_{(m), b}^{\min}$ of Igusa varieties $\Ig_{(m), b}$ over central leaves at stablizer quasi-parahoric levels are affine.
\end{corollary}

  Let us describe the boundary of $\Ig_{(m), b}^{\tor}$ in the usual formal way. Fix a $[\Phi]$. Recall that $\Ig^{\natural}_{(m), b, \Xi} \to \CE^{b, \natural}_{\Xi}$ is the pullback of $\Ig^{\natural}_{(m), b, C} \to \CE^{b, \natural}_{C}$ along $\Xi \to C$. We similarly deonte the pullbacks of $\Ig^{\natural}_{(m), b, C} \to \CE^{b, \natural}_{C}$ along $\Xi(\sigma) \to C$, $\Xi_{\sigma} \to C$, $\ovl{\Xi}_{\Sigma(\Phi)} \to C$ by $\Ig^{\natural}_{(m), b, \Xi(\sigma)}$, $\Ig^{\natural}_{(m), b, \Xi_{\sigma}}$, $\Ig^{\natural}_{(m), b, \ovl{\Xi}_{\Sigma(\Phi)}}$ respectively.

    For any $\sigma \in \Sigma(\Phi)^+$, we denote by $\Ig_{(m), b, \mathcal{Z}([\Phi, \sigma])}$ the boundary stratum of $\Ig_{(m), b}^{\tor}$ over $\CE^{b, \tor}_{\mathcal{Z}([\Phi, \sigma])}$.

  Let $(\Ig_{(m), b}^{\tor})^{\wedge}_{\mathcal{Z}([\Phi, \sigma])}$ (resp. $(\Ig_{(m), b}^{\tor})^{\wedge}_{\Phi, \sigma}$, $(\Ig_{(m), b}^{\tor})^{\wedge}_{\Phi}$) be the completion of $\Ig_{(m), b}^{\tor}$ along
    \[ \Ig_{(m), b, \mathcal{Z}([\Phi, \sigma])}, \quad (\textit{resp.}\ \bigcup_{\tau \in \Sigma(\Phi)^+,\ \Bar{\tau} \subset \Bar{\sigma}}  \Ig_{(m), b, \mathcal{Z}([\Phi, \tau])},\quad  \bigcup_{[\tau] \in \Sigma(\Phi)^+/\Gamma_{\Phi}}  \Ig_{(m), b, \mathcal{Z}([\Phi, \tau])}), \]
    and let $\Ig^{\natural}_{(m), b, \mathfrak{X}_{\sigma}}$ (resp. $\Ig^{\natural}_{(m), b, \mathfrak{X}^{\dd, \circ}_{\sigma}}$, $\Ig^{\natural}_{(m), b, \mathfrak{X}_{\Sigma(\Phi)}}$) be the completion of $\Ig^{\natural}_{(m), b, \Xi(\sigma)}$ (resp. $\Ig^{\natural}_{(m), b, \Xi(\sigma)}$, $\Ig^{\natural}_{(m), b, \ovl{\Xi}_{\Sigma(\Phi)}}$) along 
    \[ \Ig^{\natural}_{(m), b, \Xi_{\sigma}},\quad (\textit{resp.}\ \bigcup_{\tau \in \Sigma(\Phi)^+,\ \Bar{\tau} \subset \Bar{\sigma}}  \Ig^{\natural}_{(m), b, \Xi_{\tau}},\quad   \bigcup_{\tau \in \Sigma(\Phi)^+} \Ig^{\natural}_{(m), b, \Xi_{\tau}}).  \]
    Due to Proposition \ref{prop: wdt Y satisfy same prop},
    \begin{corollary}\label{cor: description of Igusa boundary}
        There is a canonical isomorphism $\Ig_{(m), b, \mathcal{Z}([\Phi, \sigma])} \cong \Ig^{\natural}_{(m), b, \Xi_{\sigma}}$ extends to isomorphisms of formal schemes:
        \[ (\Ig_{(m), b}^{\tor})^{\wedge}_{\mathcal{Z}([\Phi, \sigma])} \cong \Ig^{\natural}_{(m), b, \mathfrak{X}_{\sigma}},\quad (\Ig_{(m), b}^{\tor})^{\wedge}_{\Phi, \sigma} \cong  \Ig^{\natural}_{(m), b, \mathfrak{X}^{\dd, \circ}_{\sigma}},\quad (\Ig_{(m), b}^{\tor})^{\wedge}_{\Phi} \cong \Ig^{\natural}_{(m), b, \mathfrak{X}_{\Sigma(\Phi)}}/\Gamma_{\Phi}. \]
    \end{corollary}
\subsection{Some propositions}

\begin{proposition}\label{prop: Igusa closed embedding}
    $\Ig_{(m), b} \to \Ig_{(m), b'}$ and $\IGUSA_b \to \IGUSA_{b'}$ are closed embeddings.
\end{proposition}
\begin{proof}
    In \cite[Proposition 6.2]{hamacher2019adic}, the authors showed that $\IGUSA_b \to \IGUSA_{b'}$ is a closed embedding. For the exactly same reason, $\Ig_b \to \Ig_{b'}$ is also a universal injection. Since $\Ig_{m, b} \subset \Ig_{m, b}^{\square} = \CE^b \times_{\CE^{b'}} \Ig_{m, b'}$, given two points $(z_1, j_1), (z_2, j_2) \in \Ig_{m, b}(S)$, if they have same image in $\Ig_{m, b'}$, then $j_1 = j_2$. Since the image of $z_1$ coincides with image of $z_2$ (which we denote by $z$), there is a canonical identification $\ab_{\GG, z_1} = \ab_{\GG, z_2} = \ab_z$, and $j_1 = j_2$ gives exactly same $\lrbracket{\rho_{m, i}: \ab_{z, i}[p^m] \cong \XXp_i[p^{m}]}$. In particular, we can lift $(z_1, j_1)$ and $(z_2, j_2)$ to points $(z_1, \tilde{j}_1)$, $(z_2, \tilde{j}_2)$ in $\Ig_{b}(S)$ such that $\tilde{j}_1 = \tilde{j}_2$. Since $\Ig_b \to \Ig_{b'}$ is a universal injection, then $z_1 = z_2$, $\Ig_{m, b} \to \Ig_{m, b'}$ is a universal injection.

    Given $z \in \CE^b(k)$ with image $z' \in \CE^{b'}(k)$, let $\Spf R$ be the formal completion of $\Shum{K^{\dd}, \bZ_p}$ at $z'$, $\Spf R_G$ be the formal completion of $\Shum{K, \OO_{\breve{E}}}$ at $z$, then $\Spf R_G \to \Spf R$ is a closed embedding. Denote by $\bar{R}_G = R_G/\mathfrak{m}_{\breve{E}}$, $\bar{R} = R/pR$. Due to \cite[Proposition 5.1.1]{kim2019central}, $\Spec \hat{\OO}_{\CE^b, z}$ is the scheme with reduced subscheme structure on $\Spec \hat{\OO}_{\CE^{b'}, z'} \times_{\Spec \bar{R}} \Spec \bar{R}_G$, thus $\Spec \hat{\OO}_{\CE^b, z} \stackrel{\hat{f}}{\to} \Spec \hat{\OO}_{\CE^{b'}, z'}$ is a closed embedding. In particular, $\hat{f}$ is an unramified morphism. Due to \cite[\href{https://stacks.math.columbia.edu/tag/039I}{Tag 039I}]{stacks-project}, the morphism between the local rings $\OO_{\CE^b, z} \stackrel{f}{\to} \OO_{\CE^{b'}, z'}$ is unramified. Since $z, z'$ are randomly chosen and we are dealing with Jacobson schemes, then $\CE^b \to \CE^{b'}$ is unramified. In particular, its base change $\Ig_{(m), b}^{\square} \to \Ig_{(m), b'}$ is unramified. $\Ig_{(m), b}$ is a closed subscheme of $\Ig_{(m), b}^{\square}$, then $\Ig_{(m), b} \to \Ig_{(m), b}^{\square} \to \Ig_{(m), b'}$ is unramified.

    In summary, $\Ig_{(m), b} \to \Ig_{(m), b'}$ is universally closed (it is even finite), universally injective, unramified, thus it is a closed embedding due to \cite[\href{https://stacks.math.columbia.edu/tag/04XV}{Tag 04XV}]{stacks-project}.
\end{proof}
Due to Proposition \ref{prop: tor closed embedding}, and note that $\IGUSA_b^{\tor} \rightiso \Ig_b^{\tor}$, $\IGUSA_{b'}^{\tor} \rightiso \Ig_{b'}^{\tor}$,
\begin{corollary}\label{cor: Ig tor closed embedding}
    Let $(\Sigma, \Sigma^{\dd})$ be strictly compatible, then $\Ig_{(m), b}^{\tor} \to \Ig_{(m), b'}^{\tor}$, $\Ig_{(m), b, C}^{\natural} \to \Ig_{(m), b', C}^{\natural}$ and $\IGUSA_b^{\tor} \to \IGUSA_{b'}^{\tor}$ are closed embeddings.
\end{corollary}

We can apply \ref{eq: pink formula diagram} to Igusa varieties $\lrbracket{\Ig_{K(l^m), b}}_{m \geq 0}$, where $l$ is a prime different from $p$, $\Ig_{K(l^m), b}$ is the Igusa variety defined over the central leaf $\CE^b \subset \Shum{K(l^m), \bar{s}}$. Due to Proposition \ref{prop: pink's formula}, we have
\begin{proposition}\label{prop: pink to Igusa varieties}
    
   $ \varinjlim_m R\Gamma(\Ig_{K(l^m), b}^{\tor}, \F_l) \rightiso \varinjlim_m R\Gamma(\Ig_{K(l^m), b}, \F_l)$.
\end{proposition}

\begin{proposition}
    Keep notations in Lemma \ref{lemma: Igusa isomorphism interior}. The isomorphism $\IGUSA_{b_1} \cong \IGUSA_{b_2}$ extends to $\IGUSA_{b_1}^{\tor} \cong \IGUSA_{b_2}^{\tor}$.
\end{proposition}
\begin{proof}
    In some PEL type cases, this is included in the proof of \cite[Theorem 3.3.2, Proposition 3.3.4]{caraiani2024generic}. We give the proof in the Hodge-type case. Let $b_1', b_2'$ be the images of $b_1, b_2$ in $B(G^{\dd}, \lrbracket{\mu^{\dd}})$ respectively. Then the isogeny $\varphi: \XXp_{b_1} \to \XXp_{b_2}$ is also the isogeny $\XXp_{b_1'} \to \XXp_{b_2'}$, the isomorphism $\IGUSA_{b_1} \cong \IGUSA_{b_2}$ is the restriction of $\IGUSA_{b_1'} \cong \IGUSA_{b_2'}$ (via the closed embeddings $\IGUSA_{b_i} \hookrightarrow \IGUSA_{b_i'}$). Apply Corollary \ref{cor: Ig tor closed embedding}, the isomorphism $\IGUSA_{b_1'}^{\tor} \cong \IGUSA_{b_2'}^{\tor}$ induces the isomorphism $\IGUSA_{b_1}^{\tor} \cong \IGUSA_{b_2}^{\tor}$ we want.
\end{proof}

%% file: sections/KR_strata.tex
We follow the notations and assumptions in subsection \ref{subsec: def KR strata}. 

\subsection{Well-position}

    \begin{proposition}\label{prop: KR strata are well-positioned}\leavevmode
    \begin{enumerate}
        \item Let $w\in\Adm(\lrbracket{\mu^{-1}})_{\KKc}$, then KR strata $\KR_{\Kc, w}$, $\KR_{\Kc, \leq w}$ and their connected components are well-positioned subsets in $\Shum{\Kc, \Bar{s}}$, they are well-positioned subschemes with respect to the induced reduced subscheme structures.
        \item $\Shum{\Kc, \bar{s}}$ (resp. $\Shumm{\Kc, \bar{s}}$, $\Shumc{\Kc, \Sigma, \bar{s}}$) is a disjoint union of connected components of $\KR_{\Kc, w}$ (resp. $\KR_{\Kc, w}^{\min}$, $\KR_{\Kc, w, \Sigma}^{\tor}$) for those $w \in \Adm(\lrbracket{\mu^{-1}})_{\KKc}$.
          \item Given any KR strata $\KR_{\Kc, w_1}$, $\KR_{\Kc, w_2}$,
          \begin{align*}
               \ovl{\KR_{\Kc, w_1}} \cap \KR_{\Kc, w_2} \neq \emptyset &\Longleftrightarrow \ovl{\KR_{\Kc, w_1}^{\min}} \cap \KR_{\Kc, w_2}^{\min} \neq \emptyset \\ &\Longleftrightarrow \ovl{\KR_{\Kc, w, \Sigma}^{\tor}} \cap \KR_{\Kc, w_2, \Sigma}^{\tor} \neq \emptyset.
          \end{align*}
         \item We have closure relations:
         \begin{equation*}
     \ovl{\KR_{\Kc, w, \Sigma}^{\tor}}= \KR_{\Kc, \leq w, \Sigma}^{\tor} = \bigsqcup_{w' \leq w} \KR_{\Kc, w', \Sigma}^{\tor},\quad \ovl{\KR_{\Kc, w}^{\min}}= \KR_{\Kc, \leq w}^{\min} = \bigsqcup_{w' \leq w} \KR_{\Kc, w'}^{\min}.
\end{equation*}
    \end{enumerate}
    \end{proposition}
    \begin{proof}\leavevmode
    \begin{enumerate}
        \item  Since such defined KR strata are locally closed, and are unions of central leaves. We have shown that central leaves are well-positioned in Proposition \ref{prop: Newton central are well-posiitoned, PR}, then KR strata $\KR_{\Kc, w}$, $\KR_{\Kc, \leq w}$ are well-positioned, due to Lemma \ref{lem: union of well-positioned sets} and Proposition \ref{proposition: open-closed subschemes are well positioned}.
        \item KR strata satisfy a clean closure relation \ref{equation: closure relations of KR}, thus we apply Proposition \ref{prop: general well-position, disjoint}.
        \item Same as $(2)$.
        \item Apply Lemma \ref{lemma: closure of well-pos is well-pos} and \ref{lem: union of well-positioned sets}.
    \end{enumerate}
 
    \end{proof}

\subsection{Boundary descriptions}

In order to work with KR strata at parahoric levels on the boundary, we assume $\Kc_{\Phi, L, p} \cap G_h(\Q_p) = \Kc_{\Phi, h, p}$ for all $[\Phi] \in \Cusp_{\Kc}(G, X)$, this implies that $\Zb^{\bigsur}(\Phi)_{\Bar{s}} = \Shum{\Kc_{\Phi, h}}$ have parahoric level structures. If we are willing to consider the KR strata at quasi-parahoric levels, then this assumption is not necessary. The more important assumption is that, we need to assume $\Shum{K_{\Phi, h}}$ also supports a schematic local model diagram \ref{eq: conj, local model diagram}. This is true under the Assumption \ref{general condition}, see Theorem \ref{thm: local model diagram exists}. 

\begin{corollary}\label{cor: boundary of KR is KR, weaker}
    Assume $\Shum{K_{\Phi, h}}$ supports a schematic local model diagram \ref{eq: conj, local model diagram}. Let $Y \subset \Shum{\Kc, \Bar{s}}$ be a connected component of a KR stratum (resp. closure of a KR stratum), then its boundary strata $Y^{\natural}_{\Zb^{\bigsur}(\Phi)} \subset \Zb^{\bigsur}(\Phi)_{\Bar{s}}$ (for each $[\Phi]\in \Cusp_{\Kc}(G, X)$) is either empty or a connected component of a KR stratum (resp. closure of a KR stratum) in $\Zb^{\bigsur}(\Phi)_{\Bar{s}}$.
\end{corollary}
\begin{proof}
    As we mentioned in Remark \ref{remark: nearby cycles and KR strata}, closures of KR strata are also defined as supports of nearby cycles, we can apply Proposition \ref{prop: boundary of nby is nby}: note that KR strata are well-positioned due to Proposition \ref{prop: KR strata are well-positioned}, and $C \to \Zb^{\bigsur}$ is smooth with geometrically connected fibers due to Proposition \ref{remark: when the assumption C to Z is true}. This proves the case when $Y$ is a connected component of $\KR_{\Kc, \leq w}$. When $Y$ is a connected component of $\KR_{\Kc, w}$, since $\KR_{\Kc, w}$ is the smooth locus of $\KR_{\Kc, \leq w}$, we apply \cite[Lemma 2.2.10]{lan2018compactifications}.
\end{proof}

We have a stronger result:
\begin{proposition}\label{prop: boundary of KR is KR}
       Assume $\Shum{K_{\Phi, h}}$ supports a schematic local model diagram \ref{eq: conj, local model diagram}. Let $Y \subset \Shum{\Kc, \Bar{s}}$ be a KR stratum (resp. closure of a KR stratum), then its boundary strata $Y^{\natural}_{\Zb^{\bigsur}(\Phi)} \subset \Zb^{\bigsur}(\Phi)_{\Bar{s}}$ (for each $[\Phi]\in \Cusp_{\Kc}(G, X)$) is either empty or a KR stratum (resp. closure of a KR stratum) in $\Zb^{\bigsur}(\Phi)_{\Bar{s}}$.
\end{proposition}

First of all, the case when $Y = \KR_{\Kc, \leq w}$ comes from the case when $Y = \KR_{\Kc, w}$, due to Lemma \ref{lemma: closure of well-pos is well-pos}. We deal with the case  $Y = \KR_{\Kc, w}$. We give two proofs here, the first proof follow the spirit of the proof of Corollary \ref{cor: boundary of KR is KR, weaker} using nearby cycles, the second proof essentially uses the construction of local model diagram and plays with de Rham Hodge tensors.

\begin{proof}
  
  Since all the strata involved are locally closed subschemes in Jacobson schemes, it suffices to work with $k$-points. Fix some KR stratum $\KR_w \subset \Zb^{\bigsur}_{\Bar{s}}$, let $z_1, z_2 \in \KR_{w}(k)$. Let $Y_1, Y_2 \subset \Shum{\Kc, \bar{s}}$ be the unique KR strata such that $z_1 \in Y^{\natural}_{1, \Zb^{\bigsur}}$, $z_2 \in Y^{\natural}_{2, \Zb^{\bigsur}}$ respectively. Such $Y_1, Y_2$ exists due to Proposition \ref{prop: KR strata are well-positioned}.
  
  Due to Lemma \ref{lem: group torsor, complete local rings} and arguments in \ref{lemma: local ring and KR strata}, we have a $k$-point $y_1 \in Y_{1, \mathcal{Z}([\Phi, \sigma])} \cong Y^{\natural}_{\Xi_{\sigma}(\Phi)}$ (resp. $y_2 \in Y_{2, \mathcal{Z}([\Phi, \sigma])} \cong Y^{\natural}_{2, \Xi_{\sigma}(\Phi)}$) lifting $z_1$ (resp. $z_2$) such that $\OO_{\Zb^{\bigsur}, z_1}^{\wedge} \cong \OO_{\Zb^{\bigsur}, z_2}^{\wedge}$ implies $\OO_{\Xi(\sigma), y_1}^{\wedge} \cong \OO_{\Xi(\sigma), y_2}^{\wedge}$, then $\OO_{\Shumc{\Kc, \Sigma}, y_1}^{\wedge} \cong \OO_{\Shumc{\Kc, \Sigma}, y_2}^{\wedge}$, therefore $Y_1 = Y_2$ due to Lemma \ref{lemma: local ring and KR strata}. In particular, given a KR stratum $Y \subset \Shum{\Kc, \bar{s}}$, $Y^{\natural}_{\Zb^{\bigsur}} \cap \KR_{w}$ is non-empty implies that $\KR_w \subset Y^{\natural}_{\Zb^{\bigsur}}$. 
    
    On the other hands, let $Y \subset \Shum{\Kc, \bar{s}}$ be a KR stratum, we need to show $Y_{\Zb^{\bigsur}}^{\natural} \subset \Zb^{\bigsur}_{\bar{s}}$ is either empty or contained in a KR stratum. Let $Y_1, Y_2 \subset Y$ be two connected components of $Y$, they are well-positioned due to Proposition \ref{proposition: open-closed subschemes are well positioned}. Due to Corollary \ref{cor: boundary of KR is KR, weaker}, it suffices to prove $Y_{1, \Zb^{\bigsur}}^{\natural}$ and $Y_{2, \Zb^{\bigsur}}^{\natural}$ are contained in the same KR stratum in $\Zb^{\bigsur}_{\bar{s}}$. We fix a boundary stratum $Y_{\mathcal{Z}([\Phi, \sigma])} \cong Y_{\Xi_{\sigma}(\Phi)}^{\natural}$. Since $Y_{i, \Xi_{\sigma}(\Phi)}^{\natural}$ is Zariski locally $Y_{i, E(\sigma) \times C}^{\natural} = E(\sigma)_{\bar{s}} \times Y_{i, C}^{\natural}$ over $\Xi_{\sigma}(\Phi)_{\bar{s}}$, we can pick closed point $y_i \in Y_{i, \mathcal{Z}([\Phi, \sigma])}$, and let $x_{i} \in Y_{i, C(\Phi)}^{\natural}$ be its image, such that the isomorphism $\OO_{\Shumc{\Kc, \Sigma}, y_1}^{\wedge} \cong \OO_{\Shumc{\Kc, \Sigma}, y_2}^{\wedge}$ in Lemma \ref{lemma: local ring and KR strata} induces an isomorphism $\OO_{C, x_{1}}^{\wedge} \cong \OO_{C, x_{2}}^{\wedge}$. Let $z_1$, $z_2$ be the images of $x_{1}$, $x_{2}$ in $\Zb^{\bigsur}_{\bar{s}}$, since $C \to \Zb^{\bigsur}$ is an abelian scheme torsor,  then $\OO_{\Zb^{\bigsur}, z_1}^{\wedge} \cong \OO_{\Zb^{\bigsur}, z_2}^{\wedge}$ due to Lemma \ref{lem: group torsor, complete local rings}, then $z_1$, $z_2$ are in the same KR stratum $\KR_w \subset \Zb^{\bigsur}_{\bar{s}}$ due to Lemma \ref{lemma: local ring and KR strata}. 
\end{proof}

\begin{lemma}\label{lem: group torsor, complete local rings}
    Let $G \to S$ be a equidimensional smooth group scheme over a locally noetherian scheme $S$, $x_1$, $x_2$ be two geometric points of $G$ such that their images in $G$ and $S$ are closed points, with residue fields $k(x_1) \cong k(x_2)$, then $\OO_{S, x_1}^{\wedge}\cong \OO_{S, x_2}^{\wedge}$ if and only if $\OO_{G, x_1}^{\wedge} \cong \OO_{G, x_2}^{\wedge}$.
\end{lemma}
\begin{proof}
    Let $s_i$ be the image of $x_i$ in $S$. We first assume $s_1 = s_2 = s$ and $x_i \to S$ factors through $\Spec \kappa(s) \to S$, let $A = \OO_{S, \Tilde{s}}$ be the (\'etale) stalk of $S$ at the geometric point $\Tilde{s}$ over $s$ with maximal ideal $\mathfrak{m}$, since $x_1, x_2 \to G$ factors through $G_A$, and $\OO_{G, x_i} \cong \OO_{G_A, x_i}$, it suffices to work on $G \to \Spec A$. Moreover, let $G_{\hat{A}} = G \times_{\Spec A} \Spec \hat{A}$, where $\hat{A}$ is the $\mathfrak{m}$-completion of $A$, let $p \in G_{\hat{A}}$ be a point with image $q \in G$ lying over the closed point of $\Spec A$, then the local ring map induces an isomorphism after completion: $\OO_{G_{\hat{A}}, p}^{\wedge} \rightiso \OO_{G, q}^{\wedge}$ due to \cite[\href{https://stacks.math.columbia.edu/tag/0BG5}{Tag 0BG5}]{stacks-project}. Since $x_i \to G$ factors through $G_{\hat{A}}$, we apply this to the images of $x_i$. Since $k:= k(x_1) \cong k(x_2) \cong k(\bar{s})$, $x_1, x_2 \in G(k)$, there exists $g \in G(k)$ such that $gx_1 = x_2$. Since $G \to \Spec A$ is formally smooth, we lift $g$ to $\hat{g} \in G(\hat{A})$ over $\Spec \hat{A}$. Since $g$ induces an isomorphism $t_{g}: G \to G$ over $\Spec \hat{A}$ by left translation and maps $x_1$ to $x_2$, $t_{g}$ induces $\OO_{G_{\hat{A}}, x_1}^{\wedge} \cong \OO_{G_{\hat{A}}, x_2}^{\wedge}$, thus $\OO_{G, x_1}^{\wedge} \cong \OO_{G, x_2}^{\wedge}$ over $\Spec \hat{A}$. 
    
    By above arguements, in the remaining of the proof, we could replace $x_i$ with a point over $e(S)$, this does not change either $\OO_{S, x_i}^{\wedge}$ or $\OO_{G, x_i}^{\wedge}$. From now on, we assume $x_i$ has image in $e(S)$.

    Let $e: S \to G$ be the identity section, it is a regular locally closed embedding due to \cite[Theorem 17.12.1]{grothendieck1967elements} since $G \to S$ is smooth. Due to the previous step, we could assume $x_1$, $x_2$ are geometric points over $e(S) \subset G$. Since we are investigating local properties, replace $G$ with an open subscheme containing $x_1$ and $x_2$, we assume $e(S) \subset G$ is closed. Let $G^{\wedge}_{e(S)}$ be the formal completion of $G$ along $e(S)$. \cite[Corollary 16.9.9]{grothendieck1967elements} implies that we could find affine neighbourhood $V_i$ of $S$ at $s_i$ such that $G^{\wedge}_{e(S)}|_{V_i}$ is isomorphic to $\Spf \OO_{V_i}\llbracket T_1, \dots, T_{n_i} \rrbracket$. Since $G \to S$ is equidimensional, thus $n_1 = n_2$. 
    
    $(\Longrightarrow)$:  Since $\OO_{S, x_1}^{\wedge}\cong \OO_{S, x_2}^{\wedge}$, we could find \'etale neighbourhoods $U_i$ of $x_i \to S$ such that $U_1 \cong U_2$, see \cite[\href{https://stacks.math.columbia.edu/tag/0CAT}{Tag 0CAT}]{stacks-project}. By shrinking $V_i$ and pulling back to $U_i$, we see that $G^{\wedge}_{e(S)}|_{U_1} \cong G^{\wedge}_{e(S)}|_{U_2}$. Take further completion of $G^{\wedge}_{e(S)}|_{U_i}$ at $x_i$ (whose images are closed points), we get $\OO_{G, x_1}^{\wedge} \cong \OO_{G, x_2}^{\wedge}$.

    $(\Longleftarrow)$: Since $\OO_{G, x_1}^{\wedge} \cong \OO_{G, x_2}^{\wedge}$, then we have an isomorphism
    \[ \OO_{S, x_1}^{\wedge}\llbracket T_1, \dots, T_{n} \rrbracket \cong \OO_{S, x_2}^{\wedge}\llbracket T_1, \dots, T_{n} \rrbracket.\] Since $\OO_{S, x_i}^{\wedge}$ are Noetherian complete local ring, then $\OO_{S, x_1}^{\wedge}\cong \OO_{S, x_2}^{\wedge}$ due to \cite[Theorem 8]{ishibashi1976isomorphic} or \cite[Theorem 5]{hamann1975power}.
\end{proof}

\begin{lemma}\label{lemma: local ring and KR strata}
    Let $y_1, y_2$ be two closed points in $\Shumc{\Kc, \Sigma, \Bar{s}}$, then $y_1, y_2$ are in the same $\KR_{\Kc, w, \Sigma}^{\tor}$ if and only if $\OO_{\Shumc{\Kc, \Sigma}, y_1}^{\wedge} \cong \OO_{\Shumc{\Kc, \Sigma}, y_2}^{\wedge}$.
\end{lemma}
\begin{proof}
    Let $y \in \mathcal{Z}([\Phi, \sigma])$, due to \cite[Corollary 2.1.7]{lan2018compactifications}, there exist \'etale morphisms $(\ovl{U}, q) \to (\Shumc{\Kc, \Sigma}, y)$, $(\ovl{U}, q) \to (E(\sigma) \times_{\Spec \Z} C, (t, x))$, where $x$ is the image of $y$ in $C_{\Bar{s}}$. These \'etale morphisms restrict to the common open parts $U \to \Shum{\Kc}$, $U \to E \times_{\Spec \Z} C$. Then there exists a closed point $p \in U_{\bar{s}}$ which has image $y' \in \Shum{\Kc, \Bar{s}}$ and $(z', x) \in E_{\bar{s}} \times C_{\bar{s}}$ such that 
    \begin{equation}\label{eq: KR formal local, C}
        \OO_{\Shumc{\Kc, \Sigma}, y}^{\wedge} \cong \OO_{\ovl{U}, q}^{\wedge} \cong \OO_{E(\sigma) \times_{\Spec \Z} C, (z, x)}^{\wedge} \cong \OO_{E \times_{\Spec \Z} C, (z', x)}^{\wedge} \cong \OO_{U, p}^{\wedge} \cong \OO_{\Shum{\Kc}, y'}^{\wedge}.
    \end{equation}
    Since the partial toroidal compactifications of KR strata cover $\Shumc{\Kc, \Sigma, \Bar{s}}$ due to Proposition \ref{prop: KR strata are well-positioned}, we assume $y \in \KR^{\tor}_{\Kc, w, \Sigma}$ for some $w \in \Adm(\lrbracket{\mu^{-1}})_{\KKc}$. Since $\KR_{\Kc, w}$ is well-positioned, we can furthur take $y' \in \KR_{\Kc, w}$ due to Corollary \ref{corollay: cor 2.1.7 for well positioned Y}. Therefore, it suffices to work with $\Shum{\Kc}$.

    Recall we have a smooth morphism $\Shum{\Kc} \to [\mathbb{M}_{\GGc, \mu}/\GGc]$. There exist $z_i$ ($i = 1, 2$) in $\mathbb{M}^{\perf}_{\GGc, \mu, \bar{s}}(k)$ such that $\OO_{\mathbb{M}_{\GGc, \mu}, z_i}^{\wedge} \cong \OO_{\Shum{\Kc}, y_i}^{\wedge}$, and $y_i \in \KR_{\Kc, w_i}$ (for some $w_i \in \Adm(\lrbracket{\mu^{-1}})_{\KKc}$) if and only if $z_i \in \fl_{\GGc, w_i}^{\circ}$.
    
    $(\Longrightarrow)$: When $y_1, y_2$ are in the same KR stratum, say, $\KR_{\Kc, w}$, then $z_1, z_2 \in \fl_{\GGc, w}^{\circ}$. Note that $\fl_{\GGc, w}^{\circ}$ is a $\GGc(k)$-orbit in $\mathbb{M}_{\GGc, \mu}\otimes k$, there exists $g \in \GGc(k)$ such that $gz_1 = z_2$, we can lift $g$ to $\Tilde{g} \in \GGc(\OO_{\Breve{E}})$, the left translation $\Tilde{g}$ induces an isomorphism $\mathbb{M}_{\GGc, \mu} \otimes \OO_{\Breve{E}} \to \mathbb{M}_{\GGc, \mu} \otimes \OO_{\Breve{E}}$, thus an isomorphism $\OO_{\mathbb{M}_{\GGc, \mu}, z_1}^{\wedge} \cong \OO_{\mathbb{M}_{\GGc, \mu}, z_2}^{\wedge}$. Therefore, $\OO_{\Shum{\Kc}, y_1}^{\wedge} \cong \OO_{\Shum{\Kc}, y_2}^{\wedge}$.

    $(\Longleftarrow)$: Since $\OO_{\mathbb{M}_{\GGc, \mu}, z_1}^{\wedge} \cong \OO_{\mathbb{M}_{\GGc, \mu}, z_2}^{\wedge}$, there exist \'etale neighbourhoods $V_1 \cong V_2$ of $\mathbb{M}_{\GGc, \mu}$ at $z_1$ and $z_2$ respectively. The closures of KR strata on $\mathbb{M}_{\GGc, \mu}$ are the support of the nearby cycle of the constant sheaf over $\mathbb{M}_{\GGc, \mu}$, and the nearby cycle functor commutes with smooth morphisms, by adjusting automorphisms of $V_i$, we assume $V_1 \cong V_2$ respects the stratification on both sides, this forces $w_1 = w_2$.
\end{proof}

  Let us give another proof of Proposition \ref{prop: boundary of KR is KR}. For simplicity, we work in the setting of Kisin-Pappas integral models under the assumption \ref{general condition}.

  \begin{proof}

       We pullback the local model diagram $\Shum{K} \to [\mathbb{M}_{\GG, \mu}/\GG]$ to $W^0 \to \Shum{K}$, then we have a $\GG$-torsor $\wdt{W}^0$ over $W^0$ which trivializes de Rham cohomology with tensors. Given a geometric point $\bar{t} \to W^0_{\bar{s}}$, lift it to a discrete valuation ring $\Spec \OO_F \to W^0$, we choose a trivialization $f$ of the de Rham cohomology keeping the tensors as in equations \ref{eq: de rham Hodge tri, abelian sch} and \ref{eq: de rham Hodge tri, over B, abelian sch}. This gives a point in $(y, f) \in \wdt{W}^0(\OO_F)$ with $y \in W^0(\OO_F)$, then the $\GG$-orbit $f$ in $\mathbb{M}_{\GG, \mu}$ determines the KR stratum that the image of $y$ in $\Shum{K, \bar{s}}$ sits in, and the $\GG_h^{g}$-orbit of the $-1$-graded piece (let us denote by $f'$) in $\mathbb{M}_{\GG^g_h, \mu_h}$ determines the the boundary KR stratum that the image of $y$ in $\Shum{K_{\Phi, h}, \bar{s}}$  sits in.

      We need to show that, given two $\OO_F$-points $(y_1, f_1)$ and $(y_2, f_2)$ as above (we abbreviate as $f_1$ and $f_2$), then $f_1$ and $f_2$ are in the same $\GG$-orbit in $\mathbb{M}_{\GG, \mu}$ if and only if $f_1'$ and $f_2'$ are in the same $\GG_h$-orbit in $\mathbb{M}_{\GG_h^g, \mu_h}$. Due to the discussion above Proposition \ref{prop: boundary of Newton strata}, we furthur assume $\mu_h = \mu$ from now on.

       Note that the trivializations $f_1$, $f_2$ keeps the weight filtration, and $W_{-1} \subset F^0 \subset W_{-0}$. $f_i$ sends the Hodge filtration of the de Rham cohomology to a filtration $F^{\bullet}_{f_i}$ on $V_{\Z_{(p)}}$, this filtration determines and is determined by $0 \subset F^0_{f_i}/W_{-1}(V) \subset W_{0}(V)/W_{-1}(V)$, where $W_{\bullet}(V)$ is the fixed weight filtration on $V_{\Z_{(p)}}$ (induced from $V$). In particular, we see that $f_1$ and $f_2$ have same image in $\mathbb{M}_{\GG, \mu}$ (viewed as element in $\mathbb{M}_{\GSP, \mu}$ via the closed embedding $\mathbb{M}_{\GG, \mu} \to \mathbb{M}_{\GSP, \mu}$) if and only if $f_1'$ and $f_2'$ have same image in $\mathbb{M}_{\GG_h^g, \mu_h}$ (viewed as element in $\mathbb{M}_{\GSP_h^g, \mu_h}$ via the closed embedding $\mathbb{M}_{\GG_h^g, \mu_h} \hookrightarrow \mathbb{M}_{\GSP_h^g, \mu_h}$).

        We denote by  $x_1$ and $x_2$ the image of $f_1$ and $f_2$ in $\mathbb{M}_{\GG, \mu}$ respectively, and denote by $x_1'$ and $x_2'$ the image of $f_1'$ and $f_2'$ in $\mathbb{M}_{\GG_h^g, \mu_h}$ respectively. We use the properness of the local models. Let $\tilde{x}_i \in G/P_{\mu^{-1}}$ be the generic fiber of $x_i$, it determines and is uniquely determined by $x_i$, consider the diagram:
\[\begin{tikzcd}
	& {Q/P_{Q, \mu^{-1}}} & {G/P_{\mu^{-1}}} \\
	{G_h/P_{G_h, \mu^{-1}}} & {L/P_{L, \mu^{-1}}}
	\arrow[hook, from=1-2, to=1-3]
	\arrow[from=1-2, to=2-2]
	\arrow[hook, from=2-1, to=2-2]
\end{tikzcd}\]
        here $P_{(\ast), \mu^{-1}} = \lrbracket{g \in (\ast)| \ \lim_{t \to 0} \mu(t)^{-1}g\mu(t) \ \textit{exists}}$. Note that we have already fixed the weight filtration, thus $\tilde{x}_i \in Q/P_{Q, \mu^{-1}}$. Their projections in $L/P_{L, \mu^{-1}}$ indeed factor through $G_h/P_{G_h, \mu^{-1}}$, and coincides with $\tilde{x}_i'$.

        $(\Longrightarrow)$: If $x_1$ and $x_2$ are in the same $\GG$-orbit of $\mathbb{M}_{\GG, \mu}$, then they are in the same $\GG^g$-orbit of $\mathbb{M}_{\GG^g, \mu}$. Let $x_1 = hx_2$ for some $h \in \GG^g(\OO_{F})$. Since $h$ fixes the weight filtration, $h \in \mathcal{Q}^g(\OO_{F})$, here $\mathcal{Q}^g$ is the closure of $Q \cong Q^g$ in $\GG^g$. Consider the projection $\mathcal{Q}^g \to \LL^g$, let $l$ be the image of $h$. By the above diagram, $\tilde{x}_1' = l\tilde{x}_2' \in G_h/P_{G_h, \mu^{-1}}$, thus $x_1' = lx_2'$, $x_1'$ and $x_2'$ are in the same $\GG_h^g$-orbit, due to Lemma \ref{lem: adm injection}.

      $(\Longleftarrow)$: If $x_1' = lx_2'$ for some $l \in \GG_h^g(\OO_F)$, we can lift $l$ to some $h \in \mathcal{P}^g(\OO_F)$, here $\mathcal{Q}^g$ is the closure of $Q \cong Q^g$ in $\GG^g$. Then $\tilde{x}_1$ and $h\tilde{x}_2$ have same image in $L/P_{L, \mu^{-1}}$. Multiplying elements in $W(F)$ does not change the image of $f$ inside $\mathbb{M}_{\GG, \mu}$ since it stablizes the Hodge-filtration, therefore $x_1$ and $x_2$ are in the same $\GG$-orbit.
   \end{proof}

   \begin{lemma}\label{lem: adm injection}
       $\Adm(\lrbracket{\mu})_{\KKc_h} \to \Adm(\lrbracket{\mu})_{\KKc_L}$ is a bijection.
   \end{lemma}
   \begin{proof}
       Consider the commutative diagram:
\[\begin{tikzcd}
	{\Adm(\lrbracket{\mu})_{\KKc_h}} & {\Adm(\lrbracket{\mu})_{\KKc_L}} \\
	{\Adm(\lrbracket{\mu^{\ad}})_{\KK^{\ad, \circ}_h}} & {\Adm(\lrbracket{\mu^{\ad}})_{\KK^{\ad, \circ}_L}}
	\arrow[from=1-1, to=1-2]
	\arrow[from=1-1, to=2-1]
	\arrow[from=1-2, to=2-2]
	\arrow[from=2-1, to=2-2]
\end{tikzcd}\]
   Even though $\KKc_h \to \KK^{\ad, \circ}_h$ and $\KKc_L \to \KK^{\ad, \circ}_L$ are not surjective, with the help of \cite[Lemma 5.1.4]{shen2021ekor}, the vertical morphisms are still bijections. The bottom morphism is bijective since $L^{\ad} = G_h^{\ad} \times G_l^{\ad}$, $\mu^{\ad}$ is trivial in $G_l^{\ad}$, $\Adm(\lrbracket{\mu^{\ad}})_{\KK_L^{\circ, \ad}} = \Adm(\lrbracket{\mu^{\ad}})_{\KK_h^{\circ, \ad}} \times \Adm(\lrbracket{t^0})_{\KK_l^{\circ, \ad}}$, $\Adm(\lrbracket{\mu^{\ad}})_{\KK_l^{\circ, \ad}} = \Adm(\lrbracket{t^0}) = \lrbracket{\identity} \in \wdt{W}_{G_l^{\ad}}$.
   \end{proof}

  \subsection{Local model diagram}
   \begin{remark}
       We would like to thank Alex Youcis for pointing out a mistake in the proof of \cite[Proposition 1.0.5]{mao2024}. The proposition says that the smooth morphism $\lambda: \Shum{K} \to [\GG\backslash \lcM]$ (local model diagram constructed in \cite{kisin2018integral} and \cite{kisin2024integral}) extends to a smooth morphism $\lambda^{\tor}: \Shumc{K, \Sigma} \to [\GG\backslash \lcM]$. We will include this proposition as a part of results in a forthcoming work joint with Peihang Wu.
   \end{remark}

%% file: sections/EKOR_strata.tex
    In this section, we work under the assumption \ref{general condition} and moreover assume $G_{\Q_p}$ is unramified. Let $\Kc = K$ be a hyperspecial group, $G_{\Z_p} = \GG$ is a reductive group scheme, its special fiber $G_{\F_p}$ is also a reductive group with same root datum as $G$. We follow notations from \cite{zhang2018ekedahl}. The geometric object which serves as the index set for EO strata is the algebraic $G$-zip. In \cite{zhang2018ekedahl}, the author consturcted a smooth morphism $\zeta: \Shum{K, \bar{s}}\to G\textit{-}\ZIP^{\mu^{-1}}$. Note that we use $\mu^{-1}$ instead of $\mu$ to make the sign compatible, see \cite[Notation 3.2.3]{zhang2018ekedahl}.

    Let $w \in |G\textit{-}\ZIP^{\mu^{-1}}|$, Ekedal-Oort stratum $\EO_{K, w}$ is the fiber of $\zeta$ at $w$. Since $\zeta: \Shum{K}\to G\textit{-}\ZIP^{\mu^{-1}}$ is a smooth morphism, EO strata satisfy the closure relations
 \begin{equation}\label{eq: closure relations of EO}
   \ovl{\EO_{K, w}} = \EO_{K, \leq w}:=\sqcup_{w' \leq w} \EO_{K, w'}.
 \end{equation}

\subsection{Well-position}

\begin{proposition}\label{prop: EO sataisfy the discrete fiber property}
    EO strata satisfy the assumption \ref{assumption: fiber discreteness}. In particular, EO strata and their connected components are well-positioned subsets of $\Shum{K, \Bar{s}}$, and they are well-positioned subschemes with respect to the induced reduced subscheme structures.
\end{proposition}
\begin{proof}
    The assumption \ref{assumption: fiber discreteness} on discrete fiber property is given by \cite[Corollary 4.3.7]{goldring2019strata}. EO strata on integral model of Siegel Shimura varieties are well-positioned due to \cite[Proposition 3.5.5]{lan2018compactifications}.
\end{proof}

Repeat the proof of Proposition \ref{prop: KR strata are well-positioned}, we have
   \begin{proposition} \leavevmode
       \begin{enumerate}
           \item Let $w\in|G\textit{-}\ZIP^{\mu^{-1}}|$, then EO strata $\EO_{K, w}$, $\EO_{K, \leq w}$ and their connected components are well-positioned subsets in $\Shum{K, \bar{s}}$, they are well-positioned subschemes with respect to the induced reduced subscheme structures.
        \item $\Shum{K, \bar{s}}$ (resp. $\Shumm{K, \bar{s}}$, $\Shumc{K, \Sigma, \bar{s}}$) is a disjoint union of connected components of $\EO_{K, w}$ (resp. $\EO_{K, w}^{\min}$, $\EO_{K, w, \Sigma}^{\tor}$) for those $w \in |G\textit{-}\ZIP^{\mu^{-1}}|$.
         \item Given any EO strata $\EO_{K, w_1}$, $\EO_{K, w_2}$,
\begin{align*}
     \ovl{\EO_{K, w_1}} \cap \EO_{K, w_2} \neq \emptyset &\Longleftrightarrow \ovl{\EO_{K, w_1}^{\min}} \cap \EO_{K, w_2}^{\min} \neq \emptyset \\ &\Longleftrightarrow \ovl{\EO_{K, w, \Sigma}^{\tor}} \cap \EO_{K, w_2, \Sigma}^{\tor} \neq \emptyset.
\end{align*}
          \item 
          \begin{equation}
      \ovl{\EO_{K, w, \Sigma}^{\tor}} = \EO_{K, \leq w, \Sigma}^{\tor} = \bigsqcup_{w' \leq w} \EO_{K, w', \Sigma}^{\tor},\quad \ovl{\EO_{K, w}^{\min}} =  \EO_{K, \leq w}^{\min} = \bigsqcup_{w' \leq w} \EO_{K, w'}^{\min}. 
\end{equation}
       \end{enumerate}
   \end{proposition}

\subsection{Boundary descriptions}

    As in the proof of Proposition \cite[Proposition 3.5.5]{lan2018compactifications}, on integral model of Siegel Shimura variety, the boundary of an EO stratum is again an EO stratum, and Lemma \ref{lemma: boundary of Newton no intersection, siegel} also works for EO strata. Note that $K_{\Phi, h, p}$ are hyperspecial when $K_p$ is, and the good embeddings $\iota_{[\Phi]}$ are automatically very good embeddings due to \cite[Theorem 4.4.3, Proposition 5.3.10]{kisin2024integral}, we could consider the EO strata at boundary without furthur assumptions. Apply Proposition \ref{prop: boundaries of fiber discrete strata}.
   \begin{corollary}
       EO strata satisfy assumption \ref{assumption: fiber discreteness, enhanced}, where we take both $\lrbracket{Z_i}_{i \in I_{\Phi}}$, $\lrbracket{Z_j}_{j \in J_{\Phi^{\dd}}}$ be EO strata of $\Shum{K_{\Phi, h}, \bar{s}}$, $\Shum{K_{\Phi^{\dd}, h}^{\dd}, \bar{s}}$ respectively. Then given an EO stratum $\EO_{K, w} \subset \Shum{K, \bar{s}}$, its boundary stratum $\EO^{\natural}_{K, w, \Zb^{\bigsur}(\Phi)} \subset \Zb^{\bigsur}(\Phi)_{\bar{s}}$ (for each $[\Phi] \in \Cusp_K(G, X)$) is a union of connected components of some EO strata in $\Zb^{\bigsur}(\Phi)_{\bar{s}} = \Shum{K_{\Phi, h}, \bar{s}}$. 
   \end{corollary}
   To get deeper results, let us recall the construction of EO strata.

      We will freely use the notations from \cite{zhang2018ekedahl}. Note that $G\textit{-}\ZIP^{\mu^{-1}}$ can be viewed as the quotient stack $[E_{G, \mu^{-1}} \backslash G]$, let us introduce more details about this $E_{G, \mu^{-1}}$-torsor over $\Shum{K, \bar{s}}$. Fix a geometric point $\bar{x} \to \Shum{K, \bar{s}}$, lift $\bar{x}$ to a $W(k(\bar{x}))$-point $\tilde{x} \to \Shum{K}(W(k(\bar{x})))$, a $W(k(\bar{x}))$-point in the $E_{G, \mu^{-1}}$-torsor $\wdt{\Shum{K, \bar{s}}} \to \Shum{K, \bar{s}}$ is given by $(\bar{x}, g_{t}, (p_+, p_-))$ with a chosen $g_{t} \in G(W(k(\bar{x}))_{\Q})$ (depends on $\bar{x}$ and $t$) and $(p_+, p_-) \in E_{G, \mu^{-1}}(k(\bar{x}))$, let us explain how $g_{t}$ is chosen in \cite{zhang2018ekedahl}: 

       Let $M = H^1_{\dR}(\ab_{\GG, \tilde{x}}/W(k(\bar{x}))) \cong H^1_{\crys}(\ab_{\GG, \bar{x}}/W(k(\bar{x})))$, the absolute Frobenious on $\ab_{\GG, \bar{x}}$ induces a $\sigma$-linear map $\varphi: M \to M$ whose linearization will be denoted as $\varphi^{\mathrm{lin}}$. Consider the Hodge filtration $M^1 \subset M$,then $M^1_{k(\bar{x})}$ is the kernel of the Frobenius $\varphi$ on $H^1_{\dR}(\ab_{\GG, \tilde{x}}/W(k(\bar{x})))$. In particular, $\varphi(M^1) \subset pM$, $\varphi/p: M^1 \to M$ is well-defined. We can pick any splitting $M = M^0 \oplus M^1$ ($M^1$ is a direct summand of $M$), then the linear map
       \[ \alpha: M^{(\sigma)} = M \otimes_{W(k(\bar{x})), \sigma} W(k(\bar{x})) = M^{0(\sigma)} \oplus M^{1(\sigma)} \stackrel{\varphi^{\mathrm{lin}}|_{M^{0(\sigma)}} + (\varphi/p)^{\mathrm{lin}}|_{M^{1(\sigma)}}}{\longrightarrow} M\]
       is an isomorphism (\cite[Lemma 3.2.2]{zhang2018ekedahl}).

       Consider the trivial $G_{\Z_p} \otimes W(k(\bar{x}))$-torsor $I_{\tilde{x}}$ over $\Spec W(k(\bar{x}))$:
       
       \[ I_{\tilde{x}} := \Isom_{W(k(\bar{x}))}((V_{\Z_p}^{\vee}\otimes W(k(\bar{x})), (s_{\alpha}\otimes \identity)), (\mathcal{V}_{\tilde{x}}, (t_{\alpha, \dR, \tilde{x}}))),\] 
       here $\mathcal{V}$ is the relative de Rham cohomology $H^1_{\dR}(\ab_{\GG})$. Fix $t \in I_{\tilde{x}}(W(k(\bar{x})))$ from now.
       
       Let $\mu^{\prime}: (\Gm)_{W(k(\bar{x}))} \to G_{\Z_p} \otimes W(k(\bar{x}))$ be a cocharacter such that the inverse of $\mu^{\prime}_{\bar{x}}:= \bar{x}\mu^{\prime}\bar{x}^{-1}$ induces the Hodge filtration on $\mathcal{V}_{\tilde{x}}$. Then the Hodge filtration $\mathcal{V}_{\tilde{x}}^1 \subset \mathcal{V}_{\tilde{x}}$ is defined by $\mu^{\prime, -1}$. To be more precise, $\mathcal{V}_{\tilde{x}}^1$ is the weight $1$ subspace under $\mu^{\prime, -1}$. In the above paragraph, we can take $M:=\mathcal{V}_{\tilde{x}}$, $M^1:=\mathcal{V}_{\tilde{x}}^1$, $M^0:=\mathcal{V}_{\tilde{x}}^0$, where $\mathcal{V}_{\tilde{x}}^0$ is the weight $0$ subspace under $\mu^{\prime, -1}$. Let $\xi$ be the $W(k(\bar{x}))$-linear isomorphism $V^{\vee}_{\Z_p} \otimes W(k(\bar{x})) \rightiso (V^{\vee}_{\Z_p} \otimes W(k(\bar{x})))^{(\sigma)}$ given by $v\otimes w\mapsto v \otimes 1 \otimes w$, $\xi_t = t^{(\sigma)} \circ \xi $, then $g_{t}$ is defined as $t^{-1}\circ \alpha \circ \xi_t$, it is an element of $G_{\Z_p}(W(k(\bar{x})))$. Let $(\mathcal{V}_{\tilde{x}})_0$ (resp. $(\mathcal{V}_{\tilde{x}})_1$) be the sub $W(k(\bar{x}))$-module of $\mathcal{V}_{\tilde{x}}$ that is generated by $\varphi(\mathcal{V}_{\tilde{x}}^0)$ (resp. $\frac{\varphi}{p}(\mathcal{V}_{\tilde{x}}^1)$), then $V_{\Z_p}^{\vee} \otimes W(k(\bar{x})) = t^{-1}((\mathcal{V}_{\tilde{x}})_0) \oplus t^{-1}((\mathcal{V}_{\tilde{x}})_1)$ is the splitting induced by $\nu = g_{t} \mu^{\prime(\sigma)}g_{t}^{-1}$, to be more precise, $t^{-1}((\mathcal{V}_{\tilde{x}})_i)$ is the weight $i$ subspace under $\nu^{-1}$.

       Now we explain how $E_{G, \mu^{-1}}$ acts on such $g_{t}$. Fix a Hodge cocharacter $\mu_{W(\kappa)}$ (\cite[\S 3.2.2]{zhang2018ekedahl}), $P_+$ (resp. $P_-$) be the parabolic subgroups of $G_{\Z_p}\otimes W(\kappa)$ ($\kappa$ is the residue field of $E_v$) such that its Lie algebra is the sum of spaces with non-negative weights (resp. non-positive weights) in $\Lie G_{\Z_p}\otimes W(\kappa)$ under $\ad \circ \mu^{-1}$, $U_+$ (resp. $U_-$) be the unipotent radical of $P_+$ (resp. $P_-$), $L$ be the common Levi subgroup of $P_+$ and $P_-$. Let $C^{\bullet}$ (resp. $D_{\bullet}$) be the descending (resp. ascending) filtration on $V^{\vee}_{\Z_p} \otimes W(\kappa)$ such that $C^i$ (resp. $D_j$) is the weight $\geq i$ submodule (resp. weight $\leq j$ submodule) with respec to $\mu$ (resp. $\mu^{(\sigma)}$). Then $P_+$ (resp. $P_-$) is the stabilizer of $C^{\bullet}$ (resp. $D_{\bullet}$). Let $S$ be a scheme over the finite filed $\kappa$, Recall that
       \[  E_{G, \mu^{-1}}(S):= \lrbracket{(p_+:=lu_+, p_-:=l^{(\sigma)}u_-): l \in L(S), u_+ \in U_+(S), u_- \in U_-^{(\sigma)}(S)} \]
       here $L \to L^{(\sigma)}$ is the relative Frobenious of $L$. $E_{G, \mu^{-1}}(S)$ acts on $G_{\kappa}(S)$ from the left via $(p_+, p_-)\cdot g := p_+gp^{-1}_-$.
       
       Note that $\mu^{\prime}_t:= t\mu't^{-1}$ is $G_{\Z_p}(W(k(\bar{x})))$-conjugated to the fixed Hodge cocharacter $\mu_{W(\kappa)}$, we can take a $g_1 \in G_{\Z_p}(W(k(\bar{x})))$ such that $g_1(\mu \otimes W(k(\bar{x})))g_1^{-1} = \mu'$. We can consider the following trivial $P_+$ and $P^{\sigma}_-$ torsors respectively (\cite[Corollary 3.2.7]{zhang2018ekedahl}):
       \[ I_{\tilde{x}, +} := \Isom_{W(k(\bar{x}))}((V_{\Z_p}^{\vee}\otimes W(k(\bar{x})), (s_{\alpha}\otimes \identity), C^{\bullet}), (\mathcal{V}_{\tilde{x}}, (t_{\alpha, \dR, \tilde{x}}), \mathcal{V}_{\tilde{x}} \supset \mathcal{V}_{\tilde{x}}^1)) \subset I_{\tilde{x}},  \]
       \[ I_{\tilde{x}, -} := \Isom_{W(k(\bar{x}))}((V_{\Z_p}^{\vee}\otimes W(k(\bar{x})), (s_{\alpha}\otimes \identity), D_{\bullet}), (\mathcal{V}_{\tilde{x}}, (t_{\alpha, \dR, \tilde{x}}), (\mathcal{V}_{\tilde{x}})_0 \subset \mathcal{V}_{\tilde{x}})) \subset I_{\tilde{x}}, \]
       then $I_{\tilde{x}, +} = t \cdot g_1(P_+\otimes W(k(\bar{x})))$, $I_{\tilde{x}, -} = t \cdot g_{t}g_1^{(\sigma)}(P_-\otimes W(k(\bar{x})))^{(\sigma)}$. By replacing $t$ with $t g_1$ (that is to say, in the trivialization of $I_{\tilde{x}}$, we assume this trivialization comes from the trivialization $I_{\tilde{x}, +} \subset I_{\tilde{x}}$). Under this choice, $I_{\tilde{x}, +} = t (P_+\otimes W(k(\bar{x})))$, $I_{\tilde{x}, -} = t \cdot g_{t}(P_-\otimes W(k(\bar{x})))^{(\sigma)}$. We will always pick $t \in  I_{\tilde{x}, +}(W(k(\bar{x})))$ for each point $\bar{x}$.

       Let $I_{\bar{x}}$, $I_{\bar{x}, +}$, $I_{\bar{x}, -}$ be the special fiber of $I_{\tilde{x}}$, $I_{\tilde{x}, +}$, $I_{\tilde{x}, -}$ respectively. $(I_{\bar{x}}, I_{\bar{x}, +}, I_{\bar{x}, -}, \iota_{\bar{x}})$ determines a $G$-$\ZIP^{\mu^{-1}}$ over $\bar{x}$, here $\iota_{\bar{x}}$ is an isomorphism of $L^{(\sigma)}$-torsor $I_{\bar{x}, +}^{(\sigma)}/U_+^{(\sigma)}\rightiso I_{\bar{x}, -}/U_-^{(\sigma)}$, please refer to the concrete construction in \cite[\S 3.2.6]{zhang2018ekedahl}.

       Given $g \in G(k(\bar{x}))$, we can construct a standard $G$-$\ZIP^{\mu^{-1}}(k(\bar{x}))$: $I_{\underline{g}} = (I_g, I_{g, +}, I_{g, -}, \iota_g)$. Here $I_g = k(\bar{x}) \times_{\kappa} G_{\kappa}$, $I_{g, +} = k(\bar{x}) \times_{\kappa} P_+ \subset I_{g}$, and $I_{g, -} = g(k(\bar{x}) \times_{\kappa} P_-^{(\sigma)}) \subset I_g^{(\sigma)} \cong I_g$, and the left multiplication by $g$ induces an isomorphism of $L^{(\sigma)}$-torsors 
       \[ \iota_{g}: I_{g, +}^{(\sigma)}/U_+^{(\sigma)} = k(\bar{x}) \times_{\kappa} P^{(\sigma)}_+/U^{(\sigma)}_+ \cong k(\bar{x}) \times_{\kappa} P^{(\sigma)}_-/U^{(\sigma)}_- \rightiso g(k(\bar{x}) \times_{\kappa} P^{(\sigma)}_-/U^{(\sigma)}_-) = I_{\bar{x}, -}/U_-^{(\sigma)}.  \]

         Given $f = (p_+, p_-) \in E_{G, \mu^{-1}}(k(\bar{x}))$, it maps $I_{\underline{g}}$ to $I_{\underline{g}'}$ for some $g' \in G(k(\bar{x}))$. $f = f_+: I_{g, (+)} \to I_{g', (+)}$ is given by left action by $p_+$, and $f_-: I_{g, -} = g(k(\bar{x})\times_{\kappa} P^{(\sigma)}_-) \to I_{g', -} = g'(k(\bar{x})\times_{\kappa} P^{(\sigma)}_-)$ is given by left action $g'p_-g^{-1} = p_+$. In particular, $g'$ has to be $p_+gp_-$. And any morphism $I_{\underline{g}}$ to $I_{\underline{g}'}$ is given in this way.

       Given $f = (p_+, p_-) \in E_{G, \mu^{-1}}(k(\bar{x}))$, it acts on $I_{\bar{x}}$ and $I_{\bar{x}, +}$ by right action of $p_+$, and acts on $I_{\bar{x}, -}$ by right action of $p_-$, and it maps $\iota_{\bar{x}}$ to $\iota_{\bar{x}} \circ l^{(\sigma)}$ (this can be seen by tracing the construction of $\iota_{\bar{x}}$ in \cite[\S 3.2.6]{zhang2018ekedahl}), the following diagram commutes:
\[\begin{tikzcd}
	{I_{\bar{x}, +}^{(\sigma)}/U_+^{(\sigma)}} & {I_{\bar{x}, +}^{(\sigma)}/U_+^{(\sigma)}} \\
	{I_{\bar{x}, -}/U_-^{(\sigma)}} & {I_{\bar{x}, -}/U_-^{(\sigma)}}
	\arrow["{\cdot l^{(\sigma)}}", from=1-1, to=1-2]
	\arrow["{\iota_{\bar{x}}}"', from=1-1, to=2-1]
	\arrow["{\iota_{\bar{x}}\circ l^{(\sigma)}}"', from=1-2, to=2-2]
	\arrow["{\cdot l^{(\sigma)}}", from=2-1, to=2-2]
\end{tikzcd}\]

       Once we fix a trivialization $t$ of $I_{\bar{x}, +}$ as before, we get a standard $G$-$\ZIP^{\mu^{-1}}$ over $x$ 
       \[ (I_{\bar{x}}, I_{\bar{x}, +}, I_{\bar{x}, -}, \iota_{\bar{x}}) \cong (I_{g_{t}}, I_{g_{t}, +}, I_{g_{t}, -}, \iota_{g_{t}}) =: I_{\underline{g}_t}. \]
       Then the action of $f = (p_+, p_-) \in E_{G, \mu^{-1}}(k(\bar{x}))$ on $(I_{\bar{x}}, I_{\bar{x}, +}, I_{\bar{x}, -}, \iota_{\bar{x}})$ can be realized by its action on $I_{\underline{g}}$ introduced above.

     \begin{proposition}\label{prop: boundary of EO strata}
         Let $Y \subset \Shum{K, s}$ be an EO stratum, then its boundary stratum $Y^{\natural}_{\Zb^{\bigsur}(\Phi)} \subset \Zb^{\bigsur}(\Phi)_{s}$ (for each $[\Phi]\in \Cusp_K(G, X)$) is either empty or an EO strata on $\Zb^{\bigsur}(\Phi)_{s}$.
     \end{proposition}
     \begin{proof}
         We apply similar arguments as in the proof of Proposition \ref{prop: boundary of KR is KR}. Given a geometric point $\bar{x} \to W^0$, we lift it to a $W(k(\bar{x}))$-point $\tilde{x} \to W^0$, we can always take the section $t \in I_{\tilde{x}}$ preserving the weight filtrations, and its conjugated reduction $t\circ g$ in the $\gr^W_{-1}$-piece gives a section in $I_{\tilde{z}}$, here $\wdt{z}$ is the composition $\wdt{x} \to \Shum{K_{\Phi, h}}$, and $I_{\tilde{z}}$ is the similarly defined $G_{h, \Z_p}^g\otimes W(k(\bar{x}))$-torsor over $\Spec W(k(\bar{x}))$ using the integral model $\Shum{K_{\Phi, h}}$ in place of $\Shum{K}$. Note that on the boundary, everything needs to be conjugated by $g$ (coming from cusp label $\Phi = ((P, X_{P, +}), g)$) as in the proof of Proposition \ref{prop: boundary of Newton strata}, \ref{prop: boundary of central leaves}, \ref{prop: boundary of KR is KR}. Without lose of generality, we take $g = 1$.
         
         Recall that $g_{t} = t^{-1} \circ \alpha \circ \xi_t$, and $\alpha$ keeps the weight filtration (Frobenius action does not change the weight filtration), thus $g_{t} \in Q_{\Z_p}(W(k(\bar{x})))$. Consider its projection $g_{t}' \in L_{\Z_p}(W(k(\bar{x})))$, it falls in $G_{h, \Z_p}(W(k(\bar{x})))$. As in the proof of Corollary \ref{prop: boundary of KR is KR}, $W_{-1} \subset \mathcal{V}^1 \subset W_{0}$, when $t \in I_{\tilde{x}, +} \subset I_{\tilde{x}}$, then its reduction is in $I_{\tilde{z}, +} \subset I_{\tilde{z}}$. This pair $(g_{t}, g_{t}')$ gives the index of the EO stratum that the images of $\bar{x}$ fall into under $W^0_{\bar{s}} \to \Shum{K, \bar{s}}$ and under $W^0_{\bar{s}} \to \Shum{K_{\Phi, h}, \bar{s}}$ respectively.
         
         Given two geometric points $\bar{x}, \bar{x}' \to W^0_s$, if they are in the same EO stratum under $W^0_{\bar{s}} \to \Shum{K, \bar{s}}$, then $g_{t}$ and $g_{t'}$ are in the same $E_{G, \mu^{-1}}$-orbit. Let us write $g_{t'} = p_+g_{\bar{t}}p_-^{-1}$. Since $(I_{\bar{x}}, I_{\bar{x}, +}, I_{\bar{x}, -}, \iota_{\bar{x}})$ and $(I_{\bar{x}'}, I_{\bar{x}', +}, I_{\bar{x}', -}, \iota_{\bar{x}'})$ are isomorphic, we can adjust $(I_{\bar{x}}, I_{\bar{x}, +}, I_{\bar{x}, -}, \iota_{\bar{x}})$ by self-isomorphism (which does not change the EO stratum that $\bar{x}$ falls into) such that the trivialization $t$ maps to $t'$, and that $p_+$ as well as $p_-$ keep the weight filtration, thus $p_+, p_- \in Q_{\Z_p}(W(k(\bar{x})))$. Consider their reductions $p_+', p_-' \in L_{\Z_p}(W(k(\bar{x})))$, then $g_{t'}' = p_+'g_{t}'p_-^{\prime, -1}$. Then the images of $\bar{x}, \bar{x}'$ in $\Shum{K_{\Phi, h}, \bar{s}}$ are in the same EO stratum in $\Shum{K_{\Phi}, s}$ due to Lemma \ref{lem: E bijection}.

         Now given $\bar{x}, \bar{x}' \to W^0_{\bar{s}}$, assume their images in $\Shum{K_{\Phi, h}, \bar{s}}$ are in the same EO stratum, we show that their images in $\Shum{K, \bar{s}}$ are in the same EO stratum. We use the notations from above paragraph, assume $g_{t'}' = p_+'g_{t}'p_-^{\prime, -1}$ for some $p_+' \in P_{h, +, \Z_p}(W(k(\bar{x})))$, $p_-' \in P_{h, -, \Z_p}(W(k(\bar{x})))$ (here $P_{h, \pm}$ are defined similarly in $G_h$ as we define $P_{\pm}$ in $G$ using the cocharacter $\mu_h = \mu$). On the other hand, it follows from the construction that the action of both $g_{t}$ and $g_{t'}$ factor through graded pieces of the weight filtrations, we can furthur assume $g_{t}$ and $g_{t'}$ are in $L_{\Z_p}(W(k(\bar{x})))$. Therefore, use the natural embeddings $P_{h, \pm} \to P_{\pm}$ induced by $G_h \to G$, we have $g_{t'} = p_+g_{\bar{t}}p_-^{-1}$, $g_{\bar{t}}$ and $g_{t'}$ are in the same $E_{G, \mu^{-1}}$-orbit.

     \end{proof}

     \begin{lemma}\label{lem: E bijection}
         $|[E_{G_h, \mu^{-1}}\backslash G_h]| \to |[E_{L, \mu^{-1}}\backslash L]|$ is a bijection.
     \end{lemma}
     \begin{proof}
         Recall that $|[E\backslash G]| = \prescript{I}{}{W_G}$, where $W_G$ is the Weyl group of $G$, $I$ is the type of parabolic group $P$ determined by $\mu^{-1}$ and $W_{G, I}\subset W_G$ is the group generated by the reflections in $I$, and $\prescript{I}{}{W_G}$ is the set of all $w\in W_G$ that have minimal length in the cosets $W_{G, I}w$. $\prescript{I}{}{W_G}$ is a set of representatives of left cosets $W_{G, I}\backslash W_G$. Since $W_L = W_h \times W_l$ (here $W_l$ is the Weyl group of the quotient reductive group $G_l = L/G_h$), $W_{L, I_L} = W_{h, I_h} \times W_{l, I_l}$ and since $\mu$ is trivial in $G_l$, then $W_{l, I_l} = W_l$, thus $W_{G_h, I_h}\backslash W_{G_h} \rightiso W_{L, I_L}\backslash W_L$.
     \end{proof}

\subsection{Comparison}

Since the partial minimal compactifications of EO strata on the integral model of Siegel modular varieties are affine, see \cite[Theorem C]{boxer2015torsion}, due to Propsition \ref{prop: affineness in general} and \cite[Corollary 4.3.7]{goldring2019strata}:
  \begin{corollary}{\cite[Proposition 6.3.1]{goldring2019strata}}\label{cor: min EO are affine}
      Partial minimal compactifications of EO strata on the integral model of Hodge-type Shimura varieties are affine.
  \end{corollary}
   Recall that the smooth morphism $\zeta: \Shum{K, \bar{s}}\to G\textit{-}\ZIP^{\mu^{-1}}$ extends to a morphism $\zeta^{\tor}_{\Sigma}: \Shumc{K, \Sigma, \bar{s}}\to G\textit{-}\ZIP^{\mu^{-1}}$ due to \cite[Theorem 6.2.1]{goldring2019strata}. Denote by $(\Shumc{K, \Sigma, \bar{s}})^{w}$ the fiber of $\zeta^{\tor}_{\Sigma}$ at any $w\in |G\textit{-}\ZIP^{\mu^{-1}}|$.

   \begin{corollary}
       $(\Shumc{K, \Sigma, \bar{s}})^{w} = \EO^{\tor}_{w, \Sigma}$.
   \end{corollary}
     \begin{proof}
         Consider the commutative diagram as in the proof of \cite[Proposition 6.3.1]{goldring2019strata}:
\[\begin{tikzcd}
	& {\Shumc{K, \Sigma, \bar{s}}} && {\Shumc{K^{\dd}, \Sigma^{\dd}, \bar{s}}} \\
	{\Shum{K, \bar{s}}} && {\Shum{K^{\dd}, \bar{s}}} \\
	& {G\ZIP^{\mu^{-1}}} && {G^{\dd}\ZIP^{\mu^{\dd, -1}}}
	\arrow["\varphi", from=1-2, to=1-4]
	\arrow["{\zeta^{\tor}_{\Sigma}}"{pos=0.7}, from=1-2, to=3-2]
	\arrow["{\zeta^{\dd, \tor}_{\Sigma^{\dd}}}"', from=1-4, to=3-4]
	\arrow[from=2-1, to=1-2]
	\arrow["{\varphi^{\tor}_{\Sigma}}"{pos=0.4}, from=2-1, to=2-3]
	\arrow["\zeta", from=2-1, to=3-2]
	\arrow[from=2-3, to=1-4]
	\arrow["{\zeta^{\dd}}", from=2-3, to=3-4]
	\arrow["{\varphi^{\ZIP}}", from=3-2, to=3-4]
\end{tikzcd}\]
         \cite[Lemma 6.3.3]{goldring2019strata} shows that $(\Shumc{K^{\dd}, \Sigma^{\dd}, s})^{w'}$ coincides with $\EO^{\tor}_{w', \Sigma^{\dd}}$ for every $w' \in |G^{\dd}\ZIP^{\mu^{\dd, -1}}|$, and such defined $\zeta^{\dd, \tor}_{\Sigma^{\dd}}$ coincides with the one introduced in \cite[Corollary 3.5.8]{lan2018compactifications} (in the Siegel case). Due to \cite[Corollary 4.3.7]{goldring2019strata}, $\varphi^{\ZIP}$ had discrete fibers, thus $(\varphi^{\tor}_{\Sigma})^{-1}(\EO^{\tor}_{w', \Sigma^{\dd}})$ is a topologically disjoint union of $\EO^{\tor}_{w, \Sigma}$ for some $w \in |G\textit{-}\ZIP^{\mu^{-1}}|$ due to Lemma \ref{lemma: pullback of well-positioned is well-positioned} and Proposition \ref{proposition: open-closed subschemes are well positioned}, then $(\Shumc{K, \Sigma, s})^{w} = \EO^{\tor}_{w, \Sigma}$ follows from the above commutative diagram.
     \end{proof}

%% file: sections/closure.tex
In this part, we work with general Pappas-Rapoport integral models, and assume He-Rapoport axioms \cite{he2017stratifications} hold. Here in the index set we use $\mu$ instead of $\mu^{-1}$ to be compatible with \cite{he2017stratifications}.

Let $K$ be a parahoric subgroup (instead of using $\Kc$), and $I \subset K$ be an Iwahori subgroup. Let $\KK$ and $\Breve{I}$ be the $\bZ_p$-points of corresponding group schemes.

 \begin{remark}
     Under the setting of Kisin-Pappas integral models, and assume $\GG = \GGc$, then the He-Rapoport axioms (except for \cite[Axiom 4(c)]{he2017stratifications}) are verified in \cite{zhou2020mod}, and the remaining \cite[Axiom 4(c)]{he2017stratifications} is verified in \cite[Corollary 1.5]{gleason2023connected}.
 \end{remark}

\subsection{EKOR strata}\label{EKOR strata}

Since central leaves are well-positioned due to Proposition \ref{prop: Newton central are well-posiitoned, PR}, and EKOR strata satisfy closure relations 
\begin{equation}\label{eq: closure relations of EKOR}
   \ovl{\EKOR_{K, w}} = \EKOR_{K, \leq w}:= \sqcup_{w' \leq w} \EKOR_{K, w'}. 
 \end{equation}
due to \cite[Theorem 6.15]{he2017stratifications}, repeat the proof of Proposition \ref{prop: KR strata are well-positioned}, we have
\begin{proposition}\label{prop: EKOR are well-positioned}\leavevmode
   \begin{enumerate}
        \item Let $w \in \prescript{\KK}{}{\Adm(\lrbracket{\mu})}$, then the EKOR stratum $\EKOR_{K, w}$, $\EKOR_{K, \leq w}$ and their connected components are well-positioned subsets in $\Shum{K, \Bar{s}}$, they are well-positioned subschemes with respect to the induced reduced subscheme structures. 
        \item $\Shum{K, \bar{s}}$ (resp. $\Shumm{K, \bar{s}}$, $\Shumc{K, \Sigma, \bar{s}}$) is a disjoint union of connected components of $\EKOR_{K, w}$ (resp. $\EKOR_{K, w}^{\min}$, $\EKOR_{K, w, \Sigma}^{\tor}$) for those $w \in \prescript{\KK}{}{\Adm(\lrbracket{\mu})}$.
          \item Given any EKOR strata $\EKOR_{K, w_1}$, $\EKOR_{K, w_2}$,
          \begin{align*}
               \ovl{\EKOR_{K, w_1}} \cap \EKOR_{K, w_2} \neq \emptyset &\Longleftrightarrow \ovl{\EKOR_{K, w_1}^{\min}} \cap \EKOR_{K, w_2}^{\min} \neq \emptyset \\ &\Longleftrightarrow \ovl{\EKOR_{K, w, \Sigma}^{\tor}} \cap \EKOR_{K, w_2, \Sigma}^{\tor} \neq \emptyset.
          \end{align*}
          \item 
                    \begin{equation*}
      \ovl{\EKOR_{K, w, \Sigma}^{\tor}} = \EKOR_{K, \leq w, \Sigma}^{\tor} = \bigsqcup_{w' \leq w} \EKOR_{K, w', \Sigma}^{\tor},
      \end{equation*}
      \begin{equation*}
      \ovl{\EKOR_{K, w}^{\min}} =  \EKOR_{K, \leq w}^{\min} = \bigsqcup_{w' \leq w} \EKOR_{K, w'}^{\min}. 
\end{equation*}
    \end{enumerate}
\end{proposition}
\begin{remark}
   The first two statements only need the fact that EKOR strata are locally closed and are unions of central leaves. The later two statements furthur need the fact that EKOR strata satisfy the closure relation \ref{eq: closure relations of EKOR}. Other unrelated axioms in \cite{he2017stratifications} are not needed.
\end{remark}

    Under the setting of Kisin-Pappas integral models, assume $\GG = \GGc$ is connected, use the local definition of EKOR strata as in \cite[\S 3]{shen2021ekor}, EKOR strata can be realized as EO strata inside KR strata (using general $G$-$\ZIP$ instead of the one coming directly from the cocharacter $\mu^{-1}$). Assume the boundary morphisms are very good as in the discussions in subsection \ref{subsec: boundary descriptions, newton}, we can apply the arguments in the proof of Proposition \ref{prop: boundary of KR is KR} and \ref{prop: boundary of EO strata}, and show that given a EKOR stratum $Y \subset \Shum{K, \bar{s}}$, its boundary stratum $Y^{\natural}_{\Zb^{\bigsur}(\Phi)} \subset \Zb^{\bigsur}(\Phi)_{\bar{s}}$ (for each $[\Phi]\in \Cusp_K(G, X)$) is either empty or an EKOR strata on $\Zb^{\bigsur}(\Phi)_{\bar{s}}$. We do not spell out the details. Note that 
    $\prescript{\KKc_h}{}{\Adm(\lrbracket{\mu})} \to \prescript{\KKc_L}{}{\Adm(\lrbracket{\mu})}$ is a bijection, the proof is same as the proof of Lemma \ref{lem: adm injection}.

\subsection{Minimal EKOR strata}

    Due to \cite[Theorem 6.15]{he2017stratifications}, among EKOR strata, there is a minimal one, which is the unique closed EKOR strata associated to $\tau_{\lrbracket{\mu}}:=(\identity, \mu^{\sharp})\in W_{\aff}\rtimes\pi_1(G)_I\cong\widetilde{W}$.
     \begin{corollary}\label{cor: minimal EKOR has zero dim}
     	$\EKOR_{K, \tau_{\lrbracket{\mu}}}=\EKOR_{K, \tau_{\lrbracket{\mu}}, \Sigma}^{\tor}$.
     \end{corollary}
     \begin{proof}
     Recall that, $\tau_{\lrbracket{\mu}}$ has length $0$, and lies in $\prescript{\KK}{}{\Adm(\lrbracket{\mu})}$, and is the minimal element with respect to the partial order $\leq_{K, \sigma}$. By \cite[Theorem 6.21]{he2017stratifications}, $\EKOR_{K, \tau_{\lrbracket{\mu}}}$ has same dimension as $\KR_{I, \tau_{\lrbracket{\mu}}}$. Argue as in the proof of \cite[Lemma 4.3.7]{daniels2024igusa}, dimension of $\KR_{I, \tau_{\lrbracket{\mu}}}$ is the length of $\tau_{\lrbracket{\mu}}$. In particular, $\EKOR_{K, \tau_{\lrbracket{\mu}}}$ has dimension $0$.
     
     Since $\EKOR_{K, \tau_{\lrbracket{\mu}}}$ is closed, $\EKOR_{K, \tau_{\lrbracket{\mu}}, \Sigma}^{\tor}$ is the closure of $\EKOR_{K, \tau_{\lrbracket{\mu}}}$ in $\Shumc{K, \Sigma, \bar{s}}$. Since EKOR strata are well-positioned, due to Proposition \ref{prop: Y same as Shum}, the toroidal compactification $\EKOR_{K, \tau_{\lrbracket{\mu}}} \to \EKOR_{K, \tau_{\lrbracket{\mu}}, \Sigma}^{\tor}$ is \'etale locally an affine toroidal embedding, it has no boundary since $\EKOR_{K, \tau_{\lrbracket{\mu}}}$ has full codimension $\dim \Shum{K, \bar{s}}$.
     \end{proof}

     Since EO strata coincide with EKOR strata at hyperspecial level, this implies that
    \begin{corollary}{{\cite[Condition 6.4.2]{goldring2019strata}}}\label{cor: minimal EO has zero dim}
         The minimal EO stratum satisfies $\EO_e = \EO_{e, \Sigma}^{\tor}$.
     \end{corollary}

\subsection{Closure of irreducible components of EKOR strata}

We prove a conjecture in \cite{van2024mod} under some mild conditions. In this subsection, we need to assume the following:
\begin{assumptionaxiom}\label{assumption, res of Hodge bundle torsion}
    The restriction of the Hodge bundle on each Kottwitz-Rapoport stratum at Iwahori level is torsioned, in particular, each Kottwitz-Rapoport stratum at Iwahori level is quasi-affine. Note that, under the setting of Kisin-Pappas integral model and $\GGc = \GG$, this assumption was verified in \cite[Theorem 3.5.9]{shen2021ekor}.
\end{assumptionaxiom}

     \begin{lemma}\label{lemma: KR intersects with minimal KR}
     Assume $G^{\ad}$ is $\Q$-simple or has no $\PGL_2$-factor\footnote{This condition is only used in Proposition \ref{prop: KR are quasi-affine}, we expect this condition can be easily removed once we set up and verify the well-position of KR strata in abelian-type cases}, and assume \ref{assumption, res of Hodge bundle torsion}. Let $w\in\Adm(\lrbracket{\mu})_{\Breve{I}} = \Adm(\lrbracket{\mu})$ be any element, $Y$ be any irreducible component of the KR stratum $\KR_{I, w}$, then the closure of $Y$ in $\Shum{I, \bar{s}}$ intersects nontrivially with the minimal KR stratum $\KR_{I, \tau_{\lrbracket{\mu}}}$.
     \end{lemma}
     
     \begin{corollary}
         Keep same notations and assumptions as in Lemma \ref{lemma: KR intersects with minimal KR}. Let $\ovl{Y} \subset \KR_{I, \leq w}$ be an irreducible component, then $\ovl{Y}$ intersects nontrivially with the minimal KR stratum $\KR_{I, \tau_{\lrbracket{\mu}}}$.
     \end{corollary}
     \begin{proof}
         Since $\KR_{I, \leq w}$ is normal due to the existence of the local model diagram and \cite[Proposition 3.7]{anschutz2022p}, the concept of irreducible components coincides with connected components. If $\ovl{Y}$ is a connected component of the noetherian space $\KR_{I, \leq w}$, then $Y = \ovl{Y} \cap \KR_{I, w}$ is non-empty and open dense in $\ovl{Y}$, and is irreducible. In particular, $Y$ is a connected component of the smooth scheme $\KR_{I, w}$, then we apply Lemma \ref{lemma: KR intersects with minimal KR}. 
     \end{proof}

          \begin{corollary}{{\cite[Conjecture 4.3.1]{van2024mod}}}\label{corollary: Hoften 3.7.5 removed}
     	 Keep same notations and assumptions as in Lemma \ref{lemma: KR intersects with minimal KR}. Let $w\in\prescript{\KK}{}{\Adm(\lrbracket{\mu})}$ be any element, $Y$ be any irreducible component of $\EKOR_{K, \leq w} = \ovl{\EKOR_{K, w}}$, then $Y$ intersects nontrivially with the minimal stratum $\EKOR_{K, \tau_{\lrbracket{\mu}}}$.
     \end{corollary}
     \begin{proof}
        Under the inclusion $\prescript{\KK}{}{\Adm(\lrbracket{\mu})} \subset \Adm(\lrbracket{\mu})$, we regard $w \in \Adm(\lrbracket{\mu})$, then by \cite[Theorem 6.21]{he2017stratifications}, the proper surjection $\pi: \Shum{I} \to \Shum{K}$ restricts to a finite morphism:
     	\[  \pi_{w}: \KR_{I, w} \to \EKOR_{K, w},  \]
       and $\dim \KR_{I, w} = \dim \EKOR_{K, w}$. $\pi_{w}$ is surjective when \cite[axiom 4(c)]{he2017stratifications} holds.
     	Consider the restriction $\pi_{\leq w}$ of $\pi$, it is a closed morphism:
     	\[  \pi_{\leq w}: \KR_{I, \leq w} \hookrightarrow \pi^{-1}(\EKOR_{K, \leq w}) \to \EKOR_{K, \leq w}. \]

        Given any irreducible component $Y$ of $\EKOR_{K, \leq w}$, let $Y' = \pi_{\leq w}^{-1}(Y)$. Since $\pi_{\leq w}$ is a closed map, then $\pi_{\leq w}$ is surjective, $Y'$ is nonempty. Since the target and the source of $\pi_{\leq w}$ have the same dimension, and $\pi_{\leq w}$ is dominant, then $Y'$ is a union of irreducible components (which are also connected components) of $\KR_{I, \leq w}$. Lemma \ref{lemma: KR intersects with minimal KR} implies $Y'\cap \KR_{I, \tau_{\lrbracket{\mu}}}\neq\emptyset$, \cite[Theorem 6.21]{he2017stratifications} implies that $\pi(\KR_{I, \tau_{\lrbracket{\mu}}}) \subset \EKOR_{K, \tau_{\lrbracket{\mu}}}$, then
     	\[ Y\cap\EKOR_{K, \tau_{\lrbracket{\mu}}} \supseteq \pi_{\leq w}(Y')\cap \pi(\KR_{I, \tau_{\lrbracket{\mu}}}) \supseteq  \pi_{\leq w}(Y'\cap \KR_{I, \tau_{\lrbracket{\mu}}}) \neq \emptyset \] 
     \end{proof}

\begin{proposition}\label{prop: KR are quasi-affine}
    Keep same notations and assumptions as in Lemma \ref{lemma: KR intersects with minimal KR}. Let $w\in\Adm(\lrbracket{\mu})$. We abbreviate $\KR_{I, w}$ by $\KR_w$. Let $\KR_w^{\min, \stein}$ be the Stein factorization of $\oint_{I, \Sigma, w}:=\oint_{I, \Sigma}|_{\KR^{\tor}_{w, \Sigma}}: \KR^{\tor}_{w, \Sigma} \to \KR^{\min}_w$, then $\KR_w^{\min, \stein} \to \KR^{\min}_w$ is a finite morphism which induces a homeomorphism on the underlying topological space. Moreover, $\KR_{w}^{\min, \stein}$ is quasi-affine.
\end{proposition}
\begin{proof}
     Consider the Stein factorization $\stein: \KR_{w, \Sigma}^{\tor} \to \KR_w^{\min, \stein}$ of $\oint_{I, \Sigma, w}: \KR_{w, \Sigma}^{\tor} \to \KR_w^{\min}$. Since $\KR_w^{\min, \stein} \to \KR_w^{\min}$ is finite surjective and $\oint_{I, \Sigma, w}$ has geometric connected fibers due to Lemma \ref{lemma: exclusiveness is redundant}, then $\KR_w^{\min, \stein} \to \KR_w^{\min}$ is a universal homeomorphism. Since $\KR_w$ (resp. $\KR_{\leq w}$) is normal, then $\KR_{w, \Sigma}^{\tor}$ (resp. $\KR_{\leq w, \Sigma}^{\tor}$) is normal due to \cite[Proposition 2.3.13]{lan2018compactifications}, then $\KR_w^{\min, \stein}$ (resp. $\KR_{\leq w}^{\min, \stein}$) is normal due to \cite[Lemma 3.39]{mao2025boundary}. By assumption, each $\KR_w$ is quasi-affine. We show that $\KR_{w}^{\min, \stein}$ is quasi-affine as follows
     
    Consider the following three cases:
    \begin{enumerate}
        \item Assume $\Shum{I}$ is proper, then there is no boundary stratum, $\KR_w = \KR_w^{\min}$ is quasi-affine.
        
        \item Assume there exists a boundary stratum of $\Shumm{I}$ that has codimension $1$. Let $\Zb(\Phi) \rightiso \mathcal{Z}([\Phi])$ be such a stratum, $\Zb^{\bigsur}(\Phi) \to \Zb(\Phi)$ is finite and does not change the codimension. Since $\Zb^{\bigsur}(\Phi)$ is normal by construction, and is flat over $\OO_{E_{(v)}}$ due to \cite[Theorem 5.2.11]{pera2019toroidal}, thus the generic fiber $\Zb(\Phi)_{\eta}$ has codimension $1$ in $\shum{I}$, which forces $G^{\ad}$ to be $\PGL_2$ (here we use the assumption that $G^{\ad}$ is $\Q$-simple or has no $\PGL_2$-factor). Therefore, $\Shumm{I}$ is of relative dimension $1$ over $\OO_{E_{(v)}}$.

        Consider any $1$-dimensional KR stratum $\KR_w$, its closure $\KR_{\leq w}$ consists of $\KR_w$ and $\KR_{\tau_{\lrbracket{\mu}}}$, and the latter is a finite union of closed points and is non-empty due to the condition \cite[Axiom 5]{he2017stratifications}. Due to Lemma \ref{lemma: closure of well-pos is well-pos}, the complement of $\KR^{\min}_{w}$ in $\KR^{\min}_{\leq w}$ is a finite union of closed points. Since $\KR^{\min}_{\leq w}$ is closed in $\Shumm{I, \Bar{s}}$, thus is proper over $k$, then $\KR^{\min}_{w}$ is the complement of a nonempty union of finite points in a proper normal curve over $k$, we claim it is affine: A normal curve is regular, a regular proper curve with a nonempty union of finite points removed is affine. We apply this to the normal proper curve $\KR_w^{\min, \stein}$. Since $\KR_w^{\min, \stein} \to \KR_w^{\min}$ is finite surjective, then $\KR_w^{\min}$ has to be affine, it's a theorem of Chavelley, see \cite[\href{https://stacks.math.columbia.edu/tag/01ZT}{Tag 01ZT}]{stacks-project}.
        
        \item Assume all the boundary strata of $\Shumm{I}$ have codimension at least $2$, we first show that for each $\KR_w$, all the boundary strata of $\KR_w^{\min}$ has codimension at least $2$ (if non-empty): if not, then $\KR_w^{\min}$ has a stratum $\KR_{w, [\Phi]} \subset \mathcal{Z}([\Phi])$ which has codimension $1$ in $\KR_w^{\min}$. Consider the morphism $\pi_w: \KR_{w, \mathcal{Z}([\Phi, \sigma])} \to \KR_{w, [\Phi]}$ which is restricted from $\pi: \mathcal{Z}([\Phi, \sigma]) \to \mathcal{Z}([\Phi])$. Note that $\pi$ can be identified with $\Xi_{\sigma}(\Phi) \to C(\Phi) \to \Zb^{\bigsur}(\Phi) \to \Zb(\Phi)$, thus $\pi$ is a composition of smooth morphisms and a quasi-finite morphism. Since $\KR_{w, \mathcal{Z}([\Phi, \sigma])} = \KR_{w, \Xi_{\sigma}(\Phi)}^{\natural}$ is the pullback of $\KR_{w, \Zb^{\bigsur}(\Phi)}^{\natural}$ under $\pi$, the relative dimension of $\pi_w$ equals the relative dimension of $\pi$. The codimension of $\KR_{w, [\Phi]}$ in $\KR_w^{\min}$ is the difference between the dimension $\KR_{w, \Zb^{\bigsur}(\Phi)}^{\natural}$ and $\KR_w$. Note that $\KR_{w, \mathcal{Z}([\Phi, \sigma])}$ also has codimension at least $1$ in $\KR_{w, \Sigma}^{\tor}$, thus $\pi_w$ has relative dimension zero, $\pi$ has relative dimension zero. Note that there exists a boundary stratum $\mathcal{Z}([\Phi, \sigma])$ which has codimension $1$ in $\Shumc{I, \Sigma}$, thus $\mathcal{Z}[\Phi]$ has codimension $1$ in $\Shumm{I}$, contradicts with the assumption.
        
         Consider the Hodge bundle $\omega$ on $\Shum{I}$, it has an extension $\omega^{\tor}$ over $\Shumc{I, \Sigma}$. Its restrction (still denoted by) $\omega$ on $\KR_w$ is torsioned by assumption, i.e., $\omega^{\otimes n} = \OO_{\KR_w}$ for some $n$. Consider the restriction of $\omega^{\tor}$  on $\KR_w^{\tor}$, which we still denote by $\omega^{\tor}$. Since the boundary strata of $\KR_w^{\min}$ have codimensions at least $2$, the identity section $s \in \Gamma(\KR_w, \omega^{\otimes n})$ extends to a section $s^{\tor} \in \Gamma(\KR_{w, \Sigma}^{\tor}, \omega^{\tor, \otimes n})$ due to the Koecher's principle, see \cite[Theorem 2.5.11]{lan2018compactifications}. In \cite[Theorem 2.5.11]{lan2018compactifications}, the condition that $\OO \otimes_{\Z} \Q$ is simple is for PEL-type cases, and is not essential, the Koecher's principle also works in the Hodge-type case as well as well-positioned subschemes in it, see \cite[Remark 2.5.15]{lan2018compactifications}.

         Let $\omega^{\min, \stein}$ be the pullback of $\omega^{\min}$ under $\KR^{\min, \stein}_w \to \KR^{\min}_w$, then $\stein^*\omega^{\min, \stein} = \omega^{\tor}$. Since $\stein_{*}\OO_{\KR_{w, \Sigma}^{\tor}} = \OO_{\KR_w^{\min, \stein}}$, $\stein_*\stein^*\omega^{\min, \stein} \rightiso \omega^{\min, \stein}$ is an isomorphism due to the projective formula, the global section $s^{\tor}$ induces a global section $s^{\stein} \in \Gamma(\KR_w^{\min, \stein}, \omega^{\min, \stein, \otimes n})$.
    
         The vanishing locus of any section on a normal scheme is either nonempty or has codimension $1$. Since all the boundary strata of $\KR_w^{\min}$ (thus of $\KR_w^{\min, \stein}$) are of codimension at least $2$, and $s = s^{\stein}|_{\KR_w}$ does not vanish, thus $s^{\stein}$ does not vanish anywhere. Therefore, $\omega^{\min, \stein, \otimes n} \cong \OO_{\KR_w^{\min, \stein}}$, $\KR_w^{\min, \stein}$ is quasi-affine.
    \end{enumerate}
\end{proof}
\begin{remark}
     In \cite[Theorem 3.5.9]{shen2021ekor}, the authors showed that the hodge bundle $\omega$ on $\KR_w$ is torsioned. If moreover, for each $w$, there exists an integer $N_w$ and a section $s \in \Gamma(\KR_{\leq w}, \omega^{\otimes N_w})$ whose non-vanishing locus is exactly $\KR_w$ (compared with the Hasse invariants in \cite[Corollary I.2.3]{goldring2019strata}), then the assumption and arguments in Proposition \ref{prop: KR are quasi-affine} show that $s$ extends to $\Gamma(\KR^{\min, \stein}_{\leq w}, \omega^{\min, \stein, \otimes N_w})$ whose non-vanishing locus is exactly $\KR^{\min, \stein}_w$, this implies that each $\KR^{\min, \stein}_w$ is affine, thus each $\KR^{\min}_w$ is affine due to \cite[\href{https://stacks.math.columbia.edu/tag/01YQ}{Tag 01YQ}]{stacks-project}.
\end{remark}
     \begin{remark}
         Since $\KR_w^{\min, \stein}$ (resp. $\KR_{\leq w}^{\min, \stein}$) is normal, $\KR_w^{\min, \stein} \to \KR_w^{\min}$ (resp. $\KR_{\leq w}^{\min, \stein} \to \KR_{\leq w}^{\min}$) is a universal homeomorphism, then the notions of irreducible components and connected components of $\KR_w^{\min}$ (resp. $\KR_{\leq w}^{\min}$) coincide. 
     \end{remark}

     The proof of Lemma \ref{lemma: KR intersects with minimal KR} is based on \cite[\S 6, \S 7]{gortz2012supersingular}. It follows from following lemmas:
     \begin{lemma}\label{lemma: cap is union of connected KR}
     	Let $w, w' \in \Adm(\lrbracket{\mu})$, $Y$ be an irreducible component of $\KR_{I, w}$. If the closure $\overline{Y}$ intersects with $\KR_{I, w'}$, then the intersection is a union of connected components of $\KR_{I, w'}$.
     \end{lemma}
     \begin{proof}
     	Since $\KR_{I, \leq w}$ is normal, any two connected components $Y_1, Y_2 \subset \KR_{I, w}$ are topologically disjoint, i.e. $\overline{Y_1}\cap \overline{Y_2}=\emptyset$, thus $\ovl{Y} \subset \KR_{I, \leq w}$ is a connected component, thus is open and closed. By the closure relations of KR strata, $w'\leq w$, the closure $\KR_{I, \leq w}$ contains the entire $\KR_{I, w'}$. Thus $\ovl{Y} \cap \KR_{I, w'}$ is open and closed, thus a union of connected components.
     \end{proof}
     \begin{lemma}\label{lemma: Y neq Y closure}
     	Keep same notations and assumptions as in Lemma \ref{lemma: KR intersects with minimal KR}. Let $w\in\Adm(\lrbracket{\mu})$ such that $l(w)\neq 0$, $Y$ be an irreducible component of $\KR_{I, w}$, then $Y\neq \overline{Y}$.
     \end{lemma}
     \begin{proof}
     	We prove this lemma by contradiction. Assume $Y=\overline{Y}$. First of all, $Y$ is well-positioned due to Proposition \ref{proposition: open-closed subschemes are well positioned}, and $Y^{\min}$ is open and closed in $\KR^{\min}_w$. Consider the Stein factorization $Y^{\tor}_{\Sigma} \to Y^{\min, \stein} \to Y^{\min}$, $Y^{\min, \stein}$ is open and closed in $\KR^{\min, \stein}_w$, thus is quasi-affine due to Proposition \ref{prop: KR are quasi-affine}. On the other hand, $Y^{\min} = \ovl{Y}^{\min}$ is closed in $\Shumm{I}$, thus is proper. Therefore, $Y^{\min, \stein}$ is a proper and quasi-affine scheme of finite-type over $k$, thus is of $\dim 0$, which contradicts with the assumption $\dim Y = l(w) > 0$, here we use the dimension conjecture of KR strata stated in \cite[Remark 3.2]{he2017stratifications}, which is true due to the proof of \cite[Lemma 4.3.7]{daniels2024igusa}.
     \end{proof}
     \begin{proof}[Proof of \ref{lemma: KR intersects with minimal KR}]
     	Let $w\in\Adm(\lrbracket{\mu})$ be arbitary element with $l(w) \neq 0$, $Y$ be any irreducible component of $\KR_{I, w}$. By Lemma \ref{lemma: Y neq Y closure}, $\overline{Y}\neq Y$, there must be a $w'\in\Adm(\lrbracket{\mu})$ such that $w'\leq w$ and $\overline{Y}\cap\KR_{I, w'}\neq\emptyset$. By Lemma \ref{lemma: cap is union of connected KR}, $\overline{Y}\cap\KR_{I, w'}$ is a union of connected components (which are also irreducible components) of the smooth stratum $\KR_{I, w'}$, let $Y'$ be one of such connected component. Since $\overline{Y'} \subset \overline{Y}$ and $l(w')<l(w)$, we replace $Y$ with $Y'$ and replace $w$ with $w'$, and apply this process until $l(w')=0$. Since $\tau_{\lrbracket{\mu}}$ is the unique length $0$ element in $\Adm(\lrbracket{\mu})$, we are done.
     \end{proof}
     \begin{remark}
        Under the setting of Kisin-Pappas integral models, Proposition \ref{prop: KR are quasi-affine} and Lemma \ref{lemma: Y neq Y closure} also hold for EKOR strata at arbitary parahoric level, since in this case EKOR strata are all quasi-affine due to \cite[Theorem 3.5.9]{shen2021ekor}. However, Lemma \ref{lemma: cap is union of connected KR} may not hold in this case, as we lack prior knowledge that every stalk of $\EKOR_{K, \leq w}$ is an integral domain. Closures of KR strata $\KR_{I, \leq w}$ are normal, but clousres of EKOR (even EO) strata might not be normal, see \cite{koskivirta2018normalization}.
     \end{remark}

     We state the connectedness of EKOR strata following \cite{van2024mod}.
     \begin{definition}\label{def: connectedness of strata}
     	Let $Y$ be a locally closed subset (resp. subscheme) of $\Shum{K, \bar{s}}$, we say $Y$ is connected if the induced morphism $\pi_0(Y)\to\pi_0(\Shum{K, \bar{s}})$ is a bijection.
     \end{definition}

     Since \cite[Conjecture 4.3.1]{van2024mod} is proved in Corollary \ref{corollary: Hoften 3.7.5 removed} under the condition $G^{\ad}$ is $\Q$-simple or has no $\PGL_2$-factor, we have
     \begin{theorem}\cite[Theorem 4.5.2]{van2024mod}
         Under the setting of Kisin-Pappas integral model, assume $G^{\ad}$ is $\Q$-simple or has no $\PGL_2$-factor, $\GG = \GGc$. Then every EKOR stratum which intersects with a $\Q$-non basic Newton stratum (\cite[Definition 4.5.1]{van2024mod}) is connected.
     \end{theorem}

     Also, under the Kisin-Pappas integral model, assume $\GG = \GGc$, then in \cite[Theorem 2(2)]{van2024mod} we could replace the condition $G_{\Q_p}$ is unramified or $\Shum{K}$ is proper with the condition that $G^{\ad}$ is $\Q$-simple or has no $\PGL_2$-factor.

%% file: sections/change_of_parahorics.tex
 We work with Pappas-Rapoport integral models and assume He-Rapoport \cite{he2017stratifications} axioms hold. Here in the index set we use $\mu$ instead of $\mu^{-1}$ to be compatible with \cite{he2017stratifications}. 
 
 Let $(G, X)$ be a Hodge-type Shimura datum. Let $K_{1, p} \subset K_{2, p}$ be two quasi-parahoric subgroups of $G(\Q_p)$, choose $K_1^p \subset K_2^p$, let $K_1 = K_{1, p}K_1^p$, $K_2 = K_{2, p}K_2^p$, then we have a morphism $\pi: \Shum{K_1} \to \Shum{K_2}$ which extends the finite morphism $\shu{K_1} \to \shu{K_2}$, due to \cite[Corollary 4.3.2]{pappas2024p} and \cite[Corollary 4.1.10]{daniels2024conjecture}.
\begin{theorem}{\cite[Theorem 4.5, Proposition 4.11]{mao2025boundary}}\label{Theorem: change of parahoric}
    Let $K_{1, p} \subset K_{2, p}$ be two quasi-parahoric subgroups of $G(\Q_p)$.
    \begin{enumerate}
        \item The morphism $\pi: \Shum{K_1} \to \Shum{K_2}$ extends (uniquely) to a proper morphism $\pi^{\min}: \Shumm{K_1} \to \Shumm{K_2}$ and $\pi^{\min, -1}(\Shum{K_2}) = \Shum{K_1}$ (thus $\pi$ is proper).
        \item Let $\Sigma_1$ be induced by $\Sigma_2$, then we have a proper morphism $\shuc{K_1, \Sigma_1} \to \shuc{K_2, \Sigma_2}$ extending the projection $\shu{K_1} \to \shu{K_2}$, then there is a (unique) proper morphism $\pi^{\tor}: \Shumc{K_1, \Sigma_1} \to \Shumc{K_2, \Sigma_2}$ extending both $\shuc{K_1, \Sigma_1} \to \shuc{K_2, \Sigma_2}$ and $\Shum{K_1} \to \Shum{K_2}$, and makes the following diagram commutes:
\[\begin{tikzcd}
	{\Shumc{K_1, \Sigma_1}} & {\Shumc{K_2, \Sigma_2}} \\
	{\Shumm{K_1}} & {\Shumm{K_2}}
	\arrow["{\pi^{\tor}}", from=1-1, to=1-2]
	\arrow["{\oint_{K_1, \Sigma_1}}", from=1-1, to=2-1]
	\arrow["{\oint_{K_2, \Sigma_2}}", from=1-2, to=2-2]
	\arrow["{\pi^{\min}}", from=2-1, to=2-2]
\end{tikzcd}\]
\item The morphisms $\pi^{\tor}$, $\pi^{\min}$ satisfy the properties list in Proposition \ref{proposition: functorialities on toroidal compactifications}, with \emph{finite morphisms} replaced with \emph{proper morphisms}, or more precisely, see \cite[Proposition 4.11]{mao2025boundary} for details.
        
    \end{enumerate}
\end{theorem}

Repeat the proof of Proposition \ref{lemma: pullback of well-positioned is well-positioned}, and note that the assumption \ref{assumption: C to Z} holds in this case, due to Proposition \ref{remark: when the assumption C to Z is true}, we have following lemma:
\begin{lemma}\label{lem: pullback of well position, change of parahoric}
    Let $Y_2 \subset \Shum{K_2, T}$ be a well-positioned subset (resp. subscheme) with respect to $Y_2^{\natural} = \lrbracket{Y^{\natural}_2(\Phi_2)}_{\Phi_2 \in \Cusp_{K_2}(G, X)}$, then $Y_1:= \pi^{-1}(Y_2) \subset \Shum{K_1, T}$ is a well-positioned subset (resp. subscheme) with with respect to $Y_1^{\natural} = \lrbracket{Y^{\natural}_1(\Phi_1)}_{\Phi_1 \in \Cusp_{K_1}(G, X)}$, where $Y^{\natural}_1(\Phi_1)$ is the pullback of $Y^{\natural}_2(\Phi_2)$ under the restriction of $\pi^{\min}$: $\Zb_1(\Phi_1) \to \Zb_2(\Phi_2)$, where $\Phi_2 = \pi_*\Phi_1$. Moreover, $Y^{\tor}_{1, \Sigma_1}$ (resp. the underlying set of $Y^{\min}_1$) is the pullback (resp. preimage) of $Y^{\tor}_{2, \Sigma_2}$ (resp. the underlying set of $Y^{\min}_2$).
\end{lemma}

\begin{proposition}\label{prop: extension of pullback, well-positioned}
    Let $Y_1 \subset \Shum{K_1, T}$ and $Y_2 \subset \Shum{K_2, T}$ be well-positioned subscheme such that $\pi: \Shum{K_1} \to \Shum{K_2}$ induces a morphism $\pi_Y: Y_1 \to Y_2$, then $\pi$ extends to morphisms 
    \begin{equation}\label{eq: ext, min tor, level}
     \pi^{\min}_Y: Y^{\min}_1 \to Y^{\min}_2 \quad \pi^{\tor}_Y: Y^{\tor}_{1, \Sigma_1} \to Y^{\tor}_{2, \Sigma_2}.
    \end{equation}
    \begin{enumerate}
       \item If $\pi_Y$ is proper, then $\pi^{\min}_Y$ and $\pi^{\tor}_Y$ are proper.
        \item  If $\pi_Y$ is proper surjective, then $\pi^{\min}_Y$ and $\pi^{\tor}_Y$ are proper surjective.
        \item If $\pi_Y$ is smooth (resp. \'etale), then $\pi^{\tor}_Y$ is smooth (resp. \'etale).
        \item If $\pi_Y^{\min}$ is quasi-finite (resp. finite), then $\pi_Y^{\tor}$ is quasi-finite (resp. finite).
    \end{enumerate}
\end{proposition}
\begin{proof}
 Let $\Bar{\pi}_Y: \ovl{Y}_1 \to \pi^{-1}(\ovl{Y}_2) \to \ovl{Y}_2$ be the restriction of $\pi$, it is a proper morphism. Note that $\pi(\ovl{Y}_1) \subset \ovl{Y}_2$, and $\pi^{-1}(Y_{2, 0}) \cap \ovl{Y}_1 \subset Y_{1, 0}$. The proper morphism $\pi^{\min}$ restricts to a proper morphism:
     \begin{equation}
         \Bar{\pi}_Y^{\min}: \ovl{Y}^{\min}_1 \hookrightarrow \pi^{\min, -1}(\ovl{Y}^{\min}_2) \stackrel{\pi^{\min}}{\longrightarrow}  \ovl{Y}^{\min}_2.
     \end{equation}
   Due to Lemma \ref{lemma: closure of well-pos is well-pos}, $Y_{2, 0}$, $\ovl{Y}_1$, $Y_{1, 0}$ are well-positioned. Due to Lemma \ref{lem: pullback of well position, change of parahoric}, $\pi^{-1}(Y_{2, 0})$ is well-positioned. Due to Lemma \ref{lemma: intersection of closed well-positioned subschemes} and \ref{lemma: subset of well-positioned subschemes},
   \begin{equation}\label{eq: pi maps complements in a good way}
       \pi^{-1}(Y_{2, 0}^{\min}) \cap \ovl{Y}_1^{\min} = (\pi^{-1}(Y_{2, 0}) \cap \ovl{Y}_1)^{\min} \subset Y_{1, 0}^{\min},
   \end{equation}
      then $\pi^{\min}$ induces a map $Y^{\min}_1 \to Y^{\min}_2$, which we denote by $\pi_Y^{\min}$. Similarly, we have $\pi^{\tor}_Y$.

 $(1)$: We first claim $\Bar{\pi}_Y^{-1}(Y_2) = Y_1$ if and only if $\pi_Y$ is proper: Note that $\pi_Y: Y_1 \to Y_2$ factors as the open dense embedding $j: Y_1 \hookrightarrow \bar{\pi}_Y^{-1}(Y_2)$ and the proper morphism $\Bar{\pi}_Y: \bar{\pi}_Y^{-1}(Y_2) \to Y_2$. If $\Bar{\pi}_Y^{-1}(Y_2) = Y_1$, then $\pi_Y$ is proper. If $\pi_Y$ is proper, then $j$ has closed image due to \cite[\href{https://stacks.math.columbia.edu/tag/01W6}{Tag 01W6}]{stacks-project}. Since $j$ is open dense, then $j$ is an isomorphism, $\Bar{\pi}_Y^{-1}(Y_2) = Y_1$. 
 
 We show $\pi^{\min}_Y$ is proper under this condition: Since $\Bar{\pi}_Y^{\min}$ is a proper morphism, it suffices to show $\Bar{\pi}_Y^{\min, -1}(Y^{\min}_2) = Y^{\min}_1$. Recall $\ovl{Y}_1 = Y_1 \sqcup Y_{1, 0}$, $\ovl{Y}_2 = Y_2 \sqcup Y_{2, 0}$, $\Bar{\pi}_Y^{-1}(Y_{2, 0}) \subset Y_{1, 0}$. Since $\Bar{\pi}_Y^{-1}(Y_2) = Y_1$ by assumption, then $\Bar{\pi}_Y^{-1}(Y_{2, 0}) = Y_{1, 0}$, and $Y_{1, 0}^{\min} = \Bar{\pi}_Y^{\min, -1}(Y_{2, 0}^{\min})$ by Lemma \ref{lem: pullback of well position, change of parahoric}. By definitions, $\ovl{Y}^{\min}_1 = Y^{\min}_1 \sqcup Y_{1, 0}^{\min}$, $\ovl{Y}^{\min}_2 = Y^{\min}_2 \sqcup Y_{2, 0}^{\min}$, thus $\Bar{\pi}_Y^{\min, -1}(Y^{\min}_2) = Y^{\min}_1$.

 $(2)$:  If $\pi_Y(Y_1) = Y_2$, then $\pi(\ovl{Y}_1) = \ovl{Y}_2$, $\Bar{\pi}_Y^{\min}(\ovl{Y}^{\min}_1) = \ovl{Y}^{\min}_2$. Since $\pi_Y$ is proper, from last paragraph, we have $\Bar{\pi}_Y^{\min, -1}(Y_{2, 0}^{\min}) = Y_{1, 0}^{\min}$ and $\Bar{\pi}_Y^{\min, -1}(Y^{\min}_2) = Y^{\min}_1$, then $\Bar{\pi}_Y^{\min}(Y^{\min}_1) = Y^{\min}_2$.

 $(3)$: Let $x_1 \in Y_{\mathcal{Z}_1([\Phi_1, \sigma_1])}$, $x_2 = \pi(x_1) \in Y_{\mathcal{Z}_2([\Phi_2, \sigma_2])}$. We identify $y_i \in Y^{\natural}_{i, \Xi_{i, \sigma_i}(\Phi_i)}$ as $x_i$ under $Y^{\natural}_{i, \Xi_{i, \sigma_i}(\Phi_i)}\cong Y_{\mathcal{Z}_i([\Phi_i, \sigma_i])}$. Apply \cite[Theorem 2.3.2(7)]{lan2018compactifications}, where we use Artin approximations to $W_{Y_1}$ and $W_{Y_2}$ in a compatible way, then we have the left diagram which restricts to the right diagram, where horizontal morphisms are \'etale.
 
\[\begin{tikzcd}
	{(Y_{1, \Sigma_1}^{\tor}, x_1)} & {(\ovl{U}_{Y_1}, u_1)} & {(Y_{1, \Xi_1(\Phi_1)(\sigma_1)}^{\natural},y_1)} & {Y_1} & {U_{Y_1}} & {Y_{1, \Xi_1(\Phi_1)}^{\natural}} \\
	{(Y_{2, \Sigma_2}^{\tor}, x_2)} & {(\ovl{U}_{Y_2}, u_2)} & {(Y_{2, \Xi_2(\Phi_2)(\sigma_2)}^{\natural}, y_2)} & {Y_2} & {U_{Y_2}} & {Y_{2, \Xi_2(\Phi_2)}^{\natural}}
	\arrow["{\pi^{\tor}_Y}", from=1-1, to=2-1]
	\arrow[from=1-2, to=1-1]
	\arrow[from=1-2, to=1-3]
	\arrow[from=1-2, to=2-2]
	\arrow["{\pi^{\natural}_{\Xi(\Phi)(\sigma)}}", from=1-3, to=2-3]
	\arrow["{\pi_Y}", from=1-4, to=2-4]
	\arrow[from=1-5, to=1-4]
	\arrow[from=1-5, to=1-6]
	\arrow[from=1-5, to=2-5]
	\arrow["{\pi^{\natural}_{\Xi(\Phi)}}", from=1-6, to=2-6]
	\arrow[from=2-2, to=2-1]
	\arrow[from=2-2, to=2-3]
	\arrow[from=2-5, to=2-4]
	\arrow[from=2-5, to=2-6]
\end{tikzcd}\]

 To show $\pi^{\tor}_Y$ is smooth (resp. \'etale) at $(x_1, x_2)$, it suffices to show $\pi^{\natural}_{\Xi(\Phi)(\sigma)}$ is smooth (resp. \'etale) at $(y_1, y_2)$. Let us denote $\pi^{\natural}_{C(\Phi)}: Y_{1, C_1(\Phi_1)}^{\natural} \to Y_{2, C_2(\Phi_2)}^{\natural}$. Since $Y_{i, \Xi_i(\Phi_i)}^{\natural} \hookrightarrow Y_{i, \Xi_i(\Phi_i)(\sigma_i)}^{\natural}$ are toroidal embedding over $Y_{i, C_i(\Phi_i)}^{\natural}$ under a split torus over $\Z$, thus $Y_{i, \Xi_i(\Phi_i)}^{\natural}$ and $Y_{i, \Xi_i(\Phi_i)(\sigma_i)}^{\natural}$ are smooth over $Y_{i, \Xi_i(\Phi_i)(\sigma_i)}^{\natural}$ with same relative dimensions, it suffices to show $\pi^{\natural}_{C(\Phi)}$ is smooth. Since $\bigcup U_{Y_i} \to Y^{\natural}_{i, C_i(\Phi_i)}$ are surjective, then $\pi_Y$ is \'etale implies that $\pi^{\natural}_{C(\Phi)}$ is smooth.

 $(4)$. If $\pi_Y^{\min}$ is finite, then $\pi_Y$ thus $\pi_Y^{\tor}$ are proper. If moreover $\pi_Y^{\tor}$ is quasi-finite, then it is finite. Therefore, it suffices to prove $\pi_Y^{\min}$ is quasi-finite implies that $\pi_Y^{\tor}$ is quasi-finite. Let us denote $\pi_{(\ast)}^{(\natural)}: Y_{1, (\ast)_1}^{(\natural)} \to Y_{2, (\ast)_2}^{(\natural)}$, where $(\ast)$ is $\Zb(\Phi)$, $\Zb^{\bigsur}(\Phi)$, $C(\Phi)$, $\Xi(\Phi)$, $\Xi(\Phi)(\sigma)$, $\Xi_{\sigma}(\Phi)$, $\mathcal{Z}([\Phi])$, $\mathcal{Z}([\Phi, \sigma])$. $\pi_Y^{\min}$ is quasi-finite implies that $\pi_{\mathcal{Z}([\Phi])}$ is quasi-finite, thus $\pi_{\Zb(\Phi)}^{\natural}$ is quasi-finite under the canonical identification $Y^{\natural}_{i, \Zb_i(\Phi_i)} \cong Y_{i, \mathcal{Z}_i([\Phi_i])}$. Since $\Zb_i^{\bigsur}(\Phi_i) \to \Zb_i(\Phi_i)$ are quasi-finite, then $\pi_{\Zb^{\bigsur}(\Phi)}^{\natural}$ is quasi-finite. Due to \cite[Proposition 4.11(3)]{mao2025boundary}, $C_1(\Phi_1) \to C_2(\Phi_2)$ factors through the finite morphism $C_1(\Phi_1) \to C_2(\Phi_2)\times_{\Zb^{\bigsur}_2(\Phi_2)} \Zb^{\bigsur}_1(\Phi_1)$ and the projection $C_2(\Phi_2)\times_{\Zb^{\bigsur}_2(\Phi_2)} \Zb^{\bigsur}_1(\Phi_1) \to C_2(\Phi_2)$. Then $\pi_{C(\Phi)}^{\natural}$ is quasi-finite. Similarly, $\pi_{\Xi(\Phi)}^{\natural}$ is quasi-finite, then $\pi_{\Xi_{\sigma}(\Phi)}^{\natural}$ is quasi-finite. Under the canonical identification $Y^{\natural}_{i, \Xi_{i, \sigma_i}(\Phi_i)} \cong Y_{i, \mathcal{Z}_i([\Phi_i, \sigma_i])}$, then $\pi_{\mathcal{Z}([\Phi, \sigma])}$ is quasi-finite. Since $\pi_{\mathcal{Z}([\Phi, \sigma])}$ is quasi-finite for any $[\Phi, \sigma]$, $\pi_Y^{\tor}$ is quasi-finite.

\end{proof}
   
   Let $K_1 \subset K_2$ (instead of $\Kc_1 \subset \Kc_2$) be two parahoric subgroups with parahoric group schemes $\GG_1$ and $\GG_2$ respectively, $b_1 \in C(\GG_1, \lrbracket{\mu})$, $b_2 \in C(\GG_2, \lrbracket{\mu})$, such that $\pi: \Shum{K_1} \to \Shum{K_2}$ maps $\CE^{b_1}$ to $\CE^{b_2}$. 
   \begin{lemma}\label{lemma: finite surj, central leaves}
       $\pi$ induces a finite surjective morphism:
     \begin{equation}\label{eq: finite surj, central leaves}
        \pi_{b_1}:  \CE^{b_1} \to \CE^{b_2}.
     \end{equation}
   \end{lemma}
   \begin{proof}
       Let $\NE^{[b_1]}$, $\NE^{[b_2]}$ be the unique Newton strata that contain $\CE^{b_1}$, $\CE^{b_2}$ respectively. Due to \cite[Axiom 3]{he2017stratifications}, $[b_1] = [b_2] \in B(G)$, and $\NE^{[b_1]} = \pi^{-1}(\NE^{[b_2]})$. Since central leaves are closed in Newton strata, thus $\pi_{b_1}$ is a proper morphism. Since $\pi_{b_1}$ is quasi-finite and surjective due to \cite[Axiom 4(c)]{he2017stratifications}, thus $\pi_{b_1}$ is finite and surjective.
   \end{proof}
     \begin{corollary}\label{cor: central leaves extends, levels}
         The finite surjection $\pi_{b_1}$ extends to finite surjections:
         \begin{equation}
           \pi_{b_1}^{\min}:   \CE^{b_1, \min} \to \CE^{b_2, \min}, \quad \pi_{b_1}^{\tor}:   \CE^{b_1, \tor}_{\Sigma_1} \to \CE^{b_2, \tor}_{\Sigma_2},
         \end{equation}
         If $\pi_{b_1}$ is \'etale, then $\pi_{b_1}^{\tor}$ is also \'etale.
     \end{corollary}
     \begin{proof}
         We apply Proposition \ref{prop: extension of pullback, well-positioned}. It suffices to show $\pi_{b_1}^{\min}$ is quasi-finite. We apply \cite[Proposition 4.11]{mao2025boundary}, which gives the structures of $\pi^{\min}$ and $\pi^{\tor}$. Fix $[\Phi_2] \in \Cusp_{K_2}(G, X)$, there are only finitely many $[\Phi_1] \in \Cusp_{K_1}(G, X)$ over $[\Phi_2]$. For such $[\Phi_1]$, the proper morphism $\pi^{\min}$ induces a morphism $\Zb_{K_1}(\Phi_1) \to \Zb_{K_2}(\Phi_2)$. To show $\pi_{b_1}^{\min}$ is quasi-finite, it suffices to show each $\CE^{b_1}_{[\Phi_1]} \to \CE^{b_2}_{[\Phi_2]}$ is quasi-finite. Since $\Zb_{K_i}^{\bigsur}(\Phi_i) \to \Zb_{K_i}(\Phi_i)$ are quasi-finite for $i = 1, 2$, to show $\pi_{b_1}^{\min}$ is quasi-finite, it suffices to show the canonical morphism $\CE^{b_1, \natural}_{\Zb_{K_1}^{\bigsur}(\Phi_1)} \to \CE^{b_2, \natural}_{\Zb_{K_2}^{\bigsur}(\Phi_2)}$ induced by $\pi^{\min}$ is quasi-finite. 

        Since $\Zb_{K_1}^{\bigsur}(\Phi_1) \to \Zb_{K_2}^{\bigsur}(\Phi_2)$ is again the change-of-parahoric morphism between integral models of Hodge-type Shimura varieties, see \cite[Proposition 4.11(2)]{mao2025boundary}. Due to Proposition \ref{prop: Newton central are well-posiitoned, PR}, both $\CE^{b_i, \natural}_{\Zb_{K_i}^{\bigsur}(\Phi_i)}$ ($i = 1, 2$) are unions of topologically disjoint connected components of central leaves. Such unions are finite unions due to noetherian properties. The finiteness part of \cite[Axiom 4(c)]{he2017stratifications} implies that $\CE^{b_1, \natural}_{\Zb_{K_1}^{\bigsur}(\Phi_1)} \to \CE^{b_2, \natural}_{\Zb_{K_2}^{\bigsur}(\Phi_2)}$ is quasi-finite.

     \end{proof}

     Next, we consider the EKOR strata. Let $K_1 \subset K_2$ be two parahoric subgroups, $w \in \prescript{\KK_1}{}{ \Adm(\lrbracket{\mu})} \subset \wdt{W}$, due to \cite[Proposition 6.11, 6.21]{he2017stratifications}, there exists a subset $\Sigma_{\KK_2}(w) \subset \wdt{W}_{\KK_2}w\wdt{W}_{\KK_2} \cap \prescript{\KK_2}{}{\wdt{W}}$ (the notations here are slightly different), such that $\pi: \Shum{K_1} \to \Shum{K_2}$ induces
     \begin{equation}
         \pi(\EKOR_{K_1, w}) = \bigsqcup_{w' \in \Sigma_{\KK_2}(w)}\EKOR_{K_2, w'}.
     \end{equation}
     In particular, when $w \in \prescript{\KK_2}{}{ \Adm(\lrbracket{\mu})}$, we regard $w$ as an element in $\prescript{\KK_1}{}{ \Adm(\lrbracket{\mu})}$, then $\pi$ induces a finite surjection:
     \begin{equation}\label{eq: finite surj, KR to EKOR}
       \pi_{w}:  \EKOR_{K_1, w} \to \EKOR_{K_2, w}.
     \end{equation}
     \begin{corollary}\label{cor: EKOR extends, level}
         Let $w \in  \prescript{\KK_2}{}{\Adm(\lrbracket{\mu})}$, then the finite surjection $\pi_{w}$ extends to finite surjections:
  \begin{equation}
      \pi_{w}^{\min}:  \EKOR_{K_1, w}^{\min} \to \EKOR_{K_2, w}^{\min}, \quad \pi_{w}^{\tor}:  \EKOR_{K_1, w, \Sigma_1}^{\tor} \to \EKOR_{K_2, w, \Sigma_2}^{\tor}
  \end{equation}
       If $\pi_w$ is \'etale, then $\pi^{\tor}_w$ is \'etale.
     \end{corollary}
     \begin{proof}
        We apply Proposition \ref{prop: extension of pullback, well-positioned}. It suffices to show $\pi_w^{\min}$ is quasi-finite. Both $\EKOR_{K_1, w} = \bigsqcup_{i \in I}\CE^{b_i}$ and $\EKOR_{K_2, w} = \bigsqcup_{j \in J}\CE^{b_j'}$ are unions of central leaves and central leaves are well-positioned due to Proposition \ref{prop: Newton central are well-posiitoned, PR}, thus $\EKOR_{K_1, w}^{\min} = \bigsqcup_{i \in I}\CE^{b_i, \min}$ and $\EKOR_{K_2, w}^{\min} = \bigsqcup_{j \in J}\CE^{b_j', \min}$ due to Lemma \ref{lem: union of well-positioned sets}. Since $\EKOR_{K_1, w} \to \EKOR_{K_2, w}$ is quasi-finite, thus for each $j \in J$, $\pi_{w}^{-1}(\CE^{b'_j})= \sqcup_{i \in I_j}\CE^{b_i}$ is a finite union, i.e. $I_j \subset I$ is a finite set. In particular, $\pi_{ w}^{\min, -1}(\CE^{b'_j, \min})= \sqcup_{i \in I_j}\CE^{b_i, \min}$ due to Lemma \ref{lem: pullback of well position, change of parahoric}. Since $I_j$ is finite for each $j$, and $\CE^{b_i, \min} \to \CE^{b_j', \min}$ is quasi-finite for each $i \in I_j$ due to Corollary \ref{cor: central leaves extends, levels}, then $\pi_{ w}^{\min}$ is quasi-finite.
     \end{proof}